\documentclass{article}
\usepackage{amsfonts}
\usepackage{nebula}
\usepackage{makeidx}
\usepackage{multicol}
\usepackage{amscd}

\makeindex

\newtheorem{Theorem}{Theorem}[chapter]
\newtheorem{Lemma}[Theorem]{Lemma}
\newtheorem{Corollary}[Theorem]{Corollary}
\newtheorem{Proposition}[Theorem]{Proposition}

\newtheorem{Conjecture}[Theorem]{Conjecture}

\newtheorem{Definition}[Theorem]{Definition}
\newtheorem{Example}[Theorem]{Example}

\newtheorem{Remark}[Theorem]{Remark}
\def\thesubsection{\thesection.\Alph{subsection}}

\begin{document}

\author{Evgueni Tevelev}
\title{Projectively Dual Varieties}
\maketitle

\pagenumbering{Roman}

\chapter*{Preface}
\markboth{Preface}{Preface}

During several centuries various reincarnations of projective duality have inspired
research in algebraic and differential
geometry, classical mechanics, invariant theory, combinatorics, etc.
On the other hand, projective duality is simply the systematic
way of recovering the projective variety from the set of its tangent hyperplanes.
In this survey we have tried to collect together 
different aspects of projective duality and points of view on it.
We hope, that the exposition is quite informal and requires only a standard
knowledge of algebraic geometry and algebraic (or Lie) groups theory.
Some chapters are, however, more difficult and use the modern
intersection theory and homology algebra.
But even in these cases we have tried to give simple
examples and avoid technical difficulties.

An interesting feature of projective duality is given by the observation
that most important examples carry the natural action
of the Lie group. 
This is  especially true for projective varieties 
that have extremal properties from the point
of view of projective geometry.
We have tried to stress this phenomenon
in this survey and to discuss many variants of it.
However, one aspect is completely omitted -- we are not discussing
the dual varieties of toric varieties and the corresponding
theory of $A$-discriminants. This theory is presented
in the beautiful book~\cite{GKZ2} and we feel no need to
reproduce it.

Parts of this survey were written during my 
visits to the Erwin Shroedinger Institute in Vienna
and Mathematic Institute in Basel. I would like to thank
my hosts for the warm hospitality. I have discussed
the contents of this book with many people, including
E.~Vinberg,
V.~Popov, A.~Kuznetsov, S.~Keel, H.~Kraft, D.~Timashev,
D.~Saltman, P.~Katsylo,  and learned a lot from them. 
I am especially grateful to F.~Zak for providing a lot of
information on projective duality
and other aspects of projective geometry.

\vskip 0.2cm
\noindent Edinburgh, November 2001 \hfill {\it Evgueni Tevelev}

\tableofcontents

\newpage
\pagenumbering{arabic}
\setcounter{page}{1}  

\def\A{{\cal A}}
\def\ad{\mathop{\rm ad}}
\def\Ad{\mathop{\rm Ad}}
\def\Ann{\mathop{\rm Ann}\nolimits}
\def\Aut{\mathop{\rm Aut}}
\def\b{{\goth b}}
\def\Bbb{\mathbb}
\def\bigover#1#2{{{\textstyle #1\strut\over\strut\textstyle#2}}}
\def\Base{\mathop{\rm Base}}
\def\Bil{\mathop{\rm Bil}}
\def\Ca{{\Bbb C}\hbox{\rm a}}
\def\C{{\Bbb C}}
\def\cal{\mathcal}
\def\char{\hbox{\rm char}}
\def\Cl{\hbox{\rm Cl}}
\def\Coker{\mathop{\rm Coker}}
\def\codim{\mathop{\rm codim}\nolimits}
\def\Cone{\mathop{\rm Cone}}
\def\const{\mathop{\rm const}}
\def\contr{\mathop{\rm contr}}
\def\corank{\mathop{\rm corank}}
\def\defect{\mathop{\rm def}}
\def\D{{\cal D}}
\def\Det{\mathop{\rm Det}}
\def\Der{\mathop{\rm Der}}
\def\diag{\mathop{\rm diag}}
\def\div{\mathop{\rm div}}
\def\Div{\mathop{\rm Div}}
\def\ds{\displaystyle}
\def\dualA{{\A_{n,k}^0}^{\!\!\!\!*}}
\def\dual#1{{{#1}^{*}}}
\def\E{{\cal E}}
\def\End{\mathop{\rm End}}
\def\endproof{\qed\smallskip}
\def\eps{\varepsilon}
\def\F{{\cal F}}
\def\fchoose#1#2{\left\{\matrix{#1\cr #2}\right\}}
\def\FF{{\Bbb F}{\Bbb F}}
\def\g{\goth{g}}
\def\gl{\goth g\goth l}
\def\GL{\mathop{\rm GL}\nolimits}
\def\goth#1{{\mathfrak #1}}
\def\Gr{\mathop{\rm Gr}\nolimits}
\def\h{{\goth h}}
\def\H{{\Bbb H}}
\def\hht{\mathop{\rm ht}}
\def\Hes{\mathop{\rm Hes}}
\def\Hom{\mathop{\rm Hom}}
\def\Id{\mathop{\rm Id}}
\def\Im{\mathop{\rm Im}}
\def\k{{\goth k}}
\def\Ker{\mathop{\rm Ker}}
\def\l{{\goth l}}
\def\L{{\cal L}}
\def\la{\mathop{\leftarrow}}
\def\Lag{\mathop{\rm Lag}}
\def\Lie{\mathop{\rm Lie}}
\def\M{{\cal M}}
\def\Mat{\mathop{\rm Mat}\nolimits}
\def\mod{\,\mathrel{\rm mod}\,}
\def\n{{\goth n}}
\def\O{{\cal O}}
\def\OO{{\Bbb O}}
\def\ov{\overline}
\def\p{\goth{p}}
\def\P{\Bbb P}
\def\PGL{\mathop{\rm PGL}}
\def\Pic{\mathop{\rm Pic}}
\def\pr{\mathop{\rm pr}} 
\def\proof{\noindent{\bf Proof.\ }}
\def\PP{{\cal P}}
\def\Q{\Bbb Q}
\def\QQ{{\cal Q}}
\def\R{\Bbb R}
\def\ra{\mathop{\rightarrow}}
\def\rank{\mathop{\rm rank}}
\def\rk{\mathop{\rm rk}}
\def\res{\mathop{\rm res}}
\def\S{{\cal S}}
\def\sketch{\noindent{\bf Sketch of the proof.\ }}
\def\schoose#1#2{\left[\matrix{#1\cr #2}\right]}
\def\Sec{\mathop{\rm Sec}}
\def\Sing{\mathop{\rm Sing}}
\def\sl{{\goth s\goth l}}
\def\SL{\mathop{\rm SL}} 
\def\so{{\goth s\goth o}}
\def\SO{\mathop{\rm SO}} 
\def\sp{{\goth s\goth p}}
\def\Sp{\mathop{\rm Sp}} 
\def\Span{\mathop{\rm Span}} 
\def\Spin{\mathop{\rm Spin}} 
\def\SS{{\Bbb S}}
\def\ssl{{\goth s\goth l}}
\def\Supp{\mathop{\rm Supp}} 
\def\t{{\goth t}}
\def\Tan{\mathop{\rm Tan}}
\def\TanStar{\mathop{\rm TanStar}}
\def\text#1{\hbox{#1}}
\def\Tr{\mathop{\rm Tr}}
\def\Tran{\mathop{\rm Tran}}
\def\uu{{\goth u}}
\def\un{\underline}
\def\X{\text{{\bf X}}}
\def\z{{\goth z}}
\def\Z{\Bbb Z}

\chapter{Dual Varieties}

\section*{Preliminaries} 
The projective duality gives a remarkably simple
method to recover any projective variety from
the set of its tangent hyperplanes.
In this chapter we recall this classical notion,
give the proof of the Reflexivity Theorem and its
consequences, provide a number of examples and motivations, 
and fix the notation to be used throughout the book.
The exposition is fairly standard and classical.
In the proof of the Reflexivity Theorem we follow \cite{GKZ2}
and deduce this result from the classical theorem
of symplectic geometry saying that any conical Lagrangian
subvariety of the cotangent bundle is equal to
some conormal variety.

\section{Definitions and First Properties}

For any finite-dimensional complex vector space $V$ we denote
by $\P(V)$ its projectivization, that is, the set of $1$-dimensional
subspaces. For example, if $V=\C^{n+1}$ is a standard complex
vector space then $\P^n=\P(\C^{n+1})$ \index{standard projective space}
is a standard complex projective space. A point of $\P^n$
is defined by $(n+1)$ homogeneous coordinates $(x_0:\ldots:x_1)$, $x_i\in\C$,
which are not all equal to $0$ and are considered up to  a scalar 
multiple. \index{homogeneous coordinates}

If $U\subset V$ is a non-trivial linear subspace then $\P(U)$ is a subset
of $\P(V)$, subsets of this form are called projective subspaces.
Projective subspaces of dimension $1$, $2$, or of codimension $1$
are called lines, planes, 
and hyperplanes.\index{projective subspace}\index{hyperplane}

For any vector space $V$ we denote by $\dual V$ the dual vector 
space,  \index{dual vector space}
the vector space of linear forms on $V$. Points of 
the dual projective space $\dual{\P(V)}=\P(\dual V)$ 
correspond to hyperplanes in $\P(V)$.
Conversely, to any point $p$ of $\P(V)$, we can associate
a hyperplane in $\dual{\P(V)}$, namely the set of all hyperplanes
in $\P(V)$ passing through $p$.
Therefore, $\P(V)^{**}$ is naturally identified with $\P(V)$.
Of course, this reflects  nothing else but a usual canonical isomorphism
$V^{**}=V$.

To any vector subspace $U\subset V$ we associate its 
annihilator \index{annihilator}
$\Ann(U)\subset \dual V$. Namely, $\Ann(U)=\{f\in\dual V\,|\,f(U)=0\}$.
We have $\Ann(\Ann(U))=U$. 
This corresponds to the projective duality
between projective subspaces in $\P(V)$ and $\dual{\P(V)}$:
for any projective subspace $L\subset\P(V)$ we 
denote by $\dual L\subset\dual{\P(V)}$
its dual projective subspace, parametrizing
all hyperplanes that contain~$L$.  \index{dual projective subspace}

Remarkably, the projective duality between projective subspaces
in $\P^n$ and $\dual{\P^n}$
can be extended
to the involutive correspondence between irreducible
algebraic subvarieties in $\P^n$ and $\dual{\P^n}$.

First, suppose that $X\subset\P^n$ is a smooth irreducible
algebraic subvariety. For any $x\in X$, we denote by 
$\hat T_xX\subset\P^n$ an embedded projective tangent 
space.  \index{embedded projective tangent space}
More precisely, if $X\in\P(V)$ is any projective variety
then we define the cone $\Cone(X)\subset V$ over it 
as a conical variety formed by all lines $l$ such that $\P(l)\in X$.
If $x\in X$ is a smooth point then any non-zero point $v$
of the corresponding line is a smooth point of $\Cone(X)$
and $\hat T_x(X)$ is defined as $\P(T_v\Cone(X))$,
where $T_v\Cone(X)$ is a tangent space of $\Cone(X)$ at $v$
considered as a linear subspace of $V$ 
(it does not depend on a choice of $v$).
For any hyperplane $H\subset\P^n$, we say 
that $H$  \index{tangent hyperplane}
is tangent to $X$ at $x$ if $H$ contains $\hat T_xX$.
We define the dual variety $\dual X\subset\dual{\P^n}$
as the set of all hyperplanes tangent to $X$. \index{dual variety}

In other words, a hyperplane $H$ belongs to $\dual X$
if and only if the intersection $X\cap H$ (regarded as a scheme)
is singular (is not a smooth algebraic variety).
In most parts of this book we shall be interested
only in dual varieties of smooth varieties,
so this description of $\dual X$ will be sufficient.
However, with this definition we can not expect the
duality $X^{**}=X$ because $\dual X$ can be singular
(and in most interesting cases it is actually singular).
So we should define $\dual X$ for a singular $X$ as well.
There are in fact two possibilities: first, we can imitate
the previous definition and consider embedded tangent spaces
at all points, not necessarily smooth.
But it turns out that the dual variety defined in this fashion
does not have good properties, for example
it can be reducible. The better way is to pick only
the `main' component of the dual variety.

\begin{Definition}\rm
Let $X\subset \P^n$ be an irreducible projective variety.
A hyperplane $H\subset\P^n$ is called tangent to $X$ 
if it contains an embedded tangent space $\hat T_xX$
at some smooth point $x\in X$.
The closure of the set of all tangent hyperplanes
is called the dual variety $\dual X\subset\dual{\P^n}$.
\end{Definition}

We shall discuss the Reflexivity Theorem $X^{**}=X$
and its consequences in the next section.
First we shall establish some simple properties of dual varieties
and give further definitions.

\begin{Definition}\rm
Let $X\subset \P^n$ be an irreducible projective variety
with the smooth locus $X_{sm}$.
Consider the set $I_X^0\subset \P^n\times\dual{\P^n}$
of pairs $(x,H)$ such that $x\in X_{sm}$ and $H$ is the hyperplane
tangent to $X$ at $x$.
The Zariski closure $I_X$ of $I_X^0$ is called
the conormal variety of $X$.  \index{conormal variety}
\end{Definition}

The projection $\pr_1:\,I^0_X\to X_{sm}$ makes $I_X^0$
into a bundle over $X_{sm}$ whose fibers are projective
subspaces of dimension $n-\dim X-1$.
Therefore, $I^0_X$ and $I_X$ are irreducible
varieties of dimension $n-1$.
By definition, $\dual X$ is the image of the projection
$\pr_2:\,I_X\to\dual{\P^n}$.
Therefore, we have the following proposition:

\begin{Proposition}
$\dual X$ is an irreducible variety.
\end{Proposition}

Moreover, since $\dim I_X=n-1$, we can expect that in `typical' cases
$\dual X$ is a hypersurface. Having this in mind,
we give the following definition:

\begin{Definition}\rm \index{defect of a projective variety}
The number $\codim_{\dual{\P^n}}\dual X-1$ is called
the defect of $X$, denoted by $\defect X$.
\end{Definition}

Typically, $\defect X=0$. In this case $\dual X$ is defined
by an irreducible homogeneous polynomial $\Delta_X$.

\begin{Definition}\rm \index{discriminant}
$\Delta_X$ is called the discriminant of $X$.
\end{Definition}

If $\defect X>0$ then, for convenience, we set $\Delta_X=1$.
Clearly, $\Delta_X$ is defined only up to a scalar multiple.
Roughly speaking, the study of dual varieties and discriminants
includes the following $3$ steps:

\begin{itemize}
\item To define some nice natural class of projective varieties $X$.
\item To find all exceptional cases when $\defect X>0$.
\item In the remaining cases, to say something about $\Delta_X$ or $\dual X$.
\end{itemize}

In the last step the minimal program is to determine
the degree of $\Delta_X$, the maximal program is to determine $\Delta_X$
as a polynomial, or at least to describe its monomials.
Another interesting problem is to describe singularities of $\dual X$.

\begin{Example}\rm
The most familiar example of a discriminant $\Delta_X$ is, of course,
the discriminant of a binary form. In order to show
that it actually coincides with some $\Delta_X$ we need first
to give an equivalent definition of $\Delta_X$.
Suppose that $x_1,\ldots,x_k$ are some local coordinates
on $\Cone(X)\subset V$. Any $f\in V^*$, a linear form on $V$,
being restricted to $\Cone(X)$ becomes an algebraic function
in $x_1,\ldots,x_k$. Then $\Delta_X$ is just an irreducible
polynomial, which vanishes at $f\in V^*$ whenever the
function $f(x_1,\ldots,x_k)$ has a multiple root, that is, 
vanishes at some $v\in\Cone(X)$, $v\ne0$, together 
with all first derivatives $\partial f/\partial x_i$.

Consider now the $d$-dimensional projective space $\P^d=\P(V)$
with homogeneous coordinates $z_0,\ldots,z_d$,
and let $X\subset\P^d$ be the Veronese curve \index{Veronese curve}
$$(x^d:x^{d-1}y:x^{d-2}y^2:\ldots:xy^{d-1}:y^d),\ x,y\in\C,\ (x,y)\ne(0,0)$$
(the image of the Veronese embedding $\P^1\subset\P^d$).
Any linear form $f(z)=\sum a_iz_i$ is uniquely determined
by its restriction to the cone $\Cone(X)$, which is a binary form
$f(x,y)=\sum a_ix^{d-i}y^i$. Therefore, $f\in\Cone(\dual X)$
if and only if $f(x,y)$ vanishes at some point $(x_0,y_0)\ne (0,0)$
(so $(x_0:y_0)$ is a root of $f(x,y)$) with its first derivatives
(so $(x_0:y_0)$ is a multiple root of $f(x,y)$).
It follows that $\Delta_X$ is the classical discriminant of  a binary 
form.  \index{discriminant of a binary form}
\end{Example}

\section{Reflexivity Theorem}
In this section we shall prove the Reflexivity Theorem.
We follow the exposition in \cite{GKZ2}, which shows that
there exists a deep connection between projective duality
and symplectic geometry. Other proofs, including the
investigation of a prime characteristic case, could be found, e.g.,
in \cite{Se}, \cite{M}, \cite{Wa}.

\begin{Theorem}\label{ReflexivityTh}$\ $\index{Reflexivity Theorem}
\begin{enumerate}
\item[\rm(a)] For any irreducible projective variety $X\subset\P^n$,
we have $X^{**}=X$. 
\item[\rm(b)] More precisely,
If $z$ is a smooth point of $X$
and $H$ is a smooth point of $\dual X$, then $H$ is tangent to $X$
at $z$ if and only if $z$, regarded as a hyperplane in $\dual{\P^n}$,
is tangent to $\dual X$ at $H$.
\end{enumerate}
\end{Theorem}

The proof will be given in the next section.

\subsection{The Conormal Variety}

We shall need
some standard definitions.

\begin{Definition}\rm
If $X$ is a smooth algebraic variety, then $TX$ denotes the tangent
bundle\index{tangent bundle} of $X$. 
If $Y\subset X$ is a smooth algebraic subvariety
then $TY$ is a subbundle in $TX|_Y$. 
The quotient $TX|_Y/TY$
is called the normal bundle\index{normal bundle} 
of $Y$ in $X$, denoted by $N_YX$.
By taking dual bundles we obtain the cotangent bundle\index{cotangent bundle} 
$T^*X$ and the conormal bundle\index{conormal bundle} $N_Y^*X$.
The conormal bundle can be naturally regarded
as a subvariety of $T^*X$.
\end{Definition}

Recall that 
$I_X\subset\P^n\times\dual{\P^n}$ is the
conormal variety and the dual 
variety $\dual X$ coincides with $\pr_2(I_X)$. 
The projection $\pr_1:\,I_X^0\to X_{sm}$ is a projective bundle,
where $\pr_1$, $\pr_2$ denote the projections of $\P^n\times\dual{\P^n}$
to its factors. More precisely, $\pr_1$ identifies $I_X^0$
with the projectivization $\P(N_{X_{sm}}^*\P^n)$ of the conormal bundle
$N_{X_{sm}}^*\P^n$. 
Indeed, the choice of a hyperplane $H\subset\P^n$
tangent to $X_{sm}$ at $x$ is equivalent to the choice of a hyperplane
$T_xH$ in the tangent space $T_x\P^n$, which contains $T_xX$.
The equation of this hyperplane is an element of 
$N_{X_{sm}}^*\P^n$ at $x$.
The Reflexivity Theorem can be reformulated as follows:
\begin{equation}
I_X=I_{\dual X}.\label{incidcoinc:eq}
\end{equation}
It is more convenient to prove~(\ref{incidcoinc:eq})
by working with vector spaces instead of projective spaces.

We assume that $\P^n=\P(V)$ and $\dual{\P^n}=\P(\dual V)$.
Then we have affine cones 
$Y=\Cone(X)\subset V$ and $\dual Y=\Cone(\dual X)\subset \dual V$.
We denote by $\Lag(Y)$ the closure of the conormal bundle
$N^*_{Y_{sm}}V$ in the cotangent bundle $T^*V$.

The space $T^*V$ is canonically identified with $V\times\dual V$.
Denote by $\pr_1$, $\pr_2$ the projections of this product to its factors.
Then $\dual Y$ coincides with $\pr_2(\Lag(Y))$. Therefore,
(\ref{incidcoinc:eq}) can be reformulated as follows:
\begin{equation}
\Lag(Y)=\Lag(\dual Y),\label{incidconescoinc:eq}
\end{equation}
where we identify $T^*V$ and $T^*\dual V$ with $V\times \dual V$.

Recall that a smooth algebraic variety $M$ is called the symplectic
variety\index{symplectic variety}, 
if $M$ admits a symplectic structure\index{symplectic structure}, that is, 
a differential $2$-form $\omega$ with the following properties:

\begin{itemize}
\item $\omega$ is closed, $d\omega=0$.
\item $\omega$ is non-degenerate, for any $p\in M$
the restriction of $\omega$ on $T_pM$ is a non-degenerate skew product.
\end{itemize}

In this case $\dim M$ is necessarily even.
An irreducible closed subvariety $\Lambda\subset M$ is called 
Lagrangian\index{Lagrangian variety}
if $\dim\Lambda=\dim M/2$ and the restriction of $\omega$
to the smooth locus $\Lambda_{sm}$ vanishes as a $2$-form
(is totally isotropic).

The most important and the most typical (having in mind the Darboux Theorem)
example of a symplectic variety is the cotangent bundle $T^*X$
of a smooth algebraic variety $X$. $T^*X$ carries a canonical 
symplectic structure defined as 
follows.\index{symplectic structure on the cotangent bundle}
Let $(x_1,\ldots,x_n)$ be a local coordinate system in $X$.
Let $\xi_i$ be the fiberwise linear function on $T^*X$
given by the pairing of a $1$-form with the vector field $\partial/\partial x_i$.
Then $(x_1,\ldots,x_n,\xi_1,\ldots,\xi_n)$ forms a local
coordinate system in $T^*X$. 
The form $\omega$ is defined by
\begin{equation}
\omega=\sum_{i=1}^nd\xi_i\wedge dx_i.\label{omegaincoor:eq}
\end{equation}
It easy to give an equivalent definition of $\omega$ 
without any coordinate systems.
We shall define a canonical $1$-form $\nu$ on $T^*X$ and then
we shall take $\omega=d\nu$.
Let $\pi:\,T^*X\to X$ be the canonical projection.
Let $p\in T^*X$ and $v\in T_p(T^*X)$ be a vector tangent to $T^*X$ at $p$.
Then $\nu(v)=p(\pi_*v)$.

An important (and typical, as we shall see) example
of a Lagrangian subvariety in $T^*X$ is obtained as follows.
Let $Y\subset X$ be any irreducible subvariety with smooth locus $Y_{sm}$
and let
$$\Lag(Y)=\overline{N^*_{Y_{sm}}X}$$
be the closure of the conormal bundle of $Z_{sm}$ (in $T^*X$).
Clearly, $\Lag(Y)$ is a conical subvariety (invariant
under dilations of fibers of $T^*X$).

\begin{Theorem}\label{ConicalLagrangian}$\ $
\begin{enumerate}
\item[{\rm (a)}] $\Lag(Y)$ is a Lagrangian subvariety.
\item[{\rm (b)}] Any conical Lagrangian subvariety has the form $\Lag(Y)$
for some irreducible subvariety $Y\subset X$.
\end{enumerate}
\end{Theorem}

\proof
Let us show that $\Lag(Y)$ is Lagrangian.
Clearly, $\dim\Lag(Y)=\dim Y+(\dim X-\dim Y)=\dim T^*X/2$.
So we need only to verify that $\omega|_{\Lag(Y)}=0$.
It is sufficient to consider the smooth locus of $Y$ only.
Let $x_1,\ldots,x_n$ be a local coordinate system on $X$
such that $Y$ is locally defined by equations $x_1=\ldots=x_r=0$.
Then the fibers of the conormal bundle over points of $Y$
are generated by $1$-forms $dx_1,\ldots,dx_r$.
Hence $\xi_{r+1}=\ldots=\xi_n=0$ on $N^*_{Y_{sm}}X$ and 
by (\ref{omegaincoor:eq}) we see that $\omega=0$ on 
$N^*_{Y_{sm}}X$. Therefore, $\Lag(Y)$ is indeed Lagrangian.

Suppose now that $\Lambda\subset T^*X$ is a conical Lagrangian
subvariety. We take $Y=\pr(\Lambda)$, where $\pr:\,T^*X\to X$ is the projection,
and claim that $\Lambda=\Lag(Y)$. 
It suffices to show that $\Lambda\subset\Lag(Y)$,
because $\Lambda$ and $\Lag(Y)$ are irreducible varieties of the 
same dimension.
In turn, to prove that $\Lambda\subset\Lag(Y)$ it suffices
to check that for any smooth point $y\in Y$ the fiber
$\pr^{-1}(y)\cap\Lambda$ is contained in the conormal
space $(N^*_YX)_y$.
Let $\xi$ be any covector from 
$\pr^{-1}(y)\cap\Lambda$. Since $T^*_yX$ is a vector space,
we can regard $\xi$ as a ``vertical'' 
tangent vector to $T^*X$ at a point $y\subset X\subset T^*X$,
where we identify $X$ with the zero section of $\pr$.
Since $\Lambda$ is conical, $\xi\in T_\xi\Lambda$. Therefore,
since $\Lambda$ is Lagrangian, $\xi$ is orthogonal with respect to $\omega$
to any tangent vector from $T_\xi\Lambda$ and, hence,
to any tangent vector $v\in T_yY$.
But by (\ref{omegaincoor:eq}) it is easy to see that
this is equivalent to $\xi\in(N^*_YX)_y$.
\endproof

Now we can prove the Reflexivity Theorem.
We shall prove (\ref{incidconescoinc:eq}).
The identification $T^*V=V\times V^*=T^*V^*$
takes the canonical symplectic structure on $T^*V$
to minus the canonical symplectic structure on $T^*V^*$.
Therefore $\Lag(Y)$ regarded as a subvariety of $T^*V^*$
is still Lagrangian. Moreover, clearly $\Lag(Y)\subset V\times V^*$
is invariant under dilations of $V$ and $V^*$,
therefore, $\Lag(Y)$ is a conical Lagrangian variety of $T^*V^*$.
Therefore, by Theorem \ref{ConicalLagrangian} 
$\Lag(Y)=\Lag(Z)$, where $Z$ is the projection of $\Lag(Y)$
on $V^*$. But this projection coincides with $\dual Y$.
Therefore, $\Lag(Y)=\Lag(\dual Y)$.
The Reflexivity Theorem is proved.

\subsection{Applications of the Reflexivity Theorem}
Suppose that $X\subset\P^n$ and $\dual X\subset\dual{\P^n}$
are projectively dual varieties, 
$I_X=I_{\dual X}\subset \P^n\times\dual{\P^n}$
is the conormal variety.
We have the diagram of projections
$$X\mathop{\longleftarrow}^{\pr_1}I_X\mathop{\longrightarrow}^{\pr_2}\dual X.$$

\begin{Theorem}$\ $
\begin{enumerate}
\item[{\rm(a)}] If $X$ is smooth then $I_X$ is smooth.
\item[{\rm(b)}] If $\dual X$ is a hypersurface then $\pr_2$ is a birational 
isomorphism.
\item[{\rm(c)}] If $X$ is smooth and $\dual X$ is a hypersurface then $\pr_2$
is a resolution of singularities.\index{desingularization of a dual variety}
\end{enumerate}
\end{Theorem}

\proof
The map $\pr_1$ is a projective bundle over a smooth locus
of $X$. Therefore, if $X$ is smooth then $I_X$ is smooth.
If $\dual X$ is a hypersurface then $\dim\dual X=\dim I_X=n-1$.
Since $\pr_2$ is generically a projective bundle, it is
birational.
Finally, (c) follows from (a) and (b).
\endproof

Typically the dual variety $\dual X\subset\dual{\P^n}$ is 
a hypersurface. Namely, we shall see that if $\dual X$ is not
a hypersurface then $X$ is ruled in projective subspaces.\index{ruled variety}

\begin{Definition}\rm
We say that $X$ is ruled in projective subspaces of dimension~$r$
if for any $x\in X$ there exists a projective subspace $L$ of dimension~$r$
such that $x\in L\subset X$. By a standard closedness argument
it is sufficient to check this property only for points $x$
from some Zariski open dense subset $U\subset X$.
\end{Definition}

Recall, that the number $\defect X=\codim_{\dual{\P^n}}\dual X-1$
is called the defect of $X$.

\begin{Theorem}\label{FirstEin}
Suppose that $\defect X=r\ge1$. Then
\begin{enumerate}
\item[{\rm (1)}] $X$ is ruled in projective subspaces
of dimension $r$.
\item[{\rm (2)}] If $X$ is smooth then for any $H\in{\dual X}_{sm}$
the contact locus\index{contact locus} $\Sing(H\cap X)$ is a projective subspace of dimension $r$
and the union of these projective subspaces is dense in $X$.
\end{enumerate}
\end{Theorem}

\proof
By Reflexivity Theorem~\ref{ReflexivityTh}, 
(a) is equivalent to the following:
if $\codim X=r+1$ then $\dual X$ is ruled in projective subspaces of dimension $r$.
The condition for a hyperplane $H$ to be tangent to $X$ at a smooth
point $x$ is that $\hat T_xH$ contains $\hat T_xX$.
For a given $x$, all $H$ with this property form a projective subspace
of dimension $r$. But the set of hyperplanes of $\dual X$ tangent
to $X$ at some smooth point obviously contains 
a Zariski open subset of $\dual X$.

(b) is proved by the same argument involving Reflexivity Theorem.
\endproof

\begin{Example}\label{DualOfTheCurve}\rm\index{dual of the curve}
Suppose that $X\subset\P^n$ is a non-linear curve.
Then $\dual X$ is a hypersurface.
Indeed, $X$ obviously could not contain a projective 
subspace $\P^k$ for $k>0$.
\end{Example}

The following Theorem is also an easy corollary of Reflexivity
Theorem. It allows to find singular points of hyperplane sections
of smooth projective varieties.

\begin{Theorem}
Suppose that $X\subset\P^{n-1}$ is smooth and $\dual X\subset\dual{\P^{n-1}}$
is a hypersurface. Let $z_1,\ldots,z_n$ be homogeneous
coordinates on $\P^{n-1}$ and $a_1,\ldots,a_n$ the dual homogeneous 
coordinates on $\dual{\P^{n-1}}$.
Suppose that $f=(a_1,\ldots,a_n)$ is a smooth point
of $\dual X$. Then the hyperplane section
$\{f=0\}$ of $X$ has a unique singular point with coordinates
given by $({\partial\Delta_X\over\partial a_1}(f):\ldots:
{\partial\Delta_X\over\partial a_n}(f))$.
\end{Theorem}

\proof
Let $H\subset \P^{n-1}$ be the hyperplane corresponding
to $f$. By the Reflexivity Theorem, $H$ is tangent to $X$
at $z$ if and only if the hyperplane in $\P^{n-1}$
corresponding to $z$ is tangent to $\dual X$ at $f$.
Since $\dual X$ is smooth at $f$, such a point $z$
is unique and is given by $z_i={\partial\Delta_X\over\partial a_i}(f)$.
\endproof


\section{Dual Plane Curves}

\subsection{Parametric Representation of the Dual Plane Curve}

Perhaps the most classical example of a dual variety
is the dual curve\index{dual plane curve} 
$\dual C$ of a non-linear plane curve $C\subset\P^2$.
By definition, generic points of $\dual C$ are the tangents to $C$
at smooth points.
In this case the Reflexivity Theorem has a fairly intuitive meaning.
The tangent line $\hat T_p\in\dual{\P^2}$ at a smooth point $p\in C$
is the limit of secants $\overline{pq}$ for $q\in C$, $q\to p$.
Similarly, the point in $\P^2$ that corresponds to 
the tangent to $\dual C\subset\dual{\P^2}$
at a non-singular point $\hat T_p$ is the limit of the
intersection points of the tangents $\hat T_p$ and $\hat T_q$ as $q\to p$.
Of course, this point is $p$.

In fact, it is quite easy to write down a parametric representation
of $\dual C$ using a given parametric representation of $C$.
Let $x,y,z$ be homogeneous coordinates on $\P^2$ and
$p,q,r$ the dual homogeneous coordinates on $\dual{\P^2}$.
We choose the affine chart $\C^2=\{z\ne0\}\subset\P^2$ with 
affine coordinates $x,y$, so the third homogeneous coordinate $z$
is set to be $1$.
The dual chart $\dual{\C^2}\subset\dual{\P^2}$
with coordinates $p,q$ is obtained by setting the third
homogeneous coordinate $r$ in $\dual{\P^2}$ to be $-1$.
Then $\dual{\C^2}$ consists of lines in $\P^2$
not passing through the point $(0,0)\in\C^2\subset\P^2$.
Every such line that meets $\C^2$ is given by 
the affine equation $px+qy=1$, the line in $\P^2$
with coordinates $p=q=0$ is the line ``at infinity''.
Suppose that a local parametric equation of $C$ has the form
$x=x(t)$, $y=y(t)$, where $t$ is a local coordinate on $C$,
and $x(t)$, $y(t)$ are analytic functions.
By definition, the dual curve $\dual C$
has the parametrization $p=p(t)$, $q=q(t)$,
where $p(t)x+q(t)y=1$ is the affine equation
of the tangent line to $C$ at the point $(x(t),y(t))$.
Therefore, we have
\begin{equation}
p(t)={-y'(t)\over x'(t)y(t)-x(t)y'(t)},\quad
q(t)={x'(t)\over x'(t)y(t)-x(t)y'(t)}.\label{eq:paramrepr}
\end{equation}
Applying this formula two times we
obtain Reflexivity Theorem once again.

\begin{Example}\rm
In analysis, there is a well-known duality between
Banach spaces $L_p$ and $L_q$ for ${1\over p}+{1\over q}=1$,
see, e.g., \cite{Ru}.
Let us give an algebraic version of this duality.

Consider the curve $X$ in $\C^2$
given by
$$x^a+y^a=1, \quad a>1,\ a\in\Q.$$
This curve, and its closure in $\P^2$, is usually
called the Fermat curve\index{Fermat curve} (especially if $a\in\Z$).
If $a\not\in\Z$, this equation involves multivalued fractional
power functions, but it is possible to put this equation
into the polynomial form.
In order to find the dual curve, we can use a parametric representation
of $X$ of the form
$$x=t,\ y=\sqrt[a]{1-t^a}.$$
Using (\ref{eq:paramrepr}) we get the parametric representation
of the dual curve as follows:
$$p=t^{a-1},\ q=(1-t^a)^{a-1\over a}.$$
The relation between $p$ and $q$ has the form
$$p^b+q^b=1,\ \hbox{\rm where}\ {1\over a}+{1\over b}=1.$$
\end{Example}

\subsection{The Legendre Transformation and Caustics}
The projective duality is closely related to the Legendre
transformation of classical mechanics, and, therefore, is in some
sense analogous to the duality of the Lagrange and the Hamilton pictures
of classical mechanics. To illustrate this analogy
it will be sufficient to recall the classical
definition of the Legendre transformation of real functions
in one variable.
Details can be found in \cite{Ar}.

Suppose that $y=f(x)$ is a smooth convex real function, $f''(x)>0$.
The Legendre transformation\index{Legendre transformation} 
of the function $f$
is a new function $g$ of a new variable $p$,
which is constructed in the following way.
Consider the line $y=px$.
We take the point $x=x(p)$ at which the graph of $y=f(x)$
has a slope $p$ (so $f'(x(p))=p$), and define $g(p)$ as $g(p)=px(p)-f(x(p))$.
Equivalently, we define $x(p)$ as a unique point, where
the function $F(p,x)=px-f(x)$ has a maximum with respect to $x$
and define $g(p)=F(p,x(p))$.

The Legendre transformation is easily seen to be involutive.
To see how it is related to the dual curve, let us notice
that the straight line $y=G(x,p)=xp-g(p)$ is nothing else but
a tangent line to the graph of $f$ with slope $p$.

To link projective duality and Legendre transformation
we need a notion of a caustic curve\index{caustic}.
To introduce it, let us express the projective duality 
entirely in terms of the projective plane $\P^2$.
By definition, a tangent line to a curve $C$ at some point $x$
is the line that contains $x$ and which is infinitesimally
close to the curve $C$ near $x$.
A point of $\dual{\P^2}$ is a line $l\subset\P^2$.
A curve in $\dual{\P^2}$ is a $1$-parameter family of lines in $\P^2$.
For example, a line in $\dual{\P^2}$ is a pencil $\dual x$
of all lines in $\P^2$ passing through a given point $x\in\P^2$.
The dual curve $\dual C$ is a $1$-parameter family
of tangent lines to $C$.
Suppose now that $C'\subset\dual{\P^2}$ is some curve
(some $1$-parameter family of lines in $\P^2$).
Let us find a geometric interpretation of the dual curve $\dual{C'}\subset\P^2$.
Take some line $l\in C'$. The condition that
$\dual x$ is tangent to $C'$ at $l$ means that the line $l\in\dual{\P^2}$
is a member of a family $C'$ and other lines from $C'$ near $l$
are infinitesimally close to the pencil of lines $\dual x$.
This is usually expressed by saying that $x$ is a caustic point for 
the family of lines $C'$.
The set of all caustic points of the family of lines $C'$
is usually called the caustic (or the envelope) of $C'$.
This is nothing else but the projectively dual curve $\dual{C'}$.
Now the Reflexivity Theorem means
that any curve $C\subset\P^2$ coincides with the caustic curve
for the family of its tangent lines.
The ``dual'' form of this theorem is less intuitively
obvious, it means that any $1$-parameter family of lines in $\P^2$
consists of tangent lines to some curve $C$ and this curve
is a caustic curve for this family of lines.
For (real) families of lines the caustic could be found (locally)
with the help of the Legendre transformation:

\begin{Theorem}
Consider a family of real straight lines $y=px-g(p)$.
Then its caustic curve has the equation $y=f(x)$,
where $f$ is the Legendre transformation of $g$.
\end{Theorem}

This Theorem is, in fact, a version of the Reflexivity Theorem,
see \cite{Ar} for the proof.

\subsection{Correspondence of Branches. Pl\"ucker Formulas.}

Even if a plane curve $C\subset\P^2$ is smooth,
the dual curve $\dual C\subset\dual{\P^2}$ almost always
has singularities. We have a natural map $C\to \dual C$,
sending a point $p\in C$ to a tangent line $l$ to $C$ at $p$.
This map is clearly a resolution of singularities.
In general, curves $C$ and $\dual C$ are birationally equivalent.
Indeed, consider the conormal variety $I_C\subset \P^2\times\dual{\P^2}$, 
the closure
of the set of pairs $(p,l)$, $p\in C_{sm}$, $l\in\dual{C}_{sm}$,
$l$ is tangent to $C$ at $p$.
Then, by Reflexivity Theorem, $I$ projects birationally both
on $C$ and $\dual C$. Therefore, $C$ and $\dual C$ are birationally
equivalent, in particular, they have the same
geometric genus $g$.

A line $l$ which is tangent to $C$ in at least two points,
is a singular point of $\dual C$. It is known as 
the multiple tangent\index{multiple tangent}\index{correspondence of branches}.
If a multiple tangent $l$ has exactly two tangency points on $C$
and the intersection multiplicity at each of them is equal exactly to $2$,
then $l$ is called the bitangent\index{bitangent}. 
A bitangent corresponds
to an ordinary double point\index{ordinary double point} 
of $\dual C$.

If the tangent $l=\hat T_p$ at a non-singular point $p\in C$
intersects $C$ at $p$ with multiplicity $\ge3$, it is again
a singular point on $\dual C$. If the intersection multiplicity
is precisely $3$, and $l$ is not tangent to $C$ at any other point,
then $p$ is called an inflection point\index{inflection point} 
(or flex\index{flex}) of $C$. 
Then $l$ is a cuspidal point\index{cuspidal point} 
(or cusp\index{cusp}) of $\dual C$.

Now we may introduce a class of ``generic'' curves
with singularities, which is preserved by the projective duality.
Namely, we say that a curve $C$ is generic if both $C$ and $\dual C$
have only double points and cusps as their singularities.
Suppose that $C$ is generic in this sense. 
Let $d$, $g$, $\kappa$, $\delta$, $b$, $f$
be the degree, the geometric genus, the number of cusps, the number of double points,
the number of bitangents, and the number of flexes of $C$.
Let $\dual d$, $\dual g$, $\dual \kappa$, $\dual \delta$, $\dual b$, $\dual f$ 
be the corresponding numbers for $\dual C$
($\dual d$ is also sometimes called the class\index{class of a plane curve} 
of $C$). Then by Reflexivity Theorem we have the following

\begin{Proposition}
$g=\dual g$, $\kappa=\dual f$, $\delta=\dual b$, $b=\dual\delta$, $f=\dual\kappa$.
\end{Proposition}

It turns out that there is another remarkable set of equations
linking these numbers, that was discovered by
Pl\"ucker\index{Pl\"ucker formulas} and Clebsch.
The proof can be found in \cite{GH}.

\begin{Theorem}
$$g={1\over2}(\dual d-1)(\dual d-2)-b-f,$$
$$g={1\over2}(d-1)(d-2)-\delta-\kappa,$$
$$d=\dual d(\dual d-1)-2b-3f,$$
$$\dual d=d(d-1)-2\delta-3\kappa.$$
\end{Theorem}

\begin{Example}\rm
Let $C\subset\P^2$ be a smooth conic.
In homogeneous coordinates $x_1,x_2,x_3$, the curve $C$
is given by
$$(Ax,x)=\sum_{i,j=1}^3a_{ij}x_ix_j=0,$$
where $A=||a_{ij}||$ is a non-degenerate symmetric $3\times3$-matrix.
It is clear that the tangent line to $C$ at a point
$x_0\in C$ is given by the equation $(Ax_0,x)=0$.
Hence the point $\xi\in\dual{\P^2}$ corresponding
to this tangent line has homogeneous coordinates $Ax_0$,
which implies $(A^{-1}\xi,\xi)=0$.
Therefore, $\dual C\subset\dual{\P^2}$
is also a smooth conic defined by the inverse matrix $A^{-1}$.
\end{Example}

\begin{Example}\rm
Let $C\subset\P^2$ be a smooth cubic curve.
By Bezout Theorem, $C$ is automatically generic
and Pl\"ucker formulas are always applicable.
$C$ has no bitangents and has exactly $9$ flexes.
The dual curve $\dual C$ is a very special curve of degree $6$
with $9$ cusps and no double points.
Schl\"affli has found a beautiful determinantal formula for $\dual C$.
Let $x_1,x_2,x_3$ be homogeneous coordinates in $\P^2$
and $p_1,p_2,p_3$ the dual coordinates in $\dual{\P^2}$.
Let $f(x_1,x_2,x_3)=0$ be the homogeneous equation of $C$
and $F(p_1,p_2,p_3)=0$ be the homogeneous equation of $\dual C$.
Consider the polynomial
$$V(p,x)=\left|\matrix{
0&p_1&p_2&p_3\cr
p_1&{\partial^2f\over\partial x_1\partial x_1}&
{\partial^2f\over\partial x_1\partial x_2}&
{\partial^2f\over\partial x_1\partial x_3}\cr
p_2&{\partial^2f\over\partial x_2\partial x_1}&
{\partial^2f\over\partial x_2\partial x_2}&
{\partial^2f\over\partial x_2\partial x_3}\cr
p_3&{\partial^2f\over\partial x_3\partial x_1}&
{\partial^2f\over\partial x_3\partial x_2}&
{\partial^2f\over\partial x_3\partial x_3}}
\right|.$$
Clearly $V(p,x)$ has degree $2$ in $x$.
Then Schl\"affli's formula\index{Schl\"affli formula} is as follows:

\begin{Theorem}
$$F(p_1,p_2,p_3)=\left|\matrix{
0&p_1&p_2&p_3\cr
p_1&{\partial^2V\over\partial x_1\partial x_1}&
{\partial^2V\over\partial x_1\partial x_2}&
{\partial^2V\over\partial x_1\partial x_3}\cr
p_2&{\partial^2V\over\partial x_2\partial x_1}&
{\partial^2V\over\partial x_2\partial x_2}&
{\partial^2V\over\partial x_2\partial x_3}\cr
p_3&{\partial^2V\over\partial x_3\partial x_1}&
{\partial^2V\over\partial x_3\partial x_2}&
{\partial^2V\over\partial x_3\partial x_3}}
\right|.$$
\end{Theorem}
The proof can be found in \cite{GKZ2}.
\end{Example}

\begin{Example}\rm
Suppose that $C\subset\P^2$ is a generic smooth quartic curve.
By applying Pl\"ucker formulas
we see that $C$ has genus $3$, $24$ flexes, and $28$ bitangents.
There are two known visual descriptions of these bitangents.
The first classical approach is to realize the quartic curve
as a `shade' of a cubic surface in $\P^3$.
Namely, suppose that $S\subset\P^3$ is a generic cubic surface,
$x\in S$ is a generic point.
Consider the projective plane $\P^2$ of lines in $\P^3$
passing through $x$. Then lines $l\in\P^2$ that are tangent to $S$
form a quartic curve $C$. The $28$ bitangents then coincide
with projections of $27$ lines on $S$ plus one extra line:
the projectivization of a tangent space $T_xS$.

Another classical description is less known but it is more
in style of these notes. 
Consider a generic $3$-dimensional linear system $L$
of quadrics in $\P^3$ (an element of $\Gr(3,S^2(\C^4)^*$).
Then singular members of this linear system
form a quartic curve in $\P(L)=\P^2$.
It can be shown that any generic quartic curve 
can be obtained this way (for example by calculating
the differential at a generic point of the map
from linear systems to quartic forms).
By the Bezout theorem $L$ has $2^3=8$ base points in $\P^3$.
For any two base points $p$, $q$ let $l_{pq}\subset \P(L)$
be the line formed by all quadrics that contain not only points $p$ and $q$,
but also a line $[pq]$ connecting them.
Therefore, we have $28$ lines $l_{pq}\subset\P(L)$.
Let us show that these lines are bitangents to $C$.

Fix a basis $\{e_1,e_2,e_3,e_4\}$ in $\C^4$ such that
$e_1$ belongs to $p$ and $e_2$ belongs to $q$.
Then quadrics from $l_{pq}$ have the form
$$G=\left(\matrix{
0&A^T\cr
A&B\cr
}\right),$$
therefore,
$$\det G=(\det A)^2.$$
This exactly means that $l_{pq}$ intersects $C$
in two double points, i.e.~$l_{pq}$ is a bitangent to $C$.
\end{Example}


\section{Projections and Linear Normality}\label{LinearNormSection}

\subsection{Projections}
Let $P=\P(V)$ be an $n$-dimensional projective space
and $L\subset P$ a projective subspace of dimension $k$, 
$L=\P(U)$, $U\subset V$. The quotient projective space\index{quotient projective space} 
$P/L=\P(V/U)$
has, as points, $(k+1)$-dimensional projective subspaces in $P$
containing $L$. The projection with center $L$\index{projection from a subspace} 
is the map $\pi_L:\,P\setminus L\to P/L$ 
which takes any point $x\in P\setminus L$
to the $(k+1)$-dimensional projective subspace spanned by $x$ and $L$.
In other words, it corresponds to the projection $V\to V/U$.
Classically, both $P/L$ and $\pi_L$ can be visualized
inside the initial projective space $P$. Namely,
$P/L$ can be identified with any $(n-k-1)$-dimensional
projective subspace $H\subset P$ not intersecting $L$.
Then $\pi_L$ sends any point $x\in P\setminus L$ to the
unique intersection point of $H$ and $(k+1)$-dimensional projective
subspace spanned by $x$ and $L$.
Clearly, the dual projective space $\dual{(P/L)}$
is canonically embedded in $\dual P$ as the set
of hyperplanes containing $L$, so it coincides
with $\dual L$, the dual variety of $L$.

\begin{Theorem}\label{ProjectionDual}$\ $
\begin{enumerate}
\item[{\rm (a)}] Let $X\subset P=\P(V)$ be an algebraic subvariety
not intersecting the projective subspace $L=\P(U)$ and 
such that $\dim X<\dim P/L$. Then
$$\dual{(\pi_L(X)}\subset \dual L\cap\dual X.$$
The discriminant $\Delta_{\pi_L(X)}$ 
{\rm(}as a polynomial on $(V/U)^*\subset V^*${\rm)}
is a factor of the restriction of $\Delta_X$ to $(V/U)^*$.
\item[{\rm (b)}] 
Suppose further that $\pi_L:\,X\to\pi_L(X)$ is an isomorphism.
Then 
$$\dual{(\pi_L(X)}=\dual L\cap\dual X,\ \hbox{\rm and}\ 
\Delta_{\pi_L(X)}=\Delta_X|_{(V/U)^*}.$$
\end{enumerate}
\end{Theorem}

\proof
A hyperplane in $P/L$ is just a hyperplane in $P$ containing $L$.
It is clear that if a hyperplane $H\subset P/L$ is tangent to $\pi_L(X)$
at some smooth point $y=\pi_L(x)$, where $x\in X$ is also smooth,
then $H$ as a hyperplane in $P$ is tangent to $X$ at $x$.
This proves (a).

Suppose now that $\pi_L$ is an isomorphism.
The same argument as above also shows that if 
$X^*_0\subset\dual X$ 
(resp.~$\pi_L(X)^*_0\subset\dual{\pi_L(X)}$)
is a dense subset of hyperplanes tangent to $X$ (resp.~$\pi_L(X)$)
at some smooth point then $\pi_L(X)^*_0=\dual L\cap X^*_0$,
since $\pi_L$ induces an isomorphism of smooth loci $X_{sm}\to\pi_L(X)_{sm}$.
Moreover, it is clear that at a generic point of $\pi_L(X)^*_0$
varieties $\dual L$ and $X^*_0$ intersect transversally.
Since $\dual X=\overline{X^*_0}$ and 
$\dual{\pi_L(X)}=\overline{\pi_L(X)^*_0}$,
in order to prove (b) it suffices to show that
$\dual L\cap\dual X=\overline{\dual L\cap X^*_0}$.
Here the assumption that $\pi_L$ is an isomorphism is crucial, because
otherwise critical points of $\pi_L$ may 
produce extra components of $\dual L\cap\dual X$.

We need to show that if $H(t)$ is a $1$-parameter family of hyperplanes
such that for $t\ne0$ the hyperplane $H(t)$ is tangent to $X$
at some smooth point $x(t)$ and $H(0)$ contains $L$, then there exists another
$1$-parameter family $H'(t)$ such that $H'(0)=H(0)$ and for $t\ne0$
the hyperplane $H'(t)$ is tangent to $X$ at $x(t)$ and contains $L$. 
Since $X$ is projective, hence compact, we can suppose that
the limit $x(0)=\lim_{t\to0}x(t)$ exists. However, $x(0)$ may
be a singular point of $X$. Let $\hat T_{x(0)}X$ be an embedded
Zariski projective tangent space to $X$ at $x(0)$.
Since $\pi_L$ is an isomorphism, it induces an isomorphism
of Zariski tangent spaces. Therefore 
$\hat T_{x(0)}X$ does not intersect $L$.
Let $T_0\subset \hat T_{x(0)}X$ be a limit position of embedded
tangent Zariski spaces $T_t=\hat T_{x(t)}X$ as $t\to 0$.
Then for any $t$ the projective subspace $T_t$ does not intersect $L$.
Therefore we may consider a family $T'_t$ of $(\dim L+\dim X+1)$-dimensional
projective subspaces such that for any $t$ we have $L\subset T'_t$
and for $t\ne0$ the subspace $T'_t$ is tangent to $X$ at $x(t)$.
Namely, $T'_t$ is a projective subspace spanned by $L$ and $T_t$.
Since $\dim X<\dim P/L$, we can embed $T'_t$ into a family
of hyperplanes $H'_t$ with same properties.
\endproof

\subsection{Linear Normality}

\begin{Definition}\rm
The projective variety $X\subset\P^n$ is called non-degenerate\index{non-degenerate variety}
if $X$ is not contained in any hyperplane $H$.
\end{Definition}

The next theorem shows that we may restrict ourselves
to non-degenerate varieties while studying dual varieties and discriminants.
First we shall need the following definition.
We recall from a previous section that if 
$\P^k=L\subset\P^n$ is a projective subspace
then a dual projective space $\dual{\P^k}$
is not canonically embedded in $\dual{\P^n}$, but 
(by some abuse of notation) is canonically
isomorphic to the quotient projective space $\dual{\P^n}/\dual{L}$,
where $\dual L$ is a projective subspace of $\dual{\P^n}$
projectively dual to $L$. However, if $Y\subset\dual{\P^k}$
is any subvariety then we can canonically define a cone over $Y$
with vertex $\dual L$ as the closure $\pi^{-1}_{\dual L}(Y)$.

\begin{Theorem}\label{DegenerateDual}
Let $X\subset\P^n$ be an irreducible subvariety.
\begin{enumerate}
\item[{\rm(a)}] Assume that $X$ is contained in a hyperplane $H=\P^{n-1}$.
If $\dual{X}'$ is the dual variety of $X$, when
we consider $X$ as a subvariety of $\P^{n-1}$,
then $\dual X$ is the cone over $\dual X'$ 
with vertex $p$ corresponding to $H$.
\item[{\rm (2)}] 
Conversely, if $\dual X$ is a cone with vertex $p$,
then $X$ is contained in the corresponding hyperplane $H$.
\end{enumerate}
\end{Theorem}

\proof
If $H'\ne H$ is a tangent hyperplane of $X$
then $H\cap H'$ is a tangent hyperplane of $X$ in $\P^{n-1}$.
Conversely, if $T$ is a tangent hyperplane of $X$ in $\P^{n-1}$
then each hyperplane $H'$ in $\P^n$ containing $T$ is tangent to $X$.
Therefore, $\dual X$ is the cone over $\dual X'$. This proves (a).

By Reflexivity Theorem, we also have (b). Namely,
each hyperplane which is tangent to $\dual X$ at a smooth point
necessarily contains $p$. Therefore $X=X^{**} $ is contained
in the hyperplane corresponding to $p$.
\endproof

In terms of discriminants, Theorem~\ref{DegenerateDual} can be reformulated
as follows. Consider a surjection $\pi:\,V\to U$ of vector spaces
and consider $\P(V^*)$, $\P(U^*)$. Then we have an embedding
$i:\,\P(U^*)\hookrightarrow\P(V^*)$.
Let $X\subset \P(U^*)$. Then $\Delta_X$ is a polynomial function
on $U$. If we consider $X$ as a subvariety in $\P(V^*)$
then $\Delta_{i(X)}$ is a function on $V$.
By Theorem~\ref{DegenerateDual} these polynomial functions are related
as follows:
$$\Delta_{i(X)}(f)=\Delta_X(\pi(f)).$$
In other words, $\Delta_{i(X)}$ does not depend on some of the arguments
and forgetting these arguments gives $\Delta_X$.

Theorems \ref{ProjectionDual} and \ref{DegenerateDual}
show that in order to study dual varieties and discriminants
we need to consider only projective varieties $X$ that
are non-degenerate and not equal to a non-trivial projection.
These projective varieties are called linearly normal.\index{linearly normal variety}
To give a more intrinsic definition of linearly
normal varieties we need to recall the correspondence
between invertible sheaves, linear systems, and projective embeddings.

An invertible sheaf on an irreducible algebraic variety $X$ is simply
the sheaf of sections of some algebraic line bundle.
For example, the structure sheaf of regular functions $\O_X$
corresponds to the trivial line bundle.
Usually we shall not distinguish notationally between
invertible sheaves and line bundles. 
Invertible sheaves form the group $\Pic(X)$ with respect
to the tensor product.

A Cartier divisor\index{Cartier divisor} 
on $X$ is a family $(U_i,g_i)$, $i\in I$,
where $U_i$ are open subsets of $X$ covering $X$, and $g_i$
are rational functions on $U_i$ such that $g_i/g_j$
is regular on each intersection $U_i\cap U_j$.
The functions $g_i$ are called local equations of the 
divisor\index{local equations of a divisor}.
More precisely, a Cartier divisor is an equivalence class
of such data. Two collections $(U_i,g_i)$ and $(U_i',g_i')$
are equivalent if their union is still a divisor.
Cartier divisors can be added, by multiplying
their local equations. Thus they form a group, denoted by $\Div(X)$.

If each local equation $g_i$ is regular on $U_i$,
then we say that the divisor $D$ is effective,\index{effective divisor}
and we write $D\ge0$.
The subschemes $\{g_i=0\}$ of the $U_i$ can then be glued together
into a subscheme of $X$, also denoted by $D$.
Therefore, effective Cartier divisors can be identified with subschemes
of $X$ that are locally given by one equation.
Any non-zero rational function $f\in\C(X)$ determines a Cartier
divisor $(f)$ which is said to be principal.
Principal divisors form a subgroup of $\Div(X)$.

Let ${\cal K}_X$ denote the sheaf of rational functions on $X$,
${\cal K}_X(U)=\C(U)$. To every Cartier divisor $D=(U_i,g_i)_{i\in I}$
we can attach a subsheaf $\O_X(D)\subset{\cal K}_X$. Namely, 
on $U_i$ it is defined as $g_i^{-1}\O_{U_i}$. On the intersection
$U_i\cap U_j$ the sheaves $g_i^{-1}\O_{U_i}$ and $g_j^{-1}\O_{U_j}$
coincide, since $g_i/g_j$ is invertible.
Therefore these sheaves can be pasted together
into a sheaf $\O_X(D)\subset{\cal K}_X$.
For instance, $\O_X(0)=\O_X$ and $\O_X(D_1+D_2)=\O_X(D_1)\otimes\O_X(D_2)$.
The sheaves $\O_X(D)$ are invertible.
In fact, multiplication by $g_i$ defines an isomorphism (``trivialization'')
$\O_X(D)|_{U_i}\simeq\O_{U_i}$.
Therefore, we have a homomorphism $\Div(X)\to\Pic(X)$.
This homomorphism is surjective and its kernel
consists of principal divisors.

A non-zero section of $\O_X(D)$ is a rational function $f$
on $X$ such that the functions $fg_i$ are regular on the $U_i$,
in other words, such that the divisor $(f)+D$ is effective.
If $D$ itself is effective, the sheaf $\O_X(D)$ has a canonical section
$s_D$, which corresponds to the constant function $1$.
By contrast, the sheaf $\O_X(-D)$, for $D$ effective, is an ideal sheaf in $\O_X$.
It defines $D$ as a subscheme.
The sections of invertible sheaves define some divisors.
Let $s\in H^0(X,\L)$ be a non-trivial global section of an invertible
sheaf $\L$. After choosing some trivializations $\phi_i:\,\L_{U_i}\simeq\O_{U_i}$
on a covering $(U_i)$, we obtain an effective divisor $(U_i,\phi_i(s_i))$,
which we denote by $\div(s, \L)$.
For instance, if $D$ is effective,
the canonical section $s_D$ of the sheaf $\O(D)$ defines $D$.
Thus we have established the possibility to define any effective 
divisor by one equation $s=0$, keeping in mind
that $s$ is not a function, but a section of an invertible sheaf.
If $s'$ is another non-zero global section of $\L$
then the divisors $\div(s',\L)$ and $\div(s,\L)$
differ by the divisor of a rational function $s/s'$.
One also says that they are linearly equivalent.\index{linearly equivalent divisors}

\begin{Example}\rm
Let us recall the construction of invertible sheaves on projective spaces $\P(V)$.
All these sheaves have the form $\O(d)$, where
$\O(d)$ is the sheaf of degree~$d$ homogeneous functions on $\P(V)$.
More precisely, let $\pi:\,V\setminus\{0\}\to \P(V)$ be the canonical projection.
If $U\subset\P(V)$ is a Zariski open set, then the sections of $\O(d)$ over $U$
are, by definition, regular functions $f$ on $\pi^{-1}(U)\subset V$,
which are homogeneous of degree $d$, $f(\lambda v)=\lambda^df(v)$.
It is well-known that these sheaves are invertible and any invertible
sheaf has the form $\O(d)$ for some $d$.
For example, $\O(-1)$, as a line bundle, is the tautological
line bundle, the fiber of $\O(-1)$, as a line bundle, over a point
of $\P(V)$ represented by a $1$-dimensional subspace $l$
is $l$ itself. If $d<0$ then $H^0(\P, \O(d))=0$.
If $d\ge0$ then $H^0(\P, \O(d))=S^dV^*$, homogeneous polynomials
of degree $d$. For any non-zero $s\in H^0(\P, \O(d))$ the corresponding
effective divisor $\div(s,\O(d))$ is just the hypersurface defined by 
the polynomial $s$. In particular, hyperplanes in $\P$
correspond to global sections of $\O(1)$.
\end{Example}

Suppose now that $X$ is a projective variety in $\P^n$.
Then any sheaf $\O(d)$ can be restricted on $X$,
thus we get a sheaf $\O_X(d)$ for any $d$.
In particular, we have a restriction homomorphism
of global sections:
$$V^*=H^0(\P^n,\O(1))\mathop{\to}^{\res} H^0(X,\O_X(1)).$$
Clearly this map is not injective if and only if $X$ is degenerate.
So suppose now that $X$ is non-degenerate.
In general, map $\res$ is not surjective.
Its image is a vector subspace $W\subset H^0(X,\O_X(1))$ with the following
obvious property: for any $x\in X$ there exist some section $s\in W$
such that $s(x)\ne0$. Clearly divisors of the form $\div(s,\O_X(1))$,
$s\in W$ are just hyperplane sections $X\cap H$ for various hyperplanes $H$.

In general, if $X$ is any variety with invertible sheaf $\L$
then any family of divisors $|W|$ of the form $\div(s,\Lambda)$,
where $s\in W\subset H^0(X,\L)$, $\dim W<\infty$, is
called a linear systems of divisors.\index{linear system of divisors}
If $X$ is projective (or just complete) 
then $H^0(X,\L)$ is finite-dimensional
and we can take $W=H^0(X,\L)$. The corresponding
linear system is called complete.\index{complete linear system} 
It is denoted by $|\L|$,
or by $|D|$ if $\L=\O(D)$.
Every effective divisor that is linearly equivalent to $D$
appears in $|D|$, and exactly once.

A linear system $W$ is called base-point free\index{base-point free} 
if the intersection
of all its divisors is empty, or, equivalently,
if for any $x\in X$ there exists some $s\in W$ such that
$s(x)\ne0$. As we have seen, if $X\subset\P(V)$ is a projective
variety then the image of $\res:\,V^*\to H^0(X,\O_X(1))$
defines a base-point free linear system.
It turns out that any base-point free linear system $|W|$
defines a morphism into a projective space.
Namely, for any $x\in X$ we have a hyperplane $\Gamma(x)\subset W$
consisting of all sections vanishing at $x$.
Therefore, we have a regular map $\Gamma:\,X\to\P(W^*)$.
If $\Gamma$ is an embedding then we see that $W$
is identified with the image of a restriction map
$W\to H^0(X,\O_X(1))$ and $\O_X(1)$ is identified with $\L$.
So, we need a final definition. 

\begin{Definition}\rm
A linear system $|W|$ is 
called very ample\index{very ample} 
if it defines an embedding in a projective space.
An invertible sheaf $\L$ is called very ample if the complete
linear system $|\L|$ is very ample. $\L$ is called ample\index{ample}
if some tensor power $\L^{\otimes m}$, $m>1$, is very ample.
\end{Definition}

In particular, we see that non-degenerate projective embeddings of $X$
are in a $1-1$ correspondence with very ample linear systems
of divisors.

Now we can describe linearly normal varieties.
\begin{Theorem}
A projective variety $X\subset\P(V)$ is linearly normal if and only
if the restriction homomorphism $\res:\, V^*\to H^0(X,\O_X(1))$
is an isomorphism if and only if the embedding $X\subset\P(V)$
is given by a complete linear system corresponding to some
very ample invertible sheaf $\L$ on $X$.
\end{Theorem}

\proof
Clearly $\res$ is injective if and only if $X$ is non-degenerate.
Suppose that $\res$ is not surjective. Then the linear system $|\O_X(1)|$
gives a non-degenerate embedding of $X$ to a higher-dimensional
vector space $H^0(X,\O_X(1))$. The initial embedding
is obtained from this one by the projection with center
$\P(\Ann(\Im(\res)))$.
Conversely, if $X$ can be obtained as a non-trivial isomorphic
projection of a non-degenerate variety $\tilde X$ in a larger
projective space $\P(U)$ then we have a proper embedding $V^*\subset U^*$.
The map $U^*\to H^0(X,\O_{\tilde X}(1))$, given 
by restricting linear functions from $\P(U)$ to $\tilde X$,
is an injection, since $\tilde X$ is non-degenerate.
Thus the map $\res$ is not surjective, being the composition
$V^*\hookrightarrow U^*\to H^0(X,\O_X(1))$.
\endproof

We see that in order to study dual varieties and discriminants
we can restrict ourselves to linearly normal varieties.

\begin{Definition}\rm
A pair $(X,\L)$ of a projective variety and a very ample invertible sheaf
on it is called the polarized variety.\index{polarized variety}
\end{Definition}

Any polarized variety admits a canonical embedding in a projective
space with a linearly normal image, namely,
the embedding given by the complete linear system $|\L|$.
Therefore we may speak about dual varieties, defect, discriminants, etc.
of polarized varieties without any confusion.

\section{Dual Varieties of Smooth Divisors}
\index{defect of smooth divisors}

Suppose that $X\subset\P^N$ is a smooth projective variety.
The following theorem allows to find the defect of smooth 
hyperplane or hypersurface sections of $X$.

\begin{Theorem}[\cite{E2,HK,Ho1}]\label{DefectOfDivisorTh}$\ $
\begin{enumerate}
\item[\rm(a)] Assume that $Y=X\cap H$ 
is a smooth hyperplane section of $X$.
Then 
$$\defect Y=\max\{0,\defect X-1\}.$$
Moreover, if $\dual X$ is not a hypersurface,
then the dual variety $\dual{Y}$ 
is the cone over $\dual X$ with $H$ as a vertex.
\item[\rm(b)] Assume that $Y$ is a smooth divisor corresponding
to a section of $\O_X(d)$ for $d\ge2$.
Then $\defect Y=0$.
\end{enumerate}
\end{Theorem}

\proof
Consider the conormal variety 
$$I_X=\P(N_X^*\P^N(1))\subset\P^N\times\dual{\P^N}$$
and the following diagram of projections
$$\P^N\supset X\la^{\pi_1} I_X\ra^{\pi_2} \dual X\subset\dual{\P^N}.$$
We use notation $\O_{I_X}(a,b)$ for the line bundle
$\pi_1^*\O_X(a)\otimes\pi_2^*\O_{\dual X}(b)$. 
We will need the notion of the tautological line bundle
on a projective bundle.
Let $E$ be a vector bundle on $X$. We consider the variety
$X_E=\P(E^*)$, the projectivization of the bundle $E^*$.
There is a projection $p:\,E_X\to X$ whose fibers are projectivizations
of fibers of $E^*$, and a  projection 
$\pi:\,E^*\setminus X\to X_E$, where $X$ is embedded into the total space of $E^*$
as the zero section.
We denote by $\xi(E)$ the tautological line bundle on $X_E$
defined as follows. For open $U\subset X_E$, a section of $\xi(E)$
over $U$ is a regular function on $\pi^{-1}(U)$ which is homogeneous
of degree $1$ with respect to dilations of $E^*$.
The restriction of $\xi(E)$ to every fiber $p^{-1}(x)=\P(E_x^*)$
is the tautological line bundle $\O(1)$ of the projective space
$\P(E_x^*)$.
For example, $\O_{I_X}(0,1)$ is the tautological line bundle
of $\P(N_X^*\P^N(1))$.

Suppose that $Y$ is a smooth member of a linear system $|\O_X(d)|$.
Clearly we have $N_YX=\O_Y(d)$.
Consider the following exact sequence
$$0\to \O_Y(d-1)\to N_Y\P^N(-1)\to N_X\P^N(-1)|_Y\to 0.$$
Let $D=\P(N_X^*\P^N(1)|_Y)$. Then there is a natural embedding
of $D$ as a divisor in $I_Y=\P(N_Y^*\P^N(1))$.
We see that $\O_{I_Y}(D)\simeq\O_{I_Y}(1-d,1)$.

We shall study the following diagram

$$
\begin{CD}
\dual{\P^N} & \ \supset \  & \dual Y &\ \supset \  & \pi_2(D) & \ \subset \  & \dual X & \ \subset \  & \dual{\P^N}\\
            &         & @AA{\pi_2}A   @AA{\pi_2}A            @AA{\pi_2}A           &               \\
            &              & I_Y     &\ \supset \  & D        & \ \subset \  & I_X   &              &            \\
            &         & @VV{\pi_1}V   @VV{\pi_1}V            @VV{\pi_1}V           &               \\
      \P^N  & \ \supset \  &       Y  &=&           Y         & \ \subset \  &       X & \ \subset \  &      {\P^N}\\
\end{CD}
$$

Consider the morphism $g$ given by the complete linear system $|\O_{I_Y}(0,1)|$.
Since $\pi_2^*\O_{\dual Y}(1)=\O_{I_Y}(0,1)$, we can regard $\dual Y$
as the linear projection of $g(I_Y)$. Thus $\dim \dual Y=\dim g(I_Y)$.

(a) Assume that $d=1$. Then $\O_{I_Y}(0,1)=\O_{I_Y}(D)$.
Therefore $g(D)$ is a hyperplane section of $g(I_Y)$ and
$$\dim g(D)=\dim\pi_2(D)=\dim\dual Y-1.$$

Therefore it is sufficient to prove that 
\begin{enumerate}
\item[(A)] If $\defect X=k>0$, then $\dim\pi_2(D)=\dim\dual X$.
\item[(B)] If $\defect X=0$ then $\dim\pi_2(D)=\dim\dual X-1$.
\end{enumerate}

For any $q\in\dual X$ the preimage $\pi_2^{-1}(q)$ is isomorphic
to the contact locus of the hyperplane $q=0$ and $X$. Therefore,
$\dim\pi_2^{-1}(q)\ge\defect X$ and for generic $q\in\dual X$
we have an equality. Therefore, if $\defect X=k>0$
then $\dim D\cap\pi_2^{-1}(q)>0$ and (A) follows.
Moreover, we see that any hyperplane tangent to $X$ is also tangent to $Y$,
therefore, $\dual Y$ is a cone over $\dual X$ with the vertex~$H$.
Suppose now that $\defect X=0$. Then $\pi_2$ is birational,
and, by the same argument as above,
$\dim\pi_2(D)\le\dim\dual X-1$.
Moreover, if $Y$ passes through a point $x\in X$ such that there
exists a hyperplane $H\in{\dual X}_{sm}$ tangent to $X$ at $x$
then, in fact, 
$\dim\pi_2(D)=\dim\dual X-1$.
This proves (B) for generic hyperplane sections and, hence,
for arbitrary smooth hyperplane sections by a standard argument.

(b) Assume $d\ge2$. We claim that $g|_{I_Y\setminus D}$ is one to one,
and, therefore, $\dim \dual Y=\dim g(I_Y)=\dim I_Y=N-1$
and $\dual Y$ is a hypersurface.
Indeed, let $y_1$ and $y_2$ be two distinct points in $I_Y\setminus D$.
We need to show that the linear system $|\O_{I_Y}(0,1)|$
separates $y_1$ and $y_2$. If $\pi_1(y)=\pi_1(y_2)$
then $\pi_2(y_1)\ne\pi_2(y_2)$ because $\pi_2$ maps $\pi_1^{-1}(\pi_1(y_1))$
isomorphically to a projective subspace in $\dual Y$.
So we may assume that $\pi_1(y_1)\ne\pi_1(y_2)$.
Then we can find a hypersurface $S$ of degree $d-1$ in $\P^N$
which will contain $\pi_1(y_1)$ but not $\pi_1(y_2)$.
Recall that $\O_{I_Y}(D)\simeq\O_{I_Y}(1-d,1)$.
Therefore, $D+\pi_1^{-1}(S)\in|\O_{I_Y}(0,1)|$.
But $y_1\in D+\pi_1^{-1}(S)$ and $y_2\not\in D+\pi_1^{-1}(S)$.
\endproof

\chapter{Dual Varieties of Algebraic Group Orbits}

\section*{Preliminaries}
It is quite difficult to find the explanation of the following phenomenon:
most interesting discriminants and dual varieties appear when the
projective variety $X$ is either homogeneous or a closure of a homogeneous
variety. There are three known big classes of examples:
projective embeddings of flag varieties, algebraic groups
acting on projective spaces with finitely many orbits, and
projective embeddings of torical varieties. The last case
leads to the theory of A-discriminants studied in detail
in the beautiful book \cite{GKZ2} and omitted here.

\section{Polarized Flag Varieties}\label{DefinitionsNotationsFlags}

\subsection{Definitions and Notations}

Let $G$ be a connected simply-connected semisimple 
complex algebraic group with 
a Borel subgroup\index{Borel subgroup} $B$ and a 
maximal torus\index{maximal torus} $T\subset B$.
Let $\PP$ be the character group\index{character group} of $T$ 
(the weight lattice\index{weight lattice}).
Let $\Delta\subset \PP$ be the set of roots of $G$ relative to~$T$.
To every root $\alpha\in\Delta$ we can assign the $1$-dimensional
unipotent subgroup $U_\alpha\subset G$. 
We define the negative roots\index{negative roots} $\Delta^-$
as those roots $\alpha$ such that $U_\alpha\subset B$.
The positive roots\index{positive roots} 
are $\Delta^+=\Delta\setminus \Delta^-=-\Delta^-$.
Let $\Pi\subset\Delta^+$ be simple roots\index{simple roots}, 
$\Pi=\{\alpha_1,\ldots,\alpha_n\}$,
where $n=\rank G=\dim T$. Any root $\alpha\in\Delta$ is an integral
combination $\sum_in_i\alpha_i$ with nonnegative $n_i$
(for $\alpha\in\Delta^+$) or nonpositive $n_i$ 
(for $\alpha\in\Delta^-$). If $G$ is simple then we use
the Bourbaki numbering\index{Bourbaki numbering} \cite{Bo} of simple roots.

The weight lattice $\PP$ is generated as a $\Z$-module by the
fundamental weights\index{fundamental weights}
$\omega_1,\ldots,\omega_n$ dual to the simple roots
under the Killing form\index{Killing form}
$$\langle \omega_i|\alpha_j\rangle=
{2(\omega_i,\alpha_j)\over(\alpha_j,\alpha_j)}=\delta_{ij}.$$
A weight $\lambda=\sum_in_i\omega_i$ is called dominant\index{dominant weight}
if all $n_i\ge0$. We denote the dominant weights by $\PP^+$.
Dominant weights parametrize finite-dimensional 
irreducible $G$-modules: to any $\lambda\in\PP^+$
we assign an irreducible $G$-module $V_\lambda$ with highest weight $\lambda$.
A weight $\lambda=\sum_in_i\omega_i$ is called strictly dominant\index{strictly dominant weight}
if all $n_i>0$. We denote the strictly dominant weights by $\PP^{++}$.
There is a partial order\index{partial order on the weight lattice} 
on $\PP$: $\lambda>\mu$ if $\lambda-\mu$ is an integral combination of simple
roots with nonnegative coefficients.
If $\lambda\in\PP^+$ then we denote by $\lambda^*\in\PP^+$ the 
highest weight of the dual $G$-module $V_\lambda^*$.
Let $W$ be the Weil group of $G$ relative to~$T$. If $w_0\in W$
is the longest element then $\lambda^*=-w_0(\lambda)$.

The character group of $B$ is identified with a character group of $T$.
Therefore, for any $\lambda\in \PP$ we can assign a $1$-dimensional $B$-module $\C_\lambda$,
where $B$ acts on $\C_\lambda$ by a character $\lambda$.
Now we can define the twisted product $G\times_B\C_\lambda$ to be the 
quotient of $G\times \C_\lambda$ by the diagonal action $B$: 
$$b\cdot(g,z)=(gb^{-1},\lambda(b)z).$$
Projection onto the first factor induces the map $G\times_B\C_\lambda\to G/B$,
which realizes the twisted product as an equivariant line bundle $\L_\lambda$ 
on $G/B$ with fiber $\C_\lambda$.
It is well-known that the correspondence 
$\lambda\to \L_\lambda$ is an isomorphism
of $\PP$ and $\Pic(G/B)$. 
$\L_\lambda$ is ample if and only if 
$\L_\lambda$ is very ample if and only if
$\lambda\in \PP^{++}$.
By the Borel--Weil--Bott theorem~\cite{Bot}\index{Borel--Weyl-Bott theorem}, 
for strictly dominant $\lambda$
the vector space of global sections $H^0(G/B,\L_\lambda)$
is isomorphic as a $G$-module to $V_\lambda$, the irreducible $G$-module
with highest weight $\lambda$.
The embedding $G/B\subset \P(V_{\lambda^*})$
identifies $G/B$ with the projectivization of the cone of highest
weight vectors. The dual variety $\dual{(G/B)}$ therefore
lies in $\P(V_\lambda)$ and parametrises global sections 
$s\in H^0(G/B,\L_\lambda)$ such that the scheme of zeros $Z(s)$ is a singular divisor.

More generally, consider any flag variety of the form 
$G/P$, where $P\subset G$
is an arbitrary parabolic subgroup\index{parabolic subgroup}.
The subgroup $P\subset G$ is called parabolic if $G/P$
is a projective variety. $P$ is parabolic if and only if
it contains some Borel subgroup.
Up to conjugacy, me may assume that $P$ contains $B$.
The combinatorial description is as follows.
Let $\Pi_P\subset\Pi$ be some subset of simple roots.
Let $\Delta^+_P\subset\Delta^+$ denote the positive roots
that are linear combinations of the roots in $\Pi_P$.
Then $P$ is generated by $B$ and 
by the root groups $U_\alpha$ for $\alpha\in\Delta^+_P$.
We denote $\Pi\setminus\Pi_P$ by $\Pi_{G/P}$ and 
$\Delta^+\setminus\Delta^+_P$ by $\Delta^+_{G/P}$.
A parabolic subgroup is maximal\index{maximal parabolic subgroup} 
if and only if $\Pi_{G/P}$
is a single simple root.

The fundamental weights $\omega_{i_1},\ldots\omega_{i_k}$
dual to the simple roots in $\Pi_{G/P}$ generate the sublattice
$\PP_{G/P}$ of $\PP$. We denote $\PP^+\cap \PP_{G/P}$ by $\PP^+_{G/P}$.
The subset $\PP^{++}_{G/P}\subset\PP^+_{G/P}$ consists of all weights
$\lambda=\sum n_k\omega_{i_k}$ such that all $n_k>0$.
Any weight $\lambda\in \PP_{G/P}$ defines a character of $P$,
and therefore a line bundle $\L_\lambda$ on $G/P$.
Then the following is well-known.
\begin{itemize}
\item
The correspondence $\lambda\to \L_\lambda$ is an isomorphism
of $\PP_{G/P}$ and $\Pic(G/P)$. 
In particular, $\Pic(G/P)=\Z$ if and only if $P$ is maximal.
\item
$\L_\lambda$ is ample if and only if 
$\L_\lambda$ is very ample if and only if 
$\lambda\in \PP^{++}_{G/P}$.
\item
If $\lambda\in \PP^+_{G/P}$ then the linear system corresponding
to $\L_\lambda$ is base--point free. The corresponding map given by sections
is a factorization $G/P\to G/Q$, where $Q$ is a parabolic subgroup
such that $\Pi_Q$ is a union of $\Pi_P$ and all simple roots in $\Pi_{G/P}$
orthogonal to $\lambda$.
\item
For strictly dominant $\lambda\in \PP^{++}_{G/P}$
the vector space of global sections $H^0(G/P,\L_\lambda)$
is isomorphic as a $G$-module to $V_\lambda$, the irreducible $G$-module
with highest weight $\lambda$.
\end{itemize}

The embedding $G/P\subset \P(V_{\lambda^*})$
identifies $G/P$ with the projectivization of the cone of highest
weight vectors. The dual variety $\dual{(G/P)}$ therefore
lies in $\P(V_\lambda)$ and parametrises global sections 
$s\in H^0(G/P,\L_\lambda)$ 
such that the scheme of zeros $Z(s)$ is a singular divisor.

\subsection{Basic Examples}

\begin{Example}[Projective spaces and Grassmanians]\rm
Let $V=\C^n$ be a finite-dimensional vector space,
$G=\SL(V)$. The projective space $\P(V)$ is identified
with a flag variety $G/P$, where 
$P$ is the normalizer in $G$ of a line in $V$.
$P$ is a maximal parabolic subgroup, therefore $\Pic(\P(V))=\Z$.
All line bundles on $\P(V)$ have a form $\O(d)$ for $d\in\Z$.
The line bundle $\O(d)$ is ample if and only if it is very ample
if and only if $d>0$. In this case $H^0(\P(V),\O(d))$ is canonically
isomorphic to $S^dV^*$. The embedding 
$\P(V)\subset \P(S^dV)$ is called the Veronese embedding\index{Veronese embedding}.
It assigns to a line spanned by $v\in V$ the line in $S^dV$ spanned by $v^d$.
The zero scheme $Z(s)$ of the global section $s\in S^d(V^*)$
are just the hypersurface in $\P(V)$ of degree $d$ given by $s=0$.
Therefore, the dual variety in this case parametrises singular hypersurfaces.
In particular, the corresponding discriminant is the classical discriminant
of a homogeneous form of degree $d$ in $n$ 
variables.\index{classical discriminant}

The Grassmanian\index{Grassmanian} 
$\Gr(k,V)$ of $k$-dimensional vector subspaces of $V$
is a flag variety $G/P$, where $P$ is the normalizer in $\SL(V)$
of a $k$-dimensional subspace in $V$.
The $k$-dimensional subspaces in $V$ are in the bijective correspondence
with $(k-1)$-dimensional subspaces in $\P(V)$.
Having this in mind we shall sometimes write $\Gr(k-1,\P(V))$
instead of $\Gr(k,V)$. 
The Grassmanian $\Gr(k,V)$ carries a natural equivariant
vector bundle $\S$, called the tautological vector bundle.
The fiber of $\S$ over a $k$-dimensional linear subspace $U\subset V$
is $U$ itself. Therefore, $\S$ is a subbundle of the trivial
bundle $\Gr(k,V)\times V$ with fiber $V$ at each point.

$P$ is a maximal parabolic subgroup corresponding to the fundamental
weight $\omega_k$, therefore $\Pic(\Gr(k,V))=\Z$.
The ample generator of $\Pic(\Gr(k,V))$ is the line bundle $\Lambda^k\S^*$.
The corresponding embedding is called the Pl\"ucker embedding\index{Pl\"ucker embedding},
$\Gr(k,V)\subset\P(\Lambda^kV)$.
Namely, if $U\in\Gr(k,V)$, $u_1,\ldots,u_k$ is any basis of $U$
then $U\mapsto [u_1\wedge\ldots\wedge u_k]$.
The discriminant of the Pl\"ucker embedding was studied by Lascoux \cite{Las}.
Clearly, all minimal flag varieties of $\SL(V)$
are isomorphic to Grassmanians.

The tangent space to $\Gr(k,V)$ at the point $U$ is a vector
space with an additional structure. Namely, it is easy to 
verify the natural isomorphism
$$T_U\Gr(k,V)=\Hom(U,V/U)=U^*\otimes V/U.$$
In particular, the tangent bundle $T\Gr(k,V)$ is canonically
isomorphic to $\S^*\otimes(V/\S)$.
Notice also that $\Gr(k,V)\simeq\Gr(n-k,V^*)$. 
Namely, any subspace $U\in\Gr(k,V)$
maps to its annihilator $\Ann U\in\Gr(n-k,V^*)$.
\end{Example}

\begin{Example}[Spinor varieties]\label{SpinorExample}\rm
Let $V=\C^n$ be a vector space equipped with a non-degenerate
symmetric scalar product $Q$ or a non-degenerate symplectic form $\omega$
($n$ should be even in this case).
Let $G=\SO(V)$ or $G=\Sp(V)$, respectively.
For any $k\le n/2$ let $\Gr_Q(k,V)$ (resp.~$\Gr_\omega(k,V)$)
denote the isotropic Grassmanian\index{isotropic Grassmanian} 
of $k$-dimensional isotropic
subspaces. Recall that a subspace $U\subset V$ is called
isotropic\index{isotropic subspace} 
if for any $v_1,v_2\in U$ we have $Q(v_1,v_2)=0$
(resp.~$\omega(v_1,v_2)=0$).
Then $\Gr_Q(k,V)$ and $\Gr_\omega(k,V)$ are projective
$G$-equivariant varieties. It easily follows from the Witt theorem
that $\Gr_Q(k,V)$ and $\Gr_\omega(k,V)$ are irreducible
homogeneous varieties (hence flag varieties)
in all cases except the symmetric case when $n=2k$.
In the latter case $\Gr_Q(k,\C^{2k})$ consists of two
homogeneous components $\Gr^+_Q(k,\C^{2k})$ and $\Gr^-_Q(k,\C^{2k})$.
These two varieties are in fact isomorphic to each other
(as algebraic varieties, not as algebraic varieties with the group $G$ action).
Indeed, let $s$ be any orthogonal reflection. Then 
$s(\Gr^\pm_Q(k,\C^{2k}))=\Gr^\mp_Q(k,\C^{2k})$.
Moreover, $\Gr^\pm_Q(k,\C^{2k})$ is isomorphic to $\Gr_Q(k-1,\C^{2k-1})$.
To see this let $H=\C^{2k-1}\subset\C^{2k}$
be any hyperplane such that the restriction of $Q$ on $H$
is non-degenerate. Then we have obvious maps
$$\Gr^\pm_Q(k,\C^{2k})\to\Gr_Q(k-1,\C^{2k-1}),\quad
U\mapsto U\cap H.$$
It is easy to see that both maps are isomorphisms.
The variety $\Gr^\pm_Q(k,\C^{2k})$
is called the spinor variety\index{spinor variety}, denoted by ${\Bbb S}_k$.

The introduced varieties represent the complete list
of minimal flag varieties of groups $\SO(V)$ and $\Sp(V)$.
Namely, $\Gr_\omega(k,V)$, $k=1,\ldots,n/2$, is the minimal flag variety
of $\Sp(V)$ corresponding to the fundamental weight $\omega_k$
(for example, $\Gr_\omega(1,V)=\P(V)$).
$\Gr_Q(k,V)$, $0<k<n/2$,
is the minimal flag variety
of $\SO(V)$ corresponding to the fundamental weight $\omega_k$
(for example, $\Gr_Q(1,V)$ is a quadric hypersurface in $\P(V)$).
Finally,
$\Gr^\pm_Q(k,\C^{2k})$ 
is the minimal flag variety
of $\SO(\C^{2k})$ corresponding to the fundamental weights $\omega_{k-1}$ and
$\omega_k$.

In particular, the Picard group of all these varieties
is isomorphic to $\Z$.
Let us describe the embedding corresponding to the ample generator of $\Z$.

Consider the symplectic case $G=\Sp(V)$ first.
Then we have the Pl\"ucker embedding
$\Gr_\omega(k,V)\subset\P(\Lambda^k V)$.
However, $\Lambda^kV$ is reducible as $\Sp(V)$-module,
it contains the submodule $\omega\wedge\Lambda^{k-2}V$.
The complement $\Lambda^k_0V$ of this submodule
is an irreducible module with highest weight $\omega_k$.
The Pl\"ucker embedding in fact gives the embedding
of $\Gr_\omega(k,V)$ in $\P(\Lambda^k_0 V)$
as the projectivization of the cone of highest weight vectors.
This embedding corresponds to the ample generator of the Pickard group.

Now consider the orthogonal case $G=\SO(V)$.
Then we have the Pl\"ucker embedding
$\Gr_Q(k,V)\subset\P(\Lambda^k V)$.
Suppose first that $k<m$ if $n=2m+1$
or $k<m-1$ if $n=2m$.
Then $\Lambda^k V$ is irreducible as $G$-module
and has the highest weight $\omega_k$.
Therefore, this embedding 
corresponds to the ample generator of the Pickard group.

Finally, the embedding of the spinor variety ${\Bbb S}_m$
corresponding to the ample generator of the Pickard group
can be described as follows.
Let $W$ be the spinor representation of $\Spin(2m-1)$
or the half-spinor representation of $\Spin(2m)$
(they are isomorphic with respect to the embedding 
$\Spin(2m-1)\subset\Spin(2m)$).
Then ${\Bbb S}_m$ embeds into $\P(W)$ as the projectivization
of the highest weight vector orbit.

To give a more precise model, let us recall the definition
of the half-spinor representation\index{half-spinor representation}.
Let $V=\C^{2m}$ be an even-dimensional vector space with
a non-degenerate symmetric scalar product $Q$.
Let $\Cl(V)$ be the Clifford algebra\index{Clifford algebra} of $V$.
Recall that $\Cl(V)$ is in fact a superalgebra,
$\Cl(V)=\Cl^0(V)\oplus\Cl^1(V)$, where
$\Cl^0(V)$ (resp.~$\Cl^1(V)$) is a linear span of elements
of the form $v_1\cdot\ldots\cdot v_r$, $v_i\in V$, $r$ is even 
(resp.~$r$ is odd). Then we have
$$\Spin(V)=\{a\in\Cl^0(V)\,|\,a=v_1\cdot\ldots\cdot v_r,\ v_i\in V,\ Q(v_i,v_i)=1\}.$$
$\Spin(V)$ is a connected simply-connected
algebraic group.

Given $a\in\Cl(V)$, $a=v_1\cdot\ldots\cdot v_r$, $v_i\in V$,
let $\overline a=(-1)^rv_r\cdot\ldots\cdot v_1$.
This is a well-defined involution of $\Cl(V)$.
Using it, we may define an action $R$ 
of $\Spin(V)$ on $V$ by formula $R(a)v=a\cdot v\cdot\overline a$.
Then $R$ is a double covering $\Spin(V)\to\SO(V)$.

Let $U\subset V$ be a maximal isotropic subspace, hence $\dim U=m$.
Take also any maximal isotropic subspace $U'$
such that $U\oplus U'=V$. For any $v\in U$ we define
an operator $\rho(v)\in\End(\Lambda^*U)$ by the formula
$$\rho(v)\cdot u_1\wedge\ldots\wedge u_r=v\wedge u_1\wedge\ldots u_r.$$
For any $v\in U'$ we define an 
operator $\rho(v)\in\End(\Lambda^*U)$ by the formula
$$\rho(v)\cdot u_1\wedge\ldots\wedge u_r=\sum_{i=1}^r(-1)^{i-1}Q(v,u_i)
\wedge u_1\wedge\ldots u_{i-1}\wedge u_{i+1}\wedge\ldots\wedge u_r.$$
These operators are well-defined and
therefore by linearity we have a linear map
$\rho:\, V\to\End(\Lambda^*U)$.
It is easy to check that for any $v\in V$ we have
$\rho(v)^2=Q(v,v)\Id$.
Therefore we have a homomorphism
of associative algebras
$\Cl(V)\to\End(\Lambda^*U)$, i.e., $\Lambda^*U$
is a $\Cl(V)$-module, easily seen to be irreducible.
Therefore, $\Lambda^*U$ is also a $\Spin(V)$-module, now reducible.
However, the even part $\Lambda^{ev}U$ is an irreducible
$\Spin(V)$-module, called the half-spinor module $\S^+$.
Now let us describe the embedding of the spinor variety
${\Bbb S}_m\subset\P(\S^+)$.
This embedding will correspond to the ample generator
of $\Pic({\Bbb S}_m)$.
By definition, 
${\Bbb S}_m$ is an $\SO(V)$-orbit of a fixed maximal isotropic subspace,
say $U'$. We claim that the normalizer of $U'$ in $\Spin(V)$
(acting on $V$ via the representation $R$)
is equal to the normalizer of $1\in\Lambda^{ev}U$ in $\Spin(V)$
(acting on $\Lambda^{ev}$ via the representation $\rho$).
This will show that ${\Bbb S}_m$ is naturally isomorphic
to the projectivization of the orbit of $1$ in the half-spinor
representation.
Consider the map
$$\Psi:\,V\times\Lambda^{ev}U\to\Lambda^{odd}U,\quad (v,w)\mapsto\rho(v)w.$$
$\Spin V$ acts on $V$ via the representation $R$ and on 
$\Lambda^{ev}$ and $\Lambda^{odd}$ via the representation $\rho$.
We claim that $\Psi$ is $\Spin(V)$-equivariant.
Indeed, for any $g\in\Spin(V)$ we have
\begin{eqnarray*}
\Psi(R(g)v,\rho(g)w)=\rho(R(g)v)(\rho(g)w)=\rho(gv\overline g)(\rho(g)w)\cr
=\rho(gv\overline gg)w=\rho(g)(\rho(v)w)=\rho(g)\Psi(v,w),
\end{eqnarray*}
since $g\overline g=1$.
It is easy to check that 
$\Psi(U',1)=0$ and if for any $v\in U'$ we have $\Psi(v,w)=0$
then $w=\lambda 1$ for some $\lambda\in\C$.
Suppose now that $g\in\Spin(V)$ normalizes $1\in\Lambda^{ev}U$.
Then for any $v\in U'$ we have 
$$\Psi(R(g)v,\rho(g)1)=\rho(g)\Psi(v,1)=0,$$
therefore $g$ normalizes $U'$.
On the contrary, if $g$ normalizes $U'$ then for any $v\in U'$
$$0=\Psi(R(g)(R(g^{-1})v),\rho(g)1)=\Psi(v,\rho(g)1),$$
therefore $g$ normalizes $1\in\Lambda^{ev}U$.
\end{Example}

\begin{Example}[Severi varieties]\label{KJHVKJNGVKJNGVK}\rm
Let ${\Bbb D}_\R$ denote one of four division algebras
over $\R$: real numbers $\R$, complex numbers $\C$, quaternions $\H$,
or octonions\index{octonions}~$\OO$. 
Let ${\Bbb D}={\Bbb D}_\R\otimes_\R\C$ be
its complexification, therefore $\Bbb D$ is either $\C$, $\C\oplus \C$,
$\Mat_2(\C)$, or $\Ca$ (the algebra of complex Cayley numbers\index{Cayley numbers}).
All these algebras have the standard involution.
Let ${\cal H}_3({\Bbb D})$ denote the $3\times3$ Hermitian 
matrices over $\Bbb D$. If $x\in{\cal H}_3({\Bbb D})$
then we may write
$$x=\left(\matrix{
c_1&u_1&u_2\cr
\ov{u_1} &c_2& u_3\cr
\ov{u_2} & \ov{u_3}& c_3
}\right), \quad c_i\in\C,\ u_i\in{\Bbb D}.
$$
${\cal H}_3({\Bbb D})$ is a Jordan algebra\index{Jordan algebra}
with respect to the operation $x\circ y=\bigover{xy+yx}{2}$,
that is, for any $a$, $b$ we have 
$$a\circ b=b\circ a\quad\text{and}\quad
((a\circ a)\circ b)\circ a=(a\circ a)\circ (b\circ b),$$ 
see~\cite{J}.
In case ${\Bbb D}=\Ca$ this algebra is called the Albert algebra\index{Albert algebra}.
The group of automorphisms $G=\Aut({\cal H}_3)$
is a semisimple algebraic group of type $\SO_3$, $\SL_3$, $\Sp_6$,
or $F_4$, respectively.
The action of $G$ on ${\cal H}_3$ has two irreducible components:
scalar matrices and traceless matrices.
It is possible to enlarge the group of automorphisms
by considering the group of norm similarities\index{norm similarities}.
Namely, ${\cal H}_3({\Bbb D})$ has a cubic form $\det$ defined
by the formula
$$\det(x)={1\over6}\left((\Tr x)^3+2\Tr(x^3)-3(\Tr x)\Tr(x^2)\right).$$
The group of norm similarities $\tilde G$ is defined as a subgroup in $\SL({\cal H}_3)$
preserving $\det$. $\tilde G$ is a semisimple algebraic group
of type $\SL_3$, $\SL_3\times\SL_3$, $\SL_6$, or $E_6$, respectively.
The representation of $\tilde G$ on ${\cal H}_3$ is irreducible.
The projectivization of the highest weight vector orbit
is called the Severi variety\index{Severi variety} corresponding to $\Bbb D$.
One may think of Severi varieties as of complexifications of 
projective planes over ${\Bbb D}_\R$.

The Severi variety corresponding to $\Ca$
is the minimal flag variety of $E_6$ (corresponding to the fundamental
weight $\omega_1$) embedded into $\P^{26}$ by the ample generator
of its Picard group. The corresponding elements of the Albert algebra
are known as elements of rank $1$, see~\cite{J}.
Other Severi varieties have the following description.

If ${\Bbb D}_\R=\R$ then ${\cal H}_3$ is equal to the Jordan algebra
of symmetric $3\times 3$ complex matrices and
$\det$ is the standard determinant.
$\tilde G=\SL_3$ acts on symmetric matrices
as on quadratic forms, the Severi variety is the projectivization
of the cone of rank $1$ symmetric matrices. Therefore, this Severi
variety is isomorphic to $\P^2$ in the second Veronese embedding.

If ${\Bbb D}_\R=\C$ then ${\cal H}_3$ is equal to the Jordan algebra
of $3\times 3$ complex matrices and
$\det$ is the standard determinant.
$\tilde G=\SL_3\times\SL_3$ acts on matrices
by left and right multiplication, 
the Severi variety is the projectivization
of the cone of rank $1$ matrices. Therefore, this Severi
variety is isomorphic to $\P^2\times\P^2$ in the Segre embedding.

If ${\Bbb D}_\R=\H$ then ${\cal H}_3$ is equal to the Jordan algebra
$\Lambda^2\C^6$ with the multiplication given by 
$\alpha_1\circ\alpha_2=\alpha_1\wedge\alpha_2$, where we identify
$\Lambda^2\C^6$ and $\Lambda^4\C^6$ by means of a standard
symplectic form $\omega$ (it is easy to see that $\Lambda^4\C^6=\omega\wedge\Lambda^2\C^6$).
The function $\det$ is defined as $\det(\alpha)=\alpha\wedge\alpha\wedge\alpha$
(in the matrix form $\det$ is equal to the Pfaffian of the skew-symmetric
matrix corresponding to $\alpha$). 
$\tilde G=\SL_6$ acts on $\Lambda^2\C^6$ traditionally.
The Severi variety is the projectivization of the variety
of decomposable bivectors of the form $u\wedge v$, $u,v\in\C^6$. 
Therefore, this Severi
variety is isomorphic to $\Gr(2,\C^6)$ in the Pl\"ucker embedding.

For ${\Bbb D}_\R$ equal to $\R$, $\C$, or $\H$
it is possible to define higher dimensional Severi varieties
in the same way (isomorphic to $\P^n$ in the second Veronese embedding,
$\P^n\times\P^n$ in the Segre embedding, or $\Gr(2,\C^{2n+2})$
in the Pl\"ucker embedding). However, for ${\Bbb D}_\R=\OO$
only the case of $3\times3$ matrices makes sense.
\end{Example}

\begin{Example}[Adjoint varieties]\rm
For the adjoint representation of $\SL_n=\SL(V)$ there is a natural
notion of the discriminant\index{discriminant of a linear operator} defined as follows.
For any operator $A\in\sl(V)$ let $P_A=\det(t\Id-A)$
be the characteristic polynomial.
Then its discriminant, $D(A)=D(P_A)$ is a homogeneous
form on $\sl(V)$ of degree $n^2-n$.
Clearly $D(A)\ne0$ if and only if all eigenvalues of $A$
are distinct, that is, if $A$ is a regular semisimple operator.
This notion can be carried over any simple Lie algebra.
Before doing that, notice that $D(A)$ can be also defined as follows.
Consider the characteristic polynomial 
$$Q_A=\det(t\Id-\ad(A))=\sum_{i=0}^{n^2-1}t^iD_i(A)$$
of the adjoint operator $\ad(A)$. If $\lambda_1,\ldots,\lambda_r$
are eigenvalues of $A$ (counted with multiplicities)
then 
the set of eigenvalues of $\ad(A)$ consists of 
$n-1$ zeros and 
$\lambda_i-\lambda_j$, $i,j=1,\ldots,n$, $i\ne j$.
Therefore, $D_0(A)=\ldots=D_{n-2}(A)=0$
and $D_{n-1}(A)$ coincides with $D(A)$ up to a non-zero scalar.

Suppose now that $\g$ is a simple Lie algebra of rank $r$. Let $\dim\g=n$.
For any $x\in \g$ let 
$$Q_x=\det(t\Id-\ad(x))=\sum_{i=0}^{n}t^iD_i(x)$$
be the characteristic polynomial of the adjoint operator.
Then $D(x)=D_r(x)$ is called 
the discriminant\index{discriminant of an element in a Lie algebra} 
of $x$.
Clearly, $D$ is a homogeneous $\Ad$-invariant polynomial on $\g$
of degree $n-r$.
Since the dimension of the centralizer $\g_x$ of any element $x\in\g$
is greater or equal to $r$, it follows that $D_0=\ldots=D_{r-1}=0$,
therefore $D(x)=0$ if and only if $\ad(x)$ has the eigenvalue $0$
with multiplicity $>r$. We claim that actually $D(x)\ne0$ if and only
if $x$ is regular semisimple
(recall that $x$ is called regular if $\dim\g_x=r$). Indeed, if $x$ is semisimple then $\ad(x)$
is a semisimple operator, therefore $D(x)=0$ if and only if
the dimension of the centralizer $\g_x$ is greater than $r$,
i.e.~$x$ is not regular. If $x$ is not semisimple then we take
the Jordan decomposition $x=x_s+x_n$, where $x_s$ is semisimple,
$x_n$ is nilpotent, and $[x_s,x_n]=0$. Then $x_s$ is automatically
not regular, therefore since $Q_x=Q_{x_s}$ we have $D(x)=D(x_s)=0$
(of course, this last equality also follows from the $\Ad$-invariance
of $D$ and the fact that the $\Ad$-orbit of $x$ contains the $\Ad$-orbit
of $x_s$ in its closure).

To study $D(x)$ further, we can use the Chevalley restriction
theorem (see~\cite{PV}) $\C[\g]^G=\C[\t]^W$, where $\t\subset\g$
is any Cartan subalgebra and $W$ is the Weil group.
Let $\Delta\subset\t^*$ be the root system, $\#\Delta=n-r$.
Then, clearly, for any $x\in\t$ $D(x)=0$ if and only if $x$ is not regular
if and only if $\alpha(x)=0$ for some $\alpha\in\Delta$.
Since $\deg D=n-r$ and $D|_\t$ is $W$-invariant, it easily
follows that 
$$D|_\t=\prod_{\alpha\in\Delta}\alpha.$$
The Weyl group acts transitively on the set of roots of 
the same length. Therefore, $D|_\t$, and hence $D$,
is irreducible if and only if all roots in $\Delta$ have the same
length, i.e.~$\Delta$ is of type $A$, $D$, or $E$.
If $\Delta$ is of type $B$, $C$, $F$, or~$G$, we have
$\Delta=\Delta_s\cup\Delta_l$, where $\Delta_s$ is the set
of short roots and $\Delta_l$ is the set of long roots.
Then we have $D=D_lD_s$, where $D_l$ and $D_s$ are irreducible polynomials
and 
$$D_l|_\t=\prod_{\alpha\in\Delta_l}\alpha, \quad
D_s|_\t=\prod_{\alpha\in\Delta_s}\alpha.$$
In the $A-D-E$ case we also set $D_l=D$ to simplify notations.

We are going to show that $D_l$ is also the discriminant in 
our regular sense. The adjoint representation $\Ad:\,G\to\GL(\g)$
is irreducible. In $A-D-E$ case let $\O=\O_l$
be the $\Ad$-orbit of any root vector,
In $B-C-F-G$ case let $\O_l\subset\g$ (resp.~$\O_s\subset\g$)
be the $\Ad$-orbit of any long root vector (resp.~any short root vector).
Then $\O_l$ is the orbit of the highest weight vector.
Its projectivization is called the adjoint variety.\index{adjoint variety}
Both orbits $\O_l$ and $\O_s$ are conical.
Let $X_l=\P(\O_l)$, $X_s=\overline{\P(\O_s)}$.

\begin{Theorem}\label{KJHKJHGJKFGH}
$D_l$ is the discriminant of $X_l$.
$D_s$ is the discriminant of $X_s$.
\end{Theorem}

\proof
We identify $\g$ and $\g^*$ via the Killing form.
Let $\alpha\in\Delta$ be any root, $e_\alpha\in \g$
the corresponding root vector. Since $[\g,e_\alpha]^\perp=\g_{e_\alpha}$,
we have
$$\dual{\left(\overline{\P(\Ad(G)\cdot e_\alpha)}\right)}=
\overline{\P(\Ad(G)\cdot\g_{e_\alpha})}.$$
Let 
$$\t_\alpha=\{x\in\t\,|\,\alpha(x)=0\}, \quad
\hat\t_\alpha=\t_\alpha+e_\alpha.$$
Then, clearly, $\hat\t_\alpha\in\g_{e_\alpha}$.
Moreover, if $\alpha$ is long then 
$D_l|_{\hat\t_\alpha}=D_l|_{\t_\alpha}=0$.
If $\alpha$ is short then 
$D_s|_{\hat\t_\alpha}=D_s|_{\t_\alpha}=0$.
Therefore, in order to prove Theorem it suffices to check that
\begin{equation}
\dim \Ad(G)\cdot\hat\t_\alpha=n-1.\label{KJHGJKGFJKKL}
\end{equation}
Clearly, for generic $x\in\hat\t_\alpha$ we have
\begin{equation}
\dim \Ad(G)\cdot\hat\t_\alpha=\dim G+\dim\hat\t_\alpha-\dim \Tran(x,
\hat\t_\alpha),\label{KJHKJHGFK}
\end{equation}
where 
$$\Tran(x,\hat\t_\alpha)=\{g\in G\,|\,\Ad(g)x\in \hat\t_\alpha\}.$$
Let $x=y+e_\alpha\in\hat\t_\alpha$, where $y\in\t_\alpha$ is such that
for any $\beta\in\Delta\setminus\{\pm\alpha\}$ we have $\beta(y)\ne0$.
Suppose that $\Ad(g)x\in\hat\t_\alpha$. Since $[e_\alpha,\t_\alpha]=0$,
it easily follows that $\Ad(g)y\in\t_\alpha$ and $\Ad(g)e_\alpha=e_\alpha$.
Under our assumptions on $y$ this means that
$$\Tran(x,\hat\t_\alpha)=(N_G(\t_\alpha))_{e_\alpha},$$
where $N_G(\t_\alpha)$ is the normalizer of $t_\alpha$ in $G$.
Therefore, in order to prove (\ref{KJHGJKGFJKKL}) using
(\ref{KJHKJHGFK}), it suffices to check that
\begin{equation}
\dim (\g_{\t_\alpha})_{e_\alpha}=r,\label{JHJHGDKHDKHGDF}
\end{equation}
where $\g_{\t_\alpha}$ is the centralizer of $t_\alpha$ in $\g$.
Now, $\g_{\t_\alpha}$ is a Levi subalgebra equal to
$\t+\C e_\alpha+\C e_{-\alpha}$.
Therefore, 
$(\g_{\t_\alpha})_{e_\alpha}=\t_\alpha+\C e_\alpha$
and (\ref{JHJHGDKHDKHGDF}) follows.
\endproof
\end{Example}

 
\section{The Pyasetskii Pairing}\label{PyasetskiiSection} 

\subsection{Actions With Finitely Many Orbits}

In \cite{Py} Pyasetskii has shown that if a connected algebraic group
acts linearly on a vector space with a finite number
of orbits then the dual representation has the same 
property and, moreover, the number of orbits
is the same. In this section we deduce this result
from the projective duality and give some examples
and applications.

\begin{Theorem}\label{JHGJHGJHGFK}
Suppose that a connected algebraic group $G$ acts
on a projective space $\P^n$ with a finite number of orbits.
Then the dual action $G:\dual{\P^n}$ has the same
number of orbits. Let $\P^n=\mathop{\sqcup}\limits_{i=1}^N\O_i$
and $\dual{\P^n}=\mathop{\sqcup}\limits_{i=1}^N\O'_i$ be the orbit
decompositions. Let $\O_0=\O'_0=\emptyset$.
Then the bijection is defined as follows:
$\O_i$ corresponds to $\O'_j$ if and only if
$\overline{\O_i}$ is projectively dual to $\overline{\O'_j}$.
\end{Theorem}

\proof
Indeed, take any $G$-orbit $\O\subset\dual{\P^n}$.
If $\overline\O\ne\dual{\P^n}$ then the projectively
dual variety $\dual{\overline\O}$ is the closure of some
orbit $\O'\subset\P^n$, since it is $G$-invariant.
Therefore, by the Reflexivity Theorem 
we have $\overline\O=\dual{\overline{\O'}}$.
\endproof

\begin{Remark}\rm
If the number of orbits for the action $G:\P^n$
is infinite then, in general, there is no natural
bijection between orbits in $\P^n$ and $\dual{\P^n}$.
However, it can be shown (see e.g.~\cite{Py}) 
that the modalities of these actions
are equal (modality is equal to the maximal number of parameters
that the family of orbits can depend on).
This also easily follows from the arguments involving
projective duality.
\end{Remark}

The linear action of an algebraic group on a vector space $V$
is called conical\index{conical action} 
if any $G$-orbit $\O$ is conical\index{conical orbit}, i.e.~$\O$
is preserved by homotheties. The action is conical if and only if 
for any $v\in V$, $v$ is contained in the tangent space $\g\cdot v$
to the orbit $Gv$. For example, if there are only finitely
many orbits then the action is conical.
If the action is conical then its non-zero orbits are 
in a 1-1 correspondence with orbits of the projectivization $G:\P(V)$.

\begin{Proposition}\label{JKHGKJGHFKGF}
If the linear action $G:V$ is conical then the dual action $G:\dual V$
is conical as well.
\end{Proposition}

\proof
We can embed $V=\C^n$ and $V^*=(\C^n)^*$
in projective spaces $\P^n$ and~$\dual{\P^n}$ as the affine charts.
These embeddings can be done equivariantly.
Points of $V^*$ correspond to hyperplanes in $\P^n$
not passing through the origin $0\in V\subset\P^n$.
Suppose that $\O\subset V^*$ is a non-conical orbit.
Then the dual variety $\dual{\overline\O}\subset\P^n$
intersects with $\C^n$ non-trivially.
Therefore, $\dual{\overline\O}\subset\P^n$ is the closure
of a conical variety in $\C^n$. Therefore, its dual variety
$\dual{\dual{\overline\O}}\subset\dual{\P^n}$ does not
intersect $(\C^n)^*$. But this contradicts the Reflexivity Theorem.
\endproof

Now the Pyasetskii Theorem\index{Pyasetskii Theorem} is an easy corollary
of Theorem~\ref{JHGJHGJHGFK} and Proposition~\ref{JKHGKJGHFKGF}.

\begin{Corollary}[\cite{Py}]\label{PyasetskyTh}
Suppose that a connected algebraic group $G$ acts linearly
on a vector space $V$ with a finite number of orbits.
Then the dual action $G:V^*$ has the same
number of orbits. Let $V=\mathop{\sqcup}\limits_{i=1}^N\O_i$
and $\dual{V}=\mathop{\sqcup}\limits_{i=1}^N\O'_i$ be the orbit
decompositions. 
Then the bijection is defined as follows:
$\O_i$ corresponds to $\O'_j$ if and only if
$\P(\overline{\O_i})$ is projectively dual to $\P(\overline{\O'_j})$.
\end{Corollary}

Recall that the Reflexivity Theorem is equivalent
to the formula (\ref{incidconescoinc:eq}) that
$\Lag(\O_i)=\Lag(\O'_j)$, where for any conical variety
$Y\subset V$ we denote by $\Lag(Y)$ 
the closure of the conormal bundle
$N^*_{Y_{sm}}V$ in the cotangent bundle $T^*V=V\oplus V^*$.
The varieties $\Lag(\O_i)$, $i=1,\ldots,N$ have a nice interpretation.
Namely, the action of $\g$ (the Lie algebra of $G$) on $V$
is given by the element of 
$$\Hom(\g,\End(V))=\g^*\otimes V^*\otimes V=\Hom(V^*\otimes V,\g^*).$$
The corresponding bilinear map $\mu:\,V\oplus V^*\to\g^*$
is called the moment map\index{moment map}. 
Let $L=\mu^{-1}(0)$ be its zero fiber.
Then it is easy to see that $L=\cup_i\Lag(\O_i)$
is the decomposition of $L$ into irreducible components.

The large class of linear actions of reductive groups
with finitely many orbits is provided by graded
semisimple Lie algebras.
Suppose that $\g=\oplus_{k\in\Z}\g_k$
is a graded semi-simple Lie algebra of 
the connected semi-simple group $G$.
Then there exists a unique semisimple element $\xi\in\g_0$
such that 
$$g_k=\{x\in \g\,|\,[\xi,x]=kx\}.$$
The connected component $H$ of the centralizer $G_\xi\subset G$
is a reductive subgroup of $G$ called the Levi subgroup\index{Levi subgroup}
(because it is the Levi part of the parabolic subgroup).
$\g_0$ is the Lie algebra of $H$.
$H$ acts on each graded component $\g_k$.

\begin{Theorem}[\cite{Ri,Vi}]\label{JHGFJHGFJHG}
$H$ acts on $\g_k$ with finitely many orbits.
\end{Theorem}

\proof
Clearly all elements of $\g_k$ are nilpotent in $\g$.
Since $G$ has finitely many nilpotent orbits,
it is sufficient to prove that for any nilpotent $G$-orbit~$\O$
the intersection $\O\cap\g_k$ consists of finitely many $H$-orbits.
By the general Richardson Lemma\index{Richardson Lemma} (see~\cite{Ri,PV}),
it is enough to check that for any $x\in\O\cap\g_k$
the tangent space $[\g_0,x]$ to the $H$-orbit $\Ad(H)x$
is equal to the intersection of $\g_k$ with the tangent space
$[\g,x]$ to the $G$-orbit $\Ad(G)x$.
But this is clear:
$$[\g,x]\cap\g_k=[\oplus_i\g_i,x]\cap\g_k=\oplus_i[\g_i,x]\cap\g_k=[\g_0,x].$$
Theorem is proved.
\endproof

Notice that the dual $H$-module $(\g_k)^*$ is isomorphic to $\g_{-k}$
under the Killing form and the moment map $\g_k\times\g_{-k}\to\g_0$
is given by the Lie bracket.


\subsection{The Multisegment Duality.}

Apart from actions associated with $\Z$-graded semisimple
Lie algebras, another class of actions with finitely many orbits
is provided by the theory of representations of quivers (see~\cite{Ga}).
These two classes overlap: the representations of quivers of type $A$
give the same class of actions as the standard gradings of 
of $\SL(V)$. The Pyasetskii pairing in this case was studied
in a series of papers under the name of the multisegment
duality\index{multisegment duality}, 
or the Zelevinsky involution\index{Zelevinsky involution}.
Here we give an overview of some of these results.

We fix a positive integer $r$ and 
consider the set $S=S_r$ of pairs of integers
$(i,j)$ such that $1\le i\le j\le r$.
Let $\Z_+^S$ denote the semigroup of families
$m=(m_{ij})_{(i,j)\in S}$ of non-negative integers
indexed by $S$. The set $S$ and the semigroup $\Z_+^S$
have several interpretations. First, $S_r$ is naturally
identified with the set of positive roots of type $A_r$.
Namely, if $\alpha_1,\ldots,\alpha_r$ denote simple roots
then each $(i,j)\in S$ corresponds to a positive root
$\alpha_i+\alpha_{i+1}+\ldots+\alpha_j$.
Another useful interpretation is to regard a pair $(i,j)\in S$
as a segment $[i,j]=\{i,i+1,\ldots,j\}$ in $\Z$.
A family $m=m_{ij}\in\Z_+^S$ can be regarded as a collection
of segments, containing $m_{ij}$ copies of each $[i,j]$.
Thus, elements of $\Z_+^S$ can be called multisegments\index{multisegment}.
The weight\index{weight of a multisegment} 
$|m|$ of a multisegment $m$ is defined as a sequence 
$\gamma=\{d_1,\ldots,d_r)\in\Z_+^r$ given by
$$d_i=\sum_{i\in[k,l]}m_{kl}\quad\text{for}\ i=1,\ldots,r.$$
In other words, $|m|$ records how many segments of $m$ contain any given number
$i\in[1,r]$. For any $\gamma\in\Z_+^r$ we set 
$\Z_+^S(\gamma)=\{m\in\Z_+^S\,|\,|m|=\gamma\}$.
In the language of positive roots, $\Z_+^r$ is the semigroup
generated by simple roots. Under this identification,
the elements $m\in\Z_S^+(\gamma)$ become the partitions of $\gamma$
into a sum of positive roots. The cardinality of
$\Z_S^+(\gamma)$ is known as a Kostant partition function\index{Kostant partition function} 
(of $\gamma$).

Another important interpretation of $\Z_S^+(\gamma)$ 
is that it parametrizes isomorphism classes
of representations of quivers of type $A$.
Let $A_r$ be the quiver equal to the Dynkin diagram of type $A_r$,
where all edges are oriented from the left to the right.
Let $A_r^*$ be a dual quiver\index{dual quiver} with all orientations reversed.
The representation\index{representations of quivers} 
of $A_r$ with the dimension
vector $\gamma=\{d_1,\ldots,d_r)\in\Z_+^r$ is the collection of vector
spaces $\C^{d_1}, \ldots, \C^{d_r}$ and linear maps
$$\varphi_1:\,\C^{d_1}\to\C^{d_2}, \ldots, \varphi_{r-1}:\,\C^{d_{r-1}}\to
\C^{d_r}.$$
The representation of $A_r^*$ with the dimension
vector $\gamma=\{d_1,\ldots,d_r)\in\Z_+^r$ is the collection of vector
spaces $\C^{d_1}, \ldots, \C^{d_r}$ and linear maps
$$\varphi_1:\,\C^{d_2}\to\C^{d_1}, \ldots, \varphi_{r-1}:\,\C^{d_{r}}\to
\C^{d_{r-1}}.$$
Therefore, representations of $A_r$ with dimension vector $\gamma$ 
are parametrized
by points of a vector space 
$V(\gamma)=\mathop{\bigoplus}_{i=1}^{r-1}\Hom(\C^{d_i},\C^{d_{i+1}})$.
Representations of $A_r^*$ with dimension vector $\gamma$ 
are parametrized
by points of a vector space 
$V(\gamma)^*=\mathop{\bigoplus}_{i=1}^{r-1}\Hom(\C^{d_{i+1}},\C^{d_{i}})$.
Notice that vector spaces $V(\gamma)$ and $V(\gamma)^*$ are 
naturally dual to each other.
Moreover, $V(\gamma)$ and $V(\gamma)^*$ are dual modules
of the group $G(\gamma)=\GL_{d_1}\times\ldots\times\GL_{d_r}$
with respect to the natural action.
These actions could also be described in terms
of graded Lie algebras.
Namely, consider the Lie algebra $\g=\gl_{d_1+\ldots+d_r}$.
Elements of $\g$ can be represented as block matrices with diagonal blocks
of shapes $d_1\times d_1$, $\ldots$, $d_r\times d_r$.
This determines $\Z$-grading: block diagonal matrices have grade $0$, etc.
Then $\g_0$ is a Lie algebra of $G(\gamma)$ and the action of $G(\gamma)$
on $V(\gamma)$ and $V(\gamma)^*$ is equivalent to the 
action of $G(\gamma)$ on $\g_{-1}$ and $\g_{1}$.

The orbits of $G(\gamma)$ on $V(\gamma)$ (resp.~$V(\gamma)^*$)
correspond to  isoclasses of representations of $A_r$ (resp.~$A_r^*$)
with dimension vector $\gamma$.
These orbits are parametrized by elements of 
$\Z_S^+(\gamma)$. Elements of $S$ parametrize
indecomposable $A_r$-modules (or $A_r^*$-modules).
Namely, each $(i,j)\in S$ corresponds to an indecomposable
module $R_{ij}$ with dimension vector 
$$|(i,j)|=(0^{i-1},1^{j-i+1},0^{r-j})$$
and all maps are isomorphisms, if it is possible, or zero maps
otherwise.
Then any family $(m_{ij})\in\Z_S^+$ corresponds
to an $A_r$ (or $A_r^*$) module $\oplus_SR_{ij}^{m_{ij}}$.

By Pyasetskii Theorem~\ref{PyasetskyTh}, there is a natural bijection
of $G(\gamma)$-orbits in $V(\gamma)$ and $V(\gamma)^*$.
Therefore, there exists a natural involution $\zeta$ of $\Z_+^S(\gamma)$,
which can be extended to a weight-preserving involution
of $\Z_+^S$ called the multisegment duality.

An explicit description of $\zeta$ was found in~\cite{KZ}
using the Poljak Theorem~ \cite{Po} from the combinatorial theory of networks.
To formulate it, we need the following definition.
For any multisegment $m\in\Z_+^S$ the ranks $r_{ij}(m)$ are given by
$$r_{ij}(m)=\sum_{[i,j]\subset[k,l]}m_{kl}.$$
It is easy to see that the multisegment $m$ can be recovered
from its ranks by formula
$$m_{ij}=r_{ij}(m)-r_{i-1,j}(m)+r_{i-1,j+1}(m).$$
If the multisegment corresponds to the representation
of $A_r$ given by 
$$(\varphi_1,\ldots,\varphi_{r-1})\in V(\gamma)$$
then $r_{ij}$ is equal to the rank of the map
$\varphi_{j-1}\circ\ldots\circ\varphi_{i+1}\circ\varphi_i$.
In particular, $r_{ii}=d_i$.

For any $(i,j)\in S$ let $T_{ij}$ denote the set of all maps
$\nu:\,[1,i]\times[j,r]\to[i,j]$ such that
$\nu(k,l)\le\nu(k',l')$ whenever $k\le k'$, $l\le l'$
(in other words, $\nu$ is a morphism of partially ordered sets,
where $[1,i]\times[j,r]$ is supplied with the product order).

\begin{Theorem}[\cite{KZ}]
For every $m=(m_{ij})\in\Z_+^S$ we have
$$r_{ij}(\zeta(m))=\mathop{\min}_{\nu\in T_{ij}}
\sum_{(k,l)\in[1,i]\times[j,r]}
m_{\nu(k,l)+k-i,\nu(k,l)+l-j}.$$
\end{Theorem}

The inductive description of $\zeta$ was given in \cite{MW}.
Let $m=(m_{ij})$ be a multisegment of weight $\gamma=(d_1,\ldots,d_r)$.
We set $i_1=\min\{i\,|\,d_i\ne0\}$ and define the sequence
of indices $j_1,\ldots,j_p$ as follows:
$$j_1=\min\{j\,|\,m_{i_1j}\ne0\},\ldots,
j_{t+1}=\min\{j\,|\,j>j_t,\ m_{i_1+t,j}\ne0\},$$
where $t=1,\ldots,p-1$.
The sequence terminates when $j_{p+1}$ does not exist.
Let $i_t=i_1+t-1$ for $t=1,\ldots,p+1$.
We associate to $m$ the multisegment
$m'$ given by
$$m'=m-(i_1,j_1)-(i_2,j_2)-\ldots-(i_p,j_p)+(i_2,j_1)+(i_3,j_2)+\ldots+(i_{p+1},j_p),$$
where we use the convention that $(i,j)=0$ unless $1\le i\le j\le r$.

\begin{Theorem}[\cite{MW}]
If the multisegment $m'$ is associated to $m$ then
$$\zeta(m)=\zeta(m')+(i_1,i_p).$$
\end{Theorem}

The involution $\zeta$ can also be described in terms
of irreducible finite-dimensional 
representations of affine Hecke algebras
and in terms of canonical bases for quantum groups, see~\cite{KZ}.


\section{Parabolic Subgroups With Abelian Unipotent Radical}\label{ParabolicsAura}
Let $L$ be a simple algebraic group and $P\subset L$ a parabolic subgroup
with abelian unipotent radical. In this case $\l=\Lie L$ admits
a $\Z$-grading with only three non-zero parts:
$$\l=\l_{-1}\oplus\l_0\oplus\l_1.$$
Such a grading is said to be short\index{short grading}. 
Here $\l_0\oplus\l_1=\Lie P$
and $\exp(\l_1)$ is the abelian unipotent radical of $P$.

There exists a unique semisimple element $\xi\in\l_0$
such that 
$$l_k=\{x\in \l\,|\,[\xi,x]=kx\}.$$
We denote by $G$ the connected component 
of the centralizer $L_\xi\subset L$.
Then $\l_0=\g=\Lie G$ and by Theorem~\ref{JHGFJHGFJHG}
$G$ acts on $\l_{\pm1}$ with a finite number of orbits.
In this section we describe explicitly  these $G$-orbits and the Pyasetskii
pairing. These results are well-known,
we follow mainly \cite{Pan} and \cite{MRS}.

We denote by $T$ the maximal torus in $G$ (and hence in $L$), $\t=\Lie T$,
$\Delta$ is the root system of $(\l,\t)$, 
$$\Delta=\Delta_{-1}\cup\Delta_0\cup\Delta_1$$
is the partition corresponding to the short grading.
We fix a Borel subgroup $B\subset G$ containing $T$.
This choice determines a set of positive roots  
$$\Delta_0^+\subset\Delta_0\quad
\text{and}\quad \Delta^+=\Delta^+_0\cup\Delta^1\subset\Delta.$$
We denote by $\gamma$ the highest root in $\Delta^+$, in fact in $\Delta_1$.
$\Pi$ is the set of simple roots in $\Delta^+$, then
$\Pi_0=\Pi\cap\Delta_0$ is the set of simple roots in $\Delta_0^+$.
Let $W$ (resp.~$W_0$) be the Weil group of $\l$ (resp.~of $\g$)
with respect to $\t$.
For $\alpha\in\Delta$, we let $w_\alpha$ denote the corresponding reflection
in $W$, $\alpha^\vee=\bigover{2\alpha}{(\alpha,\alpha)}$ the corresponding coroot,
and $e_\alpha$ a non-zero root vector.

\begin{Proposition}
The representation of $\g$ on $\l_1$ is faithful and irreducible. 
\end{Proposition}

\proof
Indeed, $\k=\{x\in\g\,|\,\ad(x)\l_1=0\}$ is an ideal in $\g$.
Since $\l_1$ and $\l_{-1}$ are dual $\g$-modules,
$[\k,\l_{-1}]=0$. It follows that $\k$ is an ideal in $\l$.
Therefore, $\k=0$ since $\l$ is simple.

Suppose now that $\l_1$ is not irreducible, $\l_1=\l_1^1\oplus\l_1^2$
is the $\g$-module decomposition.
Since $\l_{-1}$ is dual to $\l_1$, we have $\l_{-1}=\l_{-1}^1\oplus\l_{-1}^2$,
where $\l_{-1}^i=(\l_1^j)^\perp$ for $i\ne j$.
Then $[\l_{-1}^i,\l_1^j]=0$ for $i\ne j$. Indeed, for any $x\in\l_{-1}^i$,
$y\in\l_0$, $z\in\l_1^j$ we have $(y,[x,z])=-([x,y],z)=0$,
therefore $[x,z]=0$.
It follows that $\l_{-1}^i\oplus[\l_{-1}^i,\l_1^i]\oplus\l_1^i$
is an ideal in $l$. But $\l$ is simple. Contradiction.
\endproof

Therefore, the center of $\g$ is one-dimensional,
spanned by $\xi$, $P$ is a maximal parabolic subgroup,
and $\#\Pi_0=\#\Pi-1$. Thus, there is a unique simple root in $\Delta_1$.
Call it $\beta$. 

\begin{Proposition}\label{KJHKHlwsjkdfghvlhK}
$\ $
\begin{enumerate}
\item[\rm(a)] $\beta$ is long.
\item[\rm(b)] The $\beta$-height of any root $\alpha_0\in\Delta_1$,
i.e.~the coefficient $n_\beta$ in the sum 
$\alpha_0=n_\beta\beta+\sum_{\alpha\in\Pi_0}n_\alpha\alpha$,
is equal to $1$. 
\item[\rm(c)] For any $\alpha,\alpha'\in\Delta_1$, $(\alpha,\alpha')\ge0$.
\item[\rm(d)] $W_0$ acts transitively on long roots in $\Delta_1$.
\end{enumerate}
\end{Proposition}

\proof
(a) Indeed, since $\l_1$ is an irreducible $\g$-module,
$\beta$ is the unique lowest weight of the $\g$-module $\l_1$.
The longest element in $W_0$ takes $\beta$ to the highest weight of $\l_1$,
i.e.~to $\gamma$. Hence $\beta$ is long. 

(b) Since $\l_1=U\g\cdot\beta$, 
the $\beta$-height of any root in $\Delta_1$
is equal to $1$. 

(c) Suppose that $\alpha,\alpha'\in\Delta_1$. Since $[e_\alpha,e_{\alpha'}]=0$,
it follows that $\alpha+\alpha'$ is not a root, 
therefore, $(\alpha,\alpha')\ge0$.

(d) Suppose that $\alpha_0\in\Delta_1$ is a long root,
$\alpha_0=\beta+\sum_{\alpha\in\Pi_0}n_\alpha\alpha$.
If we have $(\alpha_0,\alpha)<0$ for some $\alpha\in\Pi_0$ then
$w_\alpha(\alpha_0)>\alpha_0$ and we may finish by induction.
Suppose that $(\alpha_0,\alpha)\ge0$ for all
$\alpha\in\Pi_0$. Since $(\alpha_0,\beta)\ge0$ by (c),
$\alpha_0$ is a dominant weight for $\l$.
Since $\alpha_0$ is long, it follows that $\alpha_0=\delta$.
Therefore, all long roots in $\Delta_1$ are $W_0$-conjugate to $\delta$.
\endproof

On the contrary, if $\l$ is an arbitrary simple Lie algebra,
$\Delta\supset\Pi$ is the root system and the simple roots,
$\beta\in\Pi$, and the $\beta$-height of the highest root is equal to $1$,
then, clearly, the maximal parabolic subgroup corresponding to $\beta$
has an abelian unipotent radical. Using this observation,
it is easy to find all parabolic subgroups with abelian unipotent radical
in simple algebraic groups. 
If $\beta=\alpha_k$ in the Bourbaki numbering of simple roots then
we shall denote the corresponding parabolic subgroup by $P_k$.
\begin{itemize}
\item
$A_r$: any maximal parabolic subgroup $P$
has an abelian unipotent radical. 
\item
$B_r$: the only possibility is $P_1$.
\item
$C_r$: the only parabolic subgroup with abelian unipotent radical is $P_r$.
\item
$D_r$: there are three possibilities: $P_1$, $P_{r-1}$, and $P_r$
\item
Finally, there are two exceptional cases: $P_1$ (or $P_6$) of $E_6$
and $P_7$ of $E_7$.
\end{itemize}

Most results about shortly graded simple Lie algebras can be proved
by induction using the following procedure.
Let 
$$\Delta'=\{\alpha\in\Delta\,|\,(\alpha,\beta)=0\}.$$
Then $\Delta'=\Delta'_{-1}\cup\Delta'_0\cup\Delta'_1$ is the root
system of the graded semisimple subalgebra $\l'\subset\l$.

\begin{Proposition}
Suppose that $\Delta'_1\ne\emptyset$. Then
\begin{enumerate}
\item[\rm(a)] $\delta\in\Delta'_1$.
\item[\rm(b)] $\l'_0=\l_0''\oplus\l_0'''$, where $\l'_{-1}\oplus\l''_0\oplus\l'_1$
is a shortly graded simple Lie algebra.
\item[\rm(c)] $\Delta'_1$ contains a unique minimal root, say $\beta'$.
$\beta'$ is long.
\end{enumerate}
\end{Proposition}

\proof
Indeed, suppose that $(\beta,\alpha_0)=0$ for some $\alpha_0\in\Delta_1$.
Then 
$$(\beta,\delta)=(\beta,\alpha_0)+
\left(\beta,\sum_{\alpha\in\Pi_0}n_\alpha\alpha\right),$$
where all $n_\alpha\ge0$.
Since $(\beta,\alpha)\le0$ for any $\alpha\in\Pi_0$,
we get $(\beta,\delta)\le0$. Since $(\beta,\delta)\ge0$ by
Proposition~\ref{KJHKHlwsjkdfghvlhK}, we see that $\delta\in\Delta'_1$.
It follows that $\delta$ is the unique highest weight in $\l'_0$-module
$\l'_1$. Therefore, $\l'_1$ is an irreducible $\l'_0$-module and 
we obtain (b). Now we can apply Proposition~\ref{KJHKHlwsjkdfghvlhK}
to the shortly graded simple Lie algebra $\l'_{-1}\oplus\l''_0\oplus\l'_1$
and get (c).
\endproof

It is possible to define two canonical strings of 
pairwise orthogonal long roots
in $\Delta_1$. The construction is originally due to Harish Chandra.
The lower canonical string\index{lower canonical string}  $\beta_1,\ldots,\beta_r$
is the cascade up from $\beta_1=\beta$
within $\Delta_1$: at each stage $\beta_{i+1}$ is the minimal root in 
$$\Delta^{(i)}_1=\{\alpha\in\Delta_1\,|
\,(\alpha,\beta_1)=\ldots=(\alpha,\beta_i)=0\}.$$
The process terminates when $\Delta^i_1$ is empty.
The previous proposition shows that
the lower canonical string $\beta_1,\ldots,\beta_r$
is well-defined, i.e.~on each step
$\Delta_1^i$ contains a unique minimal element.
Clearly, all the roots $\beta_i$ are long.

The upper canonical string\index{upper canonical string} 
$\gamma_1,\ldots,\gamma_{r'}$
is the cascade down from $\gamma_1=\gamma$
within $\Delta_1$. At each stage, $\gamma_{i+1}$ is the unique maximal
element in 
$${}^{(i)}\Delta_1=\{\alpha\in\Delta_1\,|
\,(\alpha,\gamma_1)=\ldots=(\alpha,\gamma_i)=0\}.$$
The process terminates when ${}^{(i)}\Delta_1$ becomes empty.
Obviously, the longest element in $W_0$ takes the lower canonical string
to the upper canonical string.
Hence, the upper canonical string is well-defined and both strings
have the same cardinality.
In fact, we have the following proposition

\begin{Proposition}[\cite{RRS}]\label{KJHGKJVKJVKV}
Suppose that $\beta'_1,\ldots,\beta'_{r'}$ 
is a string of pairwise orthogonal long roots in $\Delta_1$.
Then $r'\le r$ and there exists $w\in W_0$ such that
$\beta'_i=w(\beta_i)$ for $i\le r'$.
\end{Proposition}

\proof
By Proposition~\ref{KJHKHlwsjkdfghvlhK} (d) we may suppose that
$\beta'_1=\beta_1$. Now the proof goes by induction:
both $\beta'_2,\ldots,\beta'_{r'}$ and
$\beta_2,\ldots,\beta_{r}$ belong to $\Delta'_1$ and, hence,
there exists $w\in W'_0$ (the Weil group of $\l_0'$) such that
$\beta'_i=w(\beta_i)$ for $i\le r'$.
It remains to notice that $w\beta_1=\beta_1$.
\endproof

Now we can describe the orbit decomposition and the Pyasetskii
pairing.

\begin{Theorem}\label{AURAOrbitDecomp}
Let $\alpha_1,\ldots,\alpha_r\in\Delta_1$ be any maximal sequence
of pairwise orthogonal long roots.
Set $e_k=e_{\alpha_1}+\ldots+e_{\alpha_k}$ for $k=0,\ldots,r$
{\rm(}so $e_0=0${\rm)}. 
Denote the $G$-orbit of $e_k$ by $\O_k\subset\l_1$. 
Set $f_k=e_{-\alpha_1}+\ldots+e_{-\alpha_k}$ for $k=0,\ldots,r$
{\rm(}so $f_0=0${\rm)}. 
Denote the $G$-orbit of $f_k$ by $\O'_k\subset\l_{-1}$. 
Then
\begin{enumerate}
\item[\rm(a)] $\l_1=\mathop{\sqcup}\limits_{i=0}^r\O_i$,
$\l_{-1}=\mathop{\sqcup}\limits_{i=0}^r\O'_i$.
\item[\rm(b)] $O_i\subset\overline{O_j}$ if and only if $i\le j$ 
if and only if $O'_i\subset\overline{O'_j}$.
\item[\rm(c)] $\O_k$ corresponds to $\O'_{r-k}$ in the Pyasetskii
pairing.
\item[\rm(d)] If $h_k=\alpha_1^\vee+\ldots+\alpha_k^\vee$ then
$\langle e_k,h_k,f_k\rangle$ is a homogeneous $\sl_2$-subalgebra.
\end{enumerate}
\end{Theorem}

\proof
(d) 
This is obvious.

(a) 
By Proposition~\ref{KJHGKJVKJVKV} the choice of a maximal string is irrelevant.
We take $\alpha_i=\beta_i$, where $\beta_1,\ldots,\beta_r$
is the lower canonical string.
First, we need to show that any element in $\l_1$ is $G$-conjugate
to one of $e_k$. We argue by the induction on $r$.
Let $x\in\l_1$, $x\ne0$. 
We can write $x$ as 
$$x=\sum_{\alpha\in\Delta_1}x_{\alpha}e_{\alpha}.$$
Since $\l_1$ is an irreducible $G$-module,
after conjugating $x$ we may assume that $x_{\beta}\ne0$.
Since the $G$-orbit of $x$ is conical, we may assume that $x_{\beta}=1$.
Let $\hat\Delta_0=\{\alpha\in\Delta_0\,|\,(\alpha,\beta)<0$.
Clearly, $\hat\Delta_0\subset\Delta_0^+$. Then 
$$\uu=\langle e_\alpha\,|\,\alpha\in\hat\Delta_0\rangle\subset\g$$
is a unipotent subalgebra. Let $U\subset G$ be an unipotent subgroup with $\uu=\Lie U$.
Then $U$ acts on $\l_1$. Clearly, for any roots $\alpha\in\hat\Delta_0$,
$\alpha'\in\Delta_1'$ we have $\alpha+\alpha'\not\in\Delta$ by 
Proposition~\ref{KJHKHlwsjkdfghvlhK} (c). Therefore, 
$U$ acts trivially on $\l_1'$. Consider the action of $U$ on $\l_1/\l_1'$.
Clearly, the function $x_\beta$ is $U$-invariant, therefore
each $U$-orbit in $\l_1/\l_1'$ lies on a hyperplane $x_\beta=\const$.
For any root $\alpha\in\Delta_1$ such that $(\alpha,\beta)>0$ we have
$\alpha=\beta+\hat\alpha$, where $\hat\alpha\in\hat\Delta_0$.
Therefore, the tangent space $\uu\cdot (x \mod \l'_1)$ is equal to the hyperplane
$x_\beta=0$. Since any orbit of a linear unipotent group is closed \cite{OV},
it follows that the $U$-orbit of $x \mod \l_1'$ is the affine hyperplane
$x_\beta=1$. In particular, $x$ is $U$-conjugate to an element
$e_\beta+x'$, where $x'\in\l_1'$.
By induction, $x'$ is $G'$-conjugate
to one of $e_0',\ldots,e_{r-1}'$, where $\Lie G'=\l'_0$.
Since $G'\cdot e_\beta=e_\beta$, we see that any nonzero 
element of $\l_1$ is $G$-conjugate to 
$e_\beta+e'_k=e_{k+1}$ for some $k=0,\ldots,r-1$.

To finish the proof of (a) we need to show that all orbits $\O_i$
are distinct. To prove this it suffices to check that
$\dim\g\cdot e_i=\dim\O_i<\dim\O_j=\dim\g\cdot e_j$ for $i<j$.
This will follow, in turn, from the inequality
$$\dim\Ann(\g\cdot e_i)=\dim\l_{-1}^{e_i}>\dim\l_{-1}^{e_j}=\dim\Ann(\g\cdot e_i),$$
where $\l_{-1}^{e_i}=\{x\in\l_{-1}\,|\,[x,e_i]=0\}$.
Since $[f_i,\l_{-1}]=0$, it follows from the $\sl_2$-theory
that $\l_{-1}^{e_i}=\l_{-1}^{h_i}$.
But 
$$\l_{-1}^{h_i}=\langle e_{-\alpha}\,|\,\alpha\in\Delta^{(i)}\rangle,$$
therefore, $\dim\l_{-1}^{e_i}=\#\Delta^{(i)}$.
Since $\#\Delta^{(i)}>\#\Delta^{(j)}$, we see that all orbits $\O_i$
are distinct.

It is clear from the above that 
$\Ann(\g\cdot e_i)\cap Gf_{r-i}$ is open in $\Ann(\g\cdot e_i)$,
therefore, we have (c).

Finally, to prove (b) we need only to show that $\O_i\subset\overline{\O_j}$
for $i<j$. But 
$$\mathop{\lim}_{t\to-\infty}\exp{t(\alpha_{i+1}^\vee+\ldots+\alpha_{j}^\vee)}
e_j=e_i.$$
\endproof

After these case-by-case-independent considerations it
is worthy to sort out all examples.

\begin{Example}\label{Determinant}\rm
Consider the short gradings of $\l=\sl_{n+m}$.
Then 
$$G=\{(A,B)\in \GL_n\times \GL_m\,|\,\det(A)\det(B)=1\}.$$
$\l_{1}$ can be identified with $\C^n\times\C^m$.
There are $r=\min(n,m)$ non-zero $G$-orbits.
Namely, $\O_i$, $i=1,\ldots,r$, 
is the variety of $m\times n$-matrices of rank $i$.
The projectivization of $\O_1$ is identified with $X=\P^{n-1}\times\P^{m-1}$
in the Segre embedding.
Therefore, the dual variety $\dual X$ is equal to the 
projectivization of the closure $\overline{\O_{r-1}}$, the variety
of matrices of rank less than or equal to $r-1$.
$\dual X$ is a hypersurface if and only if $n=m$, in which case
$\Delta_X$ is the ordinary determinant of a square matrix.
Another interesting case is $n=2$, $m\ge2$: we see that
the Segre embedding of $\P^1\times\P^k$ is self-dual.
\end{Example}

\begin{Example}\rm
Consider the short grading of $\l=\so_{n+2}$ that corresponds to $\beta$
being the first simple root.
Then $G=\C^*\times\SO_n$ and $l_{1}=C^n$ with a simplest action.
There are two non-zero orbits: the dense one and the
self-dual quadric hypersurface $Q\subset\C^n$ preserved by $G$.
\end{Example}

\begin{Example}\label{Pfaffian}\rm
Consider the short grading of $\l=\so_{2n}=D_{n}$ that corresponds
to $\beta=\alpha_n$ or $\beta=\alpha_{n-1}$.
Then $G=\C^*\times\SL_n$ acts naturally on $\l_1=\Lambda^2\C^n$.
There are $r=[n/2]$ non-zero orbits, where $\O_i$, $i=1,\ldots,r$,
is the variety of skew--symmetric matrices of rank $2r$.
The projectivization of $\O_1$ is identified with $X=\Gr(2,n)$
in the Pl\"ucker embedding.
Therefore, the dual variety $\dual X$ is equal to the 
projectivization of the closure $\overline{\O_{r-1}}$, the variety
of matrices of rank less than or equal to $2r-2$.
$\dual X$ is a hypersurface if and only if $n$ is even, 
in which case
$\Delta_X$ is the Pfaffian of a skew-symmetric matrix.
If $n$ is odd then $\defect X=\codim\O_{r-1}=2$. 
\end{Example}

\begin{Example}\label{KGFJHFDHGFDGF}\rm
The short grading of $E_6$ gives the following action.
$G=\C^*\times\SO_{10}$ and $\l_1$ is the half-spinor representation.
Except for the zero orbit and the dense orbit there is only one orbit $\O$.
In particular, the projectivization of $\O$ is smooth and self-dual.
It is  a spinor variety $\SS_5$.
It could also be described via Cayley numbers.
Let $\Ca$ be the algebra of split Cayley numbers
(therefore $\Ca=\OO\otimes\C$, where $\O$ is the real division algebra
of octonions). Let $u\mapsto\overline u$ be the canonical involution in $\Ca$.
Let $\C^{16}=\Ca\oplus\Ca$ have octonionic coordinates $u,v$.
Then the spinor variety $\S_5\subset\P(\C^{16})$
is defined by homogeneous equations
$$u\ov u=0,\ v\ov v=0, u\ov v=0,$$
where the last equation is equivalent to $8$ complex equations.
\end{Example}

\begin{Example}\rm
The final example appears from the short grading of $E_7$.
Here $G=\C^*\times E_6$, and $\l_1=\C^{27}$ can be identified with the
exceptional simple Jordan algebra (the Albert algebra\index{Albert algebra}),
see~\cite{J}.
Then $E_6$ is the group of norm similarities\index{group of norm similarities}
and $\C^*$ acts by homotheties. There are $3$ non-zero orbits:
the dense one, the cubic hypersurface (defined by the norm
in the Jordan algebra), and the  closed conical
variety with smooth projectivization consisting of elements of rank one,
i.e.~the exceptional Severi variety
(the model of the Cayley projective plane).
\end{Example}

\chapter{The Cayley Method for Studying Discriminants}

\section*{Preliminaries} 
The remarkable observation due to Cayley \cite{Ca} 
gives an expression of the discriminant
as the determinant of a certain complex of finite-dimensional vector spaces.
Of course, complexes were unknown at his time and he was studying
resultants rather than discriminants (though he also noticed that
resultants and discriminants are, in fact, the equivalent notions).
This approach was clarified and developed in the beautiful serie
of papers by Gelfand, Kapranov, and Zelevinsky \cite{GKZ1,GKZ2}.

\section{Jet Bundles and Koszul Complexes}\index{jet bundle}
Let $X$ be a smooth irreducible algebraic variety with 
an algebraic line bundle $\L$.
We consider the bundle $J(\L)$ of first jets of sections of $\L$.
By definition, the fiber of $J(\L)$ at a point $x\in X$
is the quotient of the space of all sections of $\L$
near $x$ by the subspace of sections which vanish at $x$
with their first derivatives. In other words, $J(\L)_x=\L/I_x^2\L$,
where $I_x$ is the ideal of functions vanishing at $x$.
Thus $J(\L)$ is a vector bundle of rank $\dim X+1$.

To any section $f$ of $\L$, we associate a section $j(f)$ of $J(\L)$,
called the first jet of $f$. \index{first jet}
Namely, the value of $j(f)$ at $x$ is the class of $f$ modulo $I_x^2\L$.
The correspondence $f\mapsto j(f)$ is $\C$-linear, but being a differential operator,
is not $\O_X$-linear.

Let us summarize here without proofs some well-known 
properties\index{properties of jet bundles} of jet bundles.
\begin{Theorem}[\cite{GKZ2,KS}]\label{JetProperties}
$\quad$
\begin{enumerate}
\item[{\rm(a)}]
For any two line bundles $\L$, $\M$ on $X$, there exists
a canonical isomorphism $J(\L\otimes\M)\simeq J(\L)\otimes\M$.
\item[{\rm(b)}]
Let $\Omega^1_X$ be the sheaf of regular differential $1$-forms on $X$.
Then there exists a canonical exact sequence of vector bundles
\begin{equation}
0\to \Omega^1_X\otimes\L\to J(\L)\to \L\to0.\label{ExactSeqForJets}
\end{equation}
\item[{\rm(c)}] Let $X=\P(V)$ be a projective space. Then 
the jet bundle $J(\O_X(1))$ is a trivial 
vector bundle naturally identified with $\O_X\otimes V^*$
\end{enumerate}
\end{Theorem}

The relevance of jets to dual varieties is as follows.
Suppose that $X\subset\P(V)$ is an irreducible projective variety.
We take $\L=\O_X(1)$. Then any $f\in V^*$ is a linear function on $V$
and, hence, can be regarded as a global section of $\L$.
The following result is an immediate consequence of definitions:

\begin{Proposition}\label{DualJetVanishing}
$f\in V^*$ represents a point in the dual variety $\dual X$
if and only if the section $j(f)$ of a jet bundle $J(\L)$
vanishes at some point $x\in X$.
\end{Proposition}

\index{Koszul complex}
Let $E$ be an algebraic vector bundle of rank $r$ on an irreducible
algebraic variety $X$. For any global section $s$ of $E$
consider the following complexes of sheaves on $X$, called Koszul complexes.
$${\cal K}_+(E,s)=\left\{0\to\O_X\ra^s E\ra^{\wedge s}\Lambda^2E
\ra^{\wedge s}\ldots\ra^{\wedge s}\Lambda^rE\ra0\right\},$$
$${\cal K}_-(E,s)=\left\{0\to\Lambda^rE^*\ra^{i_s}\ldots
\ra^{i_s} \Lambda^2E^*\ra^{i_s}E^*\ra^{i_s}\O_X\ra0\right\}.$$
The differential in ${\cal K}_+$ is given by 
the exterior multiplication with $s$
and the differential in ${\cal K}_-$ is given by the contraction with $s$,
that is, by the map $\Lambda^jE^*\to\Lambda^{j-1}E^*$ dual to the map
$\Lambda^{j-1}E\to \Lambda^jE$ given by $\wedge s$.
We fix gradings in ${\cal K}_+$ and ${\cal K}_-$ by assigning
the degree $j$ to $\Lambda^jE$ in ${\cal K}_+$ and by assigning
the degree $-j$ to $\Lambda^jE^*$ in ${\cal K}_-$.
Notice that we have an isomorphism of complexes
$${\cal K}_-(E,s)\simeq{\cal K}_+(E,s)\otimes\Lambda^rE^*[r],$$
where $r$ in brackets means the shift in the grading by $r$.
The following Theorem shows that the cohomology of Koszul
complexes `represents' the vanishing set of $s$.
\index{cohomology of Koszul complexes}

\begin{Theorem}\label{KoszulGeneral}
$\quad$
\begin{enumerate}
\item[{\rm(a)}] Koszul complexes ${\cal K}_+(E,s)$ or
${\cal K}_-(E,s)$ are exact if and only if $s$ vanishes nowhere on $X$
{\rm(}the set of zeros $Z(s)$ is empty\/{\rm)}.
\item[{\rm(b)}] Suppose that $X$ is smooth and that $s$ vanishes along a 
smooth subvariety $Z(s)$ of codimension $r=\rank E$, and, moreover,
$s$ is transverse to the zero section.
Then ${\cal K}_-(E,s)$ has only one non-trivial cohomology
sheaf, namely $\O_{Z(s)}$ {\rm(}regarded as a sheaf on $X${\rm)}
in highest degree $0$. The complex ${\cal K}_+(E,s)$
in this situation has the only non-trivial cohomology
sheaf in the highest degree $r$, and this sheaf is the restriction
$\det E=\Lambda^rE|_{Z(s)}$ regarded as a sheaf on $X$.
\end{enumerate}
\end{Theorem}

\sketch
We shall prove only (a), the proof of (b) follows the same ideas
but is more technically involved, see e.g.~\cite{Fu1,GH1}.
If $V$ is a vector space and $v\in V$ is a non-zero vector
then the differential in the exterior algebra $\Lambda^*V$
given by $\wedge v$ is exact. Indeed, we can include $v$
in some basis as a first vector, then both the kernel and the image
of the differential are spanned by monomials containing~$v$.
Therefore, the dual differential $i_v$ in $\Lambda^*V^*$
is also exact. Applying this to our situation we see
that if the section $s$ does not vanish anywhere
then both Koszul complexes are exact fiberwise and, therefore,
are exact as complexes of sheaves.
However, if $Z(s)$ is not empty then the cokernel of the 
last differential in Koszul complexes is not trivial.
Indeed, if $x\in Z(s)$ then, trivializing $E$ near $x$,
we represent $s$ as a collection of functions $(f_1,\ldots,f_r)\in\O^r$.
Then the last differential in both complexes has a form $\O^r\to\O$,
$(u_1,\ldots,u_r)\mapsto\sum u_if_i$. If $s$ vanishes at $x$
then all $f_i=0$ and this map is not surjective.
\endproof

Let now $X\subset\P(V)$ be a smooth projective variety.
Then we can apply Theorem~\ref{KoszulGeneral}
to the jet bundle $E=J(\L)$, where $\L=\O_X(1)$.
Combining Proposition~\ref{DualJetVanishing} and
Theorem~\ref{KoszulGeneral}
we get the following

\begin{Theorem}[\cite{GKZ2}]\label{JetKoszul}
For any $f\in V^*$ let $j=j(f)$ denote its first jet. Then
a vector $f\in V^*$ represents a point in $\dual X$
if and only if any of the following complexes of sheaves
on $X$ is not exact:
$${\cal K}_+(J(\L),j)=\left\{0\to\O_X\ra^{j} J(\L)\ra^{\wedge {j}}\Lambda^2J(\L)
\ra^{\wedge {j}}\ldots\ra^{\wedge {j}}\Lambda^rJ(\L)\ra0\right\},$$
$${\cal K}_-(J(\L),{j})=\left\{0\to\Lambda^rJ(\L)^*\ra^{i_{j}}\ldots
\ra^{i_{j}} \Lambda^2J(\L)^*\ra^{i_{j}}J(\L)^*\ra^{i_{j}}\O_X\ra0\right\}.$$
\end{Theorem}

\section{Cayley Determinants of Exact Complexes}\index{Cayley determinant}

The determinants of exact complexes (in the implicit form)
were first introduced by Cayley in his paper \cite{Ca} on resultants.
A systematic early treatment of this subject 
was undertaken by Fisher \cite{Fi}
whose aim was to give a rigorous proof of Cayley results.
In topology determinants of complexes were introduced
in 1935 by Reidermeister and Franz \cite{Fra}.\index{Reidermeister torsion}
They used the word `torsion' for determinant-type invariants constructed.
In this section we give only definitions that are necessary
for the formulation of results related to dual varieties.
More details could be found in \cite{GKZ2,KnMu,Del,Q,RS}.

The base field $k$ can be arbitrary.
Suppose that $V$ is a finite-dimensional vector space.
Then the top-degree component of the exterior algebra 
$\Lambda^{\dim V}V$ is called
the determinant of $V$, denoted by 
$\Det V$. \index{determinant of a vector space}
If $V=0$ then we set $\Det V=k$.
It is easy to see that for any exact triple
$0\to U\to V\to W\to 0$
we have a natural isomorphism $\Det V\simeq\Det U\otimes\Det W$.

Suppose now that $V=V_0\oplus V_1$ is a finite-dimensional supervector 
space.\index{supervector space}
Then, by definition, $\Det V$ is set to be $\Det V_0\otimes(\Det V_1)^*$.
Once again, for any exact triple of supervector spaces
$0\to U\to V\to W\to 0$
we have $\Det V\simeq\Det U\otimes
\Det W$.\index{determinant of a supervector space}
For any supervector space $V$ we denote by~$\tilde V$ the new supervector
space given by $\tilde V_0=V_1$, $\tilde V_1=V_0$.
Clearly, we have a natural isomorphism 
\begin{equation}
\Det\tilde V=(\Det V)^*.\label{PairityCayleyIsomorphism}
\end{equation}

Now let $(V,\partial)$ be a finite-dimensional supervector space
with a differential~$\partial$ such that $\partial V_0\subset V_1$,
$\partial V_1\subset V_0$, $\partial^2=0$. Then $\Ker\partial$, $\Im\partial$, and
the cohomology space
$H(V)=\Ker\partial/\Im\partial$ are again supervector spaces.
We claim that there exists a natural isomorphism 
\begin{equation}
\Det V\simeq\Det H(V).\label{CayleyIsomorphism}
\end{equation}
Indeed, from the exact sequence
$$0\to\Ker\partial\to V\ra^\partial \widetilde{\Im\partial}\to0$$
we see that
$\Det V\simeq \Det(\Ker\partial)\otimes\Det(\Im\partial)^*$.
From the exact sequence
$$0\to\Im\partial\to\Ker\partial\to H(V)\to0$$
we get that $\Det H(V)\simeq\Det(\Ker\partial)\otimes\Det(\Im\partial)^*$.
Therefore, we have (\ref{CayleyIsomorphism}).
In particular, if $\partial$ is exact, $H(V)=0$, then we have
a natural isomorphism
\begin{equation}
\Det V\simeq k.\label{ExactCayleyIsomorphism}
\end{equation}

Let us fix some bases $\{e_1,\ldots,e_{\dim V_0}\}$ in $V_0$
and $\{e_1',\ldots,e_{\dim V_1}'\}$ in $V_1$.
Let $\{f_1',\ldots,f_{\dim V_1}'\}$ be a dual basis in $V_1^*$.
Then we have a basis vector 
$$e_1\wedge\ldots\wedge e_{\dim V_0}\otimes
f_1'\wedge\ldots\wedge f_{\dim V_1}'\in\Det V.$$
Therefore, if $(V,\partial,e)$ is a based supervector space
with an exact differential then by~(\ref{ExactCayleyIsomorphism})
we get a number $\det(V,\partial,e)$ called
the Cayley determinant of a based supervector space with an exact differential.
If we fix other bases
$\{\tilde e_1,\ldots,\tilde e_{\dim V_0}\}$ in $V_0$
and $\{\tilde e_1',\ldots,\tilde e_{\dim V_1}'\}$ in $V_1$
then, clearly, 
\begin{equation}
\det(V,\partial,\tilde e)=\det A_0(\det A_1)^{-1}\det(V,\partial,e),
\label{CayleyHomogeneity}
\end{equation}
where $(A_0,A_1)\in\GL(V_0)\times\GL(V_1)$ are transition matrices
from bases $e$ to bases $\tilde e$.
One upshot of this is the fact that if bases $e$ and $\tilde e$
are equivalent over some subfield $k_0\subset k$
then the Cayley determinants with respect to these bases
are equal up to a non-zero multiple from $k_0$.
Since all isomorphisms are natural so far,
the homogeneity condition (\ref{CayleyHomogeneity})
can be reformulated as follows.
For any based supervector space
$(V,\partial,e)$ with an exact differential
and any 
$(A_0,A_1)\in\GL(V_0)\times\GL(V_1)$
we can consider a new based supervector space
$(V,A\cdot\partial,e)$ with the exact differential
$A\cdot\partial$ given by 
$$A\cdot\partial|_{V_0}=A_1\cdot \partial\cdot A_0^{-1},\quad
A\cdot\partial|_{V_1}=A_0\cdot \partial\cdot A_1^{-1}.$$
Then we have
$$\det(V,A\cdot\partial,e)=(\det A_0)^{-1}\det A_1\det(V,\partial,e).$$
In particular, it is easy to see that the following formula is valid:
\begin{equation}
\det(V,\lambda\partial, e)=\lambda^{\dim(\Ker\partial)_1-
\dim(\Ker\partial)_0}\det(V,\partial,e).\label{HomogeneityDifferentialCayley}
\end{equation}

In the matrix form the Cayley determinant can be calculated as follows.
Suppose that $(V,\partial,e)$ 
is a based supervector space with an exact differential,
$$\dim V_0=\dim\Ker\partial_0+\dim\Im\partial_0=
\dim\Im\partial_1+\dim\Ker\partial_1=\dim V_1=n.$$
Therefore, both $\partial_0$ and $\partial_1$ are represented
by $(n\times n)$-matrices $D_0$ and $D_1$.
For any subsets $I,J\subset B=\{1,\ldots,n\}$
we denote by $D_0^{IJ}$ and $D_1^{IJ}$ the submatrices
of $D_0$ and $D_1$ formed by columns indexed by $I$ and rows
indexed by $J$.
Then it is easy to see that there exist subsets $I_0,I_1\subset B$
such that submatrices $D_0^{B\setminus I_0,I_1}$, $D_1^{B\setminus I_1, I_0}$
are invertible. In particular, 
$$\# I_1=\rk\partial_0,\quad \# I_0=\rk\partial_1=n-\# I_1.$$
The formula (\ref{CayleyHomogeneity}) implies that
\begin{equation}
\det(V,\partial,e)=\det D_0^{B\setminus I_0,I_1}(\det D_1^{B\setminus I_1, I_0})^{-1}.
\label{JKHGKJGHFKHGFKLJ}
\end{equation}
This formula gives an explicit matrix description of the Cayley determinant.

\smallskip

Now we consider a finite complex of finite-dimensional vector spaces
$$\ldots\ra^{\partial_{i-1}} V^i\ra^{\partial_i} 
V^{i+1}\ra^{\partial_{i+1}}\ldots.$$
Then we can define the finite-dimensional supervector
space $V=V_0\oplus V_1$,
$$V_0=\mathop{\oplus}\limits_{i\equiv0\mod2}V^i,\quad
V_1=\mathop{\oplus}\limits_{i\equiv1\mod2}V^i,$$ 
with an induced differential $\partial$.
In particular, all previous considerations are valid.
Therefore, if the complex $(V^\bullet,\partial)$ is exact
and there are some fixed bases $\{e^i_1,\ldots,e^i_{\dim V^i}\}$ 
in each component $V^i$
then we have the corresponding Cayley determinant $\det(V^\bullet,\partial,e)\in k^*$.
For example, if $L$ and $M$ are based vector spaces
and $A:\,L\to M$ is an invertible operator then
the complex $0\to L\ra\limits^A M\to 0$ is exact
and the corresponding Cayley determinant is equal to $\det A$
(if $L$ is located in the even degree of the complex).

In fact, Cayley in \cite{Ca} has found some matrix representation of his
determinant.\index{determinant of a complex}
In order to give his formula we shall write our complex
in the explicit coordinate form
\begin{equation}
V^\bullet=\{0\to k^{B_0}\ra^{D_0} k^{B_1}\ra^{D_1}\to\ldots\ra^{D_{r-1}}k^{B_r}\to0\},
\label{GHGTYTYF}
\end{equation}
where for simplicity we assume that the non-zero terms of the complex
are located in degrees between $0$ and $r>0$.
In general, if $V^\bullet[m]$ is the same complex as $V^\bullet$
but with a grading shifted by $m$,
then by~(\ref{PairityCayleyIsomorphism}) we have 
$\det(V^\bullet[m],\partial,e)=(\det(V^\bullet[m],\partial,e))^{(-1)^m}$
In the formula (\ref{GHGTYTYF}) we may suppose
that all $B_i$ are some finite sets.
Therefore, $D_i$ is a matrix with columns indexed by $B_i$
and rows indexed by $B_{i+1}$.
For any subsets $X\subset B_i$, $Y\subset B_{i+1}$
we denote by $(D_i)_{XY}$ the submatrix in $D_i$
with columns from $X$ and rows from $Y$.

A collection of subsets $I_i\subset B_i$ is called admissible\index{admissible collection}
if $I_0=\emptyset$, $I_r=B_r$, for any $i=0,\ldots,r-1$ we have
$\#(B_i\setminus I_i)=\#I_{i+1}$, and the submatrix
$(D_i)_{B_i\setminus I_i,I_{i+1}}$ is invertible.
The following theorem easily follows from (\ref{JKHGKJGHFKHGFKLJ}).

\begin{Theorem}[\cite{GKZ2}]$\ $
\begin{enumerate}
\item[{\rm(a)}] Admissible collections exist. For any admissible collection
we have 
$$\# I_i=\sum_{j=0}^{i-1}(-1)^{i-1-j}\# B_j.$$
\item[{\rm(b)}] Let $\{I_i\}$ be an admissible collection.
Denote by $\Delta_i$ the determinant of the matrix
$(D_i)_{B_i\setminus I_i,I_{i+1}}$. Then
$$\det(V^\bullet,D,B)=\prod_{i=0}^{r-1}\Delta_i^{(-1)^{i}}.$$
\end{enumerate}
\end{Theorem}

In particular, this Theorem (or formula~(\ref{HomogeneityDifferentialCayley}))
implies 

\begin{Corollary}\label{CayleyDetDegree}
For any based exact complex $(V^\bullet,d,e)$ we have\index{based exact complex}
$$\det(V^\bullet,\lambda d,e)=\lambda^{\sum_i(-1)^{i+1}\cdot i\cdot\dim V^i}\det(V^\bullet,d,e).$$
\end{Corollary}

\section{Discriminant Complexes}\index{discriminant complex}
Let $X\subset\P(V)$ be a smooth projective variety.
We denote $\L=\O_X(1)$, so any $f\in V^*$ can be regarded
as a section of $\L$. Let $\M$ be another line bundle on $X$.
We define discriminant complexes $C^*_+(X,\M)$ and $C^*_-(X,\M)$
as complexes of global sections of Koszul complexes
${\cal K}_+(J(\L),j(f))$ and ${\cal K}_-(J(\L),j(f))$
tensored by $\M$. More precisely,
$$C^i_+(X,\M)=H^0(X,\Lambda^iJ(\L)\otimes\M),$$
$$C^i_-(X,\M)=H^0(X,\Lambda^{-i}J(\L)^*\otimes\M).$$
Thus the terms of discriminant complexes are fixed
and the differentials depend on $f\in V^*$.
We shall denote this differential by $\partial_f$.
By Theorem~\ref{JetKoszul} a 
vector $f\in V^*$ represents a point in $\dual X$
if and only if any of two Koszul complexes is not exact.
We want to get the same condition but for discriminant
complexes, which are complexes of finite-dimensional vector spaces
instead of complexes of sheaves.
However, the exactness of 
a complex of sheaves does not necessarily imply
the exactness of the corresponding complex of global sections.
The obstruction to this is given by the higher cohomology
of the sheaves of the complex.

\begin{Definition}
The discriminant complexes
are called stably twisted if all terms of the corresponding Koszul
complexes have no higher cohomology.\index{stably twisted discriminant complex}
\end{Definition}

For example, suppose that $\M$ is any ample line bundle on $X$.
Ampleness is equivalent to the fact that for any coherent sheaf $\cal F$
on $X$, the sheaves ${\cal F}\otimes\M^{\otimes n}$ have no higher
cohomology for $n\gg0$. Therefore, for sufficiently big~$n$
the discriminant complexes 
$C^*_+(X,\M^{\otimes n})$ and $C^*_-(X,\M^{\otimes n})$
will be stably twisted.
In this situation we have the following Theorem
that immediately follows from Theorem~\ref{JetKoszul}
and `abstract de Rham theorem' (see e.g.\cite{GH1}).

\begin{Theorem}[\cite{GKZ2}]
Suppose that the discriminant complex $C^*_+(X,\M)$ {\rm(}resp.\ $C^*_-(X,\M)${\rm)}
is stably twisted.
Let $f\in V^*$ be such that its projectivization does not belong
to the dual variety $\dual X\subset\P(V^*)$.
Then the complex $(C^*_+(X,\M),\partial_f)$
{\rm(}resp.~$(C^*_-(X,\M),\partial_f)${\rm)} is exact.
\end{Theorem}

Remarkably, it turns out that much more is true.
Suppose that the discriminant complex $C^*_+(X,\M)$ (or~$C^*_-(X,\M)$)
is stably twisted.
Then for a generic $f\in V^*$ the complex 
$C^*_+(X,\M)$ (or~$C^*_-(X,\M)$) is exact by the previous Theorem.
Let $f_1,\ldots,f_n$ be a basis of $V^*$.
Let $\C(x_1,\ldots,x_n)=\C(V^*)$ be a field
of rational functions on $V^*$.
Consider the complex $C^*_+(X,\M)\otimes\C(V^*)$ (or~$C^*_-(X,\M)\otimes\C(V^*)$)
with the differential given by $\partial=\sum x_i\partial_{f_i}$.
Then this complex is exact, since its generic specialization is exact.
We fix some bases over $\C$ in each component 
$C^i_+(X,\M)$ (or~$C^{-i}_-(X,\M)$)
and consider them as bases over $\C(V^*)$ in each
$C^i_+(X,\M)\otimes\C(V^*)$ (or~$C^{-i}_-(X,\M)\otimes\C(V^*)$).
Then we can calculate the Cayley determinant $\Delta^+_{X,\M}$
(or~$\Delta^-_{X,\M}$) of the complex
$C^*_+(X,\M)\otimes\C(V^*)$ 
(or~$C^*_-(X,\M)\otimes\C(V^*)$), which will be
the non-zero element of the field of rational functions $\C(V^*)$.
Clearly the Cayley determinants depend only on the discriminant
complexes (up to a scalar multiple) and do not depend on the choice
of bases. The following Theorem shows that
the Cayley determinant coincides with the discriminant.

\begin{Theorem}[\cite{GKZ2}]\label{DiscriminantCayley}
If the discriminant complex $C^*_-(X,\M)$ is stably twisted
then, up to a non-zero scalar, we have
$$\Delta^-_{X,\M}=\Delta_X,$$
where $\Delta_X$ is the discriminant of $X$.
If $C^*_+(X,\M)$ is stably twisted
then, up to a non-zero scalar, we have
$$(\Delta^+_{X,\M})^{(-1)^{\dim X+1}}=\Delta_X.$$
\end{Theorem}

The proof of this Theorem involves derived categories and Cayley determinants
for complexes over arbitrary Noetherian integral domains, see~\cite{GKZ2}.
Let us consider several examples.

\subsubsection*{Sylvester Formula.} \index{classical discriminant}
Let $X=\P^1=\P(\C^2)$ be a projective
line embedded into $\P(S^d\C^2)$ via the Veronese embedding.
The space $V^*=(S^d\C^2)^*$ is the space of binary forms
$$f(x_0,x_1)=a_0x_0^d+a_1x_0^{d-1}x_1+\ldots+a_dx_1^d$$
and the discriminant $\Delta_X(f)$ is the classical discriminant
of this binary form.\index{Sylvester formula}
It is quite easy to implement Theorem~\ref{DiscriminantCayley}
in this case. 
We take a twisting bundle $\M=\O_{\P^1}(2d-3)$.
The discriminant complex $C^*_-(X,\M)$ has only
two non-zero terms and, therefore, the Cayley determinant
is reduced to a determinant of a square matrix.
More precisely, the Cayley determinant in this case
is equal to the determinant of the linear map
$$\partial_f:\,S^{d-2}\C^2\oplus S^{d-2}\C^2\to S^{2d-3}\C^2,$$
given by 
$$\partial_f(u,v)={\partial f\over\partial x_0}u+{\partial f\over\partial x_1}v.$$
The final formula coincides with the classical Sylvester
formula for the discriminant of a binary form (up to a scalar):
$$\Delta_X(f)={(-1)^{d-1}\over d^{d-2}}D,$$
where $D$ is the determinant of the following matrix
$$
\left|\matrix{
a_1   & 2a_2   & \ldots& (d-1)a_{d-1} & da_d & \ldots & 0\cr
0     & a_1    & \ldots& (d-2)a_{d-2} & (d-1)a_{d-1}& \ldots &0\cr
\vdots&\vdots  &\ddots&\vdots&\vdots&\ddots&\vdots\cr
0     &0       &\ldots&a_1&2a_2&\ldots&da_d\cr
da_0  &(d-1)a_1&\ldots&2a_{d-2}&a_{d-1}&\ldots&0\cr
0     &    da_0&\ldots&3a_{d-3}&2a_{d-2}&\ldots&0\cr
\vdots&\vdots  &\ddots&\vdots&\vdots&\ddots&\vdots\cr
0     &0       &\ldots&da_0&(d-1)a_1&\ldots&a_{d-1}\cr
}\right|.$$

\subsubsection*{The Dual of an Algebraic Curve.}\index{dual variety of a curve}
Let $X\subset\P(V)$ be a smooth algebraic curve of degree $d$
and genus $g$. 
We take the twisting bundle $\M$ to be a generic
line bundle on $X$ of degree $2d+3g-3$.
The discriminant complex $C^*_-(X,\M)$ has only
two non-zero terms and, therefore, the Cayley determinant
is given by the determinant of a square matrix.
The size of this matrix (and hence the degree of the dual 
variety $\dual X$) is equal to $\dim H^0(X,\M)$.
The latter is equal to $2d+2g-2$ by  the Riemann-Roch theorem.

\subsubsection*{Discriminant Spectral Sequences.}
One disadvantage of Theorem~\ref{DiscriminantCayley} is the restricting
condition on the twisting line bundle $\M$. It is not always
convenient to make the discriminant complexes $C^\bullet_\pm$ stably twisted.
Another possibility is to take into account the higher cohomology as well.
The standard tool for this is the spectral sequence.
Recall \cite{GH1} that if $\F^\bullet$ 
is a finite complex of sheaves on a topological
space $X$, then one can define the hypercohomology groups ${\bf H}^i(X,\F^\bullet)$.
To define them, consider the complexes of Abelian groups $C^\bullet_j$ calculating
the cohomology of every individual $\F^j$
(for example, the \u Cech complexes with respect to an appropriate open
covering of $X$). The differentials $\F^i\to\F^{i+1}$ make this collection
of complexes into a double complex $C^{\bullet\bullet}$.
The hypercohomology groups ${\bf H}^i(X,\F^\bullet)$ are the cohomology
groups of its total complex.
In particular, we have the spectral sequence of the double complex
$C^{\bullet\bullet}$ \index{spectral sequence of the double complex}
$$E_1^{pq}=H^q(X,\F^p)\Rightarrow {\bf H}^{p+q}(X,\F^\bullet).$$
The first differential $d_1$ is induced by the differential in $\F^\bullet$.
In particular, the complex of global sections
$$\ldots\to H^0(X,\F^i)\to H^0(X,\F^{i+1})\to\ldots$$
is just the bottom row of the term $E_1$.

If the complex of sheaves $\F^\bullet$ is exact then
all hypercohomology groups ${\bf H}^i(X,\F^\bullet)$ vanish.
Indeed, the hypercohomology can be calculated using
another spectral sequence `dual' to the first:
$${}'E_2^{pq}=H^q(X,{\cal H}^p(\F^\bullet))\Rightarrow{\bf H}^{p+q}(X,\F^\bullet),$$
where ${\cal H}^p(\F^\bullet)$ is the $p$-th cohomology sheaf of $\F^\bullet$
$${\cal H}^p(\F^\bullet)={\Ker\{\F^p\ra^d\F^{p+1}\}\over
\Im\{\F^{p-1}\ra^d\F^p\}}.$$

Now let $X\subset \P(V)$ be a smooth projective variety, $\L=\O_X(1)$.
Let $\M$ be any line bundle on $X$.
We define discriminant spectral sequences $C_{r,\pm}^{pq}(X,\M,f)$
to be the spectral sequences of complexes
${\cal K}_\pm(J(\L),j(f))$ tensored by $\M$.
Then the following theorem is a consequence of the discussion above and 
Theorem~\ref{JetKoszul}.\index{discriminant spectral sequence}

\begin{Theorem}[\cite{GKZ2}]
Suppose that the projectivization of $f\in V^*$
does not belong to the dual variety $\dual X$.
Then the discriminant spectral sequences 
$C_{r,\pm}^{pq}(X,\M,f)$ are exact, i.e.~they converge to zero.
\end{Theorem}

Moreover, it is possible to define the Cayley determinants
of exact spectral sequences and to prove an analogue of 
Theorem~\ref{DiscriminantCayley}, see \cite{GKZ2} for details.

\subsubsection*{Weyman's Bigraded Complexes}
The calculation of the determinant of a spectral sequence
with many non-trivial  terms may be very involved.
J.~Weyman \cite{We} has suggested a procedure that replaces
the discriminant spectral sequence by a bigraded
complex which incorporates all higher differentials at once.
Bigraded complex is, by definition, a complex $(C^\bullet,\partial)$
of vector spaces in which each term $C^i$ is equipped
with additional grading $C^i=\oplus_{p+q=i}C^{pq}$.\index{bigraded complex}
The underlying complex $(C^\bullet,\partial)$ may then be called
a total complex of the bigraded complex $C^{\bullet\bullet}$.
The differential $\partial$ in $C^\bullet$ is decomposed into a sum
$\partial=\sum_r\partial_r$, where $\partial_r$ is bihomogeneous
of degree $(r,1-r)$. We shall consider only bigraded complexes
with $\partial_r=0$ for $r<0$.

Now let $X\subset \P(V)$ be a smooth projective variety, $\L=\O_X(1)$.
Let $\M$ be any line bundle on $X$.\index{Weyman's bigraded complex}

\begin{Theorem}[\cite{We}]
There exists a bigraded complex $(C^{\bullet\bullet},\partial_f)$ with
$$C^{pq}=H^q\left(X,\Lambda^{-p}J(\L)^*\otimes\M\right)$$
and the differential $\partial_f$ depends polynomially on $f\in V^*$.
This complex has the following properties:
\begin{itemize}
\item The differential $\partial_f$ has the form $\partial_f=\sum_{r\ge1}\partial_{r,f}$,
where $\partial_{r,f}$ has bidegree $(r,1-r)$ and its matrix elements
are homogeneous polynomials of degree $r$ in coefficients of $f$.
Moreover, $\partial_{1,f}$ is induced by the differential in the Koszul complex
${\cal K}_-(J(\L),{j(f)})$.
\item The Cayley determinant of the total complex $(C^\bullet,\partial_f)$
of $C^{\bullet\bullet}$ equals, up to a scalar factor, the discriminant $\Delta_X$.
\end{itemize}
\end{Theorem}

\chapter{Resultants and Schemes of Zeros}

\section*{Preliminaries}
Classically, discriminants were studied together
with resultants. Resultants and their generalizations
appear while looking at vector bundles: they detect
their global sections that have ``correct'' schemes of zeros. Some results in this direction
are reviewed in this chapter.

\section{Ample Vector Bundles}
\index{ample vector bundle}%

Classically, a line bundle $L$ over an algebraic variety $X$
is said to be ample if the map $X\to\P^N$ associated
to the sections of some power $L^n$ gives an imbedding of $X$.
This notion has been extended to vector bundles as well 
by Hartshorne \cite{Ha3}.
Let $E$ be a vector bundle on $X$. We consider the variety
$X_E=\P(E^*)$, the projectivization of the bundle $E^*$.
There is a projection $p:\,E_X\to X$ whose fibers are projectivizations
of fibers of $E^*$, and a natural projection 
$\pi:\,E^*\setminus X\to X_E$, where $X$ is embedded into the total space of $E^*$
as the zero section.
We denote by $\xi(E)$ the tautological line bundle on $X_E$
defined as follows. For open $U\subset X_E$, a section of $\xi(E)$
over $U$ is a regular function on $\pi^{-1}(U)$ which is homogeneous
of degree $1$ with respect to dilations of $E^*$.
The restriction of $\xi(E)$ to every fiber $p^{-1}(X)=\P(E_x^*)$
is the tautological line bundle $\O(1)$ of the projective space
$\P(E_x^*)$.

\begin{Definition}\rm
A vector bundle $E$ on $X$ is called ample if $\xi(E)$
is an ample line bundle on $X_E$.
\end{Definition}

If $E$ is itself a line bundle then $X_E=X$, $\xi(E)=E$, and this 
definition gives nothing new.
Basic properties and criteria for ampleness were obtained in \cite{Ha3}.
For example, the Grothendieck definition
of ample line bundles can be transported to vector bundles as well.
Namely, $E$ is ample on $X$ if and only if for every coherent sheaf $F$,
there is an integer $n_0>0$, such that for every $n\ge n_0$,
the sheaf $F\otimes S^n(E)$
(where $S^n(E)$ is the $n$-th symmetric power of $E$)
is generated as an $\O_X$-module by its global sections.
Using this description it is easy to see, for example, that
the direct sum of ample vector bundles is ample, etc.

The following theorem is an easy application of ample vector bundles properties.

\index{dual variety of a complete intersection}%
\begin{Theorem}[\cite{E1}]\label{DualCompleteIntersection}
Assume that $X$ is a nonlinear smooth projective variety in $\P^N$.
If $X$ is a complete intersection then $\dual X$ 
is a hypersurface.
\end{Theorem}

\proof
We may assume that $X$ is nondegenerate. Let $\dim X=n$.
Then $X$ is scheme-theoretic intersection of $N-n$ hypersurfaces
$H_1,\ldots,H_{N-n}$ such that $\deg H_k=d_k\ge2$.
Therefore the normal bundle $N_X\P^N$ is equal to
$\O_X(d_1)\oplus\ldots\oplus\O_X(d_{N-n})$.
It follows that
$$N_X\P^N(-1)=\O_X(d_1-1)\oplus\ldots\oplus\O_X(d_{N-n}-1)$$
is an ample vector bundle.
Let $I_X=\P(N_X^*\P^N(1))\subset\P^N\times(\P^N)^*$
be the conormal variety and let $\pr:\,I_X\to\dual X$
be the projection. Then $\pr^*\O_{\dual X}(1)$ is the tautological
line bundle on $\P(N_X^*\P^N(1))$. Since $N_X\P^N(-1)$
is ample, it follows that $\pr$ is a finite morphism.
Therefore, 
$$\dim\dual X=\dim I_X=N-1.$$
\endproof

Sommese \cite{S1} generalized the notion of ampleness to the following.

\begin{Definition}
A vector bundle $E$ is called $k$-ample if 
$\xi(E)^n$ is spanned by global
sections for some $n>0$ and the induced map $X_E\to\P^N$
has at most $k$-dimensional fibers\index{ampleness}.
\end{Definition}

It can be verified that the number $k$ is independent
of the integer $n$ for which
$\xi(E)^n$ is spanned and a vector bundle is $0$-ample
if and only if it is ample.
This number $k$ is called the ampleness of $E$, denoted by $a(E)$.
It can be shown that $0\le a(E)\le\dim X$.
If $\xi(E)^n$ is not spanned for any $n>0$, we define $a(E)=\dim X$.

The ampleness of homogeneous vector bundles on flag varieties
was described combinatorially by Snow \cite{Sn1} inspired
by an earlier work of Goldstein~\cite{Go}
who determined the ampleness
for tangent bundles on flag varieties.

\section{Resultants}
Classically, discriminants were studied in conjunction
with resultants.
Let $X$ be a smooth irreducible projective variety
and let $E$ be a vector bundle on $X$ of rank $k=\dim X+1$.
Set $V=H^0(X,E)$. We shall assume that $E$ is very ample\index{very ample vector bundle}, 
i.e.~$\xi(E)$ is a very ample line bundle on $X_E$.
In particular, $\xi(E)$ (and hence $E$) is generated by global sections.
Notice that $H^0(X,E)=H^0(X_E,\xi(E))$, therefore
$\xi(E)$ embeds $X_E$ in $\P(V^*)$.

\begin{Definition}\rm
The resultant variety\index{resultant variety} $\nabla\subset\P(V)$ is the set of all sections
vanishing at some point $x\in X$.
\end{Definition}

\begin{Example}\rm
Suppose that $X=\P^{k-1}=\P(\C^k)$ and 
$$E=\O(d_1)\oplus\ldots\oplus\O(d_k).$$
Then $V=S^{d_1}(\C^k)^*\oplus\ldots\oplus S^{d_k}(\C^k)^*$
and $\nabla$ is the classical resultant\index{classical resultant} variety
parametrizing $k$-tuples of homogeneous forms on $\C^k$
of degrees $d_1,\ldots,d_k$ having a common non-zero root.
\end{Example}

The assertion (b) of the following Theorem is sometimes called
the Cayley trick\index{Cayley trick} who first noticed that the resultant
can be written as a discriminant.

\begin{Theorem}[\cite{GKZ2}]\label{JHGKJFKJFGK}$\ $
\begin{enumerate}
\item[\rm(a)] $\nabla$ is an irreducible hypersurface of degree $\int_Xc_{k-1}(E)$.
\item[\rm(b)] $\nabla$ is projectively dual to $X_E$.
\end{enumerate}
\end{Theorem}

\proof
$\xi(E)$ embeds $X_E$ in $\P(V^*)$ in such a way that all fibers
become projective subspaces of dimension $k-1$.
Geometrically, $\nabla$ parametrizes hyperplanes in $\P(V^*)$
that contain at least one fiber of $X_E$.
Since the embedded tangent space to any point of $X_E$
contains the fiber through it, it follows that $\nabla\supset\dual{X_E}$.
On the contrary, suppose that $v\in\nabla$.
Then the corresponding section $s\in H^0(X,E)$ vanishes at some point $x\in X$.
Therefore, $s$ defines a linear map $ds:\,T_xX\to T_xE=T_xX\oplus E_x$.
Let $\pr:\,T_xE\to E_x$ be the projection.
Since $\dim X=\dim E-1$,  $\pr\circ ds(T_xX)$ is contained
in some hyperplane $H\subset E_x$.
An easy local calculation shows that the hyperplane in $\P(V^*)$
corresponding to $v$ contains the embedded tangent space to $X_E$
at the point corresponding to $H$.
Therefore, $\nabla=\dual{X_E}$.

Now let us prove that $\nabla$ is a hypersurface.
An easy dimension count shows that it is sufficient to prove that
a generic section $s$ of $H^0(X,E)$ vanishing at $x\in X$
does not vanish at other points. In other words,
a generic hyperplane in $P=\P(V^*)$ containing some fiber $F$ of $E_X$
does not contain other fibers.
Consider the linear projection 
$\pi_F:\,P\setminus F\to P/F$.
If $F'$ is any fiber of $X_E$ distinct from $F$ then
$\pi_F(F')$ is a projective subspace of dimension $\dim F'=k-1$.
Therefore, images of fibers of $X_E$ distinct from $F$
form at most $(k-1)$-dimensional family of projective subspaces
of dimension $k-1$.
An easy dimension count shows that there exists a hyperplane $H\subset P/F$
that does not contain any projective subspace from this family.
The preimage of this hyperplane in $P$ contains only one fiber $F$.

It remains to calculate the degree of $\nabla$.
To find it, we should take two generic sections $s_1$ and $s_2$
of $V$ and find the number of values $\lambda\in\C$ such that
$s_1+\lambda s_2$ vanishes at some $x\in X$.
In other words, we need to find the number of points (with multiplicities)
of a zero-dimensional cycle $Z\subset X$  consisting of all points $x\in X$
where $s_1$ and $s_2$ are linearly dependent.
But this is exactly the geometric definition of $c_{k-1}(E)$, see \cite{Fu1}.
\endproof

The same argument also proves the following result from the folklore.
Suppose that $X\subset\P^N$ is a scroll, i.e.~a projectivization
of a vector bundle $E$ over a smooth variety $Y$ such that all fibers
are embedded linearly. Assume that $\dim E>\dim Y$.
Then 
$$\defect X=\dim E-\dim Y-1$$ 
and $\dual X$ parametrizes
hyperplanes in $\P^N$ that contain at least one fiber of the 
projective bundle $X\to Y$.

The proof of the following theorem can be found in \cite{E2}:

\begin{Theorem}[\cite{E2}]
If $X$ is a smooth $n$-dimensional projective variety in $\P^N$
and $\defect X=k\ge n/2$, then $X$ is a projective bundle $X\simeq\P_Y(F)$,
where $F$ is a vector bundle of rank $(n+k+2)/2$ on a smooth $(n-k)/2$-dimensional
variety $Y$ and the fibers are embedded linearly.
\end{Theorem}

\section{Zeros of Generic Global Sections}

Both resultants and discriminants are related to the 
following general construction.
Suppose that $X$ is a smooth projective variety and $E$
is a vector bundle on $X$ generated by global sections.
Let $Z(s)$ denote
the scheme of zeros
of any global section $s\in H^0(X,E)$. 
One might expect that for generic~$s$
the scheme $Z(s)$ is a smooth variety of codimension $\dim E$.
Then we can define the degeneration variety $D\subset H^0(X,E)$
parametrizing all global sections~$s$ such that $Z(s)$
is not smooth of expected codimension.
For example, if $E$ is a very ample line bundle
then $D$ is a cone over the dual variety.
If $E$ is a very ample vector bundle and $\dim  E=\dim X+1$
then $D$ is the resultant variety.
If $E$ is a very ample vector bundle and $\dim  E=\dim X$
then the corresponding homogeneous polynomial
can be called Bezoutian.
Indeed, if $X=\P^n$ and $E=\O(d_1)\oplus\ldots\oplus\O(d_n)$
then $D$ parametrizes sets of homogeneous forms of degrees
$d_1,\ldots,d_n$ such that the Bezout theorem is not applicable.

In general, in order to make the theory consistent,
it is necessary to impose very strong conditions on the vector
bundle $E$. These conditions are not always satisfied
even in the case of homogeneous vector bundles
on flag varieties.
Even the question whether or not a generic zero scheme
is non-empty can be quite difficult.

Assume that $G$ is a connected reductive group, 
$T$  is a fixed maximal torus,
$B$ is a fixed Borel subgroup, $T\subset B\subset G$,
$B_-$ is the opposite Borel subgroup,
$P$ is a parabolic subgroup containing $B_-$,
$\X(T)$ is the lattice of characters of $T$,
and $\lambda\in\X(T)$  is the dominant weight.
Consider the homogeneous vector bundle $\L_\lambda=G\times_PU_\lambda$ over $G/P$, where
$U_\lambda$ is the irreducible $P$-module with highest weight $\lambda$.
See Section~\ref{DefinitionsNotationsFlags} for further details.
By the Borel--Weil--Bott theorem (see~\cite{Bot}), 
$V_\lambda=H^0(G/P, \L_\lambda)$ is an irreducible $G$-module with highest
weight~$\lambda$. 

\begin{Theorem}\label{JHGFJGFDJCHC}
Let $s\in V_\lambda$ be a generic global section. Then
\begin{enumerate}
\item[{\rm(a)}] If $\dim U_\lambda>\dim G/P$, the scheme of zeros $Z_s$ is empty.
\item[{\rm(b)}] If $\dim U_\lambda\le\dim G/P$, either $Z_s$ is empty
or $s$ intersects the zero section of $\L_\lambda$ transversally and $Z_s$ 
is a smooth unmixed subvariety of expected codimension
$\dim U_\lambda$.
\item[{\rm(c)}] 
If $\dim U_\lambda=\dim G/P$, the geometric number of points
in $Z_s$ is equal to the top Chern class of $\L_\lambda$.
\end{enumerate}
\end{Theorem}
\proof
(a) is obvious and follows by an
easy dimension count.
(b) and (c) follow from the fact that $\L_\lambda$
is generated by global sections. Let us recall this argument here.
Let $\dim U_\lambda\le\dim G/P$,
and assume that every global section has a zero.
We have to prove that a generic global section~$s$ 
intersects the zero section of $\L_\lambda$ 
transversally. This will imply, in particular, that 
$Z_s$ is a smooth unmixed subvariety of codimension $\dim U_\lambda$. 
For simplicity we suppress the index~$\lambda$. 
Consider the incidence variety 
$$Z\subset G/P\times V,\quad
Z=\{(x,s)\,|\,x\in (Z_s)_{red}\}.$$
Since $G/P$ is homogeneous and $Z$ is invariant,
it follows that $Z$  is obtained by spreading
the fiber $Z_e=\{s\in V\,|\,s(eP)=0\}$ by the group $G$.
Since $U$ is irreducible, we have $\dim Z_e=\dim V-\dim U$.
Hence, $Z$ is a smooth irreducible subvariety
of dimension $\dim V+\dim G/P-\dim U$.

Let $\pi:\,Z\to V$ be the restriction to $Z$ 
of the projection of $G/P\times V$ on the second summand.
By assumption, $\pi$ is a surjection.
By Sard's lemma for algebraic varieties (see~\cite{Mum1}),
for a generic point $s\in V$ and any point $(x,s)$ 
in $\pi^{-1}s$ the differential $d\pi_{(x,s)}$ is surjective. 
We claim that $s$ has the transversal intersection with the zero section.
Indeed, 
$Z$ can be regarded as a subbundle of the trivial bundle $G/P\times V$.
Then $\L=(G/P\times V)/Z$ and 
the zero section of $\L$ is identified with $Z\mod Z$.
The section $s$ is identified with $(G/P\times\{s\})\mod Z$.
Hence, it is sufficient to prove that $G/P\times\{s\}$ 
is transversal to $Z$,
which is equivalent to the following claim:
$d\pi_{(x,s)}$ is surjective for all $(x,s)\in Z$.
Now statement (c) of the theorem follows from the standard intersection 
theory (see~\cite{Fu1}).
\endproof

\begin{Example}\rm
It should be noted that Theorem~\ref{JHGFJGFDJCHC} 
cannot be strengthened to the point
where the non-emptiness and the irreducibility of the scheme of
zeros in Theorem~\ref{JHGFJGFDJCHC} could be established apriori, as the following
example shows.
Consider the vector bundle $S^2\SS^*$ on $\Gr(k, 2n)$.
The dimension of a fiber does not exceed the dimension of
the Grassmanian as $k\le \bigover{4n-1}{3}$, but a generic section
(that is, a non-degenerate quadratic form in $\C^{2n}$) 
has a zero (that is, a $k$-dimensional isotropic subspace)
only if $k\le n$. For $k=n$ the scheme of zeros is a reducible variety of dimension
$\bigover{n(n-1)}{2}$ with two irreducible components 
(spinor varieties) that
correspond to two families of maximal isotropic subspaces on an even-dimensional quadric.
\end{Example}

Consider the $P$-submodule
$M_\lambda=\{s\in H^0(G/P, \L_\lambda)\,|\,s(P)=0\}$
in $V_\lambda$.
Here $s(P)$ is a germ of $s$ in the point of $G/P$ corresponding to~$P$.
It is easy to see that $M_\lambda$ can be characterized
as the unique maximal proper $P$-submodule of~$V_\lambda$.
Clearly $V_\lambda/M_\lambda\simeq U_\lambda$ as $P$-modules.
Consider the map
$\Psi:\,G\times M_\lambda\to V_\lambda$, $\Psi(g, v)=gv$.
Then generic global sections of~$\L$ have zeros iff
$\Psi$ is dominant iff
the differential of~$\Psi$ at a generic point is surjective. 
We write $\g$ and $\p$ for Lie algebras of $G$ and $P$, respectively.
The natural ``orbital'' map $\g\times V_\lambda\to V_\lambda$ defines
the $P$-equivariant map
$\psi:\,\g/\p\times M_\lambda\to U_\lambda$.
Then $\Psi$ is dominant iff
$\psi(\cdot,x)$ is surjective at a generic point $x\in M_\lambda$.
Let $Z_\lambda=\psi^{-1}(0)_{red}\subset\g/\p\times M_\lambda$
be the incidence variety. 
Then the following proposition
follows from an easy dimension count:

\begin{Proposition}
Suppose that
$$\dim Z_\lambda=\dim\g/\p+\dim M_\lambda-\dim U_\lambda.$$
Then generic global sections of~$\L_\lambda$ have zeros.
\end{Proposition}

Denote by $\pi$ the restriction on $Z_\lambda$ of the projection of
$\g/\p\times M_\lambda$ to the first factor.
The variety $Z_\lambda$ could be rather complicated. 
Say, it is not irreducible in general.
We shall use the fact that fibers $\pi^{-1}(x)$ are linear spaces.
So consider the function $l_\lambda(x)=\dim\pi^{-1}(x)$, $x\in\g/\p$.

\begin{Conjecture}\label{KJFKUTYFKLUFTL}
There exists an algebraic stratification 
$\g/\p=\mathop{\sqcup}\limits_{i=1}^rX_i$
such that for any $\lambda$ the function $l_\lambda(x)$ is constant
along each $X_i$: 
$$l_\lambda(x)=l_\lambda^i,\quad x\in X_i.$$
\end{Conjecture}
If such stratification exists then
\begin{equation}
\dim Z_\lambda=\mathop{\max}\limits_i(\dim X_i+l_\lambda^i).\label{ISOTRSUBSPeqno(3)}
\end{equation}

If $P$ acts on $\g/\p$ with a finite number of orbits
then we can take these orbits as $X_i$
(since $\psi$ is $P$-equivariant).
Unfortunately, this class of parabolic subgroups is very small. 
The dual $P$-module $(\g/\p)^*$ is isomorphic
to the representation of $P$ on the unipotent radical $\p_n$ of $\p$,
hence by Pyasetskii Theorem~\ref{PyasetskyTh} the action $P:\,\g/\p$ has a finite number 
of orbits iff the action $P:\,\p_n$ satisfies this property.
Such actions were studied in \cite{PR}, the complete classification 
for classical groups was obtained
in \cite{HRo}. 
But even in these cases the problem
of the complete description of the orbital decomposition
seems very messy.
In the next section we shall introduce another class of parabolic subgroups
that satisfy the Conjecture~\ref{KJFKUTYFKLUFTL}.

The most wonderful class of parabolic subgroups 
that will fit all our needs is the class of parabolic subgroups $P$
with abelian unipotent radical $\p_n$, see Section~\ref{ParabolicsAura}.
So from now on $P$ will be a parabolic subgroup with {\it aura}.
Let $P^-$ be an opposite parabolic subgroup,
let $L=P\cap P^-$ be a Levi subgroup with Lie algebra $\l$, 
let $\p^-_n$ be a unipotent radical of~$\p^-$.
Then $\p_n$ is abelian iff the decomposition
$\g=\p^-_n\oplus\l\oplus\p_n$ is a $\Z$-grading:
$\g=\g_{-1}\oplus\g_0\oplus\g_1$. Notice that
$\g/\p$  is isomorphic to $\p^-_n$ as a $L$-module. 
It is well-known that the action of
$L$ on $\p^-_n$ has a finite number of orbits in this case.
Take the orbital decomposition $\p^-_n=\mathop{\cup}\limits_{i=1}^rLf_i$.
For its detailed description see Theorem~\ref{AURAOrbitDecomp}.
Then the equality (\ref{ISOTRSUBSPeqno(3)}) takes a form
\begin{equation}
\dim Z_\lambda=\mathop{\max}\limits_i(\dim Lf_i+l_\lambda(f_i)).\label{ISOTRSUBSPeqno(4)}
\end{equation}
It remains to calculate $l_\lambda(f_i)$. 
If $f_i=0$ then $l_\lambda(f_i)=\dim M_\lambda$.
In the opposite case we can include $f_i$ in a $\ssl_2$-triple 
$\langle f_i, h_i, e_i\rangle$.
Evidently we can assume that $h_i\in\l$, $e_i\in\p_n$.
These $h_i$'s were also written down explicitly in Theorem~\ref{AURAOrbitDecomp}.
We shall use the following easy lemma from the $\ssl_2$-theory:

\begin{Lemma}
Let $V$ be a finite-dimensional $\ssl_2$-module, 
$\ssl_2=\langle f,h,e\rangle$.
Let $\b=\langle h,e\rangle$.
Suppose that $M\subset V$ is a $\b$-submodule such that $eV\subset M$. 
Then
$$\dim\Coker(M\subset V\mathop{\to}\limits^fV\to V/M)=\dim (V/M)^h,$$ 
where $(V/M)^h=\{v\in V/M\,|\,hv=0\}$.
\end{Lemma}
Applying this Lemma 
to the representation of $\ssl_2$ in $V_\lambda$ we get
\begin{eqnarray*}
l_\lambda(f_i)=\dim\Ker\psi(f_i,\cdot)&=
\dim M_\lambda-\dim U_\lambda+\dim\Coker\psi(f_i,\cdot)=\cr
&=\dim M_\lambda-\dim U_\lambda+\dim U_\lambda^{h_i}.
\end{eqnarray*}
Notice that this formula remains valid
for $f_i=0$ if we assume that in this case $h_i=0$.
Joining this formula with~(\ref{ISOTRSUBSPeqno(4)}) and applying the Lemma above
we get the following result.

\begin{Theorem}[\cite{T3}]\label{JHGVFJHGFHFGDJH}
Suppose that for any $i$ 
$$\dim U_\lambda^{h_i}
\le\codim_{\g/\p} Lf_i.\label{ISOTRSUBSPeqno(5)}$$
Then generic global sections of the bundle $\L_\lambda$ have zeros.
\end{Theorem}

All ingredients of the formula (\ref{ISOTRSUBSPeqno(5)}) can be easily computed
in many particular cases.
We will apply Theorem~\ref{JHGVFJHGFHFGDJH} 
for the proof of the following

\begin{Theorem}[\cite{T3}]\label{LJKGDJUFGDJKDHKJ}
Let $w\in S^dV^*$, resp.~$w\in\Lambda^d V^*$, be a generic form. Then
$V$ contains a $k$-dimensional isotropic subspace
w.r.t.~$w$ if and only if
\begin{equation}
n\ge{{\ds d+k-1\choose\ds d}\over k}+k,\ \text{resp.}
\ n\ge{{\ds k\choose\ds d}\over k}+k,\label{IsotrSubspEqno(1)}
\end{equation}
with the following exceptions {\rm(}a form is supposed to be generic{\rm)}.
$V$ contains a $k$-dimensional isotropic subspace if and only if

\begin{itemize}
\item $n\ge2k$ for $w\in S^2V^*$ or $w\in \Lambda^2V^*$.
\item $k\le n-2$ for $w\in \Lambda^{n-2}V^*$, $n$ is even.
\item $k\le 4$ for $w\in \Lambda^3V^*$, $n=7$.
\end{itemize}
\end{Theorem}

It is interesting that the variety
of $4$-dimensional isotropic subspaces of a generic skew\-symmetric
$3$-form in ${\Bbb C}^7$ is a smooth $8$-dimensional
Fano variety. Moreover, this variety is a compactification
of the unique symmetric space of the simple
algebraic group~$\text{G}_2$.

\smallskip

\proof
We choose a basis $\{e_1,\ldots,e_n\}$ in $V$ and 
identify $\GL_n$ with the group of non-singular matrices.
Let $T$, $B$, $B^-$ be the subgroups of diagonal, 
upper- and lower-triangular matrices.
We fix an integer $k$. Consider the parabolic subgroup 
$$P=\left(\matrix{
A&0\cr
*&B\cr}\right),$$
where $B$ is a $k\times k$-matrix.
Then $G/P$ is the Grassmanian $\Gr(k,V)$. 
Consider the vector bundle
$\L=S^d\S^*$ (resp.~$\L=\Lambda^d\S^*$, but only in the case $k\ge d$)
on $G/P$,
where $\S$ is the tautological bundle. 
Then $\L=\L_\lambda$, where $\lambda$ is the highest weight of the
$\GL_n$-module $S^dV^*$ (resp.~$\Lambda^dV^*$). 
Let $w\in S^dV^*$ (resp.~$w\in\Lambda^dV^*$), 
let $s_w$ be the corresponding global section. 
It is easy to see that $(Z_{s_w})_{red}$ coincides with the variety
of $k$-dimensional isotropic subspaces w.r.t.~$w$. 
Notice that inequalities (\ref{IsotrSubspEqno(1)}) are equivalent to the condition
$\dim U_\lambda\le\dim G/P$.
In the sequel we suppose that these inequalities hold.

We take
$$f_i=E_{1,n}+E_{2,n-1}+\ldots+E_{i,n+1-i},\quad i=0,\ldots,r=\min(k,n-k).$$
The unipotent radical $\p^-_n$ can be identified
with matrices of shape $k\times(n-k)$ then the orbit
$Lf_i$ is identified with the variety of matrices of rank $i$, therefore
$$\codim_{\g/\p} Lf_i=(k-i)(n-k-i).$$
We have
$$h_i=(E_{n,n}-E_{1,1})+(E_{n-1,n-1}-E_{2,2})+\ldots+
(E_{n+1-i,n+1-i}-E_{i,i}),$$
hence
\begin{eqnarray*}
\dim U_\lambda^{h_i}={d+k-i-1\choose d}\ \text{if}\ V_\lambda=S^dV^*;\cr
\dim U_\lambda^{h_i}={k-i\choose d}\ \text{if}\ V_\lambda=\Lambda^dV^*.
\end{eqnarray*}
Applying Theorem\ref{JHGVFJHGFHFGDJH} we get

\begin{Proposition}\label{JKHGJGHFJGFDJGFD}$\ $
\begin{enumerate}
\item[\rm(a)]
Suppose that for any $i=0,\ldots,r$ we have
$${\ts d+k-i-1\choose\ts d}\le(k-i)(n-k-i).$$
Then a generic form $w\in S^dV^*$ has a  $k$-dimensional
isotropic subspace.
\item[\rm(b)]
Suppose that for any $i=0,\ldots,r$ we have
$${\ts k-i\choose\ts d}\le(k-i)(n-k-i).$$
Then a generic form $w\in \Lambda^dV^*$ has a $k$-dimensional
isotropic subspace.
\end{enumerate}
\end{Proposition}

It remains to clarify when 
the conditions of Proposition~\ref{JKHGJGHFJGFDJGFD} 
follow from formulas~(\ref{IsotrSubspEqno(1)}).

\subsection*{A. Symmetric Case.} 
Clearly if $i=k$ then the conditions of Proposition~\ref{JKHGJGHFJGFDJGFD}  
are satisfied. 
Therefore it suffices to find out when the inequality
$$n\ge{\ts{ d+k-i-1\choose d}\over\ts k-i}+k+i$$ 
follows from the inequality
$$n\ge{\ts{d+k-1\choose d}\over\ts k}+k$$
as $i=1,\ldots,\min(k-1, n-k)$.
If $d=1$ then Theorem~\ref{LJKGDJUFGDJKDHKJ} 
is obvious, if $d=2$ it reduces to a well-known
result about isotropic subspaces of quadratic form.
So assume that $d\ge3$.
We use a following Lemma that can be easily verified by induction:
\begin{Lemma}\label{proclaimLemma5}
Let $d\ge3$, $\alpha\ge2$. Then 
${\ts{d+\alpha-1\choose d}\over\ts\alpha}
\ge{\ts{d+\alpha-2\choose d}\over\ts\alpha-1}+1.$
\end{Lemma}
It follows from Lemma~\ref{proclaimLemma5} that
\begin{eqnarray*}
n\ge{\ts{d+k-1\choose d}\over\ts k}+k\ge
{\ts{d+k-2\choose d}\over\ts k-1}+k+1\qquad\qquad\cr
\ge
{\ts{d+k-3\choose d}\over\ts k-2}+k+2\ge\ldots\ge
{\ts{d+k-1-(k-1)\choose d}\over\ts k-(k-1)}+k+k-1.
\end{eqnarray*}
QED.

\subsection*{B. Skew-symmetric Case.}
Clearly if $i>k-d$ then the conditions of Proposition~\ref{JKHGJGHFJGFDJGFD}
are satisfied. Therefore it suffices to check when
the inequality 
$$n\ge{\ts{k-i\choose d}\over\ts k-i}+k+i$$
follows from the inequality
$$n\ge{\ts{k\choose d}\over\ts k}+k$$ 
as
$i=1,\ldots,\min(k-d, n-k)$.
If $d=1$ then Theorem~\ref{LJKGDJUFGDJKDHKJ}  
is obvious, if $d=2$ it reduces
to a well-known result about isotropic subspaces of  a $2$-form.
So assume that $d\ge3$.
We use a following Lemma that can be easily verified by induction:
\begin{Lemma}\label{proclaimLemma6}
Let $3\le d\le\alpha-2$. Then ${\ts{\alpha\choose d}\over\ts\alpha}
\ge{\ts{\alpha-1\choose d}\over\ts\alpha-1}+1.$
\end{Lemma}

It follows from Lemma~\ref{proclaimLemma6} that
\begin{eqnarray*}
n\ge{\ts{k\choose d}\over\ts k}+k\ge
{\ts{k-1\choose d}\over\ts k-1}+k+1\qquad\qquad\cr
\ge
{\ts{k-2\choose d}\over\ts k-2}+k+2\ge\ldots\ge
{\ts{k-(k-3-d)\choose d}\over\ts k-(k-3-d)}+k+(k-3-d).
\end{eqnarray*}
It remains to consider 3 cases: $i=k-d-2$, $i=k-d-1$, $i=k-d$
as $i\le n-k$.

Let $i=k-d-2$, $n\ge2k-d-2$. 
We want to deduce the inequality
$n\ge{\ts{d+2\choose d}\over\ts d+2}+2k-d-2=2k-{d-1\over2}-2$
from $n\ge{\ts{k\choose d}\over\ts k}+k$.
It suffices to check
${k\choose d}\ge k(k-{d-1\over2}-2)$.
Since $d\ge3$ and $k-d=i+2\ge3$ we have
${k\choose d}\ge{k\choose 3}\ge k(k-3)\ge k(k-{d-1\over2}-2)$.

Let $i=k-d-1$, $n\ge2k-d-1$. 
We want to deduce the inequality
$n\ge{\ts{d+1\choose d}\over\ts d+1}+2k-d-1=2k-d$
from $n\ge{\ts{d+1\choose d}\over\ts d+1}+2k-d-1=2k-d$.
It suffices to check
${k\choose d}\ge k(k-d)$. 
Since $d\ge3$ and $k-d=i+1\ge2$ we have
${k\choose d}\ge{k\choose 3}\ge k(k-3)\ge k(k-d)$
as $k-d\ge3$.
If $k-d=2$ then
${k\choose d}={(d+2)(d+1)\over2}\ge 2(d+2)\ge k(k-d)$. 

Let $i=k-d$, $n\ge2k-d$. 
We need to deduce the inequality
$n\ge{\ts{d\choose d}\over\ts d}+2k-d=2k-d+{\ts1\over\ts d}$
from
$n\ge{\ts{k\choose d}\over\ts k}+k$.
All exceptions are defined by the following
system of equalities and inequalities:
$$n=2k-d,\quad {k\choose d}\le k(k-d),\quad k-d\ge1.$$
There exist only two possibilities: 
either $k=d+1$, $n=d+2$ or $d=3$, $k=5$, $n=7$.

In the first case
we need to clarify whether a generic form
$w\in\Lambda^{n-2}V^*$ has a $(n-1)$-dimensional isotropic subspace or not.
This holds iff a generic $2$-form in $V^*$
has a non-zero kernel.
It is well-known that this is true only for odd $n$.

To complete the proof it 
remains to find out whether a generic form
$w\in\Lambda^3(\C^7)^*$ has a $5$-dimensional isotropic subspace. 
In suitable coordinates this isotropic subspace
will coincide with the linear span $\langle e_1,\ldots,e_5\rangle$.
Therefore $w=\sum\limits_{1\le i<j<k\le 7}\alpha_{ijk}x_i\wedge x_j\wedge x_k$,
where $\alpha_{ijk}=0$ as $k\le5$.
Consider the one-parameter subgroup
$$H(t)=\diag(t^2,t^2,t^2,t^2,t^2,t^{-5},t^{-5})$$
in $\SL_7$.
Clearly, $\lim\limits_{t\to0}H(t)w=0$.
Therefore $w$ belongs to the Null-cone
of the action $\SL_7:\,\Lambda^3(\C^7)^*$.
It follows that any non-constant homogeneous invariant
of this action vanishes at $w$. 
But it is known that this action has non-constant invariants (see~\cite{PV}),
for example, the discriminant is well-defined in this case.
Therefore generic forms do not admit
$5$-dimensional isotropic subspaces.

Moreover, it can be shown that a
$3$-form in the $7$-dimensional vector space 
admits a $5$-dimensional isotropic subspace
if and only if its discriminant is equal to zero.
\endproof

\section{Moore--Penrose Inverse and Applications}

The nice notion of a generalized inverse of an arbitrary matrix
(possibly singular or even non-square) has been discovered
independently by Moore \cite{Mo} and Penrose \cite{Pe}.
The following\index{Moore--Penrose inverse} 
definition belongs to Penrose (Moore's definition is different but 
equivalent):

\begin{Definition}\rm
A matrix $A^+$ is called {\it a MP-inverse\/} of a 
matrix $A$
if 
$$AA^+A=A,\quad A^+AA^+=A^+,$$ 
and $AA^+$, $A^+A$ are Hermitian matrices.
\end{Definition}

It is quite surprising but a MP-inverse always exists and is unique.
Since the definition is symmetric with respect to $A$ and $A^+$
it follows that $(A^+)^+=A$.
If $A$ is a non-singular square matrix then $A^+$
coincides with an ordinary inverse matrix $A^{-1}$.
The theory of MP-inverses and their numerous modifications becomes
now a separate subfield of Linear Algebra \cite{CM} with various applications.
Here we show that this notion
quite naturally arises in the theory of shortly graded simple Lie algebras
and give applications.
To explain this connection let us first give another definition of a MP-inverse.

Let $A\in\Mat_{n,m}(\C)$. 
Then it is easy to see that a matrix $A^+\in\Mat_{m,n}(\C)$
is a MP-inverse of $A$ if and only if there exist Hermitian matrices 
$$B_1\in\Mat_{n,n}(\C)\ \text{and}\ B_2\in\Mat_{m,m}(\C)$$
such that the following matrices form an $\ssl_2$-triple
in $\ssl_{n+m}(\C)$:
$$E=\left(\matrix{ 0&A\cr0&0\cr}\right),\quad 
H=\left(\matrix{ B_1&0\cr0&B_2\cr}\right),\quad
F=\left(\matrix{ 0&0\cr A^+&0\cr}\right).$$

By an $\ssl_2$-triple $\langle e,h,f\rangle$ in a Lie algebra $\g$
we mean a collection of (possibly zero) vectors
such that 
$$[e,f]=h,\quad [h,e]=2e,\quad [h,f]=-2f.$$
In other words, an $\ssl_2$-triple is a homomorphic image
of canonical generators of $\ssl_2$ with respect to some homomorphism 
of Lie algebras $\ssl_2\to\g$.

This definition admits an immediate generalization.
Suppose that $\g$ is a simple complex Lie algebra, $G$
is a corresponding simple simply-connected Lie group.
Suppose further that $P$ is a parabolic subgroup of $G$
with abelian unipotent radical (with {\it aura}).
Then $\g$ admits a short grading
$$\g=\g_{-1}\oplus \g_0\oplus\g_1$$
with only three nonzero parts. Here $\p=\g_0\oplus\g_1$ is a Lie
algebra of $P$ and $\exp\g_1$ is the abelian unipotent radical of $P$.
Let $\k_0$ be a compact real form of $\g_0$.
In this section
we shall permanently consider compact real forms of reductive
subalgebras of simple Lie algebras. These subalgebras
will always be Lie algebras of algebraic reductive subgroups
of a corresponding simple complex algebraic group. 
Their compact real forms will always be understood
as Lie algebras of compact real forms of corresponding
algebraic groups. 
For example, a Lie algebra of an algebraic torus has a unique
compact real form.

Suppose now that $e\in\g_1$. It is well-known that
there exists a {\it homogeneous} 
$\ssl_2$-triple\index{homogeneous $\ssl_2$-triple} $\langle e,h,f\rangle$
such that $h\in\g_0$ and $f\in\g_{-1}$.

\begin{Definition}\rm
An element $f\in\g_{-1}$ is called {\it a MP-inverse\/} of $e\in\g_1$
if there exists a homogeneous $\ssl_2$-triple $\langle e,h,f\rangle$
with $h\in i\k_0$.
\end{Definition}

MP-inverses of elements $f\in\g_{-1}$ are defined in the same
way. It is clear that if $f$ is a MP-inverse of $e$
then $e$ is a MP-inverse of $f$. 

\begin{Example}\rm
Suppose that $G=\SL_{n+m}$ and $P\subset G$ is a maximal
parabolic subgroup of block triangular matrices
of the form 
$$\left(\matrix{ B_1&A\cr 0&B_2\cr}\right),\quad\text{where}\quad
B_1\in\Mat_{n,n},\ A\in\Mat_{n,m},\ B_2\in\Mat_{m,m}.$$
The graded components of the correspondent grading
consist of matrices of the following form:
$$\g_{-1}=\left(\matrix{ 0&0\cr A'&0\cr}\right),\quad
\g_0=\left(\matrix{ B_1&0\cr 0&B_2\cr}\right),\quad
\g_1=\left(\matrix{ 0&A\cr 0&0\cr}\right),$$
where $A'\in\Mat_{m,n}$, $B_1\in\Mat_{n,n}$, $B_2\in\Mat_{m,m}$, and
$A\in\Mat_{n,m}$.
One can take $\k_0$ to be the real Lie algebra of block diagonal
skew-Hermitian matrices with zero trace. Then $i\k_0$ is a vector space
of block diagonal Hermitian matrices with zero trace.
Therefore in this case we return to a previous definition
of a Moore--Penrose inverse.
\end{Example}

\begin{Theorem}[\cite{T4}]
For any $e\in\g_1$ there exists a unique Moore--Penrose inverse $f\in\g_{-1}$.
\end{Theorem}

The proof is similar to the proof of Theorem~\ref{KJHGKJGFJ} below.

It obviously follows that
for any non-zero $f\in\g_{-1}$ there exists a unique MP-inverse $e\in\g_1$.
So taking a MP-inverse is a well-defined involutive operation.
In general, it is not equivariant with respect to a Levi subgroup $L\subset P$
with Lie algebra $\g_0$,
but only with respect to its maximal compact subgroup $K_0\subset L$.

Let us give an intrinsic description
of the Moore-Penrose inverse in all cases arising
arising from short gradings of classical simple Lie algebras.
Exceptional cases may be found in \cite{T4}.

\subsection*{Linear Maps.}
This is, of course, the classical Moore--Penrose inverse.
Let us recall its intrinsic description.
Suppose that $\C^n$ and $\C^m$ are vector spaces
equipped with standard Hermitian scalar products.
For any linear map $F:\,\C^n\to\C^m$ its Moore-Penrose
inverse is a linear map $F^+:\,\C^m\to\C^n$ defined as follows.
Let $\Ker F\subset\C^n$ and $\Im F\subset \C^m$ be
a kernel and an image of $F$.
Let  $\Ker^\perp F\subset\C^n$ and $\Im^\perp F\subset \C^m$
be their orthogonal complements with respect 
to the Hermitian scalar products.
Then $F$ defines via restriction a bijective linear map
$\tilde F:\,\Ker^\perp F\to\Im F$.
Then $F^+:\,\C^m\to\C^n$ is a unique linear map such that
$F^+|_{\Im^\perp F}=0$ and 
$F^+|_{\Im F}=\tilde F^{-1}$.
This MP-inverse corresponds to short gradings of $\ssl_{n+m}$.

\subsection*{Symmetric and Skew-symmetric Bilinear Forms.}
Let $V=\C^n$ be a vector space equipped with 
a standard Hermitian scalar product.
For any symmetric (resp.~skew-symmetric) bilinear form $\omega$ on $V$
its Moore-Penrose inverse is 
a symmetric (resp.~skew-symmetric) bilinear form $\omega^+$ on $V^*$
defined as follows.
Let $\Ker\, \omega\subset V$ be the kernel of~$\omega$.
Then $\omega$ induces a non-degenerate bilinear form
$\tilde\omega$ on $V/\Ker\,\omega$.
Let $\Ann(\Ker\, \omega)\subset V^*$ 
be an annihilator of $\Ker\, \omega$. 
Then $\Ann (\Ker\, \omega)$ is canonically isomorphic to the dual of $V/\Ker\,\omega$.
Therefore the form
$\tilde\omega^{-1}$ on 
$\Ann (\Ker\, \omega)$ is well-defined.
The form $\omega^+$ is defined as a unique form such that
its restriction on
on $\Ann(\Ker\,\omega)$ coincides with 
$\tilde\omega^{-1}$ and its kernel is $\Ann(\Ker\,\omega)^\perp$,
the orthogonal complement with respect to a standard
Hermitian scalar product on $V^*$.
This MP-inverse corresponds to the short grading of $\sp_{2n+2}$ 
(resp.~$\so_{2n+2}$).

\subsection*{Vectors in a Vector Space With the Scalar Product.}
Suppose that $V=\C^n$ is a vector space with the standard bilinear
scalar product $(\cdot,\cdot)$.
For any vector $v\in V$ its Moore-Penrose inverse
$v^\vee$ is again a vector in $V$ defined as follows:
$$v^\vee=\cases{
\bigover{2v}{(v,v)},&\hbox{if $(v,v)\ne0$}\cr
\bigover{\overline v}{(\overline v,v)},&\hbox{if $(v,v)=0$, $v\ne0$}\cr
0,&\hbox{if $v=0$.}}$$
This MP-inverse corresponds to the short grading of~$\so_{n+2}$.

It is quite natural to ask whether it is possible
to extend the notion of the Moore--Penrose inverse
from parabolic subgroups with aura to arbitrary parabolic subgroups.
It is also interesting to consider the ``non-graded'' situation.
Let us start with it. Suppose $G$ is a simple connected simply-connected
Lie group with Lie algebra $\g$. We fix a compact real form $\k\subset\g$.

\begin{Definition}\rm
A nilpotent orbit $\O\subset\g$ is called {\it a Moore--Penrose orbit\/} if for any
$e\in\O$ there exists an $\ssl_2$-triple $\langle e,h,f\rangle$
such that $h\in i\k$.
\end{Definition}

It turns out that it is quite easy to find all Moore-Penrose orbits.
Recall that the height\index{height of a nilpotent orbit} 
$\hht(\O)$ of a nilpotent orbit $\O=\Ad(G)e$
is equal to the maximal integer $k$ such that $\ad(e)^k\ne0$.
Clearly $\hht(\O)\ge2$.

\begin{Theorem}\label{KJHGKJGFJ}
$\O$ is a Moore--Penrose orbit if and only if $\hht(\O)=2$.
In this case for any $e\in\O$ there exists a unique $\ssl_2$-triple
$\langle e,h,f\rangle$ such that $h\in i\k$.
\end{Theorem}

\proof
Let $x\to\overline x$ denotes a complex conjugation in $\g$ with respect to 
the compact form $\k$. 
Therefore $x=\overline x$ iff $x\in\k$ and $x=-\overline x$ iff $x\in i\k$.
Let $B(x,y)=\Tr\, \ad(x)\ad(y)$ be the Killing form of $\g$.
Finally, 
let $H(x,y)=-B(x,\overline y)$ be a positive-definite Hermitian form on $\g$.

\begin{Lemma}\label{OLJHGKLLJ:LKGJFKGJF}
We fix a nilpotent element $e\in\g$.
Suppose that $\langle e,h,f\rangle$ is an $\ssl_2$-triple in $\g$ such that 
$h\in i\k$. Then for any other $\ssl_2$-triple 
$\langle e,h',f'\rangle$ we have $H(h,h)<H(h',h')$.
In particular,
if there exists an $\ssl_2$-triple $\langle e,h,f\rangle$
with $h\in i\k$ then the $\ssl_2$-triple with this property is unique.
\end{Lemma}

\proof
Recall that if $\langle e,h,f\rangle$ is an $\ssl_2$-triple then 
$h$ is called\index{characteristic of a nilpotent} {\it a characteristic\/} of~$e$.
Consider the subset $\cal H\subset\g$ consisting 
of all possible characteristics of~$e$. It is well-known that
$\cal H$ is an affine subspace in $\g$ such that the corresponding
linear subspace is precisely the unipotent radical $\z^u_{\g}(e)$ of the centralizer
$\z_{\g}(e)$ in~$\g$ of the element $e$.
Since $H(h',h')$ is a strongly convex function on $\cal H$, 
there exists a unique element $h_0\in\cal H$
such that $H(h_0,h_0)<H(h',h')$ for any $h'\in\cal H$, $h'\ne h_0$.
We need to show that $h_0=h$.
It is clear that an element $h_0\in\cal H$ minimizes 
$H(h,h)$ on $\cal H$ iff $H(h_0,\z^u_{\g}(e))=0$ iff
$B(\overline h_0,\z^u_{\g}(e))=0$.
If $h\in \cal H\cap i\k_0$ then $\overline h=-h$ and we have
$$B(\overline h,\z^u_{\g}(e))=
-B(h,\z^u_{\g}(e))=
-B([e,f],\z^u_{\g}(e))=
B(f,[e,\z^u_{\g}(e)])=0.$$
Therefore $h=h_0$.
\endproof

Suppose that $\langle e,h,f\rangle$ is an $\ssl_2$-triple in $\g$.
Consider the grading $\g=\mathop{\oplus}\limits_k\g_k$ such that $x\in\g_k$
iff $[h,x]=kx$. Let $\n_+=\mathop{\oplus}\limits_{k>0}\g_k$, $\n_-=\mathop{\oplus}\limits_{k<0}\g_k$.
It is well known that $\z_\g^u(e)\subset\n_+$.

\begin{Lemma}
Suppose that $\z_\g^u(e)=\n_+$. Then $\O=\Ad(G)e$ is a Moore--Penrose orbit.
\end{Lemma}

\proof
We need to prove that for any element $e'\in\O$ there exists
an $\ssl_2$-triple $\langle e',h',f'\rangle$ such that
$h'\in i\k$, where $\k$ is a {\it fixed} compact real form of $\g$.
Clearly it is sufficient to prove that for an {\it arbitrary}
compact real form $\k$ there exists an $\ssl_2$-triple
$\langle e,h,f\rangle$ with $h\in i\k$.
According to the proof of the previous Lemma
we should choose 
$h$ to be a unique characteristic such that
$B(\overline h,\z^u_{\g}(e))=0$, where
$x\to\overline x$ denotes a complex conjugation in $\g$ with respect to 
the compact form $\k$.
It remains to prove that $h\in i\k$.
Since $B$ is a non-degenerate $\ad$-invariant scalar product on $\g$
and $\z_\g^u(e)=\n_+$
it follows that $\overline h\in\p$, where $\p=\g_0\oplus\n_+$.
Let $\l$ be some ``standard'' compact real form of $\g$
such that $h\in i\l$ and $\tilde\n_{\pm}=\n_{\mp}$,
where 
$x\to\tilde x$ denotes a complex conjugation in~$\g$ with respect to 
the compact form $\l$.
Let $P\subset G$ be a parabolic subgroup of $G$
with the Lie algebra $\p$, let $H\subset P$ be its Levi subgroup
with the Lie algebra $\g_0$.
Since all compact real forms of a semisimple complex Lie algebra
are conjugated by elements of any fixed Borel subgroup
it follows that there exists $g\in P$ such that
$\Ad(g)\k=\l$.
Therefore 
$$\widetilde{\Ad(g)h}=\Ad(g)\overline h\subset\Ad(g)(\p)=\p.$$
We can express $g$ as a product $uz$, where $u\in\exp(\n_+)$, $\Ad(z)h=h$.
Then 
$\widetilde{\Ad(g)h}=\widetilde{\Ad(u)h}$.
If $u$ is not the identity element of $G$
then $\Ad(u)h=h+\xi$, where $\xi\in\n_+$ and $\xi\ne0$.
Therefore $\widetilde{\Ad(u)h}=-h+\tilde\xi$.
But $\tilde\xi\in\n_{-}$ and hence
$\widetilde{\Ad(u)h}\not\in\p$, contradiction.
Therefore $u$ is trivial and since
$\Ad(z)h=h$ we finally get
$$\overline h=\tilde h=-h.\qed$$
\endproof

Now we shall try to reverse this argument.

\begin{Lemma}
Suppose that $\O=\Ad(G)e$ is a Moore--Penrose orbit. Then $\z_\g^u(e)=\n_+$.
\end{Lemma}

\proof
We choose a standard compact real form $\l$ 
as in the proof of the previous Lemma.
Clearly, $\z_\g^u(e)$ is a graded subalgebra of $\n_+=\mathop{\oplus}\limits_{k>0}\g_k$.
Suppose, on the contrary, that $\z_\g^u(e)\ne\n_+$. Let $\xi\in\g_p$, $p>0$,
be a homogeneous element that does not belong to $\z_\g^u(e)$.
Let $u=\exp(\xi)$. Let $e'=\Ad(u)e$. We claim that all characteristics
of $e'$ do not belong to $i\l$. Indeed, all characteristics
of $e'$ have a form $\Ad(u)h+\Ad(u)x$, where $x\in\z_\g^u(e)$.
Suppose that for some $x$ we have $\Ad(u)h+\Ad(u)x\in i\l$.
Since $h\in i\l$, $\tilde\n_{\pm}=\n_{\mp}$, and
$\Ad(u)(h+x)-h\in\n_+$ it follows that 
$\Ad(u)(h+x)=h$. In $\n_+$ modulo 
$\mathop{\oplus}\limits_{k>p}\g_k$ we obtain the equation
$[\xi,h]+x=0$, but $[h,\xi]=p\xi$ and therefore $\xi\in\z_\g^u(e)$.
Contradiction.
\endproof

Now we can finish the proof of Theorem~\ref{KJHGKJGFJ}.
Combining previous lemmas we see that 
$\O$ is a Moore--Penrose orbit if and only if 
$\z_\g^u(e)=\n_+$.
It follows from the $\ssl_2$-theory that
$\dim\z_\g^u(e)=\dim\g_1+\dim\g_2$.
Therefore $\z_\g^u(e)=\n_+$ if and only if $\g_p=0$ for $p>2$.
Clearly, this is precisely equivalent to $\hht(\O)=2$.
In this case for any $e\in\O$ there exists a unique $\ssl_2$-triple
$\langle e,h,f\rangle$ such that $h\in i\k$ by Lemma~\ref{OLJHGKLLJ:LKGJFKGJF}.
\endproof

Now let us turn to the graded situation.
Suppose that $\g$ is a $\Z$-graded simple Lie algebra, 
$\g=\mathop{\oplus}\limits_{k\in\Z}\g_k$.
Let $P\subset G$ be a parabolic subgroup with Lie algebra 
$\p=\mathop{\oplus}\limits_{k\ge0}\g_k$. Let $L\subset P$ be a Levi subgroup with Lie algebra
$\g_0$. 
We choose a compact real form $\k_0$ of $\g_0$.
Suppose now that $e\in\g_k$. It is well-known that
there exists a homogeneous $\ssl_2$-triple $\langle e,h,f\rangle$
with $h\in\g_0$ and $f\in\g_{-k}$.

\begin{Definition}\rm
Take any $k>0$ and any $L$-orbit $\O\subset\g_k$.
Then $\O$ is called {\it a Moore-Penrose orbit\/} if for any $e\in\O$
there exists a homogeneous $\ssl_2$-triple $\langle e,h,f\rangle$
such that $h\in i\k_0$. In this case $f$ is called a MP-inverse of $e$.
A grading is called {\it
a Moore--Penrose grading
in degree $k>0$\/} if all $L$-orbits in $\g_k$ are Moore-Penrose. 
A grading is called {\it a Moore--Penrose grading\/} if it
is a Moore--Penrose grading in any positive degree.
A parabolic subgroup $P\subset G$ is called {\it a Moore--Penrose
parabolic subgroup\/} if there exists a Moore--Penrose grading 
$\g=\mathop{\oplus}\limits_{k\in\Z}\g_k$ such that 
$\p=\mathop{\oplus}\limits_{k\ge0}\g_k$ is a Lie algebra of $P$.
\end{Definition}

One should be careful comparing graded and non-graded situation:
if $\O\subset\g_k$ is a Moore--Penrose $L$-orbit then
$\Ad(G)\O\subset\g$ is not necessarily a Moore--Penrose
$G$-orbit.
Let us give a criterion 
for an $L$-orbit to be Moore--Penrose.
Suppose $\O=\Ad(L)e\subset\g_k$. Take any homogeneous $\ssl_2$-triple
$\langle e,h,f\rangle$. Then $h$ defines a grading 
$\g_0=\mathop{\oplus}\limits_{n\in\Z}\g_0^n$,
such that $\ad(h)|_{\g_0^n}=n\cdot\Id$.

\begin{Theorem}[\cite{T4}]\label{KJHGKHGDFKJFDLK}
$\O$ is a Moore--Penrose orbit if and only if 
$\ad(e)g_0^n=0$ for any $n>0$.
In this case for any $e'\in\O$ there exists a unique homogeneous 
$\ssl_2$-triple $\langle e',h',f'\rangle$ such that
$h'\in i\k_0$.
\end{Theorem}

The proof is similar to the proof of Theorem~\ref{KJHGKJGFJ}.

\begin{Example}\label{LKJHGKHGFDJKD}\rm
Suppose that $G$ is a simple group of type $\text{G}_2$.
We fix a root decomposition.
There are two simple roots $\alpha_1$ and $\alpha_2$ such that
$\alpha_1$ is short and $\alpha_2$ is long.
There are $3$ proper parabolic subgroups: Borel subgroup $B$
and two maximal parabolic subgroups $P_1$ and $P_2$ such that
a root vector of $\alpha_i$ belongs to a Levi subgroup of $P_i$.
Then the following is an easy application of Theorem~\ref{KJHGKHGDFKJFDLK}.
$B$ is a Moore--Penrose parabolic subgroup (actually
Borel subgroups in all simple groups are Moore--Penrose parabolic
subgroups
with respect to any grading).
$P_1$ is not Moore--Penrose, but it is a Moore--Penrose
parabolic subgroup in degree $2$ (with respect to the
natural grading of height~$2$). 
$P_2$ is a Moore--Penrose parabolic subgroup.
\end{Example}

\begin{Example}\rm
Suppose $G=\SL_n$.
We fix positive integers $d_1,\ldots,d_k$ such that $n=d_1+\ldots+d_k$.
We consider the parabolic subgroup $P(d_1,\ldots,d_k)\subset\SL_n$
that consists of all upper-triangular block matrices with sizes
of blocks equal to $d_1,\ldots,d_k$. We take a standard grading.
Then $\g_1$ is identified with the linear space of all tuples of linear maps
$\{f_1,\ldots,f_k\}$,
$$\C^{d_1}\mathop{\longleftarrow}\limits^{f_1}\C^{d_2}
\mathop{\longleftarrow}\limits^{f_2}\ldots
\mathop{\longleftarrow}\limits^{f_{k-1}}\C^{d_k},$$
$\g_{-1}$ is identified with the linear space of all tuples of linear maps
$\{g_1,\ldots,g_k\}$,
$$\C^{d_1}\mathop{\longrightarrow}\limits^{g_1}\C^{d_2}
\mathop{\longrightarrow}\limits^{g_2}\ldots
\mathop{\longrightarrow}\limits^{g_{k-1}}\C^{d_k},$$
and Levi subgroup $L(d_1,\ldots,d_k)$
is just the group  
$$(A_1,\ldots,A_k)\in\GL_{d_1}\times\ldots
\times\GL_{d_k}\quad\text{such that}\quad\det(A_1)\cdot\ldots\cdot\det(A_k)=1.$$
It acts on these spaces of linear maps in an obvious way.
The most important among $L$-orbits 
are varieties of complexes\index{variety of complexes}.
To define them, let us fix in addition non-negative integers
$m_1,\ldots,m_{k-1}$ such that $m_{i-1}+m_i\le d_i$ (we set $m_0=m_k=0$),
and consider the subvariety of all tuples $\{f_1\ldots,f_{k-1}\}$
as above such that $\rk f_i=m_i$ and $f_{i-1}\circ f_i=0$ for any $i$.
These tuples form a single $L$-orbit~$\O$ called a variety of complexes.
For any tuple $\{f_1,\ldots,f_{k-1}\}\in\O$
consider the tuple $\{f_1^+,\ldots,f_{k-1}^+\}\in\g_{-1}$,
where $f_i^+$ is a classical ``matrix'' Moore--Penrose inverse of $f_i$.
It is easy to see that this new tuple is again a complex,
moreover, this complex is a Moore--Penrose inverse
(in our latest meaning of this word)
of an original complex.
In particular, orbits of complexes are Moore--Penrose orbits.
\end{Example}

From the first glance only few parabolic subgroups are Moore--Penrose.
But this is scarcely true.
For example, we have the following Theorem:

\begin{Theorem}[\cite{T4}]\label{KHJJGHFDJGFJ}
Any parabolic subgroup in $SL_n$ is Moore--Penrose.
\end{Theorem}

To explain our interest in Moore--Penrose parabolic subgroups 
let us recall Conjecture~\ref{KJFKUTYFKLUFTL} from the previous section.
Suppose once again that $G$ is a simple connected simply-connected
Lie group, $P$ is its parabolic subgroup, $\p\subset\g$
are their Lie algebras.
We take any irreducible $G$-module $V$.
There exists a unique proper $P$-submodule $M_V$ of $V$.
We have the inclusion $i:\,M_V\to V$, the projection
$\pi:\,V\to V/M_V$ and the map $R_V:\,\g\to\End(V)$ defining the representation.
Therefore we have a linear map $\tilde R_V:\,\g\to\Hom(M_V,V/M_V)$, namely
$\tilde R_V(x)=\pi\circ R_V(x)\circ i$. Clearly $\p\subset\Ker\tilde R_V$.
Finally, we have a linear map $\Psi_V:\,\g/\p\to\Hom(M_V,V/M_V)$.
The Conjecture~\ref{KJFKUTYFKLUFTL} states that
there exists an algebraic stratification 
$\g/\p=\mathop{\sqcup}\limits_{i=1}^nX_i$
such that for any $V$ the function $\rk\,\Psi_V$ is constant
along each $X_i$.

The following theorem shows the connection of this problem
with the Moore--Penrose inverse.

\begin{Theorem}[\cite{T4}]\label{JHGJCJGFJGJ}
Suppose that 
a grading 
$\g=\mathop{\oplus}\limits_{k\in\Z}\g_k$ is a Moore--Penrose grading
in all positive degrees except at most one.
Then the Conjecture is true for the corresponding parabolic subgroup $P$.
\end{Theorem}

\proof
Suppose that $G$ is a connected reductive group with a Lie algebra $\g$.
For any elements $x_1,\ldots,x_r\in\g$ let 
$\langle x_1,\ldots,x_r\rangle_{alg}$ denote
the minimal algebraic Lie subalgebra of $\g$ that contains
$x_1,\ldots,x_r$
(algebraic subalgebras are the Lie algebras
of algebraic subgroups). 
By a theorem of Richardson \cite{Ri1}
$\langle x_1,\ldots,x_r\rangle_{alg}$ is reductive if and only
if an orbit of the $r$-tuple $(x_1,\ldots,x_r)$ in $\g^r$
is closed with respect to the diagonal action of~$G$.
Suppose now that $h_1,\ldots,h_r$ are semi-simple elements
of $\g$. Consider the closed subvariety 
$\hat\O=(\Ad(G)h_1,\ldots,\Ad(G)h_r)\subset\g^r$.
For any closed $G$-orbit $\O\subset\hat\O$ let
us denote by $G(\O)$ the conjugacy class of the reductive subalgebra
$\langle x_1,\ldots,x_r\rangle_{alg}$ for $(x_1,\ldots,x_r)\in\O$.

\begin{Lemma}
There are only finitely many conjugacy classes $G(\O)$.
\end{Lemma}

\proof
We shall use induction on $\dim\g$. Suppose that the claim of the Lemma
is true for all reductive groups $H$ with $\dim H<\dim G$.
Let $\z\subset\g$ be the center of $\g$, $\g'\subset\g$ be its
derived algebra. Consider two canonical homomorphisms
$$\g\mathop{\longrightarrow}\limits^\pi\g'\quad\text{and}\quad
\g\mathop{\longrightarrow}\limits^{\pi'}\z.$$
We take any closed $G$-orbit $\O\subset\hat\O$.
Let $(x_1,\ldots,x_r)\in\O$, $y_i=\pi(x_i)$ for $i=1,\ldots,r$.
Then 
$\langle y_1,\ldots,y_r\rangle_{alg}=
\pi(
\langle x_1,\ldots,x_r\rangle_{alg})$ and, therefore, is reductive.
Let us consider two cases.

Suppose first, that 
$\langle y_1,\ldots,y_r\rangle_{alg}=\g'$.
Then $\g'$ is a derived algebra of
$\langle x_1,\ldots,x_r\rangle_{alg}$
and, therefore, 
$\langle x_1,\ldots,x_r\rangle_{alg}=
\langle \pi'(h_1),\ldots,\pi'(h_r)\rangle\oplus\g'$.
In this case we get one conjugacy class.

Suppose now, that $\langle y_1,\ldots,y_r\rangle_{alg}\ne\g'$.
Then $\langle y_1,\ldots,y_r\rangle_{alg}$ is contained in some maximal
reductive Lie subalgebra of $\g'$. It is well-known (and not difficult
to prove) that in a semisimple Lie algebra there are
only finitely many conjugacy classes of maximal reductive
subalgebras. Let $\h'$ be one of them, $\h=\z\oplus\h'\subset\g$.
Let $H$ be a corresponding reductive subgroup of $G$.
It is sufficient to prove that for any closed $G$-orbit $\O$ of 
$\hat\O$ that meets $\h^r$ there are only finitely many 
possibilities for $G(\O)$. 
It easily follows from Richardson's Lemma \cite{Ri}
that for any $i$ the intersection $\Ad(G)h_i\cap\h$ is a union of finitely many closed
$H$-orbits, say $\Ad(H)h_i^1,\ldots,\Ad(H)h_i^{s_i}$.
It remains to prove that if for some $r$-tuple
$(x_1,\ldots,x_r)\in\Ad(H)h_1^{k_1}\times\ldots\Ad(H)h_r^{k_r}$
the corresponding subalgebra
$\langle x_1,\ldots,x_r\rangle_{alg}$ is reductive
then there are only finitely many possibilities 
for its conjugacy class. But this is precisely the
claim of Lemma for the group~H, which is true by the induction
hypothesis.
\endproof

Suppose that $\k$ is a compact real form of $\g$.

\begin{Lemma}
If $r$-tuple $(x_1,\ldots,x_r)$ belongs to $(i\k)^r$,
then its $G$-orbit is closed in $\g^r$.
\end{Lemma}

\proof
Indeed, let $B$ be a non-degenerate $\ad$-invariant
scalar product on $\g$, which is negative-definite on $\k$.
Let $H(x)=-B(\overline x, x)$ be a positive-definite
$\k$-invariant Hermitian quadratic form on $\g$,
where the complex conjugation is taken with respect to $\k$.
Let $H^r$ be a corresponding Hermitian quadratic form on $\g^r$.
More precisely, $H^r(x_1,\ldots,x_r)=
H(x_1)+\ldots+H(x_r)$.
By a Kempf--Ness criterion \cite{PV}
in order to prove that the $G$-orbit 
of $(x_1,\ldots,x_r)$ is closed it is sufficient to prove 
that the real function $H^r(\cdot)$ has a critical point
on this orbit. Let us show that $(x_1,\ldots,x_r)$
is this critical point.
Indeed, for any $g\in\g$
$$-B(\overline x_1,[g,x_1])-\ldots-B(\overline x_r,[g,x_r])=
B(x_1,[g,x_1])+\ldots+B(x_r,[g,x_r])=0.$$
\endproof

Now let $G$ be a simple simply-connected Lie group,
let $\g$ be its Lie algebra with a $\Z$-grading
$\g=\mathop{\oplus}\limits_{k\in\Z}\g_k$.
Let $r$ be a maximal integer such that $\g_r\ne0$.
We denote the non-positive part of the grading
$\mathop{\oplus}\limits_{k\le0}\g_k$ by $\p$. Let $P\subset G$
be a parabolic subgroup with the Lie algebra $\p$.
We shall identify $\g/\p$ with $\mathop{\oplus}\limits_{k>0}\g_k$.
Let $L\subset G$ be a connected reductive subgroup with Lie algebra~$\g_0$.
Let $V$ be an irreducible $G$-module.
If we choose a Cartan subalgebra $\t\subset\g_0$
then the grading of $\g$ originates from some $\Z$-grading
on $\t^*$. Therefore, there exists a 
$\Z$-grading $V=\mathop{\oplus}\limits_{k\in\Z}V_k$
such that $\g_iV_j\subset V_{i+j}$.
Let $R$ be a maximal integer such that $V_R\ne0$.
It is easy to see that $M_V=\mathop{\oplus}\limits_{k<R}V_k$
(notice that $V_R$ is an irreducible $L$-module).

Now we can prove the Theorem.
It is sufficient to prove that there exists a finite set of points
$\{x_1,\ldots,x_N\}\subset\g/\p$
such that for any $x\in\g/\p$ and for any $V$
we have $\rk\,\Psi_V(x)=\rk\,\Psi_V(x_i)$ for some $i$.
Recall that $L$ has finitely many orbits on each $\g_k$, 
see Theorem~\ref{JHGFJHGFJHG}. 
We pick some $L$-orbit $\O_i$ in each $\g_i$.
Then it is sufficient to find a finite set of points as above
only for points $x\in\g/\p$ of a form $x=x_1+\ldots+x_r$, where $x_i\in\O_i$.
For any orbit $\O_i$ let ${\cal H}_i$ denote the set of all possible
homogeneous characteristics of all elements from $\O_i$.
Clearly ${\cal H}_i$ is a closed $\Ad(L)$-orbit.
Let $\hat\O={\cal H}_1\times\ldots\times{\cal H}_r\subset\g_0^r$.
Then by the first Lemma the set of conjugacy classes of subgroups $G(\O)$
for closed $L$-orbits $\O$ in $\hat\O$ is finite.
Let us show that for any $r$-tuple $(x_1,\ldots,x_r)\in\O_1\times\ldots\times\O_r$
there exists an $r$-tuple $(h_1,\ldots,h_r)\in\hat\O$
such that $h_i$ is a homogeneous characteristic of $x_i$
and an $L$-orbit $\Ad(L)(h_1,\ldots,h_r)$ is closed.
Indeed, after simultaneous conjugation of elements $x_i$ by some 
element $g\in L$ we may suppose that any $x_i$ has a homogeneous
characteristic $h_i\in i\k_0$ (in all degrees except at most one
no conjugation is needed because of Moore--Penrose property,
for one degree this is obvious). Then by the second Lemma
an orbit $\Ad(L)(h_1,\ldots,h_r)$ is closed.
Since all functions $\rk\,\Psi_V$ are $L$-invariant,
we may restrict ourselves to the points $x=\sum_i x_i\in\g/\p$
such that $x_i\in\O_i$, any $x_i$ has a homogeneous characteristic
$h_i\in i\k_0$, and a conjugacy class of $\langle h_1,\ldots,h_r\rangle_{alg}$
is fixed. We claim that any function $\rk\,\Psi_V$
is constant along the set of these points. Moreover,
we shall prove that
\begin{equation}
\rk\,\Psi_V(x)=\dim V_R-\dim V_R^{\langle h_1,\ldots,h_r\rangle_{alg}}.
\label{KJHVJGHKGHLHG}
\end{equation}
Indeed, 
$$\rk\,\Psi_V(x)=\dim\sum_i\Im\left(\ad(x_i)|_{V_{R-i}}\right).$$
Clearly $V_R$ is $\ad(h_i)$-invariant and is killed by $\ad(e_i)$, therefore
from the $\ssl_2$-theory we get that
$V_R=\mathop{\oplus}\limits_{k\ge0}V_R^k$, where $\ad(h_i)|_{V_R^k}=k\cdot\Id$.
Moreover, 
$$\Im\left(\ad(x_i)|_{V_{R-i}}\right)=\mathop{\oplus}\limits_{k>0}V_R^k.$$
Let $H$ be a contravariant Hermitian form on $V_R$ with respect to the compact form
$\k_0$ of $\g_0$.
Since $h_i\in i\k_0$ and $H$ is a contravariant form we get that
$\mathop{\oplus}\limits_{k>0}V_R^k=(V_R^0)^\perp$.
Therefore,
$$\sum_i\Im\left(\ad(x_i)|_{V_{R-i}}\right)=\left(\cap_i V_R^{h_i}\right)^\perp=
\left(V_R^{\langle h_1,\ldots,h_r\rangle_{alg}}\right)^\perp.$$
The formula~(\ref{KJHVJGHKGHLHG}) follows.
\endproof

For example, combining Theorem~\ref{JHGJCJGFJGJ} and
Theorem~\ref{KHJJGHFDJGFJ}
we get the following corollary:

\begin{Corollary}
The conjecture is true for any parabolic subgroup in $\SL_n$.
\end{Corollary}

Though the conditions of Theorem~\ref{JHGJCJGFJGJ}
are not always satisfied,
it seems that 
one can prove the conjecture for any simple group
in this direction.

\chapter{Fulton--Hansen Theorem and Applications}

\section*{Preliminaries}
Fulton--Hansen connectedness Theorem has numerous applications
to the projective geometry including famous Zak Theorem on tangencies.
Aside from dual varieties, projective geometry also studies
secant and tangential varieties, joins, etc. In this chapter
we review some basic results related to them.

\section{Fulton--Hansen Connectedness Theorem}
\index{Fulton--Hansen connectedness theorem}

\begin{Theorem}[\cite{FH}]\label{FHTheorem1}
Let $f:\,X\to\P^N$ be a finite morphism from an irreducible complete 
variety $X$ to the projective space $\P^N$.
Then for any projective subspace $L_0\subset\P^N$ such that
$\dim X>\codim_{\P^N}L_0$
the inverse image $f^{-1}(L_0)$ is connected.
\end{Theorem}

\sketch
Let $k=\dim L$. Let $U\subset\Gr(k,\P^N)$ be the open set of all subspaces
that meet $f(X)$ properly, and let $U_0\subset\Gr(k,\P^N)$
be the set of all subspaces $L$ such that $f^{-1}(L)$ is irreducible.
Then $U_0$ is non-empty by Bertini's theorem (see~\cite{Ha1} or \cite{Wei}).
It follows from Zariski's principle of degeneration (see~\cite{Fu1} or \cite{Za})
that for all $L\in U$ the preimage 
$f^{-1}(L)$ is connected, non-empty, of codimension $\dim X+k-N$.
For each $L\in U$, there is a positive cycle $f^*(L)$ on $X$
whose support is $f^{-1}(L)$. This determines a morphism from $U$
to the Chow variety $Z$ of cycles on $X$. 
Let 
$$\Gamma\subset\Gr(k,\P^N)\times Z$$
be the closure of the graph of this morphism, 
$\pi:\,\Gamma\to\Gr(k,\P^N)$ the projection.
Since $Z$ is complete, $\pi$ is proper as well as birational.
Zariski's principle of degeneration implies that for any $\gamma\in\Gamma$
the corresponding cycle $C_\gamma$ is connected 
(i.e.~its support $|C_\gamma|$ is connected).
Zariski's main theorem implies that $\pi^{-1}(L_0)$ is connected.
It follows that the union of the $|C_\gamma|$ for $\gamma\in\pi^{-1}(L_0)$
is connected. It remains to prove that
$$f^{-1}(L_0)=\mathop{\bigcup}\limits_{\gamma\in\pi^{-1}(L_0)}|C_\gamma|.$$
The inclusion $\supset$ follows from the continuity of limit cycles.
So it remains to prove the inclusion $\subset$.
Let $x\in f^{-1}(L_0)$ and let $z=f(x)$.
Consider the subvariety $G$ of $\Gr(k,\P^N)$ consisting of all subspaces
passing through $z$ (therefore $G$ is isomorphic to $\Gr(k-1,\P^{N-1})$.
An easy dimension count shows that $G$ meets $U$.
Let $\Gamma'\subset\Gamma$ be an irreducible subvariety that maps
birationally to~$G$.
For each $L\subset G\cap U$, the cycle $f^*(L)$ contains $x$,
therefore, for any $\gamma\in\Gamma'$ the cycle $C_\gamma$ contains $x$.
In particular, there exists the limit cycle 
$C_\gamma$ for $\gamma\in\Gamma'\cap\pi^{-1}(L_0)$ that
contains $x$.
\endproof

Now we shall derive a number of corollaries.

\begin{Corollary}[\cite{FH}]
Let $f:\,X\to\P^N$ be a morphism from an irreducible complete 
variety $X$.
Then for any projective subspace $L\subset\P^N$ such that
$\dim f(X)>\codim_{\P^N}L$
the inverse image $f^{-1}(L)$ is connected.
\end{Corollary}

\proof
This follows from Theorem~\ref{FHTheorem1}
using the Stein factorization (see e.g.~\cite{Ha1}).
Indeed, for a non-finite $f$, let $f=g\circ h$ be its Stein factorization,
i.e.~$h:\,X\to X'$ has connected fibers, and $g:\,X'\to \P^N$
is finite. Then $\dim X'=\dim f(X)$,
therefore $g^{-1}(L)$ is connected by
Theorem~\ref{FHTheorem1}, therefore $f^{-1}(L)$ is also connected.
\endproof

\begin{Theorem}[\cite{FH}]\label{FHTheorem2}
Let $\P=\P^N\times\ldots\times\P^N$ {\rm(}$r$ copies{\rm)},
and let $\Delta$ be the diagonal.
If $X$ is an irreducible variety, $f:\,X\to\P$ a morphism,
with $$\dim f(X)>\codim_\P\Delta=(r-1)N,$$ 
then $f^{-1}(\Delta)$ is connected.
\end{Theorem}

\proof
Identify $\P$ as the set of ordered $r$-tuples of points of $\P^N$,
so that a point in $\P$ has $r$ factors in $\P^N$; the diagonal $\Delta$
consists of those points whose factors are all equal.

For simplicity suppose that $f$ is a finite morphism,
a general case is handled with a  help of the Stein factorization
as in the proof of the previous Corollary.
Suppose that $f^{-1}(\Delta)$ is not connected, choose
points $x$ and $x'$ in different connected components of $f^{-1}(\Delta)$.
Choose a hyperplane $H$ in $\P^N$ not containing
any of the $r$ factors of $f(x)$ or $f(x')$.
Let $\P_0$ be the open set of $\P$ consisting 
of those points none of whose factors lie in $H$.
Then $\P_0$ is a product of $r$ copies of $\C^N$,
which is identified with $\C^{rN}$ and therefore with
$\P^{rN}\setminus\{z_0=0\}$, where $z_0,\ldots,z_{rN}$
are the homogeneous coordinates on $\P^{rN}$.
Let $W\subset\P\times\P^{rN}$ be the closure of the graph
of this birational correspondence.
If coordinates are chosen on $\P^N$
such that $H$ is the hyperplane $x_0=0$,
and the homogeneous coordinates on the $k$-th
copy of $\P^N$ are $x_0^k,\ldots,x_N^k$,
then a point $(x^1)\times\ldots\times(x^r)\times(z)$
in $\P\times\P^{rN}$ belongs to $W$ if and only if
there are constants $\lambda_1,\ldots,\lambda_r$
such that
$$(z_0,z_{(k-1)N+1},\ldots,z_{kN})=\lambda_k(x_0^k,\ldots,x_N^k)$$
for all $k=1,\ldots,r$.
So $W$ is defined by the equations
$$\cases{
x_j^kz_0=x_0^kz_{(k-1)N+j} 
&\hbox{$1\le k\le r$, $1\le j\le N$}\cr
x_i^kz_{(k-1)N+j}=x_j^kz_{(k-1)N+i} 
&\hbox{$1\le k\le r$, $1\le i<j\le N$.}}$$
Let $\alpha:\,W\to\P$, $\beta:\,W\to\P^{rN}$ be the two
birational projections.
Let $L$ be the $N$-dimensional linear subspace of $\P^{rN}$
defined by the equations
$$z_{(k-1)N+j}=z_{kN+j},\quad 1\le k\le r-1,\ 1\le j\le N.$$
Then it is easy to check that
\begin{equation}
\alpha^{-1}(\Delta)\supset\beta^{-1}(L)\hbox{\rm\ and\ }
\alpha^{-1}(\Delta\cap\P_0)\subset \beta^{-1}(L). \label{JHGFJKHTDJ}
\end{equation}
Let $\tilde X$ be the irreducible component of the fiber product
$X\times_{\P}W$ that maps onto $X$.
Then there is a diagram
$$
\begin{CD}
\tilde X @>g>> W@>{\beta}>>\P^{rN}\\
@V{\pi}VV @VV{\alpha}V\\
X @>>f>\P
\end{CD}
$$
where the square commutes, and $\pi$ is surjective.
Then equations (\ref{JHGFJKHTDJ}) imply that
\begin{equation}
f^{-1}(\Delta)\supset\pi(g^{-1}\beta^{-1}(L)),\hbox{\rm \ and \ }
f^{-1}(\Delta\cap\P_0)\subset\pi(g^{-1}\beta^{-1}(L)).\label{KJHFJKHFDX}
\end{equation}
Now $\dim\beta g(\tilde X)=\dim f(X)$ since $f(X)$ meets the open set $\P_0$
which corresponds isomorphically to an open set in $\P^{rN}$.
Then by Theorem~\ref{FHTheorem1} we know that $(\beta g)^{-1}(L)$
is connected. Therefore $\pi((\beta g)^{-1}(L))$ is connected.
But by (\ref{KJHFJKHFDX}) 
$\pi((\beta g)^{-1}(L))$ is then a connected subset of $f^{-1}(\Delta)$
which contains both $x$ and $x'$.
Contradiction.
\endproof

\begin{Corollary}[\cite{FH}]$\ $
\begin{enumerate}
\item[\rm(a)] Let $Y$ be an irreducible subvariety
of the projective space $\P^N$
and let $f:\,X\to\P^N$ be a morphism from an irreducible variety $X$
to $\P^N$ with $\dim f(X)>\codim_{\P^N}Y$. Then $f^{-1}(Y)$ is connected.
\item[\rm(b)] 
If $X$ and $Y$ are irreducible subvarieties
of $\P^N$ with $\dim X+\dim Y>N$ then $X\cap Y$ is connected.
\end{enumerate}
\end{Corollary}

\proof (a)
We apply Theorem~\ref{FHTheorem2} for the map $F:\,X\times Y\to \P^N\times\P^N$
given by $F=(f,\Id)$.

(b) Follows from (a).
\endproof

\begin{Corollary}[\cite{FH}]\label{FHTheorem3} $\ $
\begin{enumerate}
\item[\rm(a)]
Let $X$ be an irreducible variety of dimension $n$,
$f:\,X\to \P^N$ an unramified morphism, with $N<2n$.
Then $f$ is a closed embedding.
\item[\rm(b)] Let $X$ be an irreducible $n$-dimensional subvariety
of $\P^N$ with $N<2n$. Then the algebraic fundamental
group of $X$ is trivial {\rm(}i.e. $X$ has no non-trivial \'etale coverings{\rm)}.
\end{enumerate}
\end{Corollary}

\proof
(a) If $f:\,X\to\P^N$ is unramified, consider
the product mapping $f\times f:\,X\times X\to\P^N\times\P^N$.
The unramified assumption is equivalent to
the assertion that the diagonal $\Delta_X$ in $X\times X$
is open as well as closed in $(f\times f)^{-1}(\Delta_{\P^N})$.
Since $f$ is automatically finite, $\dim(f\times f)(X\times X)=2n$,
so $(f\times f)^{-1}(\Delta_{\P^N})$ is connected by the 
Theorem~\ref{FHTheorem2}. Hence 
$$\Delta_X=(f\times f)^{-1}(\Delta_{\P^N}),$$
so $f$ is one-to-one, and therefore a closed embedding.

(b) Follows from (a).
\endproof

\section{Secant and Tangential Varieties}

\begin{Definition}\rm
Suppose that $X\subset\P^N$ is a smooth $n$-dimensional projective variety.
The tangential variety $\Tan(X)$ is the union of all (embedded) tangent
spaces to $X$\index{tangential variety}
$$\Tan(X)=\mathop{\cup}_{x\in X}\hat T_xX.$$\index{secant variety}
The secant variety $\Sec(X)$ is the closure of the union
of all secant lines to $X$
$$\Sec(X)=\overline{\mathop{\cup}_{x,y\in X}\P^1_{xy}},$$
where $\P^1_{xy}$ is a line connecting $x$ and $y$.
\end{Definition}

Clearly, $\Tan(X)$ and $\Sec(X)$ are closed irreducible
subvarieties of $\P^N$,
with $\Tan(X)\subset\Sec(X)$, and 
$$\dim\Tan(X)\le 2n,\quad \dim\Sec(X)\le 2n+1.$$

\begin{Example}\rm\index{Segre embedding}
Let $V=\C^{k+1}\otimes\C^{l+1}$ denote the space of $(k+1)\times(l+1)$
matrices. Let $X\subset\P(V)$ be the projectivization
of the variety of rank $1$ matrices. Then $X$ is isomorphic 
to $\P^k\times\P^l$ embedded into $\P(V)$ via the Segre embedding.
It is easy to see that $\Sec(X)$ is equal to the projectivization
of the variety of matrices of rank at most $2$.
\end{Example}

\begin{Theorem}[\cite{FH}, \cite{Z2}]\label{LKGJGDFHGFX}
Suppose that $X\subset\P^N$ is a smooth projective variety.
Then either 
$$\dim\Tan(X)=2n,\quad\dim\Sec(X)=2n+1$$ 
or
$$\Tan(X)=\Sec(X).$$
\end{Theorem}

\proof
Indeed, suppose that $\Tan(X)\ne\Sec(X)$.
Since $\Sec(X)$ is irreducible, it is sufficient
to prove that in this case $\dim\Tan(X)=2n$.
Suppose, on the contrary, that $\dim\Tan(X)<2n$.
Let $L\subset\P^N$ be a linear subspace such that $L\cap\Tan(X)=\emptyset$
and $\dim L=N-1-\dim\Tan(X)$ (since $\Tan(X)\ne\Sec(X)$ we obviously
have $\dim\Tan(X)\le N-1$).
Let $\pi$ be the linear projection $\P^N\to\P^{\dim\Tan(X)}$
with center $L$. Since $L\cap\Tan(X)=\emptyset$,
the restriction of $\pi$ on $X$ is an unramified morphism.
Since $\dim\Tan(X)<2n$, it follows by Theorem~\ref{FHTheorem3}
that the restriction of $\pi$ on $X$ is a closed embedding.
Therefore, $L\cap\Sec(X)=\emptyset$.
Since $\dim L+\dim\Sec(X)\ge N$, this is a contradiction.
\endproof

It is possible to generalize this theorem in the following direction.
Suppose that $Y\subset X\subset \P^N$ are arbitrary irreducible projective
varieties. Then we can define relative secant and tangential varieties
as follows. The relative secant variety \index{relative secant variety}
$\Sec(Y,X)$ is the closure of the union of all secants
$\P^1_{x,y}$, where $x\in X$, $y\in Y$.
The relative tangential variety $\Tan(Y,X)$ is the union \index{relative tangential variety}
of tangent stars $T^\star_yX$ for $y\in Y$, where \index{tangent star}
the tangent star $T^\star_xX$ for $x\in X$ is the union of
limit positions of secants $\P^1_{x',x''}$, where $x',x''\in X$
and $x',x''\to x$.

\begin{Example}\rm
Let $V=\C^{k+1}\otimes\C^{l+1}$ denote the space of $(k+1)\times(l+1)$
matrices. Let $X^r\subset\P(V)$ be the projectivization
of the variety of matrices of rank at most $r$. 
Then $X^r\subset X^{r+1}$ and for $r<s$ we have
$$\Sec(X^r,X^s)=\cases{
X^{r+s}&\hbox{\rm if $r+s\le\min(k+1,l+1)$}\cr
\P(V)&\hbox{\rm if $r+s>\min(k+1,l+1)$.}
}$$
\end{Example}

The proof of the following theorem is analogous to the
proof of Theorem~\ref{LKGJGDFHGFX}.

\begin{Theorem}[\cite{Z2}]\label{JKHGFKDHFXFH}
Suppose that $Y\subset X\subset \P^N$ are irreducible projective
varieties.
Then either 
$$\dim\Tan(Y,X)=\dim X+\dim Y,\quad\dim\Sec(Y,X)=\dim X+\dim Y+1$$ 
or
$$\Tan(Y,X)=\Sec(Y,X).$$
\end{Theorem}

One can go even further and define the join $S(X,Y)$ of \index{join}
any two subvarieties $X,Y\subset\P(V)$,
$$S(X,Y)=\overline{
\mathop{\cup}\limits_{x\in X,\ y\in Y}\P^1_{xy}},$$
where the closure is not necessary if the two varieties
do not intersect. 
For example, if $X$ and $Y$ are projective subspaces then 
$$S(X,Y)=\P(\Cone(X)+\Cone(Y)).$$
An important result about joins is the following
Terracini lemma.\index{Terracini lemma}

\begin{Lemma}[\cite{Te}]
Let $X,Y\subset\P(V)$ be irreducible varieties and let 
$$x\in X,\quad y\in Y,\quad z\in\P^1_{xy}.$$
Then
$$\hat T_zS(X,Y)\supset S(\hat T_xX,\hat T_yY).$$
Moreover, if $z$ is a generic point of $S(X,Y)$ then equality holds.
\end{Lemma}

If an algebraic group $G$ acts on $\P^N$
with finitely many orbits, then for any two orbits $\O_1$ and $\O_2$
the join of their closures $S(\overline{\O_1},\overline{\O_2})$
is $G$-invariant therefore coincides with some orbit closure $\overline{\O}$.
This gives an interesting associative product on the set of $G$-orbits.

\begin{Example}\label{KJHGFKGHDFK}\rm
Let $L$ be a simple algebraic group and $P$ a parabolic subgroup
with abelian unipotent radical. In this case $\l=\Lie L$ admits
a short $\Z$-grading with only three non-zero parts:
$$\l=\l_{-1}\oplus\l_0\oplus\l_1.$$
Here $\l_0\oplus\l_1=\Lie P$
and $\exp(\l_1)$ is the abelian unipotent radical of $P$.
Let $G\subset L$ be a reductive subgroup with Lie algebra
$\l_0$. 
Recall that by Theorem~\ref{AURAOrbitDecomp}
$G$ has finitely many orbits in $\l_1$ naturally 
labeled by integers
from the segment $[0,r]$.
Let $X_i=\overline{\P(\O_i)}$, $i=0,\ldots,r$.
Then it is easy to see that 
$$S(X_i,X_j)=X_{\min(i+j,r)}.$$
Indeed, by associativity and induction it suffices to prove that
$S(X_i,X_1)=X_{\min(i+1,r)}$
and this easily follows from the results of Section~\ref{ParabolicsAura}.
\end{Example}

\section{Zak Theorems}

The following theorem is called Zak Theorem on 
tangencies.\index{Zak's theorem on tangencies}

\begin{Theorem}[\cite{Z2}, \cite{FL}]\label{ZakThTangencies}$\ $
\begin{enumerate}
\item[\rm(a)] 
Suppose that $X\subset\P^N$ is a smooth nondegenerate projective variety, $\dim X=n$.
If $L$ is a $k$-plane in $\P^N$ and $k\ge n$
then $$\dim\Sing(L\cap X)\le k-n.$$
\item[\rm(b)] If $X\subset\P^N$ is a non-linear smooth projective
variety, 
then $\dim\dual X\ge\dual X$.
Moreover, if $\dual X$ is smooth, then $\dim X=\dim\dual X$.
\end{enumerate}
\end{Theorem}

\proof
(a) Indeed, let $Y=\Sing L\cap X$. Then $x\in Y$ if and only if
$\hat T_xX\subset L$. Therefore, $\Tan(Y,X)\subset L$.
In particular, $\dim\Tan(Y,X)\le k$. On the over hand,
since $X$ is non-degenerate, $X\not\subset L$. Therefore $\Sec(Y,X)\not\subset L$.
In particular, $\Tan(Y,X)\ne\Sec(Y,X)$.
By Theorem~\ref{JKHGFKDHFXFH} it follows that
$\dim\Tan(Y,X)=\dim Y+n$. Finally, we have $\dim Y\le k-n$.

(b) Let $H\subset\dual X_{sm}$. Then the contact locus $\Sing X\cap H$
is the projective subspace of dimension $\defect X$.
Therefore, by (a) we have 
$$\defect X\le \dim H-\dim X=\codim X-1.$$
It follows that $\dim X\le\dim\dual X$.
If $\dual X$ is also smooth then using the reflexivity theorem
and the same argument as above we get $\dim\dual X\le\dim X$,
therefore $\dim X=\dim\dual X$.
\endproof

To remove the condition of smoothness in 
Theorem~\ref{ZakThTangencies} (a) we need the notion of Gauss images.
Let $X\subset\P^N$ be an irreducible $n$-dimensional
projective variety. \index{Gauss image}

\begin{Definition}
For any $k\ge n$ the variety
$$\gamma_k(X)=\overline{
\cup_{x\in X_{sm}}\{L\subset\Gr(k,\P^N)\,|\,L\supset\hat T_xX\}
}\subset\Gr(k,\P^N)$$
is called the $k$-th Gauss image of $X$.
\end{Definition}
It is easy to see that Gauss images are irreducible varieties.
For example, $\gamma_{N-1}(X)\subset\Gr(N-1,\P^N)=\dual{\P^N}$ 
is equal to the dual variety $\dual X$.

\begin{Example}\rm
One might expect that usually $\dim X=\dim\gamma_n(X)$.
Indeed, suppose that $\dim X>\dim\gamma_n(X)$.
Then a generic fiber $F$ of the map $X_{sm}\to\gamma_n(X)$, $x\mapsto\hat T_xX$
has positive dimension $\dim X-\dim\gamma_n(X)$.
Therefore, if the hyperplane $H$ is tangent to $X_{sm}$ at some point $x\in F$
then $H$ is tangent to $X$ at all points of $F$.
It easily follows that $\defect X>\dim F=\dim X-\dim\gamma_n(X)$.
So $X$ has a positive defect in this case.
\end{Example}

The proof of the following theorem is analogous to the proof of
Theorem~\ref{ZakThTangencies} (a).
\begin{Theorem}[\cite{Z2}]\label{ZakThTangenciesStrong}
Suppose that $X\subset\P^N$ is a nondegenerate projective variety. 
If $L\subset\gamma_k(X)$ 
then 
$$\dim\{x\in X\,\hat T_xX\subset L\}\le k-n+(\dim\Sing X+1),$$
where we set $\dim\Sing X=-1$ if $X$ is smooth.
\end{Theorem}

The following theorem is called Zak's Theorem on 
linear normality.\index{Zak's theorem on linear normality}

\begin{Theorem}[\cite{Z4}]\label{ZakThLinNorm}
Suppose that $X\subset\P^N$ is a smooth non-degenerate
projective variety, $\dim X=n$.
If 
$$\codim X<{N+4\over 3}$$
then $\Sec(X)=\P^N$.
\end{Theorem}

\proof
By Theorem~\ref{LKGJGDFHGFX} either $\Tan(X)=\Sec(X)$ or
$\dim\Sec(X)=2n+1$. Suppose that $\Sec(X)\ne\P^N$ and
$\dim\Sec(X)=2n+1$. Then $2n+1<N$, which contradicts
$\codim X<{N+4\over 3}$.

Suppose now that $\Tan(X)=\Sec(X)$ and $\Sec(X)\ne\P^N$.
Let $z\in\Sec(X)$ be a generic point, $L=\hat T_z\Sec(X)$,
and $Q=\{x\in X\,|\,z\in\hat T_xX\}$.
Then $Q$ is closed and an easy dimension count shows that $\dim Q=2n-\dim\Sec(X)$.
For any point $x\in Q$ we have $z\in\hat T_xX\subset\Sec(X)$
therefore $\hat T_xX\subset\hat T_z\Sec(X)=L$.
It follows that $\Tan(Q,X)\subset L$.
However, $\Sec(Q,X)\not\subset L$ because $L$ is a proper subspace
and $X$ is not degenerate.
By Theorem~\ref{JKHGFKDHFXFH} we have
$$\dim\Sec(Q,X)=\dim Q+\dim X+1=3n-\dim\Sec(X)+1.$$
Since $\Sec(Q,X)\subset\Sec(X)$, we have 
$\dim\Sec(Q,X)\le\dim\Sec(X)$.
Then 
$$2\dim\Sec X\ge 3n+1$$ Hence $2N>3n+1$.
But this contradicts $\codim X<{N+4\over 3}$.
\endproof

Theorem~\ref{ZakThLinNorm} is called the theorem on linear normality
due to the following corollary first conjectured in \cite{Ha2}.
\index{Hartshorne conjecture}

\begin{Corollary}
Suppose that $X\subset\P^N$ is a smooth non-degenerate
projective variety.
If 
$$\codim X<{N+2\over 3}$$
then $X$ is linearly normal.
\end{Corollary}

\proof
Recall from Section~\ref{LinearNormSection} that
if $X$ is not linearly normal then there exists a subvariety
$X'\subset\P^{N+1}$ and a point $p\in\P^{N+1}$
such that $X$ is isomorphic to $X'$ via the linear projection
with center $p$. In particular, $\Sec(X')\ne\P^{N+1}$.
But this contradicts Theorem~\ref{ZakThLinNorm}.
\endproof

Zak \cite{Z2} has also classified the varieties in the borderline 
case\index{Zak's theorem on Severi varieties}\index{Severi variety}
$$\codim X={N+4\over 3}.$$

\begin{Theorem}[\cite{Z4}]\label{ZakThOnSeveriVarieties}
Suppose that $X\subset\P^N$ is a smooth non-degenerate
$n$-dimensional projective variety such that
$\codim X={N+4\over 3}$. Then $\Sec(X)=\P^N$ except
the $4$ following cases:
\begin{itemize}
\item $n=2$, $X=\P^2$, $X\subset\P^5$ is the Veronese embedding.
\item $n=4$, $X=\P^2\times\P^2$, $X\subset\P^8$ is the Segre embedding.
\item $n=8$, $X=\Gr(2,\C^6)$, $X\subset\P^{14}$ is the Pl\"ucker embedding.
\item $n=16$, $X\subset\P^{26}$ is the projectivization of 
the highest weight vector orbit in the $27$-dimensional irreducible 
representation of $E_6$.
\end{itemize}
\end{Theorem}

The varieties listed in Theorem~\ref{ZakThOnSeveriVarieties}
are called Severi varieties after Severi who proved 
that the unique $2$-dimensional Severi variety is the Veronese surface.
Scorza and Fujita--Roberts \cite{FR} have shown that
the unique $4$-dimensional Severi variety is the 
Segre embedding of $\P^2\times\P^2$.
It is also shown in \cite{FR} that $n\equiv 0\mod16$ as $n>8$.
Finally, Tango \cite{Tan} proved that if $n>16$ then either $n=2^a$
or $n=3\cdot2^b$, where $a\ge7$ and $b\ge5$.
The $16$-dimensional Severi variety was discovered by Lazarsfeld~\cite{La}.
In fact, all Severi varieties (1)--(4)
arise from Example~\ref{KJHGFKGHDFK}.
Another description of Severi varieties 
can be found in Example~\ref{KJHVKJNGVKJNGVK}.
\index{Severi variety}

\chapter{Dual Varieties and Projective Differential Geometry}

\section*{Preliminaries}
In sections~\ref{HJGFJHGFJHGFJHF} and \ref{KHFJGFKGTDFIYD}
we review applications of the Katz dimension formula
expressing the codimension of the dual variety
in terms of the rank of a certain Hessian matrix.
The section~\ref{JHGVJHGVJGHFKGFKGV} contains 
results of Ein about the contact locus. 
In the last section~\ref{ProjectiveIISection} 
we give a useful formalism for dealing 
with dual varieties ``under the microscope''.

\section{The Katz Dimension Formula}\label{HJGFJHGFJHGFJHF}

Let $X\subset\P^n=\P(V)$ be an irreducible 
$k$-dimensional projective variety
and $\dual X\subset\dual{\P^n}=\P(V^*)$ 
be the projectively dual variety.
The Katz dimension formula expresses the \index{Katz dimension formula}
defect $\defect(X)=\codim\dual X-1$
in terms of the rank of a certain Hessian matrix. \index{Hessian matrix}

For $x\in\P(V)$ let $x^\perp\subset V^*$ denote the annihilator of the 
line in $V$ corresponding to $x$.
Let $x_0$ be a smooth point of $X$.
Then one can choose linear functionals 
$$T_0\in V^*\setminus x_0^\perp,\quad T_1,\ldots,T_k\in x_0^\perp$$ 
so that the functions
$$t_1=T_1/T_0,\ t_2=T_2/T_0, \ldots, t_k=T_k/T_0$$ 
are local coordinates
on $X$ in the neighborhood of $x_0$.
For every $U\in x_0^\perp$ the function $u=U/T^0$ on $X$ near $x_0$
is an analytic function of $t_1,\ldots,t_k$ such that $u(0,\ldots,0)=0$.
Consider the Hessian matrix
$$\Hes(u)=\Hes(U;T_0,T_1,\ldots,T_k;x_0)={\left(
{\partial^2u\over\partial t_i\partial t_j}(0,\ldots,0)
\right)}_{i,j=1,\ldots,k}.$$

\begin{Theorem}[\cite{Ka}]\label{KatzDimension}
Let $\defect(X)=\codim\dual X-1$ be the defect of $X$. Then
\begin{enumerate}
\item[\rm(a)] We have
$$\defect(X)=\min\corank\Hes(u),$$
the minimum over all possible choices of $x_0$ and $U$.
\item[\rm(b)] $\dual X$ is a hypersurface if and only if
for some choice of $x_0$ and $u$ the Hessian matrix $\Hes(u)$ is invertible.
\end{enumerate}
\end{Theorem}

\proof
Let $I_X^0\subset\P(V)\times\P(V^*)$ be a ``smooth part''
of the conormal variety consisting of all pairs $(x,H)$
such that $x\in X$ is smooth and $H$ is tangent to $X$ at $x$.
Let $\pi:\,I_X^0\to\P(V^*)$ be the second projection.
Then $\dual X$ is the closure of $\pi(I_X^0)$.
We shall compute the Jacobian matrix of $\pi$ in appropriate
local coordinates and relate it to the Hessian matrices appearing in the Theorem.
Then we shall use an obvious fact that for any regular map $\pi:\,Z\to S$
of algebraic varieties the dimension of the image is equal
to the maximal rank of the Jacobian matrix of $\pi$
at smooth points of $Z$.

Let $(x_0,H_0)\in I_X^0$.
Then one can choose 
$$T_0\in V^*\setminus x_0^\perp,\quad T_1,\ldots,T_k\in x_0^\perp$$ 
so that the functions
$$t_1=T_1/T_0,\ t_2=T_2/T_0,\ \ldots,\ t_k=T_k/T_0$$ 
are local coordinates
on $X$ in the neighborhood of $x_0$.
Extend $T_1,\ldots,T_k$ to a basis 
$$\{T_1,\ldots,T_k,U_1,\ldots,U_{n-k}\}$$
of $x_0^\perp$. Then each of functions $u_i=U_i/T_0$ on $X$ is an analytic function
of $t_1,\ldots,t_k$ near $x_0=(0,\ldots,0)$, and we have
$$u_i(0,\ldots,0)=0\ \text{for}\ i=1,\ldots,n-k.$$
Let $x\in X$ be a point close to $x_0$ with local coordinates $(t_1,\ldots,t_k)$,
so in homogeneous coordinates 
$$x=(1:t_1:\ldots:t_k:u_1:\ldots:u_{n-k}).$$
Let $H\subset\P(V^*)$ be the hyperplane defined by an equation
$$\sum_{j=0}^k\tau_jT_j+\sum_{i=1}^{n-k}\eta_iU_i=0.$$
Then $H$ is tangent to $X$ at $x$ if and only if the function
$$\tau_0+\sum_{j=1}^k\tau_jt_j+\sum_{i=1}^{n-k}\eta_iu_i$$
vanishes at $x$ together with all its first derivatives.
Therefore, for a given~$x$, the hyperplanes tangent to $X$ at $x$
form a projective space of dimension $n-k-1$ with homogeneous
coordinates $(\eta_1:\ldots:\eta_{n-k})$, and the remaining
coordinates are given by 
$$\tau_j=-\sum_{i=1}^{n-k}\eta_i{\partial u_i\over\partial t_j},\quad j=1,\ldots,k,$$
$$\tau_0=-\left(\sum_{j=1}^k\tau_jt_j+\sum_{i=1}^{n-k}\eta_iu_i\right)=
-\sum_{i=1}^{n-k}\eta_iu_i+\sum_{j=1}^k\sum_{i=1}^{n-k}\eta_it_j{\partial u_i\over\partial t_j}.$$

Without loss of generality we may assume that the hyperplane $H_0$
has coordinates 
$$(\eta_1:\ldots:\eta_{n-k})=(0:\ldots:0:1).$$
It follows that we can set $\eta_{n-k}=1$ in above formulas, and use
$$(t_1,\ldots,t_k,\eta_1,\ldots,\eta_{n-k-1})$$ 
as local coordinates
on $I_X^0$ near $(x_0,H_0)$.
In these coordinates the projection $\pi:\,I_X^0\to\P(V^*)$ has a form
$$(t_1,\ldots,t_k,\eta_1,\ldots,\eta_{n-k-1})\mapsto
(\tau_0,\tau_1,\ldots,\tau_k,\eta_1,\ldots,\eta_{n-k-1}),$$
where $\tau_i$ are given by formulas above with $\eta_{n-k}=1$.
Therefore, the Jacobian matrix of $\pi$ at the origin is the following matrix:
$${\partial(\tau_0,\tau_1,\ldots,\tau_k,\eta_1,\ldots,\eta_{n-k-1})\over
\partial(t_1,\ldots,t_k,\eta_1,\ldots,\eta_{n-k-1})}(0,\ldots,0)=$$
$$=\left(\matrix{
0&0\cr
-\Hes(u_{n-k})&-\left({\partial(u_1,\ldots,u_{n-k})\over\partial(t_1,\ldots,t_k)}
(0,\ldots,0)\right)^t\cr
0&\Id_{n-k-1}
}\right).$$
Here vertical sizes of the blocks are $1,k,n-k-1$, and horizontal
sizes are $k$, $n-k-1$.
Clearly, this matrix has rank
$$n-k-1+\rk\Hes(u_{n-k}).$$
It follows that
$$\dim \dual X=n-k-1+\max\rk\Hes(u),$$
and the Theorem is proved.
\endproof

If $X\subset \P^n$ is a hypersurface (or a complete intersection)
then it is possible to rewrite Hessian matrices in homogeneous coordinates.
In case of hypersurfaces the corresponding result was
first formulated by B.~Segre \cite{Se}:

\begin{Theorem}\label{KJKGJFKJFK}   \index{Segre Theorem}
Let $f(x_0,\ldots,x_n)$ be an irreducible homogeneous polynomial
and let $X\subset\P^n$ be the hypersurface with the equation $f=0$.
Let $m$ be the largest number with the following property:
there exists $(m\times m)$-minor of the Hessian matrix $\left(\partial^2f/\partial x_i\partial x_j\right)$
that is not divisible by $f$. Then $\dim\dual X=m-2$.
\end{Theorem}

\begin{Example}\rm 
An irreducible surface in $\P^3$ with equation $f(x_1,\ldots,x_3)=0$
is projectively dual to a space curve if and only if the Hessian
$\det\left|\partial^2f/\partial x_i\partial x_j\right|$ is divisible by $f$.
\end{Example}

The Segre Theorem is applicable in the most common situation:
the projectively dual variety is usually a hypersurface,
therefore by taking duals of hypersurfaces we can get a considerable
amount of varieties that are not hypersurfaces.

\begin{Example}\rm
The following remark belongs to Zak (unpublished).
In notations of Theorem~\ref{KJKGJFKJFK}, suppose that
the Hessian 
$H(f)=\left|\partial^2f/\partial x_i\partial x_j\right|$
is not trivial..
Then $H(f)$ is divisible by $f^{n+1-m}$ 
(this result  can be easily verified by induction
using the well-known fact that if
any $(p\times p)$-minor of a $(p+1)\times (p+1)$-matrix $A$
over the UFD $D$ is divisible by $f^l$, where $f$ is prime,
then $\det A$ is divisible by $f^{l+1}$).
Let $\deg f=d$.
We can compare degrees of $H(f)$ and $f^{n+1-m}$ and obtain the formula
$(n+1)(d-2)\ge d(n+1-m)$, i.e.~$d\ge\bigover{2(n+1)}{m}$.
In other words, if $Y$ is a projective variety in $\P^n$ such that
$\dual Y$ is a hypersurface with the non-trivial Hessian, then
$$\deg \dual Y\ge\bigover{2(n+1)}{\dim Y+2}.$$
Zak conjectures that the non-triviality
of the Hessian can be substituted for the non-degeneracy of $Y$ 
and, moreover,
that this inequality turns to the equality
if and only if  $Y$ is the projectivization
of the variety of rank $1$ elements in a simple Jordan algebra.
\end{Example}

\section{Product Theorem and Applications}\label{KHFJGFKGTDFIYD}

\subsection{Product Theorem}  \index{Product theorem}
Let $X_1\subset\P^{n_1}$ and $X_2\subset\P^{n_2}$ be two
irreducible projective varieties. 
The product $X_1\times X_2$ is then naturally
embedded in $\P^{n_1}\times\P^{n_2}$ and the latter is embedded
in $\P^n$ via the Segre embedding\index{Segre embedding}, where 
$$n+1=(n_1+1)(n_2+1).$$
Therefore we have an embedding $X_1\times X_2\subset\P^n$.
Let $\dual{(X_1\times X_2)}\subset\dual{\P^n}$ be the dual variety.
It turns out that it is quite easy to calculate $\defect X_1\times X_2$.

\begin{Theorem}[\cite{WZ}]\label{ProductTheorem}
$$\defect X_1\times X_2=\max(0, \defect X_1-\dim X_2, \defect X_2-\dim X_1).$$
\end{Theorem}

\proof
We shall use Theorem~\ref{KatzDimension}.
Let $k_1=\dim X_1$, $k_2=\dim X_2$. Let $x_0=(x_{01},x_{02})$
be a smooth point of $X_1\times X_2$.
We choose local coordinates 
$$t_{11},t_{21},\ldots,t_{k_11}$$ 
of $X_1$ near $x_{01}$
and 
$$t_{12},t_{22},\ldots,t_{k_22}$$ 
of $X_2$ near $x_{02}$.
So for $\nu=1,2$, a point $x_\nu$ close to $x_{0\nu}$ has homogeneous coordinates
$$x_\nu=(1:t_{1\nu}:\ldots:t_{k_\nu\nu}:u_{1\nu}:\ldots:u_{n_\nu-k_\nu,\nu}),$$
where each $u_{1\nu}$ is an analytic function of 
$$t_{1\nu},\ldots,t_{k_\nu\nu}$$
vanishing at the origin.

By definition of the Segre embedding, a point $(x_1,x_2)\in X_1\times X_2$
has as homogeneous coordinates all pairwise products of homogeneous
coordinates of $x_1$ and $x_2$. Therefore, the set of homogeneous coordinates
(near $(x_{01}, x_{02})$) is given by
$$(1:t_1:t_2:u_1:u_2:tt:tu:ut:uu),$$
where symbols 
$$t_1, t_2, u_1, u_2, tt, tu, ut, uu$$ 
stand for the following
sets of variables:
$$t_1=\{t_{j1}\},\ t_2=\{t_{j2}\},\ u_1=\{u_{i1}\},\ u_2=\{u_{i2}\},$$
$$tt=\{t_{j_11}t_{j_22}\},\ tu=\{t_{j1}u_{i2}\},\ 
ut=\{u_{i1}t_{j2}\},\  uu=\{u_{i_11}u_{i_22}\}.$$
The coordinates in sets $t_1$ and $t_2$ can be chosen as local coordinates
on $X_1\times X_2$ near $(x_{01},x_{02})$.
The coordinates from remaining sets are analytic functions of these local coordinates
vanishing at the origin.

To apply Theorem~\ref{KatzDimension}, we have to consider
the Hessian matrix $\Hes(u)$, where $u$ is a linear combination
of all the coordinates from the sets
$$t_1, t_2, u_1, u_2, tt, tu, ut, uu$$ 
regarded as a function
of local coordinates from $t_1$ and $t_2$.
We can write $\Hes(u)$ as
$$\Hes(u)=\left(\matrix{
A_{11}&A_{12}\cr
A_{21}&A_{22}\cr
}\right),$$
where
$$A_{\alpha\beta}=\left({\partial u\over\partial t_{i\alpha}\partial t_{j\beta}}(0,\ldots,0)
\right)_{1\le i\le k_\alpha,\ 1\le j\le k_\beta}.$$
Then it is easy to see that $A_{11}=\Hes(u_1)$, $A_{22}=\Hes(u_2)$, and
$$\Hes(u)=\left(\matrix{
\Hes(u_1)&A\cr
A^\perp &\Hes(u_2)\cr
}\right),$$
where $A$ can be any $(k_1\times k_2)$-matrix.
Therefore, by Theorem~\ref{KatzDimension} it is sufficient 
to use the following fact from linear algebra:

\begin{Lemma}
For any integers $k_1\ge c_1\ge0$, $k_2\ge c_2\ge0$,
let $\Bil(c_1, k_1; c_2, k_2)$ denote the set of all bilinear forms
on  $\C^{k_1}\oplus \C^{k_2}$ such that their restrictions
on  $\C^{k_1}$ and $\C^{k_2}$ have coranks $c_1$ and $c_2$.
Then a generic form  $\varphi\in\Bil(c_1,k_1;c_2,k_2)$ 
has corank 
$$\corank(\varphi)=\max(0,c_1-k_2,c_2-k_1).$$
\end{Lemma}

Let us prove this lemma. 
The corank of a bilinear form is equal to the dimension
of its kernel.
Let $\varphi\in\Bil(c_1,k_1;c_2,k_2)$, 
$U=\Ker\varphi|_{\C^{k_1}}$, $\dim U=c_1$.
The form $\varphi$ induces a linear map $\psi:\,U\to (\C^{k_2})^*$.
Let $U_0=\Ker\psi$. Then it is clear that $U_0\subset \Ker\varphi$ and
$\dim U_0\ge \dim U-\dim (\C^{k_2})^*=c_1-k_2$.
This argument shows that 
$$\corank(\varphi)\ge  \max(0,c_1-k_2,c_2-k_1).$$
It remains to find a form such that
this inequality turns to an equality.
We shall proceed by induction.
Suppose first that $c_1$ and $c_2$ are positive.
Take a form $\varphi'\in\Bil(c_1-1,k_1-1;c_2-1,k_2-1)$ such that
$$\corank(\varphi')=\max(0,c_1-k_2,c_2-k_1).$$
We can define the form $\varphi$ on
$$\C^{k_1}\oplus\C^{k_2}=\left(\C^{k_1-1}\oplus\langle e\rangle\right)\oplus
\left(\C^{k_2-1}\oplus\langle f\rangle\right)$$
such that its restriction on 
$\C^{k_1-1}\oplus\C^{k_2-1}$ coincides with $\varphi'$,
subspaces $\C^{k_1-1}\oplus\C^{k_2-1}$ and $\langle e,f\rangle$
are orthogonal, and $\varphi(e,e)=\varphi(f,f)=0$, $\varphi(e,f)=1$.
Then it is clear that $\varphi\in\Bil(c_1,k_1;c_2,k_2)$ and
$$\corank(\varphi)=\corank(\varphi')=\max(0,c_1-k_2,c_2-k_1).$$

Now we may assume that $c_1=0$.
Suppose that $c_2$ and $k_1$ are both positive.
Take a form $\varphi'\in\Bil(0,k_1-1;c_2-1,k_2-1)$ such that
$$\corank(\varphi')=\max(0,c_2-k_1).$$
We define the form $\varphi$ on
$$\C^{k_1}\oplus\C^{k_2}=\left(\C^{k_1-1}\oplus\langle e\rangle\right)\oplus
\left(\C^{k_2-1}\oplus\langle f\rangle\right)$$
such that the restriction of $\varphi$ on
$\C^{k_1-1}\oplus\C^{k_2-1}$ coincides with $\varphi'$,
subspaces $\C^{k_1-1}\oplus\C^{k_2-1}$ and $\langle e,f\rangle$
are orthogonal, and $\varphi(e,e)=\varphi(e,f)=1$, $\varphi(f,f)=0$.
It is clear then that $\varphi\in\Bil(0,k_1;c_2,k_2)$ and
$$\corank(\varphi)=\corank(\varphi')=\max(0,c_2-k_1).$$

It remains to consider only cases when either $k_1=0$ or $c_1=c_2=0$.
In both cases lemma is obvious.
\endproof

It follows immediately from Theorem~\ref{ProductTheorem}
that the same result holds for the product of any number of factors.
Let $X_k\subset\P^{n_k}$, $k=1\ldots,r$ be 
irreducible projective varieties. 
The product $X=X_1\times\ldots\times X_r$ is then naturally
embedded in $\P^{n_1}\times\ldots\times P^{n_r}$ and the latter is embedded
in $\P^n$ via the Segre embedding, where $n+1=(n_1+1)\times\ldots\times(n_r+1)$.
Therefore we have an embedding $X\subset\P^n$.
Let $\dual{X}\subset\dual{\P^n}$ be the dual variety.

\begin{Theorem}\label{ManyProductTheorem}
$$\defect X+\dim X=\max(\dim X, \defect X_1+\dim X_1,\ldots, \defect X_r+\dim X_r).$$
In particular, $\dual X$ is a hypersurface if and only if
$$\defect X_k+\dim X_k\le \dim X,\quad k=1\ldots,r.$$
\end{Theorem}

\subsection{Hyperdeterminants}\label{HyperdeterminantDefs}\index{hyperdeterminants}

For $i=1,\ldots,r$ consider a projective space $\P^{l_i}=\P(\C^{l_i+1})$.
The Segre embedding $\P^{l_1}\times\ldots\times\P^{l_r}\subset\P^l$,
$\P^l=\P(\C^{l_i+1}\otimes\ldots\otimes\C^{l_r+1})$,
identifies the affine cone $\Cone(\P^{l_1}\times\ldots\times\P^{l_r})$
with a variety of decomposable tensors in 
$\C^{l_i+1}\otimes\ldots\otimes\C^{l_r+1}$.
The corresponding discriminant (if exists) is called
a hyperdeterminant for the matrix format
$(l_1+1)\times\ldots\times(l_r+1)$, cf.~Example~\ref{Determinant}.
Theorem~\ref{ProductTheorem} allows to get a simple criterion
for the existence of hyperdeterminants:

\begin{Theorem}\label{ExistenceOfDeterminantsTh}\index{hyperdeterminants, existence of}
The dual variety $\dual{(\P^{l_1}\times\ldots\times\P^{l_r})}$
is a hypersurface {\rm(}and, hence, defines a hyperdeterminant{\rm)}
if and only if 
$$2l_k\le l_1+\ldots+l_r$$
for $k=1,\ldots,r$.
\end{Theorem}

\proof
Since the dual variety of $\dual{\P^n}$ is an empty set, 
$\defect\P^n=n$.
Now everything follows from Theorem~\ref{ProductTheorem}.
\endproof

\subsection{Associated Hypersurfaces}\index{associated hypersurfaces}
Let $X\subset\P(V)$ be an irreducible subvariety, $\P^l=\P(U)$.
By Theorem~\ref{ProductTheorem}, the dual variety $\dual{(X\times\P^l)}$
is a hypersurface in $\dual{\P(V\otimes U)}$ for $\defect X\le l\le \dim X$.
Since we always have $\defect X\le\dim X$, we obtain a system
of hypersurfaces naturally associated to $X$. Their
interpretation was found in \cite{WZ}.

We can identify $\dual{\P(V\otimes U)}$ with the projectivization
of the space of linear maps $\Hom(U,V^*)$.
By some abuse of notation we shall denote
a non-zero linear map and its projectivization by the same letter.
Let $\P(V\otimes U)^*_0$ be an open subset consisting
of operators having maximal possible rank $l+1$.
For any $f\in\P(V\otimes U)^*_0$ we can associate two projective
subspaces 
$$U(f)=\P(\Im(f))\subset \P(V)^*$$ and the orthogonal subspace
$U(f)^\perp\subset\P(V)$.
Recall that
$I_X^0\subset \P(V)\times\dual{\P(V)}$ is an open part of the conormal
variety consisting 
of pairs $(x,H)$ such that $x\in X_{sm}$ and $H$ is the hyperplane
tangent to $X$ at $x$.

\begin{Theorem}[\cite{WZ}]\label{ProductPn}
The dual variety $\dual{(X\times\P^l)}$ is the closure
of the set of points $f\in\P(V\otimes U)^*_0$
such that 
$$\left(U(f)^\perp\times U(f)\right)\cap I_X^0\ne\emptyset.$$
\end{Theorem}

In two extreme cases $l=\dim X$ or $l=\defect X$ this result 
can be made more precise. To go further, we need a notion 
of an associated hypersurface.
Let $X\subset\P^n$ be a $k$-dimensional irreducible subvariety.
Consider the set ${\cal Z}(X)$ of all $(n-k-1)$-dimensional
projective subspaces in $\P^n$ that intersect $X$.
This is a subvariety in the Grassmanian $\Gr(n-k,n+1)$
parametrizing all $(n-k-1)$-dimensional projective subspaces in $\P^n$.

\begin{Theorem}
The subvariety ${\cal Z}(X)$ is an irreducible hypersurface,
called an associated hypersurface of $X$.
\end{Theorem}

\proof
Let $B(X)\subset X\times{\cal Z}(X)$ be an incidence
variety consisting of pairs $(x,L)$ such that $x\in X$,
$L\in\Gr(n-k,n+1)$, and $x\in L$.
We have two projections $p:\,B(X)\to X$ and $q:\,B(X)\to {\cal Z}(X)$.
Then $p$ is a Grassmannian fibration, any fiber of $p$ \index{Grassmanian fibration}
is isomorphic to $\Gr(n-k-1,n)$. In particular, $B(X)$ (and hence ${\cal Z}(X)$)
is irreducible.
On the other hand, $q$ is birational, because a generic $(n-k-1)$-dimensional
projective subspace intersecting $X$ meets $X$ in exactly one point.
Now an easy dimension count shows that ${\cal Z}(X)$ is a hypersurface.
\endproof

Now we can relate associated hypersurfaces and dual varieties.
We continue to use the notation of Theorem~\ref{ProductPn}.
The following theorem is sometimes called a Cayley trick:\index{Cayley trick}

\begin{Theorem}[\cite{WZ}]
Let $k=\dim X$, $\dual k=\dim\dual X$. Then
\begin{enumerate}
\item[\rm (a)] Let $l=\dim X$, and 
$p_1:\,\P(V\otimes U)^*_0\to\Gr(n-k-1,n+1)$ be the map $f\mapsto U(f)^\perp$.
Then
$$\dual{(X\times\P^l)}=\overline{p_1^{-1}({\cal Z}(X))}.$$
\item[\rm (b)] Let $l=\defect X$, and 
$p_2:\,\P(V\otimes U)^*_0\to\Gr(n-\dual k-1,n+1)$ be the map 
$f\mapsto U(f)$.
Then
$$\dual{(X\times\P^l)}=\overline{p_2^{-1}({\cal Z}(\dual X))}.$$
\end{enumerate}
\end{Theorem}

We refer the reader to \cite{GKZ2} for further results on associated varieties
and Chow forms. \index{Chow forms}

\section{Ein Theorems}\label{JHGVJHGVJGHFKGFKGV}

Suppose that $X\subset\P^N$ is a non-degenerate smooth projective variety,
$\dim X=n$. \index{contact locus}\index{Ein Theorems}
For any hyperplane $H\subset\P^N$ the contact locus $\Sing X\cap H$ is 
a subvariety of $X$ consisting of all points $x\in X$ such that
the embedded tangent space $\hat T_xX$ is contained in $H$.
One can use the Jacobian ideal of $X\cap H$ to define a scheme structure 
on the contact locus, however this scheme could be not reduced.
Clearly the contact locus is non-empty if and only if $H$ belongs
to the dual variety $\dual X$. By Theorem~\ref{FirstEin}
if $\defect X=d$ then for any $H\in{\dual X}_{sm}$
the contact locus $\Sing(H\cap X)$ is a projective subspace of dimension $d$
and a union of this projective subspaces is dense in $X$.
In this section we study further properties
of the contact locus, most results here belong to L.~Ein.

\begin{Theorem}[\cite{E1}]\label{IntermediateEinTh}
Suppose that $q$ is a generic point of $X$ and $H$
is generic tangent hyperplane of $X$ at $q$, $L=\P^d$ is
the contact locus of $H$ with $X$. Then
\begin{enumerate}
\item[\rm(a)] if $p$ is a point in $\P^d$, then the tangent cone of the 
hyperplane section $H\cap X$ at $p$ is a quadric hypersurface
of rank $n-d$ in $\hat T_p(X)$.
\item[\rm(b)] Let $s_h$ be the section of $\O_X\otimes\O_X(1)$ defining
$H\cap X$. Then $s_h$ factors through $I_L^2$, 
where $I_L\subset\O_X$ is the ideal sheaf of $L$ in $X$.
\item[\rm(c)] Let $t_h$ be the section of 
$I^2_L/I^3_L\otimes\O_L(1)\cong S^2(N^*_LX)\otimes\O_L(1)$
induced by~$s_h$. Then $t_h$ defines a nonsingular quadric hypersurface in 
$\P(N^*_LX|_p)$.
\end{enumerate}
\end{Theorem}

\proof
(a) This follows from the proof of Theorem~\ref{KatzDimension}
using Sard's Lemma for algebraic varieties.

(b) We choose a local coordinate system $\{x_1,\ldots,x_n\}$ of $X$
near $p$. We may assume that $I_L$ is generated by $x_1,\ldots,x_{n-d}$.
Since $L\subset H\cap X$, we can write the power series of $s_h$
in the following form:
$$s_h=x_1f_1+\ldots+x_{n-d}f_{n-d}+\sum_{i=1}^{n-d}\sum_{j=1}^{n-d}x_ix_jg_{ij},$$
where $f_1,\ldots,f_{n-d}$ are power series in variables
$x_{n-d+1},\ldots,x_n$ only.
But $\Sing H\cap X=L$, therefore $\left.\bigover{\partial s_h}{\partial x_i}\right|_{L}=0$
for $i=1,\ldots,n-d$. Hence $$f_1=\ldots,f_{n-d}=0$$ and
$$s_h=\sum_{i=1}^{n-d}\sum_{j=1}^{n-d}x_ix_jg_{ij}.$$
Thus $s_h$ factors through $I_L^2$.

(c) We can write $g_{ij}=a_{ij}+h_{ij}$, where $a_{ij}$'s are constants
and $h_{ij}$'s are power series without the constant term. Now
$$s_h=\sum_{i=1}^{n-d}\sum_{j=1}^{n-d}x_ix_j(a_{ij}+h_{ij}).$$
By (a), $\sum\limits_{i=1}^{n-d}\sum\limits_{j=1}^{n-d}x_ix_ja_{ij}$
is a quadratic form of rank $n-d$.
But this is also the equation for the quadric hypersurface
in $\P(I_L/I_L^2)$ induced by $s_h$.
\endproof

Recall that any vector bundle on $\P^1$ has a form $\oplus_i\O(a_i)$
for some integers~$a_i$ (see e.g.~\cite{Ha1}). 

\begin{Definition}
A vector bundle $E$ on a projective space $\P^N$\index{uniform vector bundle}
is called uniform if for any line $T\subset\P^N$ the restriction $E|_T$
is a fixed vector bundle $\oplus_i\O(a_i)$.
\end{Definition}

The following theorem is the main result of this section.

\begin{Theorem}[\cite{E1}]\label{SecondEin}
Suppose that $X$ is a nonlinear smooth projective variety in $\P^N$, 
$\defect X=d>0$, $q$ is a generic point of $X$ and $H$
is generic tangent hyperplane of $X$ at $q$, $L=\P^d$ is
the contact locus of $H$ with~$X$. Then
\begin{enumerate}
\item[\rm(a)] $N_LX\cong(N_LX)^*(1)$.
\item[\rm(b)] If $T=\P^1$ is a line in $L=\P^d$, then
$$N_LX|_T\cong\cal O_T^{\oplus(n-d)/2}\oplus
\cal O_T(1)^{\oplus(n-d)/2}$$ 
{\rm(}in particular $N_LX$ is a uniform vector bundle{\rm)} and
$$N_TX\cong\cal O_T^{\oplus(n-d)/2}\oplus
\cal O_T(1)^{\oplus(n+d-2)/2}.$$
\item[\rm(c)] There is an irreducible $\bigover{3n+d-4}{2}$-dimensional
family of lines in $X$. If $p$ is a generic point in $X$,
then there is an $\bigover{n+d-2}{2}$-dimensional family of lines in $X$
through $p$.
\end{enumerate}
\end{Theorem}

\proof
(a) We continue to use the notation of Theorem~\ref{IntermediateEinTh}.
By this theorem, the section $s_h$
gives a section of
$$I_L^2/I_L^3(1)=S^2(N_L^*X)(1)\subset\Hom(N_LX,N_L^*X(1)).$$
Therefore $s_h$ induces a map from $N_L(X)$ to $N_LX(1)$.
This map is an isomorphism by Theorem~\ref{IntermediateEinTh} (c).

(b) Consider the exact sequence
$$0\to N_LX|_T\to N_L\P^N|_T\to N_X\P^N|_T\to 0,$$
where $N_L\P^N|_T=(N-d)\O_T(1)$.
Suppose that 
$$N_LX|_T=\mathop{\oplus}\limits_{i=1}^{n-d}\O_T(a_i).$$
Then exactness implies that all $a_i\le 1$.
Using the isomorphism between $N_LX$ and $N_L^*X(1)$ we observe
that each $a_i\ge0$. Therefore,
$$N_LX|_T\cong\cal O_T^{\oplus(n-d)/2}\oplus
\cal O_T(1)^{\oplus(n-d)/2}$$ 
and $N_LX$ is a uniform vector bundle.
This implies that
$$N_TX\cong\cal O_T^{\oplus(n-d)/2}\oplus
\cal O_T(1)^{\oplus(n+d-2)/2}.$$

The statement (c) is a consequence of the standard deformation theory.
Indeed, let $p\in L$ and let $T_0$ be a line in $L$ through $p$.
Since 
$$N_{T_0}X\cong\cal O_{T_0}^{\oplus(n-d)/2}\oplus
\cal O_{T_0}(1)^{\oplus(n+d-2)/2},$$
we have 
$$\dim H^0(T_0, N_{T_0}X)={n-d\over 2}+2{n+d-2\over2}={3n+d-4\over2}$$
and 
$$H^1(T_0, N_{T_0}X)=0.$$
Therefore, the Hilbert scheme of lines in $X$ is smooth \index{Hilbert scheme}
at the point $T_0$ and there is a unique irreducible component
of the Hilbert scheme containing the point $T_0$.
This component has dimension $(3n+d-4)/2$.
The proof of the second statement of (c) is similar.
\endproof

As a quite formal consequence of Theorem~\ref{SecondEin} (b)
we get the following parity theorem that was first proved
by A.~Landman, using the Picard--Lefschetz theory (unpublished):
\index{parity theorem}

\begin{Theorem}[\cite{E1}]\label{PairityTheorem}
If $\defect X>0$ then
$\dim X\equiv\defect X\,\mod\,2$.
\end{Theorem}

\begin{Example}[\cite{GH}]\label{DualOfASurface}\rm  \index{dual of the smooth surface}
Suppose that $X\subset\P^n$ is a smooth non-linear surface.
Then $\dual X$ is a hypersurface.
Indeed, if $\defect X>0$ then $\defect X$ is even by
Theorem~\ref{PairityTheorem}.
But $X$ obviously can not contain an even-dimensional projective 
subspace $\P^{\defect X}$.
\end{Example}

\begin{Example}\label{FGGDCGFD}\rm
Suppose that $X\subset\P^N$ is a smooth projective nonlinear variety, $\dim X\ge2$.
Then $\defect X\le \dim X-2$. Indeed, it is clear that $\defect X\le\dim X-1$.
Therefore, by Theorem~\ref{PairityTheorem} we have
$\defect X\le\dim X-2$. Moreover, it can be shown \cite{E1}
that $\defect X=\dim X-2$ if and only if $X$ is a projective bundle
over a curve $C$, $X=\P_C(F)$, where $F$ is a rank $n$ vector bundle
on $C$ and all fibers are embedded into $\P^N$ linearly
(c.f. Theorem~\ref{JHGKJFKJFGK} and remarks after it).
\end{Example}

\begin{Example}\label{FHGFGHJMFMGHC}\rm
Suppose that $X\subset\P^N$ is a smooth projective nonlinear variety, 
$\dim X\ge3$, $N=2n-1$, and $\dim X=\dim\dual X$.
Since $\defect X=n-2$ in this case, by Example~\ref{FGGDCGFD}
we have that
$X$ is a projective bundle
over a curve $C$, $X=\P_C(F)$, where $F$ is a rank $n$ vector bundle
on $C$ and all fibers are embedded into $\P^N$ linearly.
Moreover, by a theorem of S.~Kleiman \cite{Kl5}
in this case $X$ is the Segre embedding of $\P^1\times\P^{n-1}$.
\end{Example}

It was first observed by Griffiths and Harris \cite{GH} that
any smooth projective variety with positive defect has the negative
Kodaira dimension, $H^0(X,K_X^{\otimes m})=0$ for
$m>0$, where $K_X$ is the canonical line bundle. \index{Kodaira dimension}
This result will have important consequences in Chapter~\ref{MoriChapter}.
Here is the detailed story:

\begin{Theorem}[\cite{E1}]\label{EINKX}
Suppose that $X$ is a nonlinear smooth projective variety in $\P^N$, 
$\defect X=d>0$,
$q$ is a generic point of $X$ and $H$
is a generic tangent hyperplane of $X$ at $q$, $L=\P^d$ is
the contact locus of $H$ with~$X$. Then
\begin{enumerate}
\item[\rm(a)] $K_X|_L=\cal O_L\left(\bigover{-n-d-2}{2}\right)$.
\item[\rm(b)] The Kodaira dimension of $X$ is negative.
\item[\rm(c)] If $K_X=\O_X(a)$, then $a=\bigover{-n-d-2}{2}$.
\item[\rm(d)] If $n>\bigover{N}{2}+1$, then 
$K_X=\O_X\left(\bigover{-n-d-2}2\right)$.
\end{enumerate}
\end{Theorem}
 
\proof
(a) By Theorem~\ref{SecondEin} (b),
if $T=\P^1$ is a line in $L=\P^d$, then
$$N_LX|_T\cong\cal O_T^{\oplus(n-d)/2}\oplus
\cal O_T(1)^{\oplus(n-d)/2}.$$ 
It follows that
$\Lambda^{n-d}N_LX|_T=\O_T\left(\bigover{n-d}{2}\right)$ and, hence,
$\Lambda^{n-d}N_LX=\O_L\left(\bigover{n-d}{2}\right)$.
Therefore
$K_X|_L=\cal O_L\left(\bigover{-n-d-2}{2}\right)$ by the adjunction formula.

(b) Since there is such a plane $L$ through a generic point $p\in X$,
the linear system $|K_X^m|$ is empty for $m\ge0$.

(c) follows from (a). 

(d) follows from (c). Indeed, if $\dim X>\bigover{N}{2}+1$, then 
$\Pic X=\Z$ is generated by $\O_X(1)$ by the Barth Theorem \cite{Ba}.
\endproof

\begin{Example}[\cite{E1}]\rm
Suppose that $X$ is a nonlinear smooth projective variety in $\P^N$, 
$\codim X=2$. Then $\dual X$ is a hypersurface,
unless $X$ is the Segre embedding of $\P^1\times\P^2$ in $\P^5$.
Indeed, assume first that $\dim X\ge4$ and $\defect X=d>0$.
Then by Theorem~\ref{EINKX} (d) the canonical class 
$K_X=\O_X(\bigover{-n-d-2}{2})$. 
However, it follows from results of \cite{BC}
that in this case $X$ is automatically a complete intersection.
This contradicts Theorem~\ref{DualCompleteIntersection}.
Now, if $\dim X=1$ or $\dim X=2$ then $\dual X$ is a hypersurface
by Example~\ref{DualOfTheCurve} and Example~\ref{DualOfASurface}.
The case $\dim X=3$ is ruled out by Example~\ref{FHGFGHJMFMGHC}.
\end{Example}

In case of not generic contact loci the following version
of Theorem~\ref{EINKX} holds:

\begin{Theorem}[\cite{BFS1}]\label{Strong-Ein-EveryPoint}
Suppose that $X$ is a smooth $n$-dimensional
projective variety in $\P^N$, $\defect X=d>0$. 
Let $H\subset\dual X$. Then $\Sing H\cap X$ is the union of projective
subspaces $\P^d$ with 
$K_X|_{\P^d}=\cal O_{\P^d}\left(\bigover{-n-d-2}{2}\right)$.
\end{Theorem}

\proof
Let $I_X\subset X\times\dual{\P^N}$ 
be the conormal variety, $\pi:\,I_X\to\dual{\P^N}$ be the projection.
Then $\pi(I_X)=\dual X$. 
Let $p:\,I_X\to X$ be the first projection, then $p$ is a projective bundle.
Clearly, $I_X$ is smooth and irreducible.
The general fibers $F$ of $\pi$ are isomorphic 
to the contact loci of $X$ with generic tangent hyperplanes, 
therefore, they are isomorphic to the projective subspaces $\P^d$
and by Theorem~\ref{EINKX} we have
$p^*(K_X+((n+d)/2+1)L)_{F}\simeq\O_F$.
It thus follows that every fiber of $\pi$ is a union of
projective subspaces 
$\P^d$ with $p^*(K_X+((n+d)/2+1)L)_{\P^d}\simeq\O_{\P^d}$.
\endproof

Let $X$ be a non-linear $n$-dimensional smooth projective variety in $\P^N$
such that $\defect X=d>0$. Let $H$ be a generic tangent hyperplane of $X$.
Then the contact locus $\Sing H\cap X$ of $X$ with $H$ is a $d$-dimensional
projective subspace $L$. Let $\tilde X$ be a blowing up of $X$ along $L$.
Let $p:\,\tilde X\to X$ be the corresponding map.
Denote by $E$ the exceptional divisor and denote by $F$ the proper transform
of $H\cap X$. Therefore, we have a following diagram
$$\matrix{
\P(N_LX) & = &            E & \subset & \tilde X & \supset & F\cr
           &   & \downarrow   &         &\downarrow&         & \downarrow\cr
           &   & L            &\subset  & X        & \supset & H\cap X.
}$$
We shall denote by $\O_{\tilde X}(a,b)$ the line bundle
$p^*\O_X(a)\otimes\O_{\tilde X}(-bE)$. For example, 
$\O_{\tilde X}(F)=\O_{\tilde X}(1,2)$.

The following vanishing theorem holds.

\begin{Theorem}[\cite{E1}]\label{EinVanishingTh}
Assume that $K_X=\O_X(b)$ for some $b\in\Z$. Then
\begin{enumerate}
\item[\rm(a)] $K_X=\O_X\left(\bigover{-n-d-2}2\right)$.
\item[\rm(b)] $K_{\tilde X}=\O_{\tilde X}\left(\bigover{-n-d-2}2, -n+d+1\right)$.
\item[\rm(c)] $H^i(\O_{\tilde X}(a,1))=0$, if $i>0$ and $a\ge\bigover{n-3d-2}{2}$.
\item[\rm(d)] $H^i(\O_{\tilde X}(a,2))=0$, if $i>0$ and $a\ge\bigover{n-3d}{2}$.
\end{enumerate}
Assume further that $\bigover{n-3d-2}{2}\le0$. Then
\begin{enumerate}
\item[\rm(e)] $H^0(N_L^*X(a))=0$ for $a\le0$.
\item[\rm(f)] $H^k(N_L^*X(a))=0$ for $a\ge-d$.
\item[\rm(g)] $H^i(N_L^*X(a))=0$ if $0<i<d$ and $a\ge\bigover{n-3d}2$.
\item[\rm(h)] $H^i(N_L^*X(a))=0$ if $0<i<d$ and $a\le\bigover{d-n}2$.
\end{enumerate}
\end{Theorem}

\proof
(a) follows from Theorem~\ref{EINKX} (a).

(b) follows from (a) and the fact that $\tilde X$ is the blowing up of $X$
along $L$.

(c) Let $f:\,\tilde X\to\P^{N-1-d}$ be the projection with center $L$.
Let $Y=f(\tilde X)$. The hyperplane section $H\cap X$ will correspond
to a hyperplane section $D$ of $Y$. Then $f^{-1}(D)=E+F$ and $f^*\O_Y(1)=\O_{\tilde X}(1,1)$.
Let $y\in Y\setminus D$, $Z=f^{-1}(y)$. Suppose that $\dim Z\ge1$.
Since $Z\cap (E\cup F)=\emptyset$, $p$ maps $Z$ isomorphically
to a variety in $X$. Therefore, $\O_{\tilde X}(0,1)|_Z$ is non-trivial.
But $\O_{\tilde X}(1,1)|_Z=f^*\O_Y(1)|_Z$ is trivial. So
$\O_{\tilde X}(0,1)$ is non-trivial. Hence $Z\cap E\ne\emptyset$,
which is a contradiction. Therefore, all positive-dimensional
fibers of $f$ belong to $E\cap F$, in particular $\dim Y=\dim X$.
Now, 
$$\O_{\tilde X}(a,1)=K_{\tilde X}\otimes f^*\O_Y(n-d)\otimes\O_{\tilde X}
\left(a-{n-3d-2\over2},0\right).$$
It follows from the Grauert--Rimenschneider vanishing theorem \cite{GR}
that $H^i(\O_{\tilde X}(a,1))=0$ if $i>0$ and $a\ge\bigover{n-3d-2}2$.

(d) The proof is similar to (c).

(e) Consider the exact sequence
$$0=H^0(\O_{\tilde X}(0,1))\to H^0(\O_E(0,1))\to H^1(\O_{\tilde X}(0,2)).$$
Now $H^1(\O_{\tilde X}(0,2))=0$ by (d). So 
$$H^0(\O_E(0,1))\simeq H^0(N_L^*X)=0.$$
Hence $H^0(N_L^*X(a))=0$ for $a\le0$.

(f) Recall that $N_LX=N_L^*X(1)$. So (f) follows from (e)
and Serre duality.

(g) Consider the exact sequence
$$H^i(\O_{\tilde X}(a,1))\to H^i(\O_E(a,1))\to H^{i+1}(\O_{\tilde X}(a,2)).$$
By (c) and (d), we conclude that
$$H^i(\O_E(a,1))\simeq H^i(N_L^*X(a))=0$$
for $a\ge\bigover{n-3d}2$.

(h) follows from (g) and Serre duality.
\endproof

\begin{Theorem}[\cite{E1}]\label{ApplicationBeilinson}
Suppose that $X$ is a nonlinear smooth projective variety in $\P^N$, 
$\dim X=n\ge3$. Assume that $K_X=\O_X(b)$ for some $b\in\Z$.
Then
\begin{enumerate}
\item[\rm(a)] $\defect X\le\bigover{n-2}{2}$.
\item[\rm(b)] If $\dim X=4m+2$ and $\defect X=2m>0$ then
$$N^*_LX=H^m\left(N_L^*X(-m)\right)\otimes\Omega_L^{m}(m),$$
and $m\le 2$.
\end{enumerate}
\end{Theorem}

\proof
Consider the Beilinson spectral sequence \cite{OSS} with\index{Beilinson spectral sequence}
$$E_1^{pq}=H^q\left(N^*_LX(p)\right)\otimes\Omega^{-p}_L(-p),$$
which converges to 
$$E^i=\cases{N^*_LX&\hbox{\rm if $i=0$}\cr
0&\hbox{\rm otherwise (i.e.~$E^{pq}_\infty=0$ if $p+q\ne0$)}.}$$

(a) If $k=\defect X\ge\bigover{n-1}2$, then
$\bigover{n-3k}2-1\le\bigover{k-n}2$.
It follows from Theorem~\ref{EinVanishingTh} that
$H^q(N_L^*X(p))=0$ for $-k\le p\le0$. It follows that $N_L^*X=0$.
This is a contradiction.

(b) In this case $H^q(N_L^*X(p))=0$ for $-2m\le p\le 0$ unless $p=-m$.
This implies that $N_L^*X=H^m(N_L^*X(-m))\otimes\Omega_L^m(m)$.
So 
$$2m+2=\rank N_L^*X\ge\rank\Omega_L^m(m)={2m\choose m}.$$
It follows that $m\le2$.
\endproof

We conclude this section with 
the following simple but very useful result first appeared in \cite{LS}
(where it was attributed to the referee). It could be called
the monotonicity theorem. \index{Monotonicity theorem}

\begin{Theorem}[\cite{LS}]\label{LanteriStruppMonotonicity}
Let $X\subset\P^N$ be a smooth projective $n$-dimensional variety.
Suppose that through its generic point there passes
a smooth subvariety~$Y$ of dimension $h$ and defect $\theta$.
Then $\defect X\ge\theta-k+h$. In other words,
$$\dim X+\defect X\ge \dim Y+\defect Y.$$
\end{Theorem}

\proof
Suppose that $x\in X$ is a generic point and $H$ is a generic hyperplane
tangent to $X$ at $x$. Let $Y$ be an $h$-dimensional submanifold
passing through~$x$ such that $\defect Y=\theta$.
Since the defect is the dimension of a generic contact locus, which
is a projective subspace, there exists a $\theta$-dimensional projective
subspace $Z\subset Y$ containing $x$ such that $H$ is tangent
to $Y$ along $Z$ (though $H$ is not a generic tangent hyperplane to $Y$).
Let $f=0$ be a local equation for $H\cap X$ at $x$.
In a neighborhood $U$ of $x$, the differential $df$ annihilates
the tangent spaces $T_xX$ and $T_zY$ for every $z\in Z\cap U$.
Hence $df$ defines on $Z\cap U$ a section of the conormal bundle
$N_Y^*X|_Z$ vanishing at $x$. Since $N_Y^*X$ is a $(n-h)$-dimensional
vector bundle, it follows that this section in fact vanishes
on a subvariety $Z'\subset Z\cap U$ of codimension less then
or equal to $n-h$. Therefore, $H$ is tangent to $X$
along $Z'$. But $\Sing X\cap H$ is $(\defect X)$-dimensional, therefore
$\defect X\ge \theta-(n-h)$.
\endproof


\section{The Projective Second Fundamental Form}\label{ProjectiveIISection}

\subsection{Gauss Map}

Suppose that $X\subset\P(V)$ is a smooth irreducible $n$-dimensional
variety, not necessarily closed. 
Therefore, for any $x\in X$ we have
an embedded tangent space $\hat T_xX$.
A natural way to keep track of the motion of $\hat T_xX$ 
(while $x$ moves in $X$) is the Gauss map.\index{Gauss map}

\begin{Definition}\rm
The Gauss map $\gamma$ is defined as follows.
$$\gamma:\, X\to \Gr(n,\P(V)),\quad x\mapsto \hat T_xX,$$
where $\Gr(n,\P(V))=\Gr(n+1,V)$ is the Grassmanian of $n$-dimensional
projective subspaces in $\P^N$.
\end{Definition}

To measure how $\hat T_xX$ moves to the first
order, one calculates the differential of $\gamma$
$$(d\gamma)_x:\, T_xX\to T_{\hat T_xX}\Gr(n,\P(V))=
\Hom(\Cone(T_xX),V/\Cone(T_xX)).$$
It is easy to see that $\Cone(x)\subset\Ker(d\gamma)_x(v)$
for any $v\in T_xX$, therefore $d(\gamma)_x$ factors to a map
$$\FF^2_x:\,T_xX\to\Hom(T_xX,N_X\P(V)|_x).$$
Moreover, $\FF^2_x$ is symmetric, essentially because the Gauss map
is already the derivative of a map
and mixed partial derivatives commute. Finally, we get the following
definition.

\begin{Definition}\rm
The section $\FF^2$ of the vector bundle $S^2T^*X\otimes N_X\P(V)$
constructed above
is called the projective second fundamental form of~$X$.
\end{Definition}
\index{fundamental form}
\index{projective second fundamental form}

If $X\subset\P(V)$ is a singular variety then $\FF^2$
is still defined over the smooth locus $X_{sm}$.

It is convenient to consider $\FF^2$ as a map $\FF^2:\,N^*_X\P(V)\to S^2T^*X$
and to set $|\FF^2|=\P(\FF^2(N^*_X\P(V)))$.
One can think of $|\FF^2|_x$ as a linear family
of quadric hypersurfaces in $\P(T_xX)$.
The space $\P(N^*_X\P(V)|_x)$ has the geometric interpretation
as the space of hyperplanes tangent to $X$ at $x$, i.e.,
the hyperplanes $H$ such that $X\cap H$ is singular at $x$.
Then $\FF^2_x$ can be viewed as the map that sends a hyperplane
to the quadratic part of the singularity of $X\cap H$ at $x$.

Let $\Base|\FF^2|_x\subset\P(T_xX)$ denote the variety of directions
tangent to the lines that osculate to order $2$ at $x$
(that is, $X$ appears to contain these lines to a second order).
More precisely,
\begin{eqnarray}
\Base|\FF^2|_x=\P\{v\in T_xX\,|\,\FF^2_x(v,v)=0\}\cr
=\{p\in\P(T_xX)\,|\,p\in Q,\ \forall Q\in |\FF^2|\}.
\end{eqnarray}

Let $\Sing|\FF^2|_x$ denote the set of tangent directions
such that the embedded tangent space does not move
to first order in these directions.
\begin{eqnarray*}
\Sing|\FF^2|_x=\P\{v\in T_xX\,|\,\FF^2(v,w)=0\ \forall w\in T_xX\}\cr
=\{p\in\P(T_xX)\,|\,p\in Q_{sing}\ \forall Q\in|\FF^2|\}.
\end{eqnarray*}

\index{smoothness principle}
The following theorem may be called the smoothness principle.
\begin{Theorem}[\cite{GH}]
Let $X\subset\P^N$ be an irreducible projective variety. 
If $X$ is smooth then $\Sing|\FF^2|_x=\emptyset$ for generic $x\in X$.
\end{Theorem}

In other words, if $X$ is smooth then the Gauss map $\gamma$
is generically finite.
In fact, much more is true.

\begin{Theorem}[\cite{Z3,E4,Ra3}]
If $X$ is smooth then the Gauss map is finite and birational.
\end{Theorem}

These theorem do not imply that
generic quadrics from the second fundamental
form are not degenerate.
Indeed, the following theorem is just a  reformulation 
of Theorem~\ref{KatzDimension}:

\begin{Theorem}\label{JHKBHBHJK}
If $X$ is any projective variety and $x\in X_{sm}$
then the projective second fundamental
form $|\FF^2|_x$ of $X$ at $x$ is a system of quadrics of rank bounded above
by $\dim X-\defect X$. Moreover, if $x\in X_{sm}$ is a generic
point and $Q\subset|\FF^2|_x$ is a generic quadric then $\rank Q=\dim X-\defect X$.
\end{Theorem}

\subsection{Moving Frames}

\index{moving frames}
The most convenient formalism for dealing with second (and higher)
fundamental forms is the language of moving frames.
We shall recall it briefly here, more details and applications
can be found in~\cite{GH}.

In $\P^N$ a frame is denoted by $\{A_0,\ldots,A_N\}$.
It is given by a basis $A_0,\ldots,A_N$
for $\C^{N+1}$. 
The set of all frames forms an algebraic variety $\F(\P^n)$
the principal homogeneous space of the general linear group $\GL_{N+1}$.
Each of the vectors $A_i$ may be viewed as a mapping
$v:\,\F(\P^N)\to\C^{N+1}$. 
Expressing the exterior derivative
$dv$ in terms of the basis $\{A_i\}$ gives
$$dA_i=\sum_j\omega_{ij}A_j.$$
The $(N+1)^2$ differential forms $\omega_{ij}$
are the Maurer--Cartan forms on the group $\GL_{N+1}$, \index{Maurer--Cartan forms}
and by taking their exterior derivatives we obtain
the Maurer--Cartan equations  \index{Maurer--Cartan equations}
$$d\omega_{ij}=
\sum_k\omega_{ik}\wedge\omega_{kj}.$$
Geometrically we may think of a frame $\{A_0,\ldots,A_N\}$
as defining a coordinate simplex in $\P^N$, and then $\omega_{ij}$
gives the rotation matrix when this coordinate simplex
is infinitesimally displaced.

There is a fibering $\pi:\,\F(\P^N)\to\P^N$
given by
$$\pi\{A_0,\ldots,A_N\}\to A_0.$$
If we set $\omega_i=\omega_{0i}$ then these $1$-forms
are horizontal for the fibering $\pi$, i.e., they vanish \index{horizontal $1$-form}
at the fibers $\pi^{-1}(p)$, $p\in\P^N$.
From $d\omega_i=\sum_j\omega_j\wedge\omega_{ji}$ we see that
the forms $\{\omega_i\}$ satisfy the Frobenius integrability \index{Frobenius integrability condition}
condition. Thus we may think of the fibration $\pi$
as defined by the foliation
$$\omega_1=\ldots=\omega_N=0.$$

The equation
$$dA_0=\sum_{i=0}^N\omega_iA_i=\sum_{i=1}^N\omega_iA_i \mod A_0$$
has the following geometric interpretation. For each choice of the frame
$\{A_0,\ldots,A_N\}$ lying over $p\in\P^N$ the horizontal $1$-forms
$\omega_1,\ldots,\omega_N$ give a basis for the cotangent space
$T^*_p\P^N$. The corresponding basis $v_1,\ldots,v_N\in T_p\P^N$
for the tangent space has the property that $v_i$ is tangent to the line
$\ov{A_0A_i}$.

Assume now that we are given a connected $n$-dimensional
smooth subvariety $M\subset\P^N$, not necessarily closed.
Associated to $M$ is the submanifold $\F(M)\subset\F(\P^N)$
of Darboux frames\index{Darboux frame}
$$\{A_0;A_1,\ldots,A_n;A_{n+1},\ldots,A_N\}$$
defined by the conditions that $A_0$ lies over some $p\in M$
and $A_0,A_1,\ldots,A_n$ span $\hat T_pM$, the embedded tangent space
of $M$ at $p$. On $\F(M)$ we have the condition 
that differential $1$-forms $\omega_1,\ldots,\omega_n$ give a basis
of $T^*_pM$ and $\omega_i=0$ for $i>n$.
Therefore it is natural to think of the $\omega_1,\ldots,\omega_n$
as homogeneous coordinates in the projectivized tangent spaces
$\P(T_pM)$. For example, a quadric in $\P(T_pM)$ is defined
by an equation $\sum_{i,j=1}^nq_{ij}\omega_i\omega_j=0$, $q_{ij}=q_{ji}$.

Since $\omega_i=0$ for $i>n$, we also have $d\omega_i=0$ on $\F(M)$.
On the over hand, by the structure equations we have
$d\omega_i=\sum_{j=1}^N\omega_j\wedge\omega_{ji}$.
Therefore, for any $i>n$ we have
$$\sum_{j=1}^n\omega_j\wedge\omega_{ji}=0.$$
Since $\omega_1,\ldots,\omega_n$ are linearly independent,
an easy calculation (called Cartan Lemma\index{Cartan Lemma}) shows that
for any $\mu>n$, for any $i=1,\ldots,n$ we have
$$\omega_{i\mu}=\sum_{j=1}^nq_{ij\mu}\omega_j, \quad q_{ij\mu}=q_{ji\mu}.$$
We set for any $\mu=n+1,\ldots,N$
$$Q_\mu=\sum_{i,j=1}^nq_{ij\mu}\omega_i\omega_j.$$
The equation $Q_\mu=0$ defines a quadric hypersurface in $\P(T_pM)$.

Then an easy local calculation (see \cite{GH1})
shows that we have
\begin{Proposition}
The linear system of quadrics in $\P(T_pM)$
spanned by~$Q_\mu$, $\mu=n+1,\ldots,N$ is the projective
second fundamental form $|\FF^2|$.
\end{Proposition}

Thus, if we think of $\FF^2$ as of the map $S^2T_pM\to (N_M\P^N)_p$
and identify $(N_M\P^N)_p$ with $\C^{N+1}/\Cone(\hat T_pM)$
then for any vector $v\in T_pM$ the vector $\FF^2(v)$
is given in coordinates by
$$\FF^2(v)=\sum_{i,j=1}^n\sum_{\mu=n+1}^N\omega_i(v)\omega_j(v)A_\mu,$$
where $v=\sum_{i=1}^n\omega_i(v)v_i\in T_pM$.

Another geometric interpretation of the second fundamental
form is via the classical Meusnier--Euler Theorem.\index{Meusnier--Euler theorem}

\begin{Definition}
If $\{p(t)\}$ is any holomorphic arc in $\P^N$
described by a vector--valued function $A_0(t)$, then the osculating
sequence is the sequence of linear spaces spanned by the\index{osculating sequence}
following collection of vectors\index{osculating space}
$$\{A_0(t)\},\ \{A_0(t),A'_0(t)\},\ 
\{A_0(t),A'_0(t),A''_0(t)\},\ \ldots$$
\end{Definition}

We then have the Meusnier--Euler Theorem

\begin{Theorem}[\cite{GH1}]\label{EulerMeusnierTh}
For a tangent vector $v\in T_pM$, the normal vector $\FF^2(v,v)\in (N_M\P^N)_p$
gives the projection in $\C^{n+1}/\Cone(\hat T_pM)$ of the second osculating
space to any curve $p(t)$, $p(0)=p$ with tangent $v$ at $t=0$
\end{Theorem}

\proof
We choose an arbitrary field of Darboux frames $\{A_i(t)\}$
along $p(t)$ and write
$${dA_0\over dt}=\sum_{i=1}^n\left({\omega_i\over dt}\right)A_i \mod A_0.$$
Recall that $\omega_i$ is a $1$-form on the frame manifold, so the expression
$\omega_i\over dt$ makes sense and is equal to the value of this $1$-form
on the tangent vector to the curve $\{A_i(t)\}$ on the frame manifold.
If we interpret $\{\omega_1,\ldots,\omega_n\}$ 
as a basis in the cotangent space $T^*_pM$
then this equation shows that $\omega_i\over dt$ is equal
to the $i$-th coordinate of $v\in T_pM$ in the dual basis $\{v_1,\ldots,v_n\}$
of $T_pM$ (in fact, $v_i=A_i \mod A_0$).

Differentiating further, we get
\begin{eqnarray*}
{d^2A_0\over dt^2}=\sum_{i=1}^n
\left({\omega_i\over dt}\right)\left({dA_i\over dt}\right) \mod A_0,\ldots,A_n\cr
=\sum_{i,j=1}^n\sum_{\mu=n+1}^N q_{ij\mu}
\left({\omega_i\over dt}\right)\left({\omega_j\over dt}\right)\mod\Cone(\hat T_pM).
\end{eqnarray*}
But we also have
$$\FF^2(v,v)=\sum_{i,j=1}^n\sum_{\mu=n+1}^N q_{ij\mu}
\left({\omega_i\over dt}\right)\left({\omega_j\over dt}\right).$$
This finishes the proof.
\endproof

Finally, let us give a useful alternative coordinate description
of the second fundamental form of $M$ in $\P^N$.
This one will be absolutely identical to the description
implicitly used in the proofs of Theorems~\ref{KatzDimension} and~\ref{IntermediateEinTh}.
For $p\in M$, we choose a homogeneous coordinate system
$X_0,\ldots,X_N$ for $\P^N$ so that $p=(1:0:\ldots:0)$
and $\hat T_pM$ is spanned by the first $n+1$ coordinate vectors.
If $x_1=X_1/X_0$, $\ldots$, $x_N=X_N/X_0$ is the corresponding
affine coordinate system for the affine chart of $\P^N$
and $(z_1,\ldots,z_n)$ is any holomorphic coordinate system for $M$
centered at $p$, then for $\mu=n+1,\ldots,N$ we have 
in the neighborhood of $p$ that
$x_\mu|_M=x_\mu(z)$ vanishes to second order at $z=0$. Therefore,
$$x_\mu(z)=\sum_{i,j=1}^nq_{ij\mu}z_iz_j+\hbox{\rm (higher order terms)}.$$
The second fundamental form is the linear system
of quadrics generated by $\sum_{i,j=1}^nq_{ij\nu}dz_i dz_j$.

\subsection{Fundamental Forms of Projective Homogeneous Spaces}

It is interesting to explicitly compute
the second fundamental form of some familiar homogeneous spaces.

\begin{Example}[\cite{GH1}]\rm  \index{Segre embedding}
Suppose that $\P^m\times\P^n\subset\P^{nm+n+m}$ is the Segre embedding.
Recall that if $\P^m=\P(V)$, $\P^n=\P(W)$ then
the image $M$ of the Segre variety is the image of the natural
inclusion
$\P(V)\times\P(W)\subset\P(V\otimes W)$.
Choosing coordinates $\{X_\alpha\}$ in $V$ and $\{Y_\mu\}$ in $W$,
this inclusion is given by mapping 
$(\{X_\alpha\},\{Y_\mu\})\to\{X_\alpha Y_\mu\}$.
We denote by $\{A_0,\ldots,A_m\}$ and $\{B_0,\ldots,B_n\}$
frames for $\P(V)$ and $\P(W)$.
We shall use the range of indices $1\le i,j\le m$, $1\le\alpha,\beta\le n$.
If $A_0$ lies over $v_0\in\P(V)$ and $B_0$ over $w_0\in\P(W)$, then the
frame
\begin{equation}
\{A_0\otimes B_0,A_0\otimes B_\alpha, A_i\otimes B_0, A_i\otimes B_\alpha\}\label{JGHVFHDHFG}
\end{equation}
in $\P(V\otimes W)$ lies over $v_0\otimes w_0$.
If we write
$$dA_0=\sum_i\phi_iA_i \mod A_0,\quad
dB_0=\sum_\alpha\psi_\alpha B_\alpha \mod B_0,$$
then
$$d(A_0\otimes B_0)=\sum_\alpha\psi_\alpha A_0\otimes B_\alpha
+\sum_i\phi_i A_i\otimes B_0 \mod A_0\otimes B_0.$$
Therefore, (\ref{JGHVFHDHFG}) is a Darboux frame for $M$
and
$$\hat T_{v_0\otimes w_0}M=\P(V\otimes w_0+v_0\otimes W),$$
so that
$$(N_M)_{v_0\otimes w_0}=(V\otimes W)/(V\otimes w_0+v_0\otimes W).$$
Differentiating further, we get
$$d^2(A_0\otimes B_0)=2\sum_{i,\alpha}\phi_i\psi_\alpha A_i\otimes B_\alpha
\mod \{A_0\otimes B_0,A_i\otimes B_0, A_0\otimes B_\alpha\}.$$
Therefore, the second fundamental form is given by
$$\FF^2(v_0\otimes w+v\otimes w_0,v_0\otimes w+v\otimes w_0)=2v\otimes w\mod
V\otimes w_0+v_0\otimes W.$$
In terms of homogeneous coordinates $\{\phi_i\}$ for $T_{v_0}\P^m$
and $\{\psi_\alpha\}$ for $T_{w_0}\P^n$, $|\FF^2|$ is the linear system
of quadrics 
$$\sum_{i,\alpha}q_{i\alpha}\phi_i\psi_\alpha$$
for any matrix $q_{i\alpha}$. In other words, matrices
of quadrics from this linear system have the form
$$\left(\matrix{0& Q\cr Q^t& 0}\right),$$
where $Q$ is an arbitrary matrix.

Summarizing, we see that
the projectivization of the tangent space at any point $x\in M$
is $\P^{n+m-1}$ and $|\FF^2|_x$ is the complete linear system
of quadrics having as a base locus the union $\P^{m-1}\cap\P^{n-1}$
of two skew linear subspaces.
\end{Example}

\begin{Example}[\cite{GH1}]\rm  \index{Veronese embedding}
Suppose that $\P^n=\P(V)\subset\P(S^dV)$ is the $d$-th Veronese
embedding. 
In terms of homogeneous coordinates $\{X_0,\ldots,X_n\}$ for $\P(V)$
the embedding is given by
$$(\ldots:X_\alpha:\ldots)\mapsto (X_0^d:X_0^{d-1}X_1:\ldots:X_0^{d-1}X_n:\ldots:F_\lambda(X):\ldots),$$
where $F_\lambda$ varies over a basis for the homogeneous forms of degree $d$.
Given a point $v_0\in\P^n$, e.g., $v_0=(1:\ldots:0)$,
$v_0$ maps to $(1:\ldots:0)$;
the projectivization of the tangent space at $v_0$
is $\P^{n-1}$ 
spanned by first $n$ coordinate vectors in $\P(S^dV)$.
The affine coordinate system in the affine chart containing $v_0$
is given by $F_\lambda(X)/X_0^d=F_\lambda(X_1/X_0,\ldots,X_n/X_0)$.
Therefore, $|\FF^2|_{v_0}$ is the linear system of all quadrics.
\end{Example}

\begin{Example}[\cite{GH1}]\rm
Let $\Gr(k,V)\subset\P(\Lambda^kV)$ be the Grassmanian\index{Pl\"ucker embedding}
in the Pl\"ucker embedding, $\dim V=n$.
Then the projectivization of the tangent space at any point $x$
is $\P(\Mat_{k,n-k})$ and an easy calculation (see \cite{GH1}) shows that
$|\FF^2|_x$ generated by $(2\times2)$-minors of this matrix.
Using this description, it is quite straightforward to find the defect
of the Grassmanian in the Pl\"ucker embedding.
All we need to do is to find the maximal possible rank
of a quadric generated by $(2\times2)$-minors of a rectangular matrix.
The answer is as follows. We may suppose that $k\le n/2$.
Then the dual variety is always a hypersurface
except the case $k=1$ (when the dual variety is empty)
and the case $k=2$, $n$ is odd (when the defect is equal to $2$).
\end{Example}

\begin{Example}[\cite{Lan4}]\rm\index{spinor variety}
Let $\SS_m\subset\P(\Lambda^{ev}\C^m)$ be the spinor variety
(c.f.~Example~\ref{SpinorExample}).
Then calculating the tangent space at the point $v_0=1$
it is easy to see that $T_{v_0}\SS_m=\Lambda^2\C^m$.
Let $\Gr(2,\C^m)\subset\P(\Lambda^2\C^m)$ be the Pl\"ucker embedding
then the linear system $|\FF^2|$ consists of all quadrics
containing $\Gr(2,\C^m)$.
\end{Example}

\begin{Example}[\cite{Lan4}]\rm   \index{Severi variety}
Consider the Severi variety, c.f.~Example~\ref{KJHVKJNGVKJNGVK}.
Let ${\Bbb D}_\R$ denote one of four division algebras
over $\R$: real numbers $\R$, complex numbers~$\C$, quaternions $\H$,
or octonions $\OO$. Let ${\Bbb D}={\Bbb D}_\R\otimes_\R\C$ denotes
its complexification, therefore $\Bbb D$ is either $\C$, $\C\oplus \C$,
$\Mat_2(\C)$, or $\Ca$ (the algebra of Cayley numbers).
All these algebras have the standard involution.
Let ${\cal H}_3({\Bbb D})$ denote the $3\times3$ Hermitian 
matrices over $\Bbb D$. If $x\in{\cal H}_3({\Bbb D})$
then we may write
$$x=\left(\matrix{
c_1&u_1&u_2\cr
\ov{u_1} &c_2& u_3\cr
\ov{u_2} & \ov{u_3}& c_3
}\right), \quad c_i\in\C,\ u_i\in{\Bbb D}.
$$
The projectivization of the highest weight vector orbit with respect to the group
of norm similarities is called the Severi variety corresponding to $\D$.
One may think of Severi varieties as of complexifications of 
projective planes over~${\Bbb D}_\R$.
Take a point $p$ of the Severi variety of the form
$$p=\left(\matrix{
1&0&0\cr 0&0&0\cr 0&0&0}\right).
$$
Choose the affine coordinates $u_1,u_2,u_3\in{\Bbb D}$, $c_2,c_3\in\C$
of ${\cal H}_3({\Bbb D})$ centered at~$p$.
Then the tangent space at $p$ is $\{u_1,u_2\}$
and the second fundamental form is generated by quadrics
$$u_1\ov{u_1}=0,\ u_2\ov{u_2}=0,\ \ov{u_2}u_1=0.$$
For example, if ${\Bbb D}=\Ca$ then the base locus of this linear
system is the spinor variety $\SS_5$ (c.f.~Example~\ref{KGFJHFDHGFDGF}).
\end{Example}

The higher fundamental forms\index{higher fundamental forms} 
are defined similar to the 
second fundamental form $\FF^2$.
The most geometric way is to utilize the
Euler--Meusnier Theorem~\ref{EulerMeusnierTh} as the definition.

Recall that if $\{p(t)\}$ is any holomorphic arc in $\P^N$
described by a vector--valued function $A_0(t)$, then the osculating
sequence at $p=p(0)$ is the sequence of linear spaces spanned by the
following collection of vectors  
$$\{A_0(t)\},\ \{A_0(t),A'_0(t)\},\ 
\{A_0(t),A'_0(t),A''_0(t)\},\ \ldots$$
Suppose now that $M\subset\P^N=\P(V)$ is a smooth variety, 
$p\in M$, $\tilde T_pM\subset V$ is the cone over the embedded
tangent space $\hat T_pM$.
Then the sequence of osculated spaces $T^{(0)}_pM=\C p$, \index{osculating space}
$T^{(1)}_pM=\tilde T_pM$, $T^{(2)}_pM$, $\ldots$
is defined as follows. The osculating space $T^{(k)}_pM$
is the linear span of $k$-th osculating spaces for smooth curves
$p(t)\subset M$, $p(0)=p$.
The $k$-th normal space $N^{(k)}_pM$ is defined as $V/T^{(k)}_pM$, so
we have $N^{(1)}_pM=(N_M\P^N)_p$.
Suppose that $v\in T_pM=T^{(1)}_pM/T^{(0)}_pM$ is any tangent vector,
$\{p(t)\}$ is any holomorphic arc in $\P^N$ with tangent vector $v$ at $t=0$.
Let $A(t)\subset V$ be a curve projected to $p(t)$.
Then it is easy to see that the vector 
$$\FF^k(v)={d^kA(t)\over dt^k}\mod T^{(k-1)}_pM\in N^{(k-1)}_pM$$
depends only on $v$. This construction gives a map
$$\FF^k:\,S^kT_pM\to N^{(k-1)}_pM$$
called the $k$-th projective fundamental form.
Of course, it can be also described via moving frames,
in coordinates, by differentiating Gauss map, etc.,
see e.g.~\cite{GH1,Lan4}.
The sequence of dimensions of osculating spaces $T^{(k)}_pM$
can, of course, depend on $p\in M$.
This will not be very essential, first, because we can restrict
ourselves to an open subset of $M$ where these dimensions are fixed.
Secondly, we shall be interested only in projective
homogeneous varieties in equivariant embeddings,
in which case all points have the same rights.

\begin{Example}\rm
It is easy to see that if $M\subset\P^N$ is a hypersurface
of degree~$d$, $p\in M$ is a smooth point, 
and $l\subset\P(T^{(d)}_pM)$
is a line passing through $p$ then $l\subset M$.
More generally, if $M\subset\P^N$ is a projective variety
whose homogeneous ideal is generated by elements
of degree at most $d$, $p\in M$ is a smooth point, 
and $l\subset\P(T^{(d)}_pM)$
is a line passing through $p$ then $l\subset M$.
If $V$ is an irreducible module of a simple group $G$
and $M=G/P\subset\P(V)$ is the projectivization of the orbit of the highest
weight then its ideal is generated in degree two (see e.g.  
\cite{Li}).   \index{lines on flag varieties}
In particular, for any $p\in M$ $\Base|\FF^2|_p$
is the set of tangent directions to lines on $M$ through $p$.
Using this idea, in \cite{LM} the thorough study
of varieties of lines and linear subspaces on minimal equivariant projective
embeddings of minimal flag varieties was undertaken.
\end{Example}

Let $G$ be a connected semisimple complex Lie group, 
$\g$ its Lie algebra, $U\g$ the universal enveloping algebra of $\g$. 
$U\g$ has a natural filtration $U\g^{(k)}$ such that
the associated graded algebra is the symmetric algebra $S\g$ of $\g$. 

Let $V_\lambda$
be an irreducible $G$-module with highest weight $\lambda$, and
$v_\lambda\in V_\lambda$ a highest weight vector. 
The action of $U\g$ on $V_\lambda$ 
induces a  filtration of $V_{\lambda}$ whose $k$-th term is  
$V_{\lambda}^{(k)}=U\g^{(k)}v_{\lambda}.$

Let $x=[v_{\lambda}]\in\P(V_\lambda)$ be the line in $V_{\lambda}$ generated by 
$v_{\lambda}$, and let $X=G/P\subset \P(V_{\lambda})$ be 
the projectivization of the orbit of the highest weight vector.
Here $P\subset G$ is a parabolic subgroup, namely,
the stabilizer of $x$. Let $\p$ be the Lie algebra of $P$.
The tangent bundle $TX$ is a homogeneous vector bundle 
and we can canonically identify $T_xX$ with 
$\g/\p$ as a $P$-module. 
Then it is easy to verify that the osculating spaces and 
the fundamental forms
of $X$ have a following simple representation-theoretic interpretation. 

\begin{Proposition} \index{fundamental forms of flag varieties}
Let $X=G/P \subset\P(V_\lambda)$ be a polarized flag variety 
and let $x=eP\in X$. Then the $k$-th osculating space
$T^{(k)}_{x}X = V_{\lambda}^{(k)}$, therefore
$N^{(k)}_xX=V_\lambda/V_{\lambda}^{(k)}$.
Moreover, there is a commutative diagram   
$$\begin{CD}
&S^k\g @=   U\g^{(k)}/U\g^{(k-1)}&\\  
&@V{\pi_1}VV   @VV{\pi_2}V&\\
S^k(T_xX)=&S^k(\g/\p) @>>{\FF^k}> V_\lambda/V_\lambda^{(k)}&=N^{(k-1)}_xX\\
\end{CD}
$$
where $\pi_1$ is induced by the projection $\g\to\g/\p$
and $\pi_2$ is the map given by $u\mapsto u\cdot v_\lambda$.
\end{Proposition}

\chapter{The Degree of a Dual Variety}

\section*{Preliminaries}

Let $X\subset\P(V)$ be a smooth projective variety, $\L=\O_X(1)$.
In this chapter we review several approaches for finding
the degree of the dual variety $\dual X$.
In Section~\ref{KJHGKJGFKJGFKJGF} we give several applications
of the Katz--Kleiman Holme formula.
In Section~\ref{KJHGFJHGFDUYDF} we collect various calculations
of degrees of discriminants of polarized flag varieties.
In final sections we study the degree of the discriminant of a polarized
variety as a function at polarization and express
the degree in terms of Hilbert polynomials.

\section{Katz--Kleiman--Holme Formula}\label{KJHGKJGFKJGFKJGF}

For any vector bundle $E$ on $X$ we denote by $c_i(E)$ its $i$-th Chern
class (see e.g.~\cite{Fu1}). \index{Chern class}
For example, $c_1(\L)=H$, the hyperplane
section divisor of $X$ in $\P(V)$.
For any $0$-dimensional cycle $Z$ on $X$ we denote its degree\index{degree of a cycle} 
(`the number of points' in $Z$) by $\int_XZ$.
Let $J(\L)$ be the first jet bundle of $\L$.
Let $\dim X=n$.

\begin{Theorem}[\cite{BFS1}]\label{MJGHVJGVJGHCVVC}
$\ $
\begin{enumerate}
\item[{\rm(a)}] $\deg\Delta_X=\int_Xc_n(J(\L))$. In particular,
$\dual X$ is a hypersurface if and only if $c_n(J(\L))\ne0$.
\item[{\rm(b)}] Moreover, $\defect X=k$ if and only if 
$c_r(J(\L))=0$ for $r\ge n-k+1$,
and \hbox{$c_{n-k}(J(\L))\ne0$}. In this case
$$\deg \dual X=\int_Xc_{n-k}(J(\L))\cdot H^k.$$
\end{enumerate}
\end{Theorem}

\sketch
Suppose that $\defect X=k$.
Notice that jets $j(f)$, $f\in V^*$, span $J(\L)$.
We choose generic elements $f_1,\ldots,f_v\in V^*$, where $v=\dim V^*$.
Using the classical representation of a Chern class of a spanned
vector bundle (see e.g.~\cite{Fu1}), we have that $c_{n-k+1}(J(\L))$
is represented by the set
$$\{x\in X\,|\,j(f_1)(x),\ldots,j(f_{k+1})(x)\ \hbox{\rm are linearly dependent}\}.$$
Since $\codim_{\P(V^*)}\dual X=k+1$, generic $k$-dimensional
projective subspace in $\P(V^*)$ does not meet $\dual X$.
Therefore, any non-trivial linear combination of 
$f_1,\ldots,f_{k+1}$ does not belong to the $\Cone(\dual X)$.
But this is equivalent to
$$\dim\Span\{j(f_1)(x),\ldots,j(f_{k+1})(x)\}=k+1$$
for any $x\in X$.
This implies by the above representation that $c_{n-k+1}(J(\L)=0$
and moreover that $c_r(J(\L))=0$ for $r\ge n-k+2$.
Now, the Chern class $c_{n-k}(J(\L))$ is represented by the set
$$\{x\in X\,|\,j(f_1)(x),\ldots,j(f_{k+2})(x)\ \hbox{\rm are linearly dependent}\}.$$
This set is a union of the singular loci of the hyperplane sections
$H\cap X$, where $H\in\dual X\cap L$ for a generic $(k+1)$-dimensional
subspace $L\subset\P(V^*)$. Notice that the number of these hyperplane
sections is equal to the degree of the dual variety $\deg\dual X$.
Moreover, by dimension reasons we may suppose that
$\dual X\cap L\subset\dual X_{sm}$.
By Theorem~\ref{FirstEin} these singular loci are projective subspaces
of dimension $k$. Therefore $c_{n-k}(J(\L))$ is represented
by $\deg X^*$ $k$-dimensional projective subspaces.
In particular, $c_{n-k}(J(\L))\ne0$ and
$\deg\dual X=\int_Xc_{n-k}(J(\L))\cdot H^k$.
\endproof

Using the exact sequence~(\ref{ExactSeqForJets})
it is possible to substitute Chern classes of the jet bundle
by Chern classes of the cotangent bundle $\Omega^1_X=T^*_X$.
The resulting formula was found by N.~Katz and S.~Kleiman
\cite{Ka,Kl1} in the case where $\dual X$ is a hypersurface
and by A.~Holme \cite{Ho1,Ho2} in the general case.\index{Katz--Kleiman--Holme formula}
Consider the Chern polynomial of $X$ with respect to the given projective
embedding  \index{Chern polynomial}
$$c_X(q)=\sum_{i=0}^nq^{i+1}\int_Xc_{n-i}(\Omega^1_X)\cdot H^i.$$

\begin{Theorem}\label{KleimanFormula}$\ $
\begin{enumerate}
\item[{\rm(a)}] $\deg\Delta_X=c_X'(1)=\sum_{i=0}^n(i+1)\int_Xc_{n-i}(\Omega^1_X)\cdot H^i$.
\item[{\rm(b)}] The codimension of $\dual X$ equals the order of the zero at $q=1$
of the polynomial $c_X(q)-c_X(1)$.
If this order is $\mu$ then $\deg \dual X=c_X^{(\mu)}(1)/\mu!$.
\end{enumerate}
\end{Theorem}

The formula in (a) can be rewritten in a numerous number of ways,
sometimes more convenient for calculations. The proofs of the 
following formulas can be found in \cite{Kl1}.

\begin{Example}[\cite{Kl1}]\rm
$$\deg\Delta_X=\int_Xs_n(N_X\P^N(-1)),$$
where $s_n$ is the $n$-th Segre class\index{Segre class} (see e.g.~\cite{Fu1}).
\end{Example}

\begin{Example}[\cite{Kl1}]\label{JHGFJHGFJGHFJH}\rm  \index{dual variety of a smooth complete intersection}
$$\deg\Delta_X=\int_X{1\over c(N^*_X\P^N(1))}.$$
For example, if $X$ is a hypersurface of degree $d$ in $\P^N$,
then $N^*_X\P^N(1)=\O_X(1-d)$
and 
$$\deg\Delta_X=\int_X{1\over 1-(d-1)H}=\int_X(d-1)^nH^n=d(d-1)^n.$$
In particular, $\dual X$ is a hypersurface (if $X$ is not a hyperplane).
In the same way one can show that $\dual X$ is a hypersurface
if $X$ is a non-degenerate smooth complete intersection (c.f.~Theorem~\ref{DualCompleteIntersection}).
\end{Example}

\begin{Example}[\cite{Kl1,Lam}]\label{KJFKJHCMHFCMHFD}\rm
The following non-trivial consequence of Theorem~\ref{KleimanFormula}
is called the class formula\index{class formula}.
$$\deg\Delta_X=(-1)^n\left[\chi(X)-2\chi(X\cap H)+\chi(X\cap H\cap H')\right],$$
where $H$ and $H'$ are generic hyperplanes (such that $X\cap H$ and $X\cap H\cap H'$
are smooth) and $\chi$ is a topological Euler characteristic.
This formula is the main ingredient in the proof of the following
Landman formula (see~\cite{Lan,Kl3}):\index{Landman formula}
\begin{eqnarray*}
\deg\Delta_X=\bigl(b_n(X)-b_{n-2}(X)\bigr)
+2\bigl(b_{n-1}(X\cap H)-b_{n-1}(X)\bigr)\\
+\bigl(b_{n-2}(X\cap H\cap H')-b_{n-2}(X)\bigr),
\end{eqnarray*}
where $b_i$ is the $i$-th Betti number.
All three summands are nonnegative numbers due to the
strong and weak Lefschetz theorems. Hence the characteristic condition
for the positive defect is
$$b_n(X)=b_{n-2}(X),\ b_{n-1}(X\cap H)=b_{n-1}(X),
\ b_{n-2}(X\cap H\cap H')=b_{n-2}(X).$$
\end{Example}

\begin{Example}[\cite{Kl1}]\label{KJHFJHFDJHFDJKDFJ}\rm
Let $X\subset \P^2$ be a smooth plane curve of degree $d$.
Consider the Veronese embedding $\P^2\subset\P^N$, $N+1={r+2\choose 2}$,
and the corresponding embedding $X\subset\P^N$.
Then the dual variety $\dual X\subset\dual{\P^N}$
parametrizes plane curves of degree $r$ tangent to $X$ at some point.
It is easy to see that
$$\deg\Delta_X=d(d+2r-3).$$
For example, the variety of conics tangent to a given conic is 
a hypersurface in $\P^5$ of degree $6$.
\end{Example}

\begin{Example}\label{DegDiscriminantOfACurve}\rm
Suppose that $X\subset\P^N$ is a smooth projective curve of genus $g$ and degree $d$.
Let $H\subset X$ be a hyperplane section. 
Then $\dual X$ is a hypersurface by
Example~\ref{DualOfTheCurve}, and its degree is given by
$$\deg\Delta_X=\deg K_X+2\deg H=2g-2+2d,$$
where $K_X$ denotes the canonical divisor.
\end{Example}

\begin{Example}\label{DegDiscriminantOfASurface}\rm
Let $X\subset\P^N$ be a smooth projective surface, $H\subset X$ be
a smooth hyperplane section. Then $\dual X$ is a hypersurface by
Example~\ref{DualOfASurface}, and its degree is given by
$$\deg\Delta_X=\deg c_2(\Omega^1_X)+2K_X\cdot H+3H\cdot H=\chi(X)-2\chi(H)+H\cdot H,$$
where $\chi$ is the topological Euler characteristic.
\end{Example}

\begin{Example}\label{BooleFormula}\rm
Suppose that $X=\P^n=\P(\C^{n+1})$ is embedded into $\P(S^d\C^{n+1})$ 
via the Veronese embedding. 
Then elements of $\dual X$ parametrize hypersurfaces
with singularities of degree $d$ in $\P^n$, therefore $\Delta_X$
is the classical discriminant of homogeneous forms of degree $d$
\index{classical discriminant}
in $n+1$ variables. Let us find its degree.
The Chow ring (or the cohomology ring) of $X$
is equal to $\Z(t)/t^{n+1}$, 
the element $t$ being the class of the hyperplane.
The first Chern class of the line bundle corresponding
to the Veronese embedding is equal to $dt$.
The total Chern class $\sum_ic_i(\Omega^1_{\P^n})$
equals $(1-t)^{n+1}$, 
therefore $c_i(\Omega^1_{\P^n})=(-1)^i{n+1\choose i}t^i$.
Clearly $\int_Xt^n=1$.
So we have
$$c_X(q)=\sum_{i=0}^n q^{i+1}(-1)^{n-i}{n+1\choose i+1}d^i=
{(qd-1)^{n+1}-(-1)^{n+1}\over d}.$$
Therefore, $\deg\Delta_X=c_X'(1)=(n+1)(d-1)^n$.
This formula is known as a Boole formula.
\end{Example}

A little bit more explicit form of Theorem~\ref{KleimanFormula}
can be given as follows.

\begin{Definition}
Let $X$ be a $n$-dimensional smooth projective variety
in $\P^N$, $\L=\O_X(1)$. Then the $s$-th rank of $X$\index{ranks of a projective variety}
with respect to the projective embedding is defined as
$$\delta_s=\sum_{i=s}^n{i+1\choose s+1}e_{n-i},$$
where $e_j$ denotes the degree of $c_j(T^*X)$ with
respect to the given projective embedding.
\end{Definition}

Then Theorem~\ref{KleimanFormula} can be reformulated
as follows

\begin{Corollary}\label{Kleiman-HolmeTh}
The defect $\defect X$ is the smallest integer $r\ge0$ such that
$\delta_r\ne0$, and for $r=\defect X$ the degree
of the dual variety $\dual X$ equals $\delta_r$.
\end{Corollary}

\section{Degrees of Discriminants of Polarized Flag Varieties}\label{KJHGFJHGFDUYDF}

\subsection{The Degree of the Dual Variety to $G/B$}
We follow the notation of Section~\ref{DefinitionsNotationsFlags}.
Let $G$ be a simple algebraic group of rank $r$ with Borel subgroup $B$.
The minimal equivariant projective embedding of $G/B$ corresponds
to the line bundle $\L_\rho$, where $\rho=\omega_1+\ldots+\omega_r$
is a half sum of positive roots.
We shall follow~\cite{CW} here and calculate the degree of the corresponding
discriminant (in fact, in most cases the dual variety will
be a hypersurface, see precise statements in Section~\ref{JKHGFJKHGDKH}).

Let $\{\beta_1,\ldots,\beta_N\}=\Delta_+$ be the collection
of all positive roots and let $\{\beta_1^\vee,\ldots,\beta_N^\vee\}=\Delta^\vee_+$
be the collection of all positive coroots.
We consider the matrix
$$M=\left({(\beta_i^\vee,\beta_j)\over (\beta_i^\vee,\rho)}\right)_{i,j=1,\ldots,N}.$$
Let $P_s(M)$ be the sum of all permanents of $s\times s$ submatrices of $M$.
Recall that the permanent of a $n\times n$ matrix $(a_{ij})$ is 
equal to $\sum\limits_{\sigma\in S_n}\prod\limits_{i=1}^na_{i,\sigma(i)}$.\index{permanent}

\begin{Theorem}[\cite{CW}]
$$\deg\Delta_{(G/B,\L_\rho)}=\sum_{s=0}^N(s+1)!P_{N-s}(M).$$
\end{Theorem}

\proof
Recall that by the Kleiman formula~\ref{KleimanFormula}
we have
$$\deg\Delta_{(G/B,\L_\rho)}=
\sum_{i=0}^N(i+1)\int_{G/B}c_{n-i}(\Omega^1_{G/B})\cdot c_1(\L_\rho)^i.$$
Therefore we need to calculate these integrals.
The correspondence $\lambda\mapsto\L_\lambda$ identifies
$\PP_\Q$ (the character lattice tensored by $\Q$)
with $\Pic(G/B)\otimes\Q$. Therefore,
we have a homomorphism of commutative algebras 
$$c:\,S^\bullet(\PP_\Q)\to A^\bullet(G/B),$$ where
$S^\bullet(\PP_\Q)$ is the symmetric algebra of $\PP_\Q$
and $A^\bullet(G/B)$ is the rational Chow ring.
It is well-known (see e.g.~\cite{BGG}) that
this homomorphism is surjective and its kernel
is generated as an ideal by $W$-invariants 
$(S^\bullet(\PP_\Q))^W_+$ of positive degree,
where $W$ is the Weyl group.
For example, 
$$c_1(\L_\rho)=c(\rho).$$
Using an appropriate filtration of the cotangent bundle $\Omega^1_{G/B}$
and `the splitting principle' it is easy to see
that 
$$c_{n-i}(\Omega^1_{G/B})=c(X_i),\quad
X_i=\sum_{\matrix{\scriptstyle\Gamma\subset\Delta_+\cr \scriptstyle\#\Gamma=n-i\cr}}
\prod_{\alpha\in\Gamma}\alpha.$$
Therefore, it is sufficient to compute
$c(\rho^iX_i)$.

For any coroot $\alpha^\vee\in\PP_\Q^*$ we define a linear
function 
$$D_{\alpha^\vee}(\lambda)=\langle\alpha^\vee,\lambda\rangle.$$
This linear function extends by a Leibniz rule to the differential
operator on $S^\bullet\PP_\Q$. We consider now
the differential operator $D$ on $S^\bullet\PP_\Q$
given by
$$D={\prod_{\alpha>0}D_{\alpha^\vee}\over\prod_{\alpha>0}(\alpha^\vee,\rho)}.$$
The operator $D$ decreases the degree by $N$,
so in particular we get a linear form
$$D:\,S^N\PP_\Q\to\Q.$$
One can show (see~\cite{CW}) that for any $x\in S^N\PP_\Q$ we have
$$
D(x)=\int_{G/B}c(x).
$$
Now the claim of the theorem follows by an easy calculation.
\endproof

\subsection{Degrees of Hyperdeterminants}
Consider the flag variety $\P^{k_1}\times\ldots\times\P^{k_r}$
of the group $\SL_{k_1}\times\ldots\times\SL_{k_r}$.
The projectively dual variety of its `minimal' equivariant 
projective embedding corresponds to a nice theory
of hyperdeterminants that was initiated by Cayley and Schl\"affli 
\cite{Ca1,Ca2,Schl}.
\index{hyperdeterminants}

Let $r\ge2$ be an integer, and 
$A=(a_{i_1\ldots i_r})$, $0\le i_j\le k_j$
be an $r$-dimensional complex matrix of format $(k_1+1)\times\ldots\times(k_r+1)$.
The hyperdeterminant of $A$ is defined as follows, c.f.~\ref{HyperdeterminantDefs}.
Consider the product $X=\P^{k_1}\times\ldots\times\P^{k_r}$
of several projective spaces embedded into the projective
space $\P^{(k_1+1)\times\ldots\times(k_r+1)-1}$ via the Segre embedding.
Let $\dual X$ be the projectively dual variety.
If $\dual X$ is a hypersurface then it is defined
by a corresponding discriminant $\Delta_X$, which
in this case is called the hyperdeterminant (of format
$(k_1+1)\times\ldots\times(k_r+1)$), denoted by $\Det$.
Clearly $\Det(A)$ is a polynomial function in matrix entries
of $A$ invariant under the action of he group $\SL_{k_1+1}\times\ldots\times\SL_{k_r+1}$.
If $\dual X$ is not a hypersurface then we set $\Det=1$.
Theorem~\ref{ExistenceOfDeterminantsTh} shows that
$\dual X$ is a hypersurface (and, hence, defines a hyperdeterminant)
if and only if $2k_j\le k_1+\ldots+k_r$ for $j=1,\ldots,r$.
If for one of $j$ we have an equality $2k_j=k_1+\ldots+k_r$
then the format is called boundary. \index{boundary format}%
Let $N(k_1,\ldots,k_r)$ be the degree of the hyperdeterminant
of format $(k_1+1)\times\ldots\times(k_r+1)$.
The proof of the following theorem can be found in \cite{GKZ2} or \cite{GKZ3}.

\begin{Theorem}$\ $
\begin{enumerate}
\item[{\rm(a)}] The generating function for the degrees $N(k_1,\ldots,k_r)$
is given by
$$\sum_{k_1,\ldots,k_r\ge0}N(k_1,\ldots,k_r)z_1^{k_1}\ldots z_r^{k_r}=
{1\over
\left(
1-\sum\limits_{i=2}^r(i-1)e_i(z_1,\ldots,z_r)
\right)^2
},$$
where $e_i(z_1,\ldots,z_r)$ is the $i$-th elementary
symmetric polynomial.
\item[{\rm(b)}] The degree $N(k_1,\ldots,k_r)$ of the boundary format
is {\rm(}assuming that $k_1=k_2+\ldots+k_r${\rm)}
$$N(k_2+\ldots+k_r,k_2,\ldots,k_r)={(k_2+\ldots+k_r+1)!\over
k_2!\ldots k_r!}.$$
\item[{\rm (c)}]
The degree of the hyperdeterminant of the cubic format is given by
$$N(k,k,k)=\sum_{0\le j\le k/2}{(j+k+1)!\over (j!)^3(k-2j)!}\cdot 2^{k-2j}.$$
\item[{\rm(d)}]
The exponential generating function for the degree $N_r$
of the hyperdeterminant of format $2\times2\times\ldots\times2$ {\rm(}$r$ times{\rm)}
is given by
$$\sum_{r\ge0}N_r{z^r\over r!}={e^{-2z}\over(1-z)^2}.$$
\end{enumerate}
\end{Theorem}

\subsection{One Calculation}
In most circumstances, known formulas for the degree of the discriminant
depend on the certain set of discrete parameters.
There are two possibilities for these parameters.
First, we may fix a projective variety and vary its
polarizations. The degree of the discriminant
as a function in polarization is studied in the next section.
Second possibility is to change a variety and to fix
a polarization (in a certain sense).
For example, consider the irreducible representation
of $\SL_{n_0}$ with the highest weight $\lambda$
Then this weight can be considered as a highest weight
of $\SL_n$ for any $n\ge n_0$ with respect to the natural
embedding $\SL_n\subset\SL_{n+1}\subset\ldots$.
As a result, we shall obtain a tower of flag varieties
with ``the same'' polarization.
The degree of the corresponding discriminants will
be a function in $n$.
This function can be very complicated.
For example, in \cite{Las} a very involved combinatorial
formula for this degree was found in the case $\lambda$
is a fundamental weight
(i.e.~for the degree of a dual variety
of the Grassmanian $\Gr(k,n)$ in the Pl\"ucker embedding).
However, sometimes this function has a closed expression.
For example, the Boole formula~\ref{BooleFormula}
is a formula of this kind.
Another simple example is the formula~\ref{Pfaffian}
for the degree of the dual variety to $\Gr(2,n)$
in the Pl\"ucker embedding.
The following theorem is the mixture of these two cases.
We are not aware of any other similar formulas.

\begin{Theorem}[\cite{T1}]\label{JHGFJGHFJ}
Let $V$ be an irreducible $\SL_n$-module with the highest
weight $(a-1)\varphi_1+\varphi_2$, $a\ge2$. Then the variety
$\dual X\subset\P(V^*)$ projectively dual to the projectivization $X$
of the orbit of the highest vector is a hypersurface of degree
$${(n^2-n)a^{n+1}-(n^2+n)a^{n-1}-2n(-1)^n\over (a+1)^2}.$$
\end{Theorem}

\proof
We have to calculate the degree of~$\Delta_X$.
We use Kleiman's formula~\ref{KleimanFormula}
for the degree of the dual variety.
In the present case $X=G/P$, where $G=\SL_n$ and
$P\subset G$ is the parabolic subgroup of matrices
\begin{equation}
\left(\matrix{
*&0&0&\ldots&0\cr
*&*&0&\ldots&0\cr
*&*&*&\ldots&*\cr
\vdots&\vdots&\vdots&\ddots&\vdots\cr
*&*&*&\ldots&*\cr
}\right).\label{eqno(9)}
\end{equation}
Assume that $T\subset G$ is the diagonal torus,
$B$ is the Borel subgroup of lower-triangular matrices,
$x_1,\ldots,x_n$ are the weights of the tautological representation,
$X(T)$ is the lattice of characters of $T$, 
$S$ is the symmetric algebra of $X(T)$ (over~$\Bbb Q$),
$W\simeq S_n$ is the Weil group of $G$,
and $W_P\simeq S_{n-2}$ is the Weil group of~$P$.
It is well known (see~\cite{BGG}) that the map $c:\,X(T)\to\Pic(G/B)$
that assigns to~$\lambda$ the first Chern class of the invertible sheaf
${\cal L}_\lambda$ can be extended
to a surjective homomorphism
$c:\,S\to A^*(G/B)$ in the (rational) Chow ring,
and its kernel coincides with $S^W_+S$.
The projection $\alpha:\,G/B\to G/P$ induces an embedding 
$\alpha^*:\,A^*(G/P)\to A^*(G/B)$.
The image coincides with the subalgebra
of $W_P$-invariants.
Hence, $A^*(G/P)=S^{W_P}/S^{W}_+S^{W_P}$. 
We denote the homomorphism
$S^{W_P}\to A^*(G/P)$ by the same letter~$c$.

To apply Kleiman's formula we need
$c_1(\cal L)$ (which is equal to $c(ax_1+x_2)$) and the total Chern class of
$\Omega_Z^1$, which is equal to
$$c(\Omega_Z^1)=c\left((1-x_1+x_2)\prod_{i=3,\ldots,n}(1-x_1+x_i)
\prod_{i=3,\ldots,n}(1-x_2+x_i)\right).$$
(This can be shown by standard arguments using the filtration of~$\Omega_Z^1$
and splitting principle.)

Let $\alpha_1,\ldots,\alpha_{n-2}$ be the elementary symmetric polynomials in
$x_3,\ldots,x_{n}$. Then
$S^{W_P}=\Bbb Q[x_1,x_2,\alpha_1,\ldots,\alpha_{n-2}]$,
and the ideal $S^W_+S^{W_P}$ is generated by
$$\matrix{
&x_1+x_2+\alpha_1,\quad x_1x_2+x_1\alpha_1+x_2\alpha_1+\alpha_2,\cr
&x_1x_2\alpha_1+x_1\alpha_2+x_2\alpha_2+\alpha_3,\ \ldots,\ 
x_1x_2\alpha_{n-4}+x_1\alpha_{n-3}+x_2\alpha_{n-3}+\alpha_{n-2},\cr
&x_1x_2\alpha_{n-3}+x_1\alpha_{n-2}+x_2\alpha_{n-2},\quad x_1x_2\alpha_{n-2}.
}$$
Hence,
$\alpha_i=(-1)^i(x_1^i+x_1^{i-1}x_2+\ldots+x_2^i) \mod S^W_+S^{W_P}$, and
$A^*(G/P)$  is isomorphic to the quotient ring
$\Bbb Q[x_1, x_2]/\langle f_1, f_2\rangle$, where
$f_1=x_1^{n-1}+x_1^{n-2}x_2+\ldots+x_2^{n-1}$ and
$f_2=x_1^n$. Note that $f_1,f_2$ is the Gr\"obner basis of the ideal
$\langle f_1, f_2\rangle$ with respect to the ordering
$x_2>x_1$ (see~\cite{Berg}).
Hence, the set of $X^iY^j$, where $X=x_1 \mod \langle f_1, f_2\rangle$, 
$Y=x_2 \mod \langle f_1, f_2\rangle$, $i=1,\ldots,n-1$, $j=1,\ldots,n-2$,
is a basis of the quotient algebra.

To calculate the degree of the discriminant by Kleiman's formula,
we have to calculate
$\int_Zc(X^{n-1}Y^{n-2})$. 
Let $\tilde w_0$ be the longest element in $W$, and let
$w_0$ be the shortest element in $\tilde w_0W_P$ 
with the reduced factorization
\begin{eqnarray*}
w_0=\left(\matrix{
1&  2&3&4&\ldots&n\cr
n&n-1&1&2&\ldots&n-2\cr
}\right)\cr
=(n-1,n)(n-2,n-1)\ldots(12)\cdot(n-1,n)(n-2,n-1)\ldots(23).
\end{eqnarray*}

Let 
$$A_{w_0}=A_{(n-1,n)}A_{(n-2,n-1)}\ldots A_{(12)}
A_{(n-1,n)}A_{(n-2,n-1)}\ldots A_{(23)}$$ be the corresponding endomorphism of degree
$-(2n-3)$ in $S$, where
$$A_{(ij)}={\displaystyle id-s_{(ij)} \over\displaystyle x_i-x_j},$$ and
$s_{(ij)}$ is the reflection that transposes
$x_i$ and $x_j$.
Then
$$\int_Zc(X^{n-1}Y^{n-2})=A_{w_0}(x_1^{n-1}x_2^{n-2})$$ 
(see~\cite{BGG}).
It is obvious that
\begin{eqnarray*}
A_{w_0}(x_1^{n-1}x_2^{n-2})=
A_{\rho_1}A_{\rho_2}(x_1^{n-1}x_2^{n-2})\cr
=
A_{\rho_1}(x_1^{n-1}A_{\rho_2}(x_2^{n-2}))=
A_{\rho_1}(x_1)A_{\rho_2}(x_2),
\end{eqnarray*}
where $\rho_k=(n-1,n)(n-2,n-1)\ldots(k,k+1)$.
These factors are both equal to~$1$, since they are equal to
$\int_{\Bbb P^{n-1}}c_1({\cal O}(1))^{n-1}$ and
$\int_{\Bbb P^{n-2}}c_1({\cal O}(1))^{n-2}$, respectively,
and it is obvious that these integrals are equal to~1.
We finally obtain that $\int_Zc(X^{n-1}Y^{n-2})=1$.

It remains to calculate the polynomial
\begin{equation}
\sum_{i=0}^{2n-3}(i+1)c_{2n-3-i}(aX+Y)^i,\label{eqno(10)}
\end{equation}
in the ring
$\Bbb Q[X,Y]$ (with the basis $X^iY^j$, $i=1,\ldots,n-1$, $j=1,\ldots,n-2$,
and relations $X^{n-1}+X^{n-2}Y+\ldots+Y^{n-1}=0$, $Y^n=0$, and
$X^n=0$, which follows from the preceding relations), where
$c_k$ is the $k$th homogeneous component of the polynomial
$$(1-X+Y)\prod_{i=3}^n(1-X+x_i)
\prod_{i=3}^n(1-Y+x_i),$$
in which the $i$th symmetric function of
$x_3,\ldots,x_n$ must be replaced by 
$(-1)^i(X^i+X^{i-1}Y+\ldots+Y^i)$.
By the previous discussion, the result of this calculation will be
$\deg(\Delta_X)X^{n-1}Y^{n-2}$.

Note that the polynomial~\ref{eqno(10)} is equal to
$$F'(T)|_{T=2X+Y}=(F_1F_2F_3F_4)'(T)|_{T=2X+Y},$$
where
\begin{eqnarray*}
F_1=T,\quad F_2=T-X+Y,\cr
F_3=\prod_{i=3}^n(T-X+x_i)=\sum_{i=0}^{n-2}(T-X)^{n-2-i}(-1)^i(X^i+\ldots+Y^i),\cr
F_4=\prod_{i=3}^n(T-Y+x_i)=\sum_{i=0}^{n-2}(T-Y)^{n-2-i}(-1)^i(X^i+\ldots+Y^i).
\end{eqnarray*}
We have, further,
\begin{eqnarray*}
F_3(T)=
\left(\sum_{i=0}^{n-2}(T-X)^{n-2-i}(-1)^i
{X^{i+1}-Y^{i+1}\over X-Y}\right)=\cr
{X(T-X)^{n-2}\over X-Y}\left(\sum_{i=0}^{n-2}(T-X)^{-i}(-1)^iX^{i}\right)
-\qquad\qquad{}\cr
{Y(T-X)^{n-2}\over X-Y}
\left(\sum_{i=0}^{n-2}(T-X)^{-i}(-1)^iY^{i}\right)=\cr
{X\left((T-X)^{n-1}-(-X)^{n-1}\right)\over T(X-Y)}-
{Y\left((T-X)^{n-1}-(-Y)^{n-1}\right)\over (T-X+Y)(X-Y)}.
\end{eqnarray*}
We obtain, likewise, that
\begin{eqnarray*}
F_4(T)=
\left(\sum_{i=0}^{n-2}(T-Y)^{n-2-i}(-1)^i
{X^{i+1}-Y^{i+1}\over X-Y}\right)=\cr
{X(T-Y)^{n-2}\over X-Y}
\left(\sum_{i=0}^{n-2}(T-Y)^{-i}(-1)^iX^{i}\right)
-\qquad\qquad{}\cr
{Y(T-Y)^{n-2}\over X-Y}
\left(\sum_{i=0}^{n-2}(T-Y)^{-i}(-1)^iY^{i}\right)=\cr
{X\left((T-Y)^{n-1}-(-X)^{n-1}\right)\over (T+X-Y)(X-Y)}-
{Y\left((T-Y)^{n-1}-(-Y)^{n-1}\right)\over T(X-Y)}.
\end{eqnarray*}
We deduce from the latest formula that
\begin{eqnarray*}
F'(T)=
{n(X-Y)(T-X)^{n-1}+(-X)^n-(-Y)^n\over X-Y}\times
\qquad\qquad\cr
{(X-Y)(T-Y)^n+T(-X)^n-(T+X-Y)(-Y)^n\over T(T+X-Y)(X-Y)}+\cr
{(X-Y)(T-X)^n+(T-X+Y)(-X)^n-T(-Y)^n\over T(X-Y)}\times
\qquad\qquad\cr
{\biggl(n(T+X-Y)-(T-Y)\biggr)(T-Y)^{n-1}+(-X)^n 
\over (T+X-Y)^2}-\cr
{(X-Y)(T-X)^n+(T-X+Y)(-X)^n-T(-Y)^n\over T(X-Y)}\times
\qquad\qquad\cr
{(X-Y)(T-Y)^n+T(-X)^n-(T+X-Y)(-Y)^n\over T(T+X-Y)(X-Y)}.
\end{eqnarray*}
Using the elementary formulae
\begin{eqnarray*}
{\alpha^n(\gamma-\beta)+\beta^n(\alpha-\gamma)+\gamma^n(\beta-\alpha)\over
(\alpha-\beta)(\beta-\gamma)(\gamma-\alpha)}
=\sum_{i+j+k=n-2}\alpha^i\beta^j\gamma^k,\cr
{\alpha^n(\gamma-\beta)+\beta^n(\alpha-\gamma)+\gamma^n(\beta-\alpha)\over
(\alpha-\beta)(\gamma-\beta)}
=\sum_{i=0}^{n-2}\beta^i(\alpha^{n-1-i}-\gamma^{n-1-i}),\cr
{(n(\alpha-\beta)-\alpha)\alpha^{n-1}+\beta^{n}\over
(\alpha-\beta)^2}=\sum_{i=0}^{n-2}(n-1-i)\alpha^{n-2-i}\beta^i,
\end{eqnarray*}
we obtain
\begin{eqnarray*}
F'(T)=\biggl(n(T-X)^{n-1}+(-1)^n\sum_{i=0}^{n-1}X^{n-1-i}Y^i\biggr)\qquad\qquad\cr\times
\biggl(\sum_{i+j+k=n-2}(T-Y)^i(-X)^j(-Y)^k\biggr)+\cr
\biggl(\sum_{i=0}^{n-2}(-X)^i((T-X)^{n-1-i}-(-Y)^{n-1-i})\biggr)\qquad\qquad\cr\times
\biggl(\sum_{i=0}^{n-2}(n-1-i)(T-Y)^{n-2-i}(-X)^i\biggr)-\cr
\biggl(\sum_{i=0}^{n-2}(-X)^i((T-X)^{n-1-i}-(-Y)^{n-1-i})\biggr)\qquad\qquad\cr\times
\biggl(\sum_{i+j+k=n-2}(T-Y)^i(-X)^j(-Y)^k\biggr).
\end{eqnarray*}
Putting $T=aX+bY$ in the latest formula, we obtain

\begin{eqnarray*}
\biggl(\sum_{i=0}^{n-1}
\bigl(1+(-1)^nn(a-1)^{n-1-i}b^i{n-1\choose i}\bigr)X^{n-1-i}Y^i\biggr)\qquad\qquad\cr\times
\biggl(\sum_{i=0}^{n-2}\schoose{n-2}{i}_{a,b-1}X^iY^{n-2-i}\biggr)-\cr
\biggl(\sum_{i=0}^{n-1}\bigl(\schoose{n-1}{i}_{a-1, b}-
\schoose{n-2}{i-1}_{a-1, b}-1\bigr)X^{n-1-i}Y^i\biggr)\qquad\qquad\cr\times
\biggl(\sum_{i=0}^{n-2}
\bigl(\fchoose{n-2}{i}_{a, b-1}-\fchoose{n-3}{i}_{a, b-1}\bigr)
X^{n-2-i}Y^i\biggr)+\cr
\biggl(\sum_{i=0}^{n-1}\bigl(\schoose{n-1}{i}_{a-1, b}-
\schoose{n-2}{i-1}_{a-1, b}-1\bigr)X^iY^{n-1-i}\biggr)\qquad\qquad\cr\times
\biggl(\sum_{i=0}^{n-2}\schoose{n-2}{i}_{a, b-1}X^iY^{n-2-i}\biggr),
\end{eqnarray*}
where
$${\schoose nk}_{xy}=
\sum_{p=0}^k\sum_{q=0}^{n-k}(-1)^{p+q}x^py^q{p+q\choose p},$$
$${\fchoose nk}_{xy}=
\sum_{p=0}^k\sum_{q=0}^{n-k}(p+q+1)(-1)^{p+q}x^py^q{p+q\choose p}.
$$

It remains to calculate the degree of the discriminant, that is,
the difference between the coefficients of
$X^{n-1}Y^{n-2}$ and $X^{n-2}Y^{n-1}$ in the expression
for $F'(aX+bY)$. After some transformations we obtain
\begin{eqnarray*}
\sum_{i=0}^{n-1}
\left((-1)^nn(a-1)^{n-1-i}b^i{n-1\choose i}
+W_i\right)\times
\qquad\qquad\qquad\cr
\left(\schoose{n-2}{i}_{a,b-1}-\schoose{n-2}{i-1}_{a,b-1}\right)-\cr
\sum_{i=0}^{n-1}
\left(W_i-1\right)
\times\qquad\qquad\qquad\qquad\qquad\qquad
\qquad\qquad\qquad\cr
\left(\fchoose{n-2}{i}_{a, b-1}-\fchoose{n-3}{i}_{a, b-1}
-\fchoose{n-2}{i-1}_{a, b-1}+\fchoose{n-3}{i-1}_{a, b-1}\right),
\end{eqnarray*}
where
$$W_i=\schoose{n-1}{n-1-i}_{a-1, b}-\schoose{n-2}{n-1-i}_{a-1, b}
.$$
This is the degree of the discriminant of the irreducible
$\SL_n$-module with the highest weight $(a-b)\varphi_1+b\varphi_2$.

Substituting $b=1$ in the last formula, we obtain
$${(n^2-n)a^{n+1}-(n^2+n)a^{n-1}-2n(-1)^n\over (a+1)^2}.$$
\endproof

\section{Degree of the Discriminant as a Function in $\L$}
Suppose that $X$ is a smooth projective variety.
If $\L$ is a very ample line bundle on $X$ then we call
the pair $(X,\L)$ a polarized variety. \index{polarized variety}
Any polarized variety admits a canonical embedding
$X\subset\P(H^0(X,\L)^*)$.
Let $\dual{(X,\L)}\subset\P(H^0(X,\L))$ be the corresponding
dual variety, which parametrises singular divisors
in the linear system $|\L|$.
If $\dual{(X,L)}$ is a hypersurface then the discriminant $\Delta_{(X,\L)}$
is a unique (up to a scalar) polynomial defining $\dual{(X,\L)}$.
If $\defect(X,\L)>0$ then we set $\Delta_{(X,\L)}=1$.

\subsection{General Positivity Theorem}

We start with the following simple but very useful observation
\begin{Theorem}\label{DefectInVeronese}.
Suppose that $\L$, $\M$ are very ample line bundles on $X$.
Then $\defect(X,\L\otimes\M)=0$.
\end{Theorem}

\proof
Indeed, we have
$\deg\Delta_{(X,\L\otimes\M)}=
\int_Xc_n(J(\L\otimes\M)$ by Theorem~\ref{MJGHVJGVJGHCVVC}.
Therefore,
$$\deg\Delta_{(X,\L\otimes\M)}=
\int_Xc_n(J(\L)\otimes\M)=
\sum_{i=0}^{n} (i+1)\int_X c_1(\M)^i\cdot c_{n-i}(J(\L)).$$
Since $J(\L)$ is spanned, all summands are non-negative.
Since $\int_Xc_1(\M)^n>0$ (being equal to the degree of $X$
in the embedding determined by $\M$),
the whole sum is positive.
\endproof

For example, if $X$ is embedded in $\P(V)$
and then re-embedded in $\P(S^dV)$ via the Veronese embedding
then the dual variety $\dual X$ of this re-embedding
is a hypersurface.

Suppose now that $\L_1,\ldots,\L_r$ are 
line bundles on $X$ such that the corresponding
linear systems $|\L_i|$ have no base points
and such that for each $i$ there is a representative
of the linear system $|\L_i|$ which is a smooth divisor on $X$.
Suppose further that any line bundle $\L$ of the form
$\L=\L_1^{\otimes n_1}\otimes\ldots\otimes\L_r^{\otimes n_r}$
is very ample for positive $n_i$.
Then the degree $\deg\Delta_{(X,\L)}$ is a function in $n_1,\ldots,n_r$.
Let us introduce new non-negative integers $m_i=n_i-1$.
Then $\deg\Delta_{(X,\L)}=f(m_1,\ldots,m_r)$.

\begin{Theorem}[\cite{CW}]\label{GeneralCWPositivity}$\ $
\begin{enumerate}
\item[{\rm(a)}]
The function $f$ is a non-trivial polynomial
with non-negative coefficients.
\item[{\rm(b)}] If each $n_i\ge2$ then $\dual{(X,\L)}$ is a hypersurface.
\end{enumerate}
\end{Theorem}

\proof
Let $n=\dim X$.
Recall that by Kleiman--Katz--Holme formula~\ref{KleimanFormula}
we have
$$\deg\Delta_{(X,\L)}=\sum_{i=0}^n(i+1)\int_Xc_{n-i}(\Omega^1_X)\cdot \L^i=P'(\L),$$
where 
$$P(t)=\sum_{i=0}^nc_{n-i}(\Omega^1_X)t^{i+1}.$$
This clearly implies that $f$ is a polynomial in $n_i$ and, hence, in $m_i$.
This polynomial is non-trivial by Theorem~\ref{DefectInVeronese}.
Therefore, we need to check positivity only.
Using the Taylor expansion, it suffices to check
the following claim. 
For each $j=1,\ldots,n$ and for each collection of nonnegative
integers $u_1,\ldots,u_r$ such that $u_1+\ldots+u_r=j-1$, we have
$$\int_X P^{(j)}(\L)\cdot \L_1^{u_1}\ldots\L_r^{u_r}\ge0.$$
The conditions imposed on $\L_i$ imply that any intersection
of the form $\L_1^{u_1}\ldots\L_r^{u_r}$ (as an element of a Chow ring)
is represented by a smooth subvariety of codimension $j-1$.
Therefore, it suffices to check that for any such
subvariety $Y\subset X$ of codimension $j-1$
we have
$$\int_X P^{(j)}(\L)\cdot [Y]\ge0.$$
Using exact sequence~\ref{ExactSeqForJets}
it is easy to see that
$$P^{(j)}(\L)=j!c_{n+1-j}(J(\L)).$$
Therefore it remains to show that
$$\int_X c_{n+1-j}(J(\L))\cdot [Y]\ge0.$$
This follows from the fact that $J(\L)$ is spanned
by its global sections (see~\cite{Fu1}).
\endproof

\subsection{Applications to Polarized Flag Varieties}
The conditions of Theorem~\ref{GeneralCWPositivity}
are perfectly well satisfied if $X$ is a flag variety $G/P$.
We shall use notations from~\ref{DefinitionsNotationsFlags}.
The fundamental weights $\omega_1,\ldots\omega_k$
dual to the simple roots in $\Pi_{G/P}$ generate the lattice
$\PP_{G/P}$ isomorphic to $\Pic(G/P)$.
We shall denote the line bundle $L_{\omega_i}$ by $\L_i$.
Then any line bundle on $G/P$ have a form 
$\L=\L_1^{\otimes n_1}\otimes\ldots\otimes\L_r^{\otimes n_r}$.
$\L$ is very ample if and only if all $n_i>0$.
Then the degree $\deg\Delta_{(X,\L)}$ is a function in $n_1,\ldots,n_r$.
Let us introduce new non-negative integers $m_i=n_i-1$.
Then $\deg\Delta_{(X,\L)}=f(m_1,\ldots,m_r)$.
All conditions of Theorem~\ref{GeneralCWPositivity}
are clearly satisfied, therefore $f$ is a non-trivial
polynomial with non-negative coefficients.

\begin{Example}\rm
Let $V=\C^r$, $G=\SL(V)$ and let $P$ be the stabiliser of the line in $V$.
Then $G/P=\P(V)$. Take $\lambda=n\omega_1$, $n>0$. Then $\L_1=\O(n)$.
The degree of $\Delta_{(G/P,\lambda)}$ is the degree of the classical
discriminant of a homogeneous form of degree $n$ in $r$ variables.
Therefore by Example~\ref{BooleFormula}  \index{classical discriminant}%
we have $\deg\Delta_{(G/P,\lambda)}=r(n-1)^{r-1}$.
Therefore the polynomial $f$ in this case has a form $f=rm^{r-1}$,
so it is a polynomial with nonnegative coefficients.
\end{Example}

\begin{Example}\rm
Let $G=SL_3$ and let $P=B$ be the Borel subgroup.
The weight $\lambda$ can be written as $\lambda=m_1\omega_1+m_2\omega_2$.
One can check using the Kleiman formula that
$$\deg\Delta_{(G/B,\lambda)}=12(m_1^2m_2+m_1m_2^2)-6(m_1^2+4m_1m_2+m_2^2)
+12(m_1+m_2)-6.$$
Substituting $n_1=m_1-1$, $n_2=m_2-1$ we see that
$$\deg\Delta_{(G/B,\lambda)}=12(n_1^2n_2+n_1n_2^2)+6(n_1^2+4n_1n_2+n_2^2)
+12(n_1+n_2)+6.$$
Again, we get a polynomial with nonnegative coefficients.
\end{Example}

\begin{Example}\rm
Assume that the weight $\lambda$ satisfies the conditions
$m_i\ge2$ (is very strictly dominant). Then the results of this section
clearly imply that the dual variety $\dual{(G/P,\lambda)}$ 
is a hypersurface.
\end{Example}

\section{Gelfand--Kapranov Formula}

Theorem~\ref{DiscriminantCayley} identifies $\Delta_X$
with the determinant of any of the discriminant complexes
$C_\pm^*(X,\M)$, provided the complex is stably twisted.
By Corollary~\ref{CayleyDetDegree} 
this implies a formula for $\deg\Delta_X$.

\begin{Theorem}[\cite{GKZ2}]\label{DegreeDiscrimComplex}
If the complex $C_-^*(X,\M)$ is stably twisted then
$$\deg\Delta_X=\sum_{i=-\dim X-1}^0(-1)^i\cdot i\cdot\dim C^i_-(X,\M).$$
If the complex $C_+^*(X,\M)$ is stably twisted then
$$\deg\Delta_X=(-1)^{\dim X+1}\sum_{i=0}^{\dim X+1}(-1)^i\cdot i\cdot\dim C^i_+(X,\M).$$
\end{Theorem}

Thus $\deg\Delta_X$ can be expressed through the dimensions of vector
spaces
$$H^0\left(X,\Lambda^{-i}J(\L)^*\otimes\M\right)\ \hbox{\rm or}
\ H^0\left(X,\Lambda^{i}J(\L)\otimes\M\right).$$
It is possible to rewrite these dimensions
using simpler quantities associated with $X$, at least
in case $\M=\O_X(l)$, $l\gg0$.
For any coherent sheaf $\F$ on $X$, we write
$$\F(l)=\F\otimes\O_X(l),\ h^i(\F)=\dim H^i(X,\F),\ \chi(\F)=\sum(-1)^ih^i(\F).$$
The number $\chi(\F)$ is called the Euler characteristic of $\F$.\index{Euler characteristic}
It follows from the Riemann-Roch-Hirzebruch theorem \cite{Hir}
that $\chi(\F(l))$ is a polynomial in $l$.
This polynomial is called the Hilbert polynomial of $\F$\index{Hilbert polynomial}
and denoted by $h_\F(l)$. For $l\gg0$, the higher cohomology of $\F(l)$
vanish and we have 
$$h_\F(l)=\chi(\F(l)=h^0(\F(l)).$$
This is a more usual definition of Hilbert polynomials, see \cite{Ha1}.
Simple calculations show that
Theorem~\ref{DegreeDiscrimComplex} can be reformulated as follows.
Recall that $T_X$ is the tangent bundle of $X$ and
$\Omega^i_X=\Lambda^iT^*_X$
is the bundle of differential $i$-forms on $X$.

\begin{Theorem}[\cite{GKZ2}]
For any $l\in\Z$ we have
$$\deg\Delta_X=\sum_{i=-\dim X-1}^0(-1)^i\cdot i\cdot\left(
h_{\Lambda^{-i}T_X}(i+l)+h_{\Lambda^{-i-1}T_X}(i+l)\right),$$
$$\deg\Delta_X=(-1)^{\dim X+1}\sum_{i=0}^{\dim X+1}(-1)^i\cdot i\cdot\left(
h_{\Omega^i_X}(i+l)+h_{\Omega^{i-1}_X}(i+l)\right).$$
\end{Theorem}

This theorem was generalized in \cite{GK} to handle the
case of the positive $\defect X$.
For any $l\in\Z$ consider the polynomial $f_l(q)$ in a formal variable $q$
$$f_l(q)=\sum_{i=0}^{\dim X+1}
(-1)^i\left(h_{\Lambda^iT_X}(i+l)+h_{\Lambda^{i-1}T_X}(i+l)\right)q^i,$$
where all polynomials that have no sense are supposed to be equal to $0$.
Clearly $f_l(q)$ is also a polynomial in $l$.

\begin{Theorem}[\cite{GK,GKZ2}]
Let $l\in\Z_+$ be any non-negative integer.
The codimension of $\dual X$ equals the order of the zero
of $f_l(q)$ at $q=1$. Let this order be $\mu$
and let $f_l(q)=a_\mu(l)(1-q)^\mu+O((1-q)^{\mu+1}).$
Then $a_\mu(l)={\mu+l-1\choose \mu-1}\cdot\deg\dual X$.
\end{Theorem}

\chapter{Milnor Classes and Multiplicities of Discriminants}\label{MilnorClassesSection}

\section*{Preliminaries}
The notion of the Milnor class was introduced in \cite{Al}
generalizing Milnor numbers of isolated singularities \cite{Mi}
and Parusi\`nski $\mu$-numbers \cite{Pa1}.
We review some applications of this technique to dual varieties,
in particular for the description of their singularities.
We also give other results in this direction.
In this chapter we assume that the reader
is more familiar with modern intersection theory than
is necessary for other chapters of this book. But
to illustrate how the topology can be used for
study of dual varieties, we give
an elementary proof of the class formula in the first section
(cf.~Example~\ref{KJFKJHCMHFCMHFD}).

\section{Class Formula}

Suppose that $X\subset\P^N$ is a smooth projective $n$-dimensional
variety with the dual variety $X^*\subset(\P^N)^*$. 
Consider a generic line $L\subset(\P^N)^*$ such that
the intersection $L\cap X^*$ consists of $d$ smooth on $X^*$ points
if $\defect X=0$, where $d$ is the degree of $X^*$.
If $\defect X>0$ then this intersection is empty
and we set $d=0$.

$L$ can be considered as a pencil (one-dimensional linear system)
of divisors on $X$, therefore
it defines a rational map $F:\,X\to\P^1$.
More precisely, if $L$ is spanned
by hyperplanes $H_1$ and $H_2$  with equations $f_1=0$ and $f_2=0$
then $F$ has a form $x\mapsto (f_1(x):f_2(x))$ for any $x\in X$.
Thus $F$ is not defined along the subvariety
$X\cap H_1\cap H_2$ of codimension $2$. 
We shall blow it up and get the variety $\tilde X$ 
and the regular morphism $\tilde F:\,\tilde X\to\P^1$.
Let $D\subset\tilde X$ be a preimage of $X\cap H_1\cap H_2$.

Now let us calculate 
the topological Euler characteristic
$\chi(\tilde X)$ in two ways.
First, 
$$\begin{array}{r}
\chi(\tilde X)=\chi(\tilde X\setminus D)+\chi(D)=
\chi(X\setminus X\cap H_1\cap H_2)+\chi(D)=\cr
\chi(X\setminus X\cap H_1\cap H_2)+2\chi(X\cap H_1\cap H_2)=
\chi(X)+\chi(X\cap H_1\cap H_2).
\end{array} $$
Here we use the fact that $D$ is a $\P^1$-bundle over
$X\cap H_1\cap H_2$, therefore 
$\chi(D)=2\chi(X\cap H_1\cap H_2)$.
Notice that the blowing up is an isomorphism
on the complement of the exceptional divisor, hence
$\chi(\tilde X\setminus D)=\chi(X\setminus X\cap H_1\cap H_2)$.

On the other hand, we may use $\tilde F$ to calculate $\chi(\tilde X)$.
For $x\in\P^1$, we have $\chi(\tilde F^{-1}(x))=\chi(X\cap H_x)$,
where $H_x$ is a hyperplane corresponding to $x$.
Let $x_1,\ldots,x_d\in\P^1$ be points that correspond
to the intersection of $L$ with $\dual X$.
Then for any $x\in\P^1\setminus\{x_1,\ldots,x_d\}$,
$X\cap H_x$ is a smooth divisor.
For any $x\in\{x_1,\ldots,x_d\}$,
$X\cap H_x$ has a simple quadratic singularity
(see Theorem~\ref{IntermediateEinTh}).
It is clear then that $\chi(X\cap H_x)=\chi(X\cap H_y)=\chi(X\cap H)$,
where $x,y\in\P^1\setminus\{x_1,\ldots,x_d\}$ and $H$ 
is a generic hyperplane.
For $x\in\{x_1,\ldots,x_d\}$, near the
simple quadratic singularity of $X\cap H_x$,
the family of divisors $X\cap H_y$, $y\to x$, looks
like the family of smooth  $(n-1)$-dimensional quadrics $Q_\eps$
with equations $T_1^2+\ldots+T_{n-1}^2=\eps T_n^2$ near 
the unique singular point $(0:\ldots:0:1)$ of the quadric 
$Q_0$ with the equation $T_1^2+\ldots+T_{n-1}^2=0$.
Therefore, 
$$\chi(\tilde X)=\chi(\P^1)\chi(X\cap H)+d(\chi(Q_0)-\chi(Q_1))=
2\chi(X\cap H)+d(\chi(Q_0)-\chi(Q_1)).$$

Combining two formulas for $\chi(\tilde X)$, we get
$$\chi(X)+\chi(X\cap H_1\cap H_2)=2\chi(X\cap H)+d(\chi(Q_0)-\chi(Q_1)).$$
The Euler charachteristic
of $k$-dimensional smooth quadric is equal to $k+2$ for $k$ even
and to $k+1$ for $k$ odd.
The Euler charachteristic of a $k$-dimensional quadric
with a unique singularity is equal to $k+1$ for $k$ even and to $k+2$ for $k$ odd,
since it is a cone over $(k-1)$-dimensional quadric.
Finally, we obtain the class formula   \index{class formula}
$$d=(-1)^n\left[\chi(X)-2\chi(X\cap H)+\chi(X\cap H\cap H')\right].$$
This formula has a large number of applications, see
Example~\ref{KJFKJHCMHFCMHFD}, \ref{KJHFJHFDJHFDJKDFJ},
\ref{DegDiscriminantOfACurve}, \ref{DegDiscriminantOfASurface}, etc.
In this chapter we consider
various attempts to generalise this formula
in different directions.

\section{Milnor Classes}

Suppose that $M$ is a smooth $n$-dimensional algebraic variety, $\L$
is a line bundle on $M$, and $X$ is the zero-scheme of a section of $\L$.
The singular scheme of $X$, $\Sing X$, is the scheme supported
on the singular locus of $X$, and defined locally by the Jacobian ideal
$\left({\partial F\over x_1},\ldots,{\partial F\over x_n}\right)$,
where $x_1,\ldots,x_n$ are local parameters for~$M$, and $F$ is the section
of $\L$ defining $X$.\index{singular scheme}\index{Jacobian ideal}

\begin{Definition}
Let $Y$ be the singular scheme of a section of a line bundle
$\L$ on a smooth variety $M$. The Milnor class of $Y$ with respect to $\L$
is the class 
$$\mu_\L(Y)=c(T^*M\otimes\L)\cap s(Y,M)$$
in the Chow group $A_*Y$ of $Y$.
\end{Definition}\index{Milnor class}

Here $c$ denotes the total Chern class and $s$ is the Segre class
(in the sense of \cite{Fu1}).
The following theorem summarizes the main properties
of the Milnor class. The reader may consult \cite{Al}
for proofs and further details.

\begin{Theorem}[\cite{Al}]\label{MilnorClassPropertiesTh}$\ $
\begin{enumerate}
\item[\rm(a)] The Milnor class $\mu_\L(Y)$ only depends
on $Y$ and $\L|_Y$. That is, the definition of the Milnor class
is independent on the choice of a smooth ambient variety $M$
in which $Y$ is realized as the singular scheme of a section
of a line bundle restricting to $\L|_Y$.
\item[\rm(b)] Suppose that $\L$ is generated by global sections
and let $X_g$ be a generic section of $\L$. Then
$$\mu_\L(Y\cap X_g)=c_1(\L)\cap\mu_\L(Y).$$
\item[\rm(c)] Suppose that $Y$ is non-singular. Then
$$\mu_\L(Y)=c(T^*Y\otimes \L)\cap[Y].$$
\end{enumerate}
\end{Theorem}

\begin{Example}\label{NilnorClassPropeerties}\rm
If $Y$ is supported on a point $P$, then $\mu_\L(Y)=m[P]$,
where $m$ is the Milnor number of $X$ at $P$. \index{Milnor number}
Indeed, in this case $\mu_\L(Y)=s(Y,M)$, so $m$ is the coefficient
of $[P]$ in $s(Y,M)$, which agrees with the definition
of the Milnor number (see \cite{Fu1}).
\end{Example}

\begin{Example}\rm
Even if $\dim Y>0$, let $|Y|$ be the support variety of $Y$, assume
that $Y$ is complete. Then in \cite{Pa1} there is introduced an invariant,
the $\mu$-number of $X$ at $|Y|$, agreeing with the ordinary \index{$\mu$-number}
Milnor number in case $|Y|$ is a point. We claim that in fact
this number equals the degree of the $0$-dimensional component 
of $\mu_\L(Y)$. Namely, the differential of the local equations
of $X$ in $M$ determines a section of $(T^*M\otimes\L)|_X$.
This can be extended to a holomorphic section $s_X$ of $T^*M\otimes\L$
over the whole of $M$. Then the $\mu$-number is defined in \cite{Pa1}
as the contribution of $|Y|$ to the intersection 
of $s_X$ with the zero section of $T^*M\otimes\L$.
In the neighborhood of $|Y|$ we have the fiber diagram
$$\begin{CD}
Y@>>> M\\
@VVV @VV{s_X}V\\
M@>{s_0}>> T^*M\otimes\L
\end{CD}$$
that is, $Y$ is the scheme-theoretic intersection
of the two sections, so the contribution of $|Y|$
to the intersection number of the sections is equal (see \cite{Fu1})
the degree of the $0$-dimensional component of
$$c(T^*M\otimes\L)\cap s(Y,M),$$
which is the claim.
In fact, in view of Theorem~\ref{NilnorClassPropeerties},
the numerical information carried by the Milnor class of $Y$ (when $Y$ is complete)
is essentially equivalent to $\mu$-numbers of $|Y|$
and its generic sections $|Y\cap X_g|$, $|X\cap X_{g_1}\cap X_{g_2}|$, $\ldots$.
\end{Example}

\section{Applications to Multiplicities of Discriminants}
We shall now apply the machinery of the previous section
to our usual situation. Suppose that $X\subset\P^N$
is a smooth projective variety, $\L=\O_X(1)$, $\dual X\subset\dual{\P^N}$
is the dual variety parametrizing hyperplanes $H$ such that
the singular scheme (or the contact scheme in this case)
$\Sing X\cap H$ is not empty.
Then according to the previous section we have the Milnor class
$\mu_\L(\Sing X\cap H)\in A_*(\Sing X\cap H)$.
The problem that we address here is how to calculate
the multiplicity of $\dual X$ at $H$ in terms of the Milnor class.
\index{multiplicities of the discriminant}

\begin{Theorem}[\cite{Al,AC,Di1,N,Pa2}]\label{MilnorDefectTh}
Let $X\cap H$ be any singular hyperplane section of $X$. Then
\begin{enumerate}
\item[\rm(a)] The defect of the dual variety $\dual X$
is the smallest integer $r\ge0$ such that
$$\int c_1(\L)^rc(\L)\cap\mu_\L(\Sing X\cap H)\ne0;$$
and for $r=\defect X$ this number equals the multiplicity $m_H\dual X$
of $\dual X$ at~$H$.
\item[\rm(b)] If $\defect X=r$ then $m_H\dual X$ is equal to
$$\int (-1)^{r}\mu_\L(\Sing X\cap H)+c_1(\L)^{r+1}\cap\mu_\L(\Sing X\cap H).$$
\end{enumerate}
\end{Theorem}

\proof
(a)
Induction on $r$. If $r=0$ then the integral is equal to
\begin{eqnarray*}
\int c(\L)\cap\mu_\L(\Sing X\cap H)=\int c(\L)c(T^*X\otimes\L)\cap s(\Sing X\cap H,X)\\
=\int c((T^*X\oplus\O)\otimes\L)\cap s(\Sing X\cap H,X).\\
\end{eqnarray*}
Thinking of $H$ as a point of $\dual X$, we have
$s(H,\dual X)=(m_H\dual X)X^*$, 
thus 
$$m_H\dual X=\deg s(H,\dual X).$$
The dual variety is the projection from the conormal variety $I_X$,
$\dual X=\pi(I_X)$. Since $r=0$, the degree of this projection
is equal to $1$. The fiber $Y$ of $\pi$ over $H$ is identified with
the contact locus $\Sing H\cap X$ (as a scheme). Therefore,
we have $s(H,\dual X)=\pi_*(Y, I_X)$. Taking degrees, we get
$$m_H\dual X=\deg s(Y,I_X).$$
Now using the second projection $I_X\to X$, which is the projective
bundle, it is quite easy to see that
$$s(Y,I_X)=c(J(\L))\cap s(\Sing X\cap H,X),$$
where $J(\L)$ is the jet bundle of $\L$.
It remains to notice that $c(J(\L))=c((T^*X\oplus\O)\otimes\L)$
by the exact sequence (\ref{ExactSeqForJets}).

Now let $r>0$. 
By the similar argument as above, we have $\int c(\L)\cap\mu_\L(\Sing X\cap H)=0$.
By Theorem~\ref{DefectOfDivisorTh}, if $\dual X$ is not a hypersurface
then the dual variety $\dual{X'}$ of a generic hyperplane section
$X'=X\cap H_g$ is the cone over $\dual X$ with $H_g$ as a vertex.
The multiplicity of $\dual X$ at $H$ then evidently
equals the multiplicity of $\dual{X'}$  at $H$; but $\defect X'=\defect X-1$,
so the statement follows by induction because
\begin{eqnarray*}
\int c_1(\L)^rc(\L)\cap\mu_\L(\Sing X\cap H)\\
=\int c_1(\L)^{r-1}c(\L)\cap(c_1(\L)\cap\mu_\L(\Sing X\cap H))\\
=\int c_1(\L)^{r-1}c(\L)\cap\mu_\L(\Sing X'\cap H),
\end{eqnarray*}
where we used Theorem~\ref{MilnorClassPropertiesTh}
in the last equality.

(b) By (a) we see that $\defect X\ge r$ if and only if the components
of dimension $i$, $0\le i<r$ of $c(\L)\cap\mu_\L(\Sing X\cap H)$ vanish,
for all hyperplane sections $X\cap H$.
That is, if $\defect X\ge r$, then for all hyperplane sections $X\cap H$
$$c(\L)\cap\mu_\L(\Sing X\cap H)=A_{\dim\Sing X\cap H}+\ldots+A_r$$
with $A_j$ a class in dimension $j$, $j=r,\ldots,\dim\Sing X\cap H$
(depending on $H$).
Therefore,
$$\mu_\L(\Sing X\cap H)=c(\L)^{-1}\cap(A_{\dim\Sing X\cap H}+\ldots+A_r)$$
(with $A_r\ne0$), and 
$$\int (-1)^r\mu_\L(\Sing X\cap H)=\sum_{i=r}^{\dim\Sing X\cap H}(-1)^{i-r}
c_1(\L)^i\cap A_i$$
while
$$\int c_1(\L)^{r+1}\mu_\L(\Sing X\cap H)=
\sum_{i=r+1}^{\dim\Sing X\cap H}(-1)^{i-r+1}
c_1(\L)^i\cap A_i.$$
Hence,
$$\int (-1)^r\mu_\L(\Sing X\cap H)+c(\L)^{r+1}\cap\mu_\L(\Sing X\cap H)=
\int c_1(\L)^r\cap A_r,$$
which coincides with the formula in (a).
\endproof

In case $\defect X=0$ and $\Sing X\cap H$ is finite,
the theorem takes the most simple form:
\begin{Corollary}[\cite{Di1,N}]\label{KJGVHGFDXHGDK}
If $\defect X=0$, $H\in\dual X$, and $\Sing X\cap H$ is finite we have
$$m_H\dual X=\sum_{p\in\Sing X\cap H}\mu(X\cap H,p),$$
the sum over ordinary Milnor numbers. In particular, if $\Sing X\cap H$
is finite then $H$ is smooth on $\dual X$ if and only if
$\Sing X\cap H=\{p_0\}$ and 
the Hessian of the function $f$ defining $X\cap H$ 
is non-zero at $p_0$.
\end{Corollary}

\begin{Remark}\rm
An easy Chern class computation shows that ranks 
(see Section~\ref{KJHGKJGFKJGFKJGF})
could be expressed
via Milnor classes 
as follows:
$$\delta_i=\int_Xc_1(\L)^ic(\L)\cap\mu_\L(X).$$
Moreover, Theorem~\ref{Kleiman-HolmeTh} can be viewed as
a particular case of Theorem~\ref{MilnorDefectTh}
applied to the zero-section of $H^0(X,\L)$.
\end{Remark}

\section{Multiplicities of the Dual Variety of a Surface}
Let $X$ be a smooth projective algebraic surface, $\L$ a very ample
line bundle on $X$, $\dual X\subset|\L|$ the dual variety,
consisting of singular members of $|\L|$. By Example~\ref{DualOfASurface},
$\dual X$ is a hypersurface. We want to determine, for each $D\in|\L|$,
the multiplicity $m_D\dual X$ of the hypersurface $\dual X$ at the point $D$.
Suppose that the divisor
$$D=\sum_{i=1}^rn_iD_i,$$
where any divisor $D_i$ is reduced and irreducible.
Then we have $D_{red}=\sum D_i$.

\begin{Theorem}[\cite{Al}]\label{MultiplicityDualOfSurfaceThe}
With notation as above, we have
$$m_D\dual X=(D-D_{red})\cdot (K_X+2D+D_{red})+\sum_{p\in D_{red}}\mu(D_{red},p),$$
where $K_X$ is the canonical divisor of $X$ and $\mu$ stands for the ordinary Milnor number.
\end{Theorem}

\proof
Take a general $H\in|\L|$ (i.e.~a $H$ intersecting $D_{red}$ transversally)
and denote by $L\subset|\L|$ the pencil containing $D$ and $H$. Then
\begin{equation}
m_D\dual X=\deg\dual X-s,\label{JHKHBKHBKUY}
\end{equation}
where $s$ is the number of singular members of $L$ different from $D$
(each of this singular members has one node as a singular locus).
In order to determine $s$ we shall blow up $X$ to construct 
a family parametrized by $L$ and use Lefschetz's formula \cite{GH1}.

For each $i$ ($1\le i\le r$) denote by $p_{ij}$ ($1\le j\le H\cdot D_i$)
the points of intersection of $H$ and $D_i$.
Let $\hat X$ denote the surface obtained from $X$ by blowing up 
each point $p_{ij}$ $n_i$ times (in the direction of $H$).
The induced pencil $\hat L$ on $\hat X$ is base point free and gives
a map $f:\,\hat X\to\P^1$. If $E_{ij}^k$, $k=1,\ldots,n_i$ are the exceptional 
divisors at $p_{ij}$ then the fiber of $f$ at the point $0$
(corresponding to $D$) is
$$f^*(0)=D+\sum_{i,j}\sum_{1\le k\le n_i-1}(n_i-k)E_{ij}^k.$$
In other words, the special fiber is isomorphic to $D$
with strings of $\P^1$'s (each $\P^1$ with a certain multiplicity)
attached at the points $p_{ij}$; each string has $n_i-1$ components.

We now denote
$$D'=(f^*(0))_{red}=D_{red}+\sum_{i,j}T_{ij}$$
the reduced fiber of $f$ at $0$, where $T_{ij}=\sum_{1\le k\le n_i-1}E_{ij}^k$
is the reduced string attached at $p_{ij}$.
Then the standard argument involving Lefschetz's formula (see \cite{GH1})
gives
$$\chi(X)=2\chi(H)-H\cdot H+(\chi(D')-\chi(H))+\sum_{1\le\lambda\le s}(\chi(H_\lambda)-\chi(H)),$$
where $\chi$ denotes the topological Euler characteristic
and $H_\lambda$ are the singular fibers of $f$ for $\lambda\ne0$.
Since $\chi(H_\lambda)-\chi(H)=1$ (see (\ref{HGFJGFJGF}) below)
and $\deg\dual X=\chi(X)-2\chi(H)+H\cdot H$ (Example~\ref{DegDiscriminantOfASurface}),
we can rewrite (\ref{JHKHBKHBKUY}) as
\begin{equation}
m_D\dual X=\chi(D')-\chi(H).\label{Hogmanay1}
\end{equation}
In order to compute $\chi(D')$, denoting $T=\cap_{i,j}T_{ij}$ we have
\begin{equation}
\chi(D')=\chi(D_{red}\cap T)=\chi(D_{red})+\chi(T)-\chi(D_{red}\cap T)\label{Hogmanay2}
\end{equation}
\begin{eqnarray*}
=\chi(D_{red})+\sum_{i,j}\chi(T_{ij})-\sum_{i,j}\chi(\{p_{ij}\})\\
=\chi(D_{red})+\sum_{i,j}n_i-\sum_{i,j}1=\chi(D_{red})+D\cdot(D-D_{red}).
\end{eqnarray*}
Now we compute $\chi(D_{red})$. Let $Z=\sum_{1\le i\le r}Z_i$
be a reduced connected curve with normalization
$$\rho:\,\tilde Z=\mathop{\sqcup}_{1\le i\le r}\tilde Z_i\to Z.$$
If $p\in Z$ is a singular point, denote $B(p)=\rho^{-1}(p)$
the set of branches at $p$. Topologically,
$Z$ is obtained from the smooth real surface $\tilde Z$
by identifying each of the sets $B(p)$ to a point $p$.
Therefore, we have
\begin{equation}
\chi(\tilde Z)-\chi(Z)=\sum_{p\in Z}(b(Z,p)-1),\label{HGFJGFJGF}
\end{equation}
where $b(Z,p)$ is the number of branches of $Z$ at $p$.
Also, from the exact sequence of sheaves
$$0\to\O_Z\to\rho_*\O_{\tilde Z}\to \rho_*\O_{\tilde Z}/\O_Z\to 0$$
we obtain
$$1-r+\sum_{p\in Z}\delta(Z,p)-h^1(\O_Z)+h^1(\O_{\tilde Z})=0,$$
where we set $\delta(Z,p)=\hbox{\rm length}_p(\rho_*\O_{\tilde Z}/\O_Z)$,
Combining this with~(\ref{HGFJGFJGF}) we obtain
\begin{equation}
\chi(Z)=2-2p_a(Z)+\sum_{p\in Z}\mu(Z,p),\label{Hogmanay3}
\end{equation}
where
$$\mu(Z,p)=2\delta(Z,p)-b(Z,p)+1$$
is the Milnor number of $(Z,p)$.
Now writing together (\ref{Hogmanay1}), (\ref{Hogmanay2}), and (\ref{Hogmanay3})
we get
\begin{eqnarray*}
m_D\dual X=\chi(D_{red})-\chi(H)+D\times(D-D_{red})\\
=(2-2p_a(D_{red}))-(2-2p_a(H)+D\cdot(D-D_{red})+\sum_{p\in D_{red}}\mu(D_{red},p)\\
=-(K_X+D_{red})\cdot D_{red} +(K_X+H)\cdot H\\
+D\cdot (D-D_{red})+\sum_{p\in D_{red}}\mu(D_{red},p).
\end{eqnarray*}
After rearranging we finally obtain the claim of the Theorem.
\endproof

\begin{Example}\rm
Suppose that $X=\P^2$, $\L=\O_X(d)$.
Then $\dual X$ parametrizes plane singular curves of degree $d$.  
If $C$ is such curve then applying Theorem~\ref{MultiplicityDualOfSurfaceThe}
we get
$$m_C\dual X=\left(3(d-1)-d'\right)d'+\mu,$$
where $d'$ is the degree of $C-C_{red}$ and $\mu$
is the sum of the Milnor numbers 
of the singularities of $C_{red}$.
For example, the multiplicity of the discriminant of plane conics
at a double line is $2$.
\end{Example}

\section{Further Applications to Dual Varieties}

Theorem~\ref{PairityTheorem} shows that if $\defect X>0$
then the defect and the dimension have the same parity.
Since the defect equals the dimension of a contact
locus with a generic tangent hyperplane, 
this can be reformulated by saying that if the generic
contact locus is positive-dimensional then it has the same
parity as $\dim X$.
It is possible to generalize this result for not necessarily
generic contact loci using the machinery of Milnor classes. 
The proof of the following
theorem can be found in \cite{Al}.
\index{parity theorem}

\begin{Theorem}[\cite{Al}]
Let $X\subset\P^N$ be a smooth projective variety.
Let $H\in\dual X$ be such that the contact locus $Y=\Sing X\cap H$
is pure-dimensional and non-singular in a neighborhood
of some complete curve $C\subset Y$ of odd degree.
Then 
$$\dim Y\equiv\dim X \mod 2.$$
\end{Theorem}

For example, if $Y$ is a generic contact locus then $Y=\P^{\defect X}$.
Therefore, if $Y$ is not a point then it contains a line $C$
and we get Theorem~\ref{PairityTheorem}.

The following result follows quite immediately from 
Theorem~\ref{MilnorDefectTh}:

\begin{Theorem}[\cite{Al}]\label{GJFHGFHGFGH}
Let $X\subset\P^N$ be a smooth projective variety.
Suppose that the contact scheme $\Sing X\cap H$ of a hyperplane $H$
with $X$ is a projective space $\P^r$. Then the defect $\defect X$ is equal to $r$
and the dual variety $\dual X$ is smooth at $H$.
\end{Theorem}

\index{contact scheme}
Recall that by the contact scheme we mean the contact locus
endowed of the scheme structure defined by the Jacobian ideal.
The surprising element of this result is that we are not assuming
the hyperplane to be generic. The dimension of the contact
scheme of the hyperplane need not be equal the dimension
of the contact locus of the general hyperplane.
The claim, however, is that it is so if the contact scheme
is a projective subspace as a scheme.
Moreover, it can be shown (also using Theorem~\ref{MilnorDefectTh})
that if the contact scheme of a variety with a given hyperplane is smooth
then the defect of the variety equals the defect of the contact scheme.
This `explains' the previous theorem, because the defect
of a projective subspace equals its dimension.

In case $\defect X=0$ Theorem~\ref{GJFHGFHGFGH} can be formulated
a little bit more explicitly.
Let $X\subset\P^n$ be a smooth variety with dual variety $\dual X$
being a hypersurface. Let $I_X\subset\P^n\times\dual{\P^n}$
be an conormal variety. We denote by $\pi$ the projection $I_X\to\dual X$.
Let $U\subset\dual X$ be the set of all points,
where the projection $\pi$ is unramified.
Clearly, $U$ is non-empty and open in $\dual X$.
In other words, $H\subset U$ if the contact locus $\Sing X\cap H$
consists of a single point $p$ scheme-theoretically
(the Hessian of the function $f$ defining $X\cap H$ 
is non-zero at $p$). Then we have the following

\begin{Theorem}[\cite{Ka}]
Suppose that $X\subset\P^N$ is smooth and $\dual X$ is a hypersurface.
$U$ consists of smooth points of $\dual X$ and is the biggest
open set in $\dual X$ for which the projection $\pi:\,\pi^{-1}(U)\to U$
is an isomorphism.
\end{Theorem}

In fact, one may expect that $U=X_{sm}$. For example,
this is so if $\dual X\setminus U$ has codimension $1$.

\begin{Proposition}
Suppose that $X\subset\P^N$ is smooth and $\dual X$ is a hypersurface.
If $Z=\dual X\setminus U$ is a hypersurface in $\dual X$
(possibly reducible) then $Z={\dual X}_{sing}$.
\end{Proposition}

\proof
By a previous theorem ${\dual X}_{sing}\subset Z$.
We need to prove an opposite inclusion.
Let $Y\subset\dual X$ be the variety of all points $x\in\dual X$
such that $\dim\pi^{-1}(x)>0$. Then $Y\cap U=\emptyset$,
therefore $\codim_{\dual X}Y\ge2$.
Therefore, if $Z_0\subset Z$ is an irreducible component
then for generic point $H\in Z_0$ the contact locus
$\Sing H\cap X$ consists of finitely many points.
Therefore, $H\in{\dual X}_{sing}$ by Corollary~\ref{KJGVHGFDXHGDK}.
\endproof

\index{singular locus of the dual variety}
The description of the singular locus of the dual variety
$\dual X$ is known only in some special cases.
The most detailed study was performed for hyperdeterminants in~\cite{WZ}.
To formulate this result, we'll need the following notations.
Suppose first that $X\subset\P^N$ is a smooth variety.
We shall define two subvarieties ${\dual X}_{cusp}$ and ${\dual X}_{node}$
in $\dual X$. The variety ${\dual X}_{node}$ is the closure
of the set of hyperplanes $H$ such that $H$ is tangent to $X$ at two distinct points.
In other words, ${\dual X}_{node}=\ov{\pr_1(I_X^2)}$, where
$$I_X^2=\{(H,x,y)\in\dual X\times X\times X\,|\,x\ne y,\ \hat T_xX,\hat T_yX\in H\}.$$
The variety ${\dual X}_{cusp}$ consists of all hyperplanes $H$
such that there exists a point $x\in \Sing X\cap H$ that
is not a simple quadratic singularity. That is, 
the Hessian of the function $f$ defining $X\cap H$ 
is equal to zero at $x$.
Suppose now that $X=\P^{k_1}\times\ldots\P^{k_r}$
in the Segre embedding. 
By Theorem~\ref{ExistenceOfDeterminantsTh} 
$\dual X$ is a hypersurface if and only if
$$k_j\le\sum_{i\ne j}k_i,\quad j=1,\ldots,r.$$
In this case the defining equation of $\dual X$
is called the hyperdeterminant of the matrix format
$(k_1+1)\times\ldots\times(k_r+1)$.
Without loss of generality we may suppose that $k_1\ge k_2\ge\ldots\ge k_r$.
Then in fact we have only one inequality, namely $k_1\le k_2+\ldots+k_r$.
We call the format $(k_1+1)\times\ldots\times(k_r+1)$
boundary if $k_1=k_2+\ldots+k_r$ and interior if
$k_1<k_2+\ldots+k_r$.
For the hyperdeterminants, the variety ${\dual X}_{node}$
can be decomposed further.
Namely, for any subset $J\subset\{1,\ldots,r\}$
we set
${\dual X}_{node}(J)=\ov{\pr_1(I_X^2(J))}$, where
$$I_X^2(J)=\{(H,x,y)\in I_X^2\,|\,x^{(j)}=y^{(j)}\hbox{\rm\ if and only if\ }
j\in J\}.$$

\index{singular locus of hyperdeterminants}
\index{hyperdeterminants}
\index{matrix format}

The proof of the following theorem can be found in \cite{WZ}.
\begin{Theorem}[\cite{WZ}]
Suppose that $k_1\le k_2+\ldots+k_r$, $k_1\ge\ldots\ge k_r$.
\begin{enumerate}
\item[\rm(a)] If the format is boundary then ${\dual X}_{sing}$ is an irreducible
hypersurface in $\dual X$. Furthermore, we have 
${\dual X}_{sing}={\dual X}_{node}(\emptyset)$ if the format is different from
$4\times2\times2$. In this exceptional case ${\dual X}_{sing}={\dual X}_{node}(1)$. 
\item[\rm(b)] If the format is interior and does not belong to the following
list of exceptions, then ${\dual X}_{sing}$ has two irreducible
components ${\dual X}_{cusp}$ and ${\dual X}_{node}(\emptyset)$,
both having codimension $1$ in $\dual X$. The list of exceptional cases
consists of the following $3$- and $4$-dimensional formats:
\begin{enumerate}
\item[\rm(1)] For $2\times2\times2$ matrices ${\dual X}_{sing}={\dual X}_{cusp}$
and has three irreducible components, all of codimension $2$ in $\dual X$.
\item[\rm(2)] For $3\times2\times2$ matrices ${\dual X}_{sing}=
{\dual X}_{cusp}$ and is irreducible.
\item[\rm(3)] For $3\times3\times3$ matrices ${\dual X}_{sing}$
has five irreducible components
${\dual X}_{cusp}$, ${\dual X}_{node}(\emptyset)$, 
${\dual X}_{node}(1)$, ${\dual X}_{node}(2)$, ${\dual X}_{node}(3)$.
\item[\rm(4)] For $m\times m\times3$ matrices with $m>3$
${\dual X}_{sing}$
has three irreducible components
${\dual X}_{cusp}$, ${\dual X}_{node}(\emptyset)$, and ${\dual X}_{node}(3)$.
\item[\rm(5)] For $2\times2\times2\times2$ matrices
${\dual X}_{sing}$
has eight irreducible components
${\dual X}_{cusp}$, ${\dual X}_{node}(\emptyset)$, and ${\dual X}_{node}(i,j)$
for $1\le i<j\le 4$.
\item[\rm(6)] For $m\times m\times2\times2$ matrices with $m>2$
${\dual X}_{sing}$
has three irreducible components
${\dual X}_{cusp}$, ${\dual X}_{node}(\emptyset)$, and ${\dual X}_{node}(3,4)$. 
\end{enumerate}
In each of the cases {\rm(2)} to {\rm(6)}
all irreducible components of ${\dual X}_{sing}$
have codimension $1$.
\end{enumerate}
\end{Theorem}

\chapter{Mori Theory and Dual Varieties}\label{MoriChapter}

\section*{Preliminaries}

One of the basic ideas of the minimal model program 
(see e.g.~\cite{CKM,KMM,Ko,Wi})
is to study
projective varieties and their morphisms by understanding the
behaviour of rational curves on these varieties.
Since projective varieties with positive defect
are covered by projective lines,
it is not very surprising that the machinery of the minimal
model program can be used to study them.
This approach was implemented in a remarkable serie of papers by several
authors, starting from \cite{BSW}.
In this chapter we give a survey of these results.
In section \ref{TheNefValueSection} we give all necessary definitions
and facts from the minimal model program.
In section \ref{TheNefDefSection} we explain how to apply
this machinery to study projective varieties with positive defect.
Most results of this section are due to \cite{BSW,BFS1,Sn}.
We also give a list
of smooth projective varieties $X$ with positive defect such that
$\dim X\le10$. The study of these varieties was initiated in
\cite{E1,E2}, continued in \cite{LS} and finished in \cite{BFS1},
with many contributions of other people.
Though the natural framework for this theory is the class 
of smooth projective varieties, recently \cite{LPS}
many of these results where extended to projective
varieties with mild singularities, say with smooth normalization.

Polarized flag varieties with positive defect were first
classified in \cite{KM} using Theorem~\ref{KatzDimension}.
Later a much more quicker approach was suggested in \cite{Sn}.
The idea is to use the machinery of the minimal model program.
In Section \ref{JKHGFJKHGDKH} 
we give the classification of polarized flag varieties
with positive defect.


\section{Some Results From Mori Theory}\label{TheNefValueSection}

\subsubsection{Divisors}
Let $X$ be a smooth projective variety.
If $H$ is a very ample Cartier divisor
on $X$ (hyperplane section of some projective embedding) then
a pair $(X,H)$ will be called a polarized variety.\index{polarized variety}
A divisor $D$ on $X$ is called semi-ample if the linear system $|mD|$
is base-point-free for some $m\gg0$.\index{semi-ample divisor}
Let $K_X$ denote a canonical divisor.

We denote by $\Div(X)$ a free abelian group of Cartier divisors on $X$.
If two divisors are linearly equivalent (differ by a principal divisor)
it will be denoted by $D_1\sim D_2$.
$\Pic(X)$ is the group of algebraic linear bundles on $X$.
It is naturally isomorphic to a factorgroup of $\Div(X)$ modulo
principal divisors. 
For any $D\in\Div(X)$ the corresponding line bundle
will be denoted by $\O(D)$.

\subsubsection{$1$-cycles}
A $1$-cycle on $X$ is a finite sum 
$$Z=\sum n_iC_i$$ \index{$1$-cycle}
of irreducible curves $C_i\subset X$ with integer multiplicities $n_i$.
A $1$-cycle $Z$ is called effective if all $n_i\ge0$.\index{effective $1$-cycle}

\subsubsection{The intersection product}
If $D$ is a Cartier divisor on $X$ and $C\subset X$ is an irreducible
curve then we define an intersection product $D\cdot C\in\Z$\index{intersection product}
as $\deg g^*\O(D)|_C$, where $g:\,\bar C\to C$ is a normalisation.
This product can be extended by linearity to a $\Z$-valued pairing
of free abelian groups $\Div(X)$ and a group of $1$-cycles.
This pairing is degenerate. A $1$-cycle $C$ (resp.~a Cartier divisor $D$)
is called numerically trivial if $D\cdot C=0$ for any Cartier\index{numerically trivial}
divisor $D$ (resp.~for any $1$-cycle $C$).
Two $1$-cycles are numerically equivalent, $Z_1\equiv Z_2$, if $Z_1-Z_2$
is numerically trivial. Dually, 
two Cartier divisors are numerically equivalent, $D_1\equiv D_2$, 
if $D_1-D_2$ is numerically trivial. \index{numerically equivalent}
The group $N^1_\Z(X)$ of Cartier divisors modulo numerical equivalence
is a factorgroup of the Neron--Severi group $NS(X)$ of Cartier divisors
modulo algebraic equivalence, which is in turn\index{Neron--Severi group}
isomorphic to a subgroup of $H^2(X,\Z)$ by a standard argument
involving the exponential exact sequence of sheaves
$$0\to\Z\to\O_X\to\O^*_X\to 0.$$
Since $H^2(X,\Z)$ is a finitely generated abelian group
it follows that $N^1_\Z(X)$  
is a free abelian group of finite rank  $\rho(X)$, this rank
is called a Picard number of $X$.\index{Picard number}
Therefore, the group $N_1^\Z(X)$ of $1$-cycles
modulo numerical equivalence
is also a free abelian group of rank  $\rho(X)$.

\subsubsection{The cone of effective $1$-cycles}
We consider two real vector spaces of dimension $\rho(X)$:
$$N_1(X)=N_1^\Z(X)\otimes\R,\quad N^1(X)=N^1_\Z(X)\otimes\R.$$
These vector spaces are naturally dual to each other 
by means of the intersection product.
The main combinatorial object of the minimal model program
is defined as follows:

\index{cone of effective $1$-cycles}
\begin{Definition}$\ $
\begin{itemize}
\item The cone of effective $1$-cycles $NE_1(X)\subset N_1(X)$ is a cone generated
by classes of effective $1$-cycles. 
\item $\overline{NE}_1(X)$ is its closure in $N_1(X)$. 
\item The positive part $NE_+(X)$ of $\overline{NE}_1(X)$ is defined as
$$NE_+(X)=\{Z\in \overline{NE}_1(X)\,|\,K_X\cdot Z\ge0\}.$$
\end{itemize}
\end{Definition}

\subsubsection{Nef divisors} \index{nef divisors}
For $D\in N^1(X)$, we say that $D$ is nef (or numerically effective) 
if $D\cdot Z\ge0$ for any $Z\in NE_1(X)$ (or equivalently for any $Z\in\overline{NE}_1(X)$).
Therefore, nef vectors form a closed cone $NE^1(X)$ in $N^1(X)$, dual
to the cone $\overline{NE}_1(X)$. 
For example, if $D$ is an ample divisor then for any curve $C\subset X$
the intersection number $C\cdot D>0$. Therefore, $D$ is nef.
In fact, we have the following Kleiman criterion:

\index{Kleiman criterion}
\begin{Theorem}[\cite{Kl6}]
$D$ is ample iff
$D\cdot Z>0$ for any $Z\in\overline{NE}_1(X)$, $Z\ne0$.
\end{Theorem}

Therefore, ample divisors correspond to the interior of the cone $NE^1(X)$.

Clearly, all introduced cones
are non-degenerate (are not contained in any hyperplane).

\subsubsection{$\Q$-divisors}  \index{$\Q$-divisors}
We define $\Q$-Cartier divisors as elements of a $\Q$-vector space
$\Div(X)\otimes\Q$.
The $\Q$-Cartier divisor is called nef if the corresponding element in
$N^1(X)$ is nef.
The $\Q$-Cartier divisor is called ample if some multiple of it
is ample in a usual sense. In the sequel we shall call $\Q$-Cartier divisors
simply divisors. This will not lead to any confusion.

\subsubsection{Extremal rays}  \index{extremal ray}
An extremal ray is a half line $R$ in $\overline{NE}_1(X)$ such that
\begin{itemize}
\item $K_X\cdot R<0$ ($K_X\cdot Z<0$ for any non-zero $Z\in R$);
\item $R$ is a face of $\overline{NE}_1(X)$, i.e.~for any 
$Z_1,Z_2\in \overline{NE}_1(X)$ if $Z_1+Z_2\in R$ 
then $Z_1$ and $Z_2$ are in $R$.
\end{itemize}
An extremal rational curve $C$ is a rational curve on $X$  \index{extremal rational curve}
such that $\R_+[C]$ is an extremal ray and 
$-K_X\cdot C\le \dim X+1$.
The following is a fundamental result of Mori Theory, usually
referred to as the Mori Cone Theorem:  \index{Mori Cone Theorem}

\begin{Theorem}[\cite{Mor}]
The cone $\overline{NE}_1(X)$ is the smallest convex
cone containing $NE_+(X)$ and all the extremal rays.
For any open convex cone $U$ containing $NE_+(X)\setminus\{0\}$
there are finitely many extremal rays that do not
lie in $U\cup\{0\}$. Every extremal ray is spanned
by an extremal rational curve.
\end{Theorem}

Another fundamental result in Mori theory is the Kawamata--Shokurov
Contraction Theorem:

\index{Contraction Theorem}
\begin{Theorem}[\cite{Kaw}]
Let $D$ be a numerically effective divisor on $X$.
If, for some $a\ge1$, the divisor $aD-K_X$ is nef and big
{\rm(}i.e.~$(aD-K_X)^n>0${\rm)}, then $D$ is semi-ample.
\end{Theorem}

One of the most useful forms of the Contraction Theorem is the following

\begin{Theorem}
Let $R$ be an extremal ray in $\overline{NE}_1(X)$.
Then there exists a normal projective variety $Y$
and a surjective morphism with connected fibers $\contr_R:\,X\to Y$
with connected fibers such that for any curve $C$ on $X$,
$\contr_R$ contracts $C$ to a point if and only if $[C]\in R$.
\end{Theorem}

\proof
Indeed, let $D\in{NE}^1(X)$ be such that $D\cdot R=0$ and $D\cdot Z>0$
for any $Z\in\overline{NE}_1(X)\setminus R$.
It follows from the Mori Cone Theorem 
that we may assume $D$ is a divisor.
(In fact, usually this theorem is proved as a part of the proof of the
Mori Cone Theorem. Then the fact that $D$ can be chosen to be a divisor
is deduced from the Rationality Theorem, see below.)
The Kleiman Criterion applies that for $a\gg0$ the divisor $aD-K_X$
is ample. Therefore, $D$ is semi-ample by the Contraction Theorem.
\endproof

\subsubsection{The nef value and the nef morphism} 
Suppose that $(X,H)$ is a smooth polarized projective variety. 
Assume that $K_X$ is not nef.

\begin{Definition}\index{nef value}
$$\tau=\min\{t\in\R\,|\,K_X+tH\ \hbox{\rm is nef}\}$$
is called a nef value of $(X,H)$.
\end{Definition}

If $X$ is fixed then $\tau$ will be called a nef value of $H$,
or a nef value of $L$, where $L$ is a line bundle corresponding to $H$.
Clearly $0<t<\infty$.  Notice also that $\tau$ is a nef value of $(X,H)$
if and only if $K_X+\tau H$ is nef, but not ample.
The following theorem is known as the Kawamata Rationality Theorem:

\index{Rationality Theorem}
\begin{Theorem}[\cite{Kaw}]
The nef value $\tau$ is a rational number.
\end{Theorem}

The next theorem is a particular case of the Kawamata--Shokurov
Contraction Theorem

\begin{Theorem}
Suppose that $\tau$ is a nef value of $(X,H)$.
The linear system $|m(K_X+\tau H)|$ is base point free for $m\gg 0$.
It defines a regular morphism $\Phi:\, X\to Y$ onto a normal variety $Y$
with connected fibers. $\Phi$ is called a nef value morphism.
\end{Theorem}
\index{nef value morphism}

\subsubsection{Length of extremal rays}\index{length of extremal rays}
If $R$ is an extremal ray, then its length $l(R)$ is defined as follows
$$l(R)=\min\{-K_X\cdot C\,|\,C\ \text{is a rational curve with}\ [C]\in R\}.$$
The motivation for this definition is to give an estimate
of the dimension of locus of extremal rational curves in an extremal ray.
The idea of the proof of the following theorem is, basically, due to Mori
\cite{Mor1}. In the sequel we shall need the similar argument
at least 4 times (!), but the details will be left to the reader.

\begin{Theorem}[\cite{Wis}]\label{KJHBVJVMNBVMNBVM}
Let $C$ be a rational curve with its class in an extremal ray $R$,
such that $-K_X\cdot C=l(R)$. If $Q$ is a smooth point on $C$
then the dimension of locus of points of curves that belong to $R$
and pass through the point $Q$ is at least $l(R)-1$.
\end{Theorem}

\proof
By the local deformation theory the space of deformations
of the morphism $f:\,\P^1\to X$, $f(\P^1)=C$, has dimension at least
$$h^0(C,f^*T_X)-h^1(C,f^*T_X)=-C\cdot K_X+(1-g(C))\dim X=\dim X-C\cdot K_X.$$
Since we want to fix the point $Q$, we have $\dim X$ more restrictions.
Therefore, the space of deformations of the morphism $f:\P^1\to X$
fixing the point $Q$ has dimension at least
$$-K_X\cdot C=l(R).$$
The group of automorphisms of $\P^1$ fixing a point $x\in \P^1$
has dimension $2$, therefore, there exists a family
of rational curves on $X$ passing through $Q$ of dimension
at least $l(R)-2$. All these curves are deformations of $C$,
therefore their numerical classes belong to $R$.
To finish the proof it remains to verify that for a generic 
point $P$ in the locus of curves from our family, there
exists only finitely many curves from the family 
passing through $P$.
Suppose, on the contrary, that there exists
a family of rational curves passing through $P$ and $Q$,
parametrized by some affine curve $D$ with smooth compactification $\overline D$. 
Passing to infinite points $\overline{D}\setminus D$,
we get limit positions of these curves, that are
irreducible and reduced, since otherwise we shall obtain a contradiction
with minimality of  $-K_X\cdot C$ for rational curves
from the ray $R$.

Therefore, there exists a ruled surface $S$
and a morphism $F:\,S\to X$ such that $C_1=F^{-1}(P)$ and $C_2=F^{-1}(Q)$
are one-dimensional sections of this ruled surface over its base. 
Moreover, $\dim F(S)=2$, otherwise
all curves in our pencil represent the same $1$-cycle.
It follows that 
\begin{equation}
C_1^2<0, \ C_2^2<0,\ \text{and}\ C_1\cdot C_2=0,\label{KJHGKJGHFK}
\end{equation}
the last equality is satisfied because $C_1$ and $C_2$ do not intersect.
But this is impossible, (\ref{KJHGKJGHFK}) can not hold on a ruled
surface, see e.g.~\cite{Ha1}.
\endproof

\subsubsection{Fano varieties}\index{Fano variety}
A smooth variety $X$ is called a Fano variety
if its anti-canonical divisor $-K_X$ is ample.
Its index is defined as follows:\index{index of a Fano variety}
$$r(X)=\max\{m>0\,|\,mH\sim-K_X\ \text{for some}\ H\in\Pic(X)\}.$$
It is well-known that $r(X)\le\dim X+1$.
Indeed, suppose that $K_X=mH$, where $-H$ is ample, let $\dim X=n$.
The Poincare polynomial $\chi(vH)$ has at most $n$ zeros,
therefore $\chi(vH)\ne0$ for some $1\le v\le n+1$.
But $\chi(vH)=\pm h^n(X,vH)$ by the Kodaira vanishing theorem
and $h^n(X,vH)=h^0(X,K_X-vH)$ by the Serre duality.
Therefore, $m\le n+1$.
The following theorem was conjectured by Mukai
and proved in \cite{Wis1} using the Mori theory.

\begin{Theorem}\cite{Wis1}\label{KJHFGKHMGDFMKJ}
Let $X$ be a $n$-dimensional Fano variety of index $r(X)$.
If $r(X)>{1\over2}\dim X+1$ then $\Pic(X)=\Z$.
\end{Theorem}


\section{Mori Theory and Dual Varieties}\label{TheNefDefSection}

\subsection{The Nef Value and the Defect}

Suppose that $(X,H)$ is a smooth polarized variety. 
Then the linear system $|H|$
defines an embedding $X\subset\P^N$. Recall that $X^*\subset(\P^N)^*$
is the dual variety and 
$$\defect(X,H)=\codim X^*-1$$ 
is the defect of $(X,H)$.

\begin{Theorem}
If $\defect(X,H)>0$ then $K_X$ is not nef.
\end{Theorem}

\proof
By Theorem~\ref{FirstEin} $X$ is covered by lines. It is well-known
that if a smooth projective variety is covered by rational curves
then $K_X$ is not nef. 
In our case it also follows directly from Theorem~\ref{EINKX} (a).
\endproof

It follows that if $\defect(X,H)$ is positive, we can apply
the technique of the previous section. It turns out that 
the nef value and the defect of $(X,H)$ are connected
in a very simple way.

\begin{Theorem}[\cite{BSW,BFS1}]\label{FirstNefDef}
Assume that $\defect(X)>0$. Then 
\begin{enumerate}
\item[\rm(a)] The nef value $\tau$ of $X$ is equal to
$$\tau={\dim X+\defect(X,H)\over2}+1.$$
\item[\rm(b)] If $\Sing X\cap H$, $H\in\dual X$, is a generic contact locus
{\rm(}hence a projective subspace{\rm)}
and $l\in\Sing X\cap H$ is any line then $R=\R_+[l]$ is an extremal
ray in $\overline{NE}_1(X)$ and $\contr_R$ coincides with the nef value
morphism $\Phi$. In particular, $\Phi$ contracts 
$\Sing X\cap H$ to a point.
\item[\rm(c)] If $F$ is a generic fiber of the nef value 
morphism $\Phi$ then $\Pic F=\Z$.
\item[\rm(d)] 
$$\dim\Phi(X)\le{\dim X-\defect(X,H)\over2}.$$
\end{enumerate}
\end{Theorem}

Notice that since 
$$\dim X\equiv\defect(X,H)\mod2$$ 
by Theorem~\ref{PairityTheorem},
the nef value in this case is an integer.
Theorem \ref{FirstNefDef} can be used to show that the defect
of a polarized variety is equal to $0$.
One needs to calculate the nef value first
and then to obtain a contradiction using this Theorem.
Some examples will be discussed in next sections.

\medskip
\proof 
We denote $n=\dim X$ and $k=\defect(X)$.

\subsubsection{Step 1}
Let $\P^k=\Sing X\cap H$, $H\in\dual X$, be a generic contact locus
and $l\in\Sing X\cap H$ be any line.
Since $(-K_X)\cdot l=\bigover{n+k}2+1$ by Theorem~\ref{EINKX} (a),
it follows that 
\begin{equation}
\tau\ge\bigover{n+k}2+1.\label{KJGFKHGFDKHDGKH}
\end{equation}

\subsubsection{Step 2}
Let $\Phi$ be the nef value morphism. Since $\Phi$ is not ample,
it follows from the Kleiman Criterion and the Mori Cone Theorem 
that there exists an extremal rational curve $C$ contracted by $\Phi$,
$(K_X+\tau H)\cdot C=0$.
We claim that $C\cdot H=1$, i.e.~$C$ is a projective line. Indeed, if $C\cdot H\ge2$ then
$$\tau\le{1\over2}\tau H\cdot C={1\over2}(-K_X)\cdot C\le{1\over2}(n+1).$$
This contradicts (\ref{KJGFKHGFDKHDGKH}).
Therefore,
$$\tau=(-K_X)\cdot C=l(R),$$
where $R=\R_+[C]$ is an extremal ray.
Applying Theorem~\ref{KJHBVJVMNBVMNBVM}, wee see that if $F$
is a positive-dimensional fiber of $\Phi$ then
\begin{equation}
\dim F\ge\tau-1.\label{KJHGLJGHFCKHHGFCKHJG}
\end{equation}

\subsubsection{Step 3}
Suppose that $\tau>\bigover{n+k}2+1$.
Let $F$ be a positive dimensional fiber of $\Phi$, $x\in F$.
By Theorem~\ref{Strong-Ein-EveryPoint} there exists a line $l\subset X$
passing through $x$ and such that $(-K_X)\cdot l=\bigover{n+k}2+1$.
Therefore, $\Phi$ does not contract $l$, hence $l\not\subset F$.
Let $y\in l$, $y\not\in F$.
Arguing as in the proof of Theorem~\ref{KJHBVJVMNBVMNBVM}, we see
that there exists a $\bigover{n+k}2$-dimensional subvariety $Y\subset X$
covered by rational curves passing through $y$ that are deformations of $l$.
Applying the similar argument with ruled surfaces once again, it is easy to see
that $\dim Y\cap F=0$. However, using (\ref{KJHGLJGHFCKHHGFCKHJG})
we get that
$$\dim Y\cap F\ge \dim Y+\dim F-\dim X\ge 
\bigover{n+k}2+\bigover{n+k}2-n>0.$$
Contradiction.
This proves (a).

\subsubsection{Step 4}
Clearly, $\Phi$ contracts any generic contact locus $\Sing X\cap H$, $H\in\dual X$.
Moreover, arguing as above (or using Theorem~\ref{SecondEin} (d)),
we see that any fiber of $\Phi$ is at least $\bigover{n+k}2$-dimensional,
therefore $\dim Y\le n-\bigover{n+k}2=\bigover{n-k}2$.
This proves (d). Let $F$ be a generic fiber.
Then $K_F=K_X|_F=-\tau H|_F$. Therefore, $F$ is a Fano variety
of index at least 
$$\tau={n+k\over2}+1>{\dim F\over2}+1.$$
Therefore, $\Pic F=\Z$ by Theorem~\ref{KJHFGKHMGDFMKJ}.
This proves (c).

\subsubsection{Step 5}
It remains to prove that $\Phi$ is a contraction of 
an extremal ray $\R_+[l]$.
Suppose, on the contrary, that there exists an extremal ray $R$
such that $l\not\in R$ and $\Phi$ contracts any curve from $R$.
Let $C$ be an extremal rational curve from~$R$.
Arguing as in Step 2, we see that $l(R)=(-K_X)\cdot C=\tau$.
Consider the morphism $\Psi=\contr_R$.
By Theorem~\ref{KJHBVJVMNBVMNBVM}, there exists a fiber $F$ of $\Psi$
such that $\dim F\ge \tau-1$, let $x\in F$.
By Theorem~\ref{Strong-Ein-EveryPoint} there exists a line $l\subset X$
passing through~$x$ and such that $(-K_X)\cdot l=\tau$.
This line is either a line in a generic contact locus or its
limit position, therefore $\Psi$ does not contract $l$.
Hence $l\not\subset F$. Let $y\in l$, $y\not\in F$.
Arguing as in the proof of Theorem~\ref{KJHBVJVMNBVMNBVM}, we see
that there exists a $(\tau-1)$-dimensional subvariety $Y\subset X$
covered by rational curves passing through $y$ that are deformations of $l$.
Applying the similar argument once again, it is easy to see
that $\dim Y\cap F=0$. But
$$\dim Y\cap F\ge \dim Y+\dim F-\dim X\ge 
2\tau-2-n=k>0.$$
Contradiction.
\endproof

\subsection{The Defect of Fibers of the Nef Value Morphism}

Suppose now that $(X,H)$ is a smooth polarized variety and 
$K_X$ is not nef. Let $\Phi:\,X\to Y$ denote
a nef value morphism. 
Then its generic fiber $F$ is a smooth irreducible
variety. It has a canonical polarization $(F,H_F)$, 
where $H_F$ is a hyperplane section of $F$ in the embedding
$F\subset X\subset \P^N$.

\begin{Theorem}[\cite{BFS1}]\label{SecondNefDef}
The following conditions are equivalent:
\begin{enumerate}
\item[\rm(i)] $\defect(X,H)>0$;
\item[\rm(ii)] $\defect(F,H_F)>\dim Y$;
\item[\rm(iii)] $\defect(X,H)=\defect(F,H_F)-\dim Y>0$.
\end{enumerate}
\end{Theorem}

\proof It follows from Theorem~\ref{LanteriStruppMonotonicity}
that
$$\dim X+\defect(X,H)\ge \dim F+\defect(F,H_F).$$
Therefore,
(i) follows from (ii). Obviously
(ii) follows from (iii). Therefore,
it is sufficient to prove that (iii) follows from (i).

Since $\defect(X,H)>0$ we have 
$$\tau(X,H)={\dim X+\defect(X,H)\over 2}+1$$
by Theorem \ref{FirstNefDef}.
Since $F$ is a generic fiber of $\Phi$ we have
$$\O(K_F+\tau H_F)=\O(K_X+\tau H)|_F$$
is a trivial bundle.
Therefore, $K_F$ is not nef and $\tau(F,H_F)=\tau(X,H)$.
If $\defect(F,H_F)>0$ then again by Theorem \ref{FirstNefDef}
$$\defect(F,H_F)={\dim F+\defect(F,H_F)\over2}+1$$
and, therefore, 
$$0<\defect(X,H)=-\dim X+\dim F+\defect(F,H_F)=\defect(F,H_F)-\dim Y.$$

It remains to show
that if $\defect(X,H)>0$ then $\defect(F,H_F)>0$ as well.
Let $L=\O(H)$, $L_F=\O(H_F)$. Recall that for any line bundle $L$
on a smooth variety $X$
we denote by $J_1(X,L)$ the vector bundle of $1$-jets of $L$;
for any vector bundle $V$ we denote by $c_r(V)$ its $r$ Chern class.
By Theorem~\ref{MJGHVJGVJGHCVVC} we need to show that
$c_{\dim F}(J_1(F, L_F))=0$.
From the exact sequence
$$0\to T_F^*(X)\otimes L_F\to J_1(X,L)|_F\to J_1(F,L_F)\to0$$
we see that the total Chern class
$$c(J_1(X,L)|_F)=c(T_F^*(X)\otimes L_F)\cdot c(J_1(F,L_F)).$$
Since $F$ is a generic fiber of $\Phi$,
$T_F^*(X)\otimes L_F=L_F+\ldots+L_F$ ($\dim Y$ copies).
Therefore, $c(T_F^*(X)\otimes L_F)=(1+H_F)^{\dim Y}$.
In particular,
$$\begin{array}{l}
\displaystyle c_{\dim F}(J_1(X,L)|_F)=\cr
\displaystyle{}\qquad c_{\dim F}(J_1(F,L_F))+\sum_{i=1}^{\dim Y}{\dim Y\choose i}c_{\dim F-i}(J_1(F,L_F))\cdot H_F^i.
\end{array}$$
All summands on the right are nonnegative since $H_F$ is very ample
and $J_1(F,L_F)$ is spanned.
Therefore in order to prove that
$c_{\dim F}(J_1(F, L_F))=0$ it remains to prove that
$c_{\dim F}(J_1(X, L)|_F)=0$.

It will be easier to prove a more general statement
that 
$c_k(J_1(X, L)|_F)=0$ for $k\ge \dim F-\defect(X,H)+1$.
By the classical representation of a Chern class of a spanned
vector bundle \cite{Fu1} we see that $c_k(J_1(X, L))=0$ is represented
by a set
$C_k\subset X$ of all points where $\dim X-k+2$ generic sections
of $H^0(X, J_1(X,L))$ are not linearly independent.
We may suppose that these sections have a form
$j_1(s_1),\ldots,j_1(s_{\dim X-k+2})$, where
$s_1,\ldots,s_{\dim X-k+2}$ are generic sections of 
$H^0(X, L)$ and $j_1$ is a natural map sending a germ of a section
to its $1$-jet.
It suffices to prove that 
for $k\ge\dim F-\defect(X,H)+1$
we have $C_k\cap F=\emptyset$. 
Since $C_k\supset C_{k+1}$, we only need to prove that $C\cap F=\emptyset$,
where $C=C_{\dim F-\defect(X,H)+1}$.
Suppose that this intersection is non-empty.
Therefore, the restriction map $\Phi_C:\, C\to Y$ is surjective.
Since $\dim C=\dim Y+\defect(X,H)-1$, it follows that
a generic fiber of $\Phi_C$ has dimension $\defect(X,H)-1$.

Suppose that $x\in C\cap F$.
There exist $\lambda_i$ for $1\le i\le \dim Y+\defect(X,H)-1$,
not all zero, such that
$\sum_i\lambda_ij_1(s_i(x))=0$. Therefore 
$s=\sum_i\lambda_is_i$ is a non-trivial global section of $L$
vanishing at $x$ with its $1$-jet.
Hence the divisor $H_s$ is singular at $x$.
By Theorem~\ref{Strong-Ein-EveryPoint}
there exists a linear subspace $\P^{\defect(X,H)}\subset X$ containing $x$,
contained in $\Sing H_s$, and such that
$\O(K_X+\tau(X,H)H)|_{\P^{\defect(X,H)}}$ is a trivial bundle.
Since this is a trivial bundle, $\P^{\defect(X,H)}$ is contracted
by $\Phi$ to a point and, therefore, $\P^{\defect(X,H)}\subset F$.
On the other hand, since $\P^{\defect(X,H)}$ is contained in 
$\Sing H_s$ it is also contained in $C$.
Therefore, a generic fiber of $\Phi_C$ has dimension
greater or equal to $\defect(X,H)$. Contradiction.
\endproof


\subsection{Varieties With Small Dual Varieties}
In this section we describe without proofs smooth polarized varieties
$(X,H)$ with positive defect for $\dim X\le10$.
First we need a couple of definitions.
\index{varieties with small dual varieties}

\begin{Definition} \index{$\P^d$-bundle}
A smooth polarized variety $(X,H)$ is called a $\P^d$-bundle over a normal
variety $Y$ if there exists a surjective morphism $p:\,X\to Y$
such that all the fibers $F$ of $p$ are projective subspaces $\P^d$
of $\P^N$, where $X\subset\P^N$ is an embedding corresponding to $H$.
\end{Definition}

\index{scroll}\index{quadric fibration}\index{Del Pezzo fibration}
\index{Mukai fibration}
\begin{Definition}
A smooth polarized variety $(X,H)$ is called a scroll
(resp.~a quadric, Del Pezzo, or Mukai fibration) over a normal variety $Y$
if there exists a surjective morphism $p:\,X\to Y$ with connected
fibers such that $K_X+(\dim X-\dim Y+1)H\sim p^*H'$
(resp.~$K_X+(\dim X-\dim Y)H\sim p^*H'$, 
$K_X+(\dim X-\dim Y-1)H\sim p^*H'$, or 
$K_X+(\dim X-\dim Y-2)H\sim p^*H'$) for some ample Cartier divisor $H'$ on $Y$.
\end{Definition}

\index{Kobayashi--Ochiai characterisation}
To adjust this definition we need to note that it easily
follows from the Kobayashi--Ochiai characterisation of
projective spaces and quadrics \cite{KO}
that if $(X,H)$ is a scroll over $Y$ then a generic fiber $F$ of $p$
is a projective subspace $\P^d$ of $\P^N$, 
where $X\subset\P^N$ is an embedding corresponding to $H$.
If $(X,H)$ is a quadric fibration then a generic fiber $F$
is isomorphic to a quadric hypersurface $Q$ in some projective
space and the line bundle corresponding to this
embedding is $\O(H)|_F$.
In fact, quadric fibrations will not appear in the classification because
if $(X,H)$ is a quadric fibration then $\defect(X,H)=0$ by~\cite{BFS1}.

\subsubsection{$\dim X\le2$}
If $\dim X\le2$ then $\defect(X,H)=0$.

\subsubsection{$\dim X=3$}
Suppose that $\dim X=3$. Then $\defect(X,H)>0$ if an only if 
$(X,H)$ is a $\P^2$-bundle over a smooth curve. In this case $\defect(X,H)=1$.

\subsubsection{$\dim X=4$}
Suppose that $\dim X=4$. Then $\defect(X,H)>0$ if an only if 
$(X,H)$ is a $\P^3$-bundle over a smooth curve. In this case $\defect(X,H)=2$.

\subsubsection{$\dim X=5$}
Let $\dim X=5$ and $\defect(X,H)>0$. 
Then we have the following possibilities:
\begin{enumerate}
\item[\rm(a)] $X$ is a hyperplane section of the Pl\"ucker embedding
of the Grassmanian $\Gr(2,5)$ in $\P^9$, $\defect(X,H)=1$.
\item[\rm(b)] $(X,H)$ is a $\P^4$-bundle over a smooth curve, $\defect(X,H)=3$.
\item[\rm(c)] $(X,H)$ is a $\P^3$-bundle over a smooth surface, $\defect(X,H)=1$.
\end{enumerate}

\subsubsection{$\dim X=6$}
Let $\dim X=6$ and $\defect(X,H)>0$. 
Then either:
\begin{enumerate}
\item[\rm(a)] $(X,H)$ is the Pl\"ucker embedding of the Grassmanian
$\Gr(2,5)$ in $\P^9$, $\defect(X,H)=2$.
\item[\rm(b)] $(X,H)$ is a $\P^5$-bundle over a smooth curve, $\defect(X,H)=4$.
\item[\rm(c)] $(X,H)$ is a $\P^4$-bundle over a smooth surface, $\defect(X,H)=2$.
\end{enumerate}

\subsubsection{$\dim X=7$}
Let $\dim X=7$ and $\defect(X,H)>0$. 
Then we have the following possibilities:
\begin{enumerate}
\item[\rm(a)] $(X,H)$ is a $\P^6$-bundle over a smooth curve, $\defect(X,H)=5$.
\item[\rm(b)] $(X,H)$ is a $\P^5$-bundle over a smooth surface, $\defect(X,H)=3$.
\item[\rm(c)] $X$ is a $7$-dimensional Fano variety with 
$\Pic X=\Z$ and $K_X\sim -5H$. In this case $\defect(X,H)=1$.
An example of such variety is the section of $10$-dimensional
spinor variety $S$ in $\P^{15}$ by three generic hyperplanes \cite{Mu,LS}.
\item[\rm(d)] $(X,H)$ is a Del Pezzo fibration over a smooth curve such that 
the general fiber $F$ is isomorphic to $\Gr(2,5)$ and 
the embedding corresponding to $\O(H)|_F$ is the Pl\"ucker
embedding in $\P^9$. In this case $\defect(X,H)=1$.
\item[\rm(e)] $(X,H)$ is a scroll over a $3$-dimensional variety, $\defect(X,H)=1$.
\end{enumerate}

\subsubsection{$\dim X=8$}
Let $\dim X=8$ and $\defect(X,H)>0$. 
Then either:
\begin{enumerate}
\item[\rm(a)] $(X,H)$ is a $\P^7$-bundle over a smooth curve, $\defect(X,H)=6$.
\item[\rm(b)] $(X,H)$ is a $\P^6$-bundle over a smooth surface, $\defect(X,H)=4$.
\item[\rm(c)] $X$ is 
is a $8$-dimensional Fano variety with 
$\Pic X=\Z$ and $K_X\sim -6H$. In this case $\defect(X,H)=2$.
\item[\rm(d)] $(X,H)$ is a scroll over a $3$-dimensional variety, $\defect(X,H)=2$.
\end{enumerate}

\subsubsection{$\dim X=9$}
Let $\dim X=9$ and $\defect(X,H)>0$. 
Then either:
\begin{enumerate}
\item[\rm(a)] $(X,H)$ is a $\P^8$-bundle over a smooth curve, $\defect(X,H)=7$.
\item[\rm(b)] $(X,H)$ is a $\P^7$-bundle over a smooth surface, $\defect(X,H)=5$.
\item[\rm(c)] 
$X$ is a $9$-dimensional Fano variety with 
$\Pic X=\Z$ and $K_X\sim -7H$. In this case $\defect(X,H)=3$.
\item[\rm(d)] $(X,H)$ is a scroll over a $3$-dimensional variety, $\defect(X,H)=3$.
\item[\rm(e)] 
$X$ is a Fano $9$-dimensional variety with $K_X\sim -6H$ and
and $\Pic X=\Z$, $\defect(X,H)=1$.
\item[\rm(f)] $(X,H)$ is a Mukai fibration over a smooth curve with 
$\Pic F=\Z$ for a generic fiber $F$, $\defect(X,H)=1$.
\item[\rm(g)] $(X,H)$ is a scroll over a $4$-dimensional variety, 
$\defect(X,H)=1$.
\end{enumerate}

\subsubsection{$\dim X=10$}
Let $\dim X=10$ and $\defect(X,H)>0$. 
Then we have the following possibilities:
\begin{enumerate}
\item[\rm(a)] $(X,H)$ is a $\P^9$-bundle over a smooth curve, $\defect(X,H)=8$.
\item[\rm(b)] $(X,H)$ is a $\P^8$-bundle over a smooth surface, $\defect(X,H)=6$.
\item[\rm(c)] 
$X$ is a $10$-dimensional Fano variety with 
$\Pic X=\Z$ and $K_X\sim -8H$. In this case $\defect(X,H)=4$.
\item[\rm(d)] $(X,H)$ is a scroll over a $3$-dimensional variety, $\defect(X,H)=4$.
\item[\rm(e)] 
$X$ is a Fano $10$-dimensional variety with $K_X\sim -7H$
and $\Pic X=\Z$, $\defect(X,H)=2$;
\item[\rm(f)] $(X,H)$ is a Mukai fibration over a smooth curve 
with $\Pic F=\Z$ for a generic fiber $F$, $\defect(X,H)=2$.
\item[\rm(g)] $(X,H)$ is a scroll over a $4$-dimensional variety, 
$\defect(X,H)=2$.
\end{enumerate}

\section{Polarized Flag Varieties With Positive Defect}\label{JKHGFJKHGDKH}

\subsection{Nef Value of Polarized Flag Varieties}

We continue to use the notation from Section~\ref{DefinitionsNotationsFlags}.
Let $G$ be a connected simply-connected semisimple 
complex algebraic group with a Borel subgroup $B$ and a maximal torus $T\subset B$.
Let $\PP$ be a character group of $T$ (the weight lattice).
Let $\Delta\subset \PP$ be a set of roots of $G$ relative to $T$.
To every root $\alpha\in\Delta$ we can assign the $1$-dimensional
unipotent subgroup $U_\alpha\subset G$. We define the negative roots $\Delta^-$
as those roots $\alpha$ such that $U_\alpha\subset B$.
Let $\Pi\subset\Delta^+$ be simple roots, $\Pi=\{\alpha_1,\ldots,\alpha_n\}$,
where $n=\rank G=\dim T$. 
The root lattice $\QQ\subset \PP$ is a sublattice spanned by $\Delta$. \index{root lattice}%
For any simple root $\alpha_i\in\Pi$ we can associate a smooth
curve $C_{\alpha_i}\subset G/B$, namely 
$C_{\alpha_i}=\overline{\pi(U_{\alpha_i})}$,
where $\pi:\,G\to G/B$ is a factorization map.
By linearity to any element $\alpha=\sum n_i\alpha_i\in \QQ$
we assign a $1$-cycle $C_\alpha=\sum n_iC_{\alpha_i}$ in $G/B$.
Clearly this cycle is effective if $\alpha\in \QQ^+$,
the semigroup generated by simple roots.

$\PP$ is generated as a $\Z$-module by fundamental weights 
$\omega_1,\ldots,\omega_n$ dual to the simple roots
under the Killing form. \index{fundamental weights}%
We denote the dominant weights by $\PP^+$ and
the strictly dominant weights by $\PP^{++}$.

The character group of $B$ is identified with a character group of $T$.
Therefore, for any $\lambda\in \PP$ we can assign a $1$-dimensional $B$-module $\C_\lambda$,
where $B$ acts on $\C_\lambda$ by a character $\lambda$.
Now we can define the twisted product $G\times_B\C_\lambda$ to be the 
quotient of $G\times E$ by the diagonal action $B$: 
$$b\cdot(g,z)=(gb^{-1},b\cdot z).$$
Projection onto the first factor induces a map $G\times_B\C_\lambda\to G/B$,
which realizes the twisted product as an equivariant line bundle $\L_\lambda$ 
on $G/B$ with fiber $\C_\lambda$.

The following theorem is well-known:

\begin{Theorem}$\ $
\begin{enumerate}
\item[\rm (1)] The correspondence $\alpha\to C_\alpha$ is an isomorphism
of $\QQ$ and $N_1^\Z(G/B)$, the group of $1$-cycles
modulo numerical equivalence.
\item[\rm (2)] The correspondence $\lambda\to L_\lambda$ is an isomorphism
of $\PP$ and $\Pic(G/B)$. The latter group is, in turn,
isomorphic to $N^1_Z(G/B)$, the group of Cartier divisors
modulo numerical equivalence.
\item[\rm (3)] The pairing $\PP\times \QQ\to\Z$ given by duality
corresponds to the intersection product on $N^1_Z(G/B)\times N_1^Z(G/B)$.
\item[\rm (4)] $L_\lambda$ is nef if and only if $\lambda\in \PP^+$.
$L_\lambda$ is ample if and only if 
$L_\lambda$ is very ample if and only if
$\lambda\in \PP^{++}$.
\end{enumerate}
\end{Theorem}

From this theorem it is easy to derive the corresponding
description for flag varieties $G/P$, where $P\subset G$
is an arbitrary parabolic subgroup.
Let $\Pi_P\subset\Pi$ be some subset of simple roots.
Let $\Delta^+_P\subset\Delta^+$ denote the positive roots
that are linear combinations of the roots in $\Pi_P$.
Then we may suppose that $P$ is generated by $B$ and 
by the root groups $U_\alpha$ for $\alpha\in\Delta^+_P$.
We denote $\Pi\setminus\Pi_P$ by $\Pi_{G/P}$ and 
$\Delta^+\setminus\Delta^+_P$ by $\Delta^+_{G/P}$.
A parabolic subgroup is maximal if and only if $\Pi_{G/P}$
is a single simple root.
We denote by $\QQ_{G/P}$ the sublattice of $\QQ$ generated by $\Pi_{G/P}$
and by $\QQ^+_{G/P}$ the corresponding subsemigroup.
For any $\alpha_i\in\Pi_{G/P}$
we can associate a smooth
curve $C_{\alpha_i}\subset G/P$, namely 
$C_{\alpha_i}=\overline{\pi(U_{\alpha_i})}$,
where $\pi:\,G\to G/P$ is a factorization map.
By linearity to any element $\alpha=\sum n_i\alpha_i\in \QQ_{G/P}$
we assign a $1$-cycle $C_\alpha=\sum n_iC_{\alpha_i}$ in $G/P$.
Clearly this cycle is effective if $\alpha\in \QQ^+_{G/P}$.

The fundamental weights $\omega_{i_1},\ldots\omega_{i_k}$
dual to the simple roots in $\Pi_{G/P}$ generate the sublattice
$\PP_{G/P}$ of $\PP$. We denote $\PP^+\cap \PP_{G/P}$ by $\PP^+_{G/P}$.
The subset $\PP^{++}_{G/P}\subset\PP^+_{G/P}$ consists of all weights
$\lambda=\sum n_k\omega_{i_k}$ such that all $n_k>0$.
Any weight $\lambda\in \PP_{G/P}$ defines a character of $P$,
and, therefore, a line bundle $L_\lambda$ on $G/P$.

\begin{Theorem}\label{GP-NEdescription}$\ $
\begin{enumerate}
\item[\rm (1)] The correspondence $\alpha\to C_\alpha$ is an isomorphism
of $\QQ_{G/P}$ and $N_1^\Z(G/P)$.
\item[\rm (2)] The correspondence $\lambda\to L_\alpha$ is an isomorphism
of $\PP_{G/P}$ and $\Pic(G/P)$. The latter group is, in turn,
isomorphic to $N^1_Z(G/P)$.
In particular, $\Pic(G/P)=\Z$ if and only if $P$ is maximal.
\item[\rm (3)] The pairing $\PP_{G/P}\times \QQ_{G/P}\to\Z$ given by duality
corresponds to the intersection product on $N^1_Z(G/P)\times N_1^Z(G/P)$.
\item[\rm (4)] $L_\lambda$ is nef if and only if $\lambda\in \PP^+_{G/P}$.
$L_\lambda$ is ample if and only if 
$L_\lambda$ is very ample if and only if 
$\lambda\in \PP^{++}_{G/P}$.
\item[\rm (5)] If $\lambda\in \PP^+_{G/P}$ then the linear system corresponding
to $L_\lambda$ is base--point free. The corresponding map given by sections
is a factorization $G/P\to G/Q$, where $Q$ is a parabolic subgroup
such that $\Pi_Q$ is a union of $\Pi_P$ and all simple roots in $\Pi_{G/P}$
orthogonal to $\lambda$.
\end{enumerate}
\end{Theorem}

Now it is quite straightforward to find nef values of ample line bundles.

\begin{Theorem}\label{NefValuesMorphismsFlags}$\ $
\begin{enumerate}
\item[\rm (1)] The anticanonical bundle $K^*_{G/P}$ 
corresponds to the dominant weight
$\rho_{G/P}=\sum_{\alpha\in\Delta^+_{G/P}}\alpha$,
$\rho_{G/P}\in\PP_{G/P}^{++}$. Canonical line bundle
$K_{G/P}$ is not nef.
\item[\rm (2)] The nef value of an ample line bundle $L_\lambda$
given by the weight $\lambda\in\PP^{++}_{G/P}$
is equal to
$$\tau=\max_{\alpha\in\Pi_{G/P}}{(\rho_{G/P},\alpha)\over(\lambda,\alpha)}.$$
\item[\rm (3)] The nef value morphism $\Phi$ is a homogeneous
fiber bundle $\Phi:\,G/P\to G/Q$, where $Q$ is the parabolic
subgroup defined by $\Pi_Q=\Pi_P\cup\Pi'$,
where $\Pi'\subset \Pi_{G/P}$ is a set of simple roots
for which the above maximum occurs.
\end{enumerate}
\end{Theorem}

\proof
The tangent bundle $T(G/P)$ is an equivariant vector bundle such that
its fiber over the identity coset is isomorphic to the quotient
of Lie algebras $\Lie G/\Lie P$ on which $P$ acts
via an adjoint representation on $\Lie G$ projected to the quotient.
The $T$-weights of this representation consist of $\Delta^+_{G/P}$.
Therefore the anticanonical bundle $K^*_{G/P}=\Lambda^{\dim G/P}T(G/P)$
corresponds to the weight $\rho_{G/P}$.
Notice that for any $\alpha\in\Pi_{G/P}$
$$(\rho_{G/P},\alpha)=\left(\sum_{\beta\in\Delta^+}\beta,\alpha\right)-
\left(\sum_{\beta\in\Delta^+_P}\beta,\alpha\right).$$
The first product is equal to $2$, the second one is non-positive,
because any root from $\Delta_P^+$ is a linear combination
of simple roots from $\Pi_P$ with non-negative coefficients
and the scalar product of two distinct simple roots is non-positive.
Therefore, $\rho_{G/P}\in\PP^{++}_{G/P}$
and $K_{G/P}$ is not nef.

Let 
$$\rho_{G/P}=\sum_{\alpha_i\in\Pi_{G/P}}n_i\omega_i,\quad
\lambda_{G/P}=\sum_{\alpha_i\in\Pi_{G/P}}m_i\omega_i.$$
Then, clearly, the nef value $\tau$ of $L$
is given by
$$\tau=\max_{\alpha_i\in\Pi_{G/P}}{n_i\over m_i}=
\max_{\alpha\in\Pi_{G/P}}{(\rho_{G/P},\alpha)\over(\lambda,\alpha)}.$$

Part (3) follows from Theorem \ref{GP-NEdescription}, part (5).
\endproof

Using this theorem, it is straightforward to calculate nef values
in particular cases. Suppose that $G$ is a simple group,
$P$ is a maximal parabolic subgroup such that $\Pi_{G/P}=\{\alpha_i\}$ is
a single simple root. Suppose further that $L$ is a generator
of $\Pic G/P$, that is, $L$ corresponds
to the fundamental weight $\omega_i$.
In the following tables we give for each $G$ and $i$ 
(in the Bourbaki numbering of
fundamental weights) the dimension $\dim G/P$ and 
the nef value $\tau$ of $L$.

\bigskip

\begin{tabular}{cl}
$G=A_l$&
\begin{tabular}{|c|c|}
\hline
$i$ & $1,\ldots,l$\\
\hline
$\dim G/P$ & $i(l+1-i)$\\
\hline
$\tau$ & $l+1$\\
\hline
\end{tabular}
\\
\noalign{\smallskip}
$G=B_l$&
\begin{tabular}{|c|c|c|}
\hline
$i$ & $1,\ldots,l-1$ & $l$\\
\hline
$\dim G/P$ & $i(4l+1-3i)/2$ & $l(l+1)/2$\\
\hline
$\tau$ & $2l-i$ & $2l$\\
\hline
\end{tabular}
\\
\noalign{\smallskip}
$G=C_l$&
\begin{tabular}{|c|c|}
\hline
$i$ & $1,\ldots,l$ \\
\hline
$\dim G/P$ & $i(4l+1-3i)/2$ \\
\hline
$\tau$ & $2l-i+1$ \\
\hline
\end{tabular}
\\
\noalign{\smallskip}
$G=D_l$&
\begin{tabular}{|c|c|c|}
\hline
$i$ & $1,\ldots,l-2$ & $l-1,l$ \\
\hline
$\dim G/P$ & $i(4l-1-3i)/2$ & $l(l-1)/2$ \\
\hline
$\tau$ & $2l-i-1$ & $2l-2$\\
\hline
\end{tabular}
\\
\noalign{\smallskip}
$G=E_6$&
\begin{tabular}{|c|c|c|c|c|}
\hline
$i$ & $1,5$ & $2,4$ & $3$ & $6$ \\
\hline
$\dim G/P$ & $16$ & $25$ & $29$ & $21$\\
\hline
$\tau$ & $12$ & $9$ & $7$ & $11$\\
\hline
\end{tabular}
\\
\noalign{\smallskip}
$G=E_7$&
\begin{tabular}{|c|c|c|c|c|c|c|c|}
\hline
$i$ & $1$ & $2$ & $3$ & $4$ & $5$ & $6$ & $7$ \\
\hline
$\dim G/P$ & $27$ & $42$ & $50$ & $53$ & $47$ & $33$ & $42$\\
\hline
$\tau$ & $18$ & $13$ & $10$ & $8$ & $11$ & $17$ & $14$\\
\hline
\end{tabular}
\\
\noalign{\smallskip}
$G=E_8$&
\begin{tabular}{|c|c|c|c|c|c|c|c|c|}
\hline
$i$ & $1$ & $2$ & $3$ & $4$ & $5$ & $6$ & $7$  & $8$\\
\hline
$\dim G/P$ & $57$ & $83$ & $97$ & $104$ & $106$ & $98$ & $78$ & $92$\\
\hline
$\tau$ & $29$ & $19$ & $14$ & $11$ & $9$ & $13$ & $23$ & $17$\\
\hline
\end{tabular}
\\
\noalign{\smallskip}
$G=F_4$&
\begin{tabular}{|c|c|c|c|c|}
\hline
$i$ & $1$ & $2$ & $3$ & $4$ \\
\hline
$\dim G/P$ & $15$ & $20$ & $20$ & $15$ \\
\hline
$\tau$ & $8$ & $5$ & $7$ & $11$ \\
\hline
\end{tabular}
\\
\noalign{\smallskip}
$G=G_2$&
\begin{tabular}{|c|c|c|}
\hline
$i$ & $1$ & $2$\\
\hline
$\dim G/P$ & $5$ & $5$\\
\hline
$\tau$ & $5$ & $3$\\
\hline
\end{tabular}
\\
\end{tabular}

\subsection{Classification}

In this section we use the calculations from the previous section
and Theorems \ref{FirstNefDef}, \ref{SecondNefDef}
to classify all polarized flag varieties with positive defect.
We continue to use notations from the previous section.

Suppose first that $G$ is a simple group, $P\subset G$
is a maximal parabolic subgroup corresponding to a simple root $\alpha_i$. 
Then all very ample linear bundles on $G/P$
have a form $L^{\otimes k}$ for $k>0$, where
$L$ corresponds to the fundamental weight $\omega_i$.
If $k>1$ then by Theorem~\ref{DefectInVeronese}
we have $\defect(G/P, L^{\otimes k})=0$.
Therefore it suffices to find $\defect(G/P, L)$.
By Theorem \ref{FirstNefDef} if $\defect(G/P, L)>0$
the $\defect(G/P, L)=2\tau(G/P,L)-2-\dim G/P$.
Therefore we may look through the table of the previous
section and pick all entries such that
$d_P=2\tau(G/P,L)-2-\dim G/P>0$.
The following table contains all these entries.
Therefore all polarized flag varieties with positive
defect (such that $P$ is a maximal parabolic subgroup)
are contained in this table.

\medskip

\noindent
\begin{tabular}{|c|c|c|c|c|c|}
\hline
No. & $G$ & $i$ & $\tau(G/P, L)$ & $\dim G/P$ & $d_P$\\
\hline
\hline
1  & $A_l$                  & $1,l$   & $l+1$ & $l$    & $l$\\
\hline
2  & $A_l$, $l$ even & $2,l-1$ & $l+1$ & $2l-2$ & $2$\\
\hline
3  & $C_l$                  & $1$     & $2l$  & $2l-1$ & $2l-1$\\
4  & $B_2$                  & $2$     & $4$   & $3$    & $3$\\
\hline
5  & $B_4$                  & $4$     & $8$   & $10$   & $4$\\
6  & $D_5$                  & $4,5$   & $8$   & $10$   & $4$\\
\hline
\hline
7  & $A_l$, $l$ odd  & $2,l-1$ & $l+1$ & $2l-2$ & $2$\\
\hline
8  & $A_5$                  & $3$     & $6$   & $9$    & $1$\\
\hline
9  & $B_l$                  & $1$     & $2l-1$& $2l-1$ & $2l-3$\\
\hline
10  & $B_3$                  & $3$     & $6$   & $6$    & $4$\\
11  & $D_4$                  & $3,4$   & $6$   & $6$    & $4$\\
\hline
12  & $B_5$                  & $5$     & $10$  & $15$   & $3$\\  
13  & $D_6$                  & $5,6$   & $10$  & $15$   & $3$\\
\hline
14  & $B_6$                  & $6$     & $12$  & $21$   & $1$\\
15  & $D_7$                  & $6,7$   & $12$  & $21$   & $1$\\
\hline
16  & $C_l$                  & $2$     & $2l-1$& $4l-5$ & $1$\\
\hline
17  & $D_l$                  & $1$     & $2l-2$& $2l-2$ & $2l-4$\\
\hline
18  & $E_6$                  & $1,5$   & $12$  & $16$   & $6$\\
\hline
19  & $E_7$                  & $1$     & $18$  & $27$   & $7$\\
\hline
20  & $F_4$                  & $4$     & $11$  & $15$   & $5$\\
\hline
21  & $G_2$                  & $1$     & $5$   & $5$    & $3$\\
\hline
\end{tabular}

\medskip

In fact, only first $6$ entries of this table
correspond to polarized flag varieties with positive defect,
and entry $4$ is a particular case of entry $3$.
This can be done case-by-case (see~\cite{Sn}).
Therefore, we have the following theorem

\begin{Theorem}\label{FirstClassificationFlag}
Suppose that $P_i$ is a maximal parabolic subgroup
of a simple algebraic group corresponding to the simple root $\alpha_i$.
Let $L$ be a very ample line bundle on $G/P_i$
corresponding to the fundamental weight $\omega_i$.
Then
\begin{enumerate}
\item[{\rm (1)}] $\defect(G/P_i, L^{\otimes k})=0$ for $k>1$.
\item[{\rm (2)}] $\defect(G/P_i, L)>0$ if and only if $G/P_i$ is one 
of the following:
\begin{itemize}
\item $A_l/P_1$, $A_l/P_l$. $\defect(G/P)=\dim(G/P)=l$.
\item $A_l/P_2$, $A_l/P_{l-1}$, $l\ge4$ even. 
$\defect(G/P)=2$, $\dim(G/P)=2l-2$. 
\item $C_l/P_1$. $\defect(G/P)=\dim(G/P)=2l-1$.
\item $B_4/P_4$, $D_5/P_4$, $D_5/B_5$.
$\defect(G/P)=4$, $\dim(G/P)=10$. 
\end{itemize}
\end{enumerate}
\end{Theorem}

Notice that we do not follow the usual tradition and 
consider isomorphic polarized varieties with different
group actions as distinct varieties.
The next step is to consider flag varieties corresponding
to arbitrary parabolic subgroups in simple group.

\begin{Theorem}\label{SecondClassificationFlag}
Suppose that $(G/P,L)$ is a polarized flag variety
of a simple algebraic group $G$, where $P\subset G$
is a non-maximal parabolic subgroup.
Then $\defect(G/P,L)=0$.
\end{Theorem}

\proof
Suppose that $\defect(G/P,L)>0$. By Theorem \ref{NefValuesMorphismsFlags}
the nef value morphism is a homogeneous fiber bundle $\Phi:\,G/P\to G/Q$,
where $Q\subset G$ is another parabolic subgroup. The fiber $F$
of $\Phi$ is isomorphic to $Q/P$.

By Theorem \ref{FirstNefDef} we have $\Pic F=\Z$, therefore
$F$ is isomorphic to a quotient of a simple group by a maximal
parabolic subgroup and $Q$ is a proper subgroup of $G$.
By Theorem \ref{SecondNefDef}
$\defect(F,L)>\dim G/Q>0$. Therefore, the polarized
flag variety $(F,L)$ is one of the varieties from 
Theorem \ref{FirstClassificationFlag}.
Let us consider 3 cases.

If $F$ is equal to $A_l/P_2$ or $A_l/P_{l-1}$, $l\ge4$ even,
then $\defect(F,L)=2$. Therefore, $\dim G/Q=1$.
It follows that $G=\SL_2$, $Q$ is a Borel subgroup,
and there are no possibilities for $P$.

If $F$ is equal to
$B_4/P_4$, $D_5/P_4$, or $D_5/B_5$
then $\defect(F,L)=4$. Therefore, $\dim G/Q$ is equal to $1$, $2$, or $3$.
It follows that $G$ is equal to $\SL_2$, $\SL_3$, $\SL_4$, or $\Sp_4$.
Again none of these groups contain $B_4$ or $D_5$ as a simple
component of a Levi subgroup.

Finally, let $F$ be 
$A_l/P_1$, $A_l/P_l$, or $C_{l+1\over2}/P_1$,
so $F$ is isomorphic to $\P^l$.
We want to show that $\dim G/Q>k$, which will be a contradiction.
The set of simple roots $\Pi_Q$ contains a subset $\Psi$
with Dynkin diagram
of type $A_l$ or $C_{l+1\over2}$. Moreover, there exists
a simple root 
$\alpha\in\Pi$
such that $\alpha\not\in\Pi_Q$ and the Dynkin diagram of $\Psi'=\Psi\cup\{\alpha\}$
is connected. Then $\Psi'$ is a set of simple roots
of a simple subgroup $G'\subset G$ and $\dim G/Q\ge\dim G'/Q'$,
where $Q'=G'\cap Q$ is a maximal parabolic subgroup defined by $\Psi$.
It is easy to check that in all arising cases $\dim G'/Q'>k$.
\endproof

Now we can handle the case of an arbitrary semisimple group
and its arbitrary parabolic subgroup.

\begin{Theorem}\label{ThirdClassificationFlag}
Suppose that $(G/P,L)$ is a polarized flag variety
of a complex semisimple algebraic group $G$.
Then $\defect(G/P,L)>0$ if and only if either $(G/P,L)$
is one of polarized varieties described in 
Theorem~\ref{FirstClassificationFlag} or
the following conditions hold:
\begin{itemize}
\item $G=G_1\times G_2$, $P=P_1\times P_2$, 
where $P_1\subset G_1$, $P_2\subset G_2$.
\item $L=\pr_1^*L_1\otimes\pr_2^*L_2$, where $\pr_i:\,G/P\to G_i/P_i$
are projections and $L_i$ are very ample line bundles on $G_i/P_i$.
\item $(G_1/P_1,L_1)$ is one of 
polarized varieties described in 
Theorem~\ref{FirstClassificationFlag}.
\item $\defect(G_1/P_1,L_1)>\dim G_2/P_2$.
\end{itemize}
In this case $\defect(G/P,L)=\defect(G_1/P_1,L_1)-\dim G_2/P_2$.
\end{Theorem}

\proof
This Theorem easily follows from Theorem~\ref{ProductTheorem}.
\endproof

\chapter{Some Applications of the Duality} 

\section*{Preliminaries}

In this chapter we collect several applications, variations, and
illustrations of the  projective duality.

\section{Discriminants and Automorphisms}

\subsection{Matsumura--Monsky Theorem}

Let $D\subset\P^n$ be a smooth hypersurface of degree $d$.
It was first proved in~\cite{MM} that the group of projective
automorphisms preserving $D$ is finite if $d>2$.
In fact, it was also proved that the group of biregular
automorphisms of $D$ is finite if $d>2$ (except the cases
$d=3$, $n=2$ and $d=4$, $n=3$). Though this generalization
looks much stronger, actually it is an easy consequence of the
``projective'' version and
the Bart Theorem~\cite{Ba}.

More generally, let $G/P$ be a flag variety of a simple
Lie group and $D\subset G/P$ be a smooth ample divisor.
Let $L_\lambda=\O(D)$ be the corresponding ample line bundle,
where $\lambda\in\PP^+$ is a dominant weight.
Then one might expect that the normalizer $N_G(D)$ of $D$ in $G$ is finite
if $\lambda$ is big enough. Indeed, if $N_G(D)$
(or actually any linear algebraic group of transformations of $D$)
contains a one-parameter subgroup of automorphisms of $D$
then $D$ is covered by rational curves.
However, if $\lambda$ is big enough then the canonical
class $K_D$ is nef by the adjunction formula, 
and, therefore, $D$ can not be covered
by rational curves. Unfortunately, this transparent approach
does not give strong estimates on $\lambda$.
Much better estimates can be obtained using an original
proof of Matsumura-Monsky theorem.\index{Matsumura-Monsky theorem}

The problem can be reformulated as follows.
Suppose that $V_\lambda$ is an irreducible $G$-module
with highest weight $\lambda$. Let $\D\subset V_\lambda$
be the discriminant variety (the dual variety
to the orbit of the highest weight vector in $V_\lambda^*$,
see Section~\ref{JKHGFJKHGDKH} for further details).
We shall show that if $\lambda$ is big enough then any
point $x\in V_\lambda\setminus\D$ has a finite stabilizer $G_x$
and the orbit of $x$ is closed, $Gx=\overline{Gx}$
(therefore $x$ is a stable point of $V_\lambda$\index{stable point}
in the sense of Geometric Invariant Theory, see~\cite{MFK}).
This result can be compared with the results of 
\cite{AVE}, where all irreducible modules of simple algebraic groups
with infinite stabilizers of generic points were found
(this classification was extended later in \cite{El1} and \cite{El2}
to handle irreducible representations of semisimple groups
and any representations of simple groups as well).

If $\D$ is not a hypersurface then an easy inspection
using Theorem~\ref{ThirdClassificationFlag} shows that
the stabilizer of any point is infinite.
So from now on we shall assume that $\D$ is a hypersurface
defined by vanishing of the discriminant $\Delta$.

Recall that a dominant weight $\lambda$ is called self-dual if
$V_\lambda$ is isomorphic to $V_\lambda^*$ as a $G$-module.
Let $\PP^+_S\subset\PP^+$ be the subcone of self-dual
dominant weights. Let $\gamma$ be the highest root.

\begin{Theorem}
Let $V_\lambda$ be an irreducible representation of a simple algebraic group $G$
with highest weight $\lambda$ such that $\D$ is a hypersurface.
Suppose that $(\lambda-\gamma,\mu)>0$ for any $\mu\in\PP^+_S$.
Let $x\in V_\lambda\setminus\D$. Then $G_x$ is finite and $Gx=\overline{Gx}$.
Moreover, $G_{[x]}$ is also finite, where $[x]$ is the line spanned by $x$.
\end{Theorem}

\proof
If $G_x$ is finite and $G_{[x]}$ is infinite then 
$[x]\setminus\{0\}\subset Gx$. 
Therefore, $0\in\overline{Gx}$ and $\Delta(x)=\Delta(0)=0$,
hence, $x\in\D$.
The same argument shows that if $x\in V_\lambda\setminus\D$
then $\overline{Gx}\subset V_\lambda\setminus\D$.
If $Gx$ is not closed then $G_y$ is infinite for any point
$y\in \overline{Gx}\setminus Gx$. Therefore,
in order to prove the Theorem it suffices to prove that for any
$x\in V_\lambda\setminus\D$ the stabilizer $G_x$ is finite.

Suppose that $G_x$ is infinite. Then $G_y$ is infinite and reductive
(see e.g.~\cite{PV}) for any point $y$ from the closed orbit
in $\overline{Gx}$. Therefore, it suffices to prove that
if $S\subset G$ is a one-dimensional torus and $Sx=x$ then $x\in\D$.

Without loss of generality we may assume that $S\subset T$,
where $T$ is the fixed maximal torus, 
and $\t=\Lie T$ is the Cartan subalgebra.
Let $x=\sum_{\pi\in\PP}x_\pi$ be the weight decomposition of $x$.
Let $\Supp(x)=\{\pi\in\PP\,|\,x_\pi\ne0\}$ be the support of $x$.
For any $\mu\in\PP$ let $H_\mu$ denote the hyperplane of weights
perpendicular to $\mu$.
Then there exists $\mu\in\PP$ such that $\Supp(x)\subset H_\mu$.
Using the action of the Weil group we may assume that $\mu\in\PP^+$.

For any $\lambda\in\PP^+$ let $\lambda^*$ be the highest weight of the dual
module $V_{\lambda}^*$. Then $\lambda^*=-w_0(\lambda)$, where
$w_0$ is the longest element of the Weil group.

Suppose that $(\mu,\lambda^*-\gamma)>0$. Then for any positive root $\alpha$
we have $(\mu,\lambda^*-\alpha)>0$. Therefore,
$x$ is perpendicular to $[\g,v_{\lambda^*}]$, where $v_{\lambda^*}$
is the highest weight vector of $V_\lambda^*$.
It follows that $x\in\D$.

Suppose that $(\mu^*,\lambda^*-\gamma)>0$. Then for any positive root $\alpha$
we have $(\mu^*,\lambda^*-\alpha)>0$. Therefore,
$w_0(x)$ is perpendicular to $[\g,v_{\lambda^*}]$.
It follows that $w_0(x)\in\D$. Therefore,
It follows that $x\in\D$ as well.

Suppose now that
$(\mu,\lambda^*-\gamma)\le0$ and
$(\mu^*,\lambda^*-\gamma)\le0$.
Then $$(\mu+\mu^*,\lambda^*-\gamma)\le0.$$
But $\mu+\mu^*$ is a self-dual weight, hence this contradicts
assumptions of the Theorem.
\endproof

\begin{Example}\rm
If $G=\SL_n$ and $\lambda=\sum n_i\omega_i$, where
$\omega_1,\ldots,\omega_n$ are the fundamental weights,
then the assumptions of the Theorem are satisfied
if and only if $\sum n_i>2$, for example if $\lambda=n\omega_1$, $n>2$.
In particular, we recover the original Matsumura--Monsky Theorem.
\end{Example}

\subsection{Quasiderivations of Commutative Algebras}

\subsubsection*{Commutative algebras without identities.}

Let $V=\C^n$. Consider the vector space $\A=S^2V^*\otimes V$
parametrizing bilinear commutative multiplications in $V$. 
In the sequel we identify points of $\A$ with the
corresponding commutative algebras.

\begin{Definition}\rm
Let $A\in\A$. A non-zero element $v\in A$ is called 
a quadratic nilpotent if $v^2=0$.\index{quadratic nilpotent}
Let $\D_1\subset\A$ be a subset of all algebras containing
quadratic nilpotents.
A one-dimensional subalgebra $U\subset A$ 
is called singular if there exists linear independent\index{singular subalgebra}
vectors $u\in U$ and $v\in A$ such that
$$u^2=\alpha u,\quad uv={\alpha\over 2}v,\quad \text{where}\ \alpha\in\C.$$
Let $\D_2\subset\A$ be a subset of all algebras
containing singular subalgebras.
\end{Definition}

Then the following theorem holds.
\begin{Theorem}[\cite{T6}]\label{JHGVJHGFJHGDKHC}$\ $
\begin{enumerate}
\item[{\rm(1)}] $\D_1$ and $\D_2$ are irreducible hypersurfaces.
\item[{\rm(2)}] Let $A\in\A$. Then $A$ contains a one-dimensional subalgebra.
\item[{\rm(3)}] Let $A\in\A\setminus(\D_1\cup\D_2)$. Then
$A$ contains exactly $2^n-1$ one-dimensional subalgebras;
all these subalgebras are spanned by idempotents.
\end{enumerate}
\end{Theorem}

\proof
Clearly $\A$ can be identified with a set of
at most  $n$-dimensional linear systems of quadrics
in $V$. Namely, any linear function $f\in V^*$ 
defines a homogeneous quadratic function
$$v\to f(v^2).$$
Then $\D_1$ corresponds to the set of linear systems
of quadrics with zero resultant.
This proves that $\D_1$ is an irreducible hypersurface.

There exists, however, another useful identification. 
Any algebra $A$ determines the $n$-dimensional
linear system of {\it affine\/} quadrics
in $V$.
Namely, any linear form $f\in V^*$ 
defines an affine quadratic form
$$v\to f(v-v^2).$$
We can embed $V$ into a projective space $\P$ as an affine chart.
Then this linear system is naturally identified
with a $n$-dimensional linear system of quadrics in $\P$.
Clearly, the base points of this linear system that
do not lie on the infinity coincide with idempotents of $A$.
Infinite base points are the projectivizations of lines
spanned by quadratic nilpotents.
It is easy to see that quadrics from our linear system
intersect transversally at $0$.
Moreover, quadrics intersect non-transversally at some point
$v\ne0$ if and only if the subalgebra spanned by $v$
(if $v$ is finite) or the subalgebra with the projectivization $v$
(if $v$ is infinite) is singular. 
Therefore, (3) follows from the Bezout theorem.

$\D_2$ is irreducible since $\D_2=\GL_n\cdot\D_2'$,
where $D_2'\subset\A$ is the linear subset of all algebras
that have a fixed singular subalgebra
and a fixed line spanned by the vector $v$
from the definition of a singular subalgebra.
It is also quite easy to check that $D_2$ is 
actually a closed hypersurface.

Since quadrics of our linear system intersect
transversally at $0$, it follows that there exist
other base points, i.e.~there exists at least one $1$-dimensional subalgebra.
\endproof

The algebras $A\in \A$ that do not belong to discriminant
varieties $\D_1$ and $\D_2$ will be called regular.\index{regular algebra}
Of course, both hypersurfaces $\D_1$ and $\D_2$ can be interpreted
as ordinary discriminants.
First, we can enlarge the symmetry group and consider
$S^2(\C^n)^*\otimes \C^n$ as $\SL_n\times\SL_n$-module.
Then this module is irreducible and its discriminant variety
(the dual variety of the projectivization of the highest weight vector orbit)
coincides with $\D_1$.
Now, consider $\A=S^2V^*\otimes V$ as $\SL(V)$-module.
Then this module is reducible, $\A=\A_0+\tilde\A$,
where $\A_0$ is a set of algebras with zero trace
and $\tilde\A$ is isomorphic to $V^*$ as $\SL(V)$-module.
Consider the discriminant of $\A_0$ as a function on $\A$
(forgetting other coordinates).
Then the corresponding hypersurface is exactly $\D_2$.
If we consider the set of linear operators $\Hom(V,V)=V^*\otimes V$
instead of $\A$, then this construction will give us
the determinant and the discriminant of the linear operator
(see Example~\ref{Determinant}, Theorem~\ref{KJHKJHGJKFGH}).

\subsubsection*{Quasiderivations.}

Let $\g$ be a Lie algebra with representation $\rho:\,\g\to\End(V)$.
Consider any $v\in V$. Then subalgebra
$$\g_v=\{g\in\g\,|\,\rho(g)v=0\}$$ 
is called the annihilator of $v$. The subset\index{annihilator}
$$Q\g_v=\{g\in\g\,|\,\rho(g)^2v=0\}$$
is called the quasi-annihilator of\index{quasi-annihilator}
$v$. Clearly, $\g_v\subset Q\g_v$.
Of course, the quasi-annihilator is not a linear subspace in general,
However, we have the following version of a Jordan decomposition:\index{Jordan decomposition}

\begin{Lemma}\label{JHKGFJKDKHFDFHDJ}
Suppose that $\rho$ is the differential
of the representation of an algebraic group.
Let $g\in Q\g_v$. Consider the Jordan decomposition in
$\g$, $g=g_s+g_n$, where $g_s$ is semisimple and $g_n$ is nilpotent, $[g_s,g_n]=0$.
Then $g_s\in\g_v$ and $g_n\in Q\g_v$.
\end{Lemma}

The proof is obvious.

For example, suppose that $\g=\gl_n$
and $\rho$ is the natural representation in the vector space
$V^*\otimes V^*\otimes V$ that parametrizes bilinear multiplications in $V$.
Let $A$ be any algebra. Then the annihilator $\g_A$ 
is identified with the Lie algebra of derivations $\Der(A)$.
Operators $D\in Q\g_A$ are called quasiderivations.\index{quasiderivation}
Of course, it is possible to write down explicit
equations that determine  $Q\Der(A)=Q\g_A$ in $\End(A)$,
but this formula is quite useless (see~\cite{Vi1}).
We shall use its particular case that is quite easy to verify:

\begin{Lemma}
Let $A$ be an algebra, $D\in\End(A)$, $D^2=0$. 
Then $D\in Q\Der(A)$ if and only if for any
$x,y\in A$ we have
\begin{equation}
D(x)D(y)=D(D(x)y)+D(xD(y)).\label{quasidereq1}
\end{equation}
\end{Lemma}

\begin{Example}\rm
It was conjectured in \cite{Vi1} that
all quasiderivations of the algebra of matrices
$\Mat_n$ have the form  $D(x)=ax+xb$, where
$(a+b)^2=[a,b]$. Let us give an example of a different quasiderivation.
Consider the linear operator
$D(x)=exe$, where $e^2=0$, $e\ne0$. 
Then it is to check using (\ref{quasidereq1})
that $D$ is a quasiderivation.
However, $D$ of course can not be written in the form $ax+xb$.
\end{Example}

Quasiderivations can be used to define naive deformations.
Namely, suppose that $A$ is any algebra
with the multiplication $u\cdot v$ and $D$ is its quasiderivation.
Consider the algebra $A_D$ with the same underlying
vector space and with the multiplication
given by $u\star v=u\cdot D(v)+D(u)\cdot v-D(u\cdot v)$.
Then if $A$ satisfies any polynomial identity
then $A_D$ satisfies this identity as well.
More generally, let $\rho:\,\g\to\End(V)$
be the differential of a representation of an algebraic group $G$,
$v\in V$, $D\in Q\g_v$.
Suppose that $H\subset V$ is a closed conical $G$-equivariant hypersurface
(the vanishing set of a homogeneous $G$-semiinvariant)
and $v\in H$. Then $\rho(D)v$ also belongs to $H$.
Indeed, since $H$ is equivariant, $\exp(\lambda\rho(D))v$ belongs to $H$
for any $\lambda\in \C$. Since $D$ is a quasiderivation,
$\exp(\lambda\rho(D))v=v+\lambda\rho(D)v$.
Since $H$ is conical, $v/\lambda+\rho(D)v$ belongs to $H$.
Since $H$ is closed, $\rho(D)v$ also belongs to $H$.

Now suppose that $A$ is a regular commutative algebra.
Then we claim that $Q\Der(A)=0$.

\begin{Theorem}
Let $A\in\A$, $A\not\in\D_1\cup\D_2$. 
Then $Q\Der(A)=0$.
\end{Theorem}

\proof
By Lemma~\ref{JHKGFJKDKHFDFHDJ}
it is sufficient to check that
$A$ has no semisimple derivations
and no nilpotent quasiderivations.
Suppose that $\Der(A)$ contains a non-zero semisimple element.
Then the group $\Aut(A)$ contains a one-dimensional
algebraic torus $T$. Let $\t$ be its Lie algebra, $H\in\t$, $H\ne0$.
We way assume that the spectrum of $H$ in $A$
is integer-valued. 
Let $A_n\subset A$ be a weight space of weight $n$. 
Then $A=\oplus A_n$ is a $\Z$-grading.
Let $v\in A$ be a homogeneous element of a maximal positive
(or minimal negative) degree. 
Then $v^2=0$, hence $A\in\D_1$.

Suppose now that $E$ is a non-zero nilpotent quasiderivation.
We can embed $E$ in an $\ssl_2$-triple $\langle F,H,E\rangle\subset\End(A)$.
Let $$A=\oplus m_dR_d$$ be an $\ssl_2$-module decomposition, 
here $R_d$ is an irreducible $(d+1)$-dimensional module.
Consider also the weight decomposition
$$A=\oplus A^n\quad\text{and}\quad\A=\oplus \A^n$$ 
with respect to $H$. 
Let $$J^n=A^n+A^{n+1}+A^{n+2}+\ldots.$$
To avoid the abuse of notations denote by
$\alpha\in\A$ the point corresponding to the algebra $A$.
Let $\Supp\,\alpha$ be the support of $\alpha$ 
(i.e.~all $n\in\Z$ such that $\alpha_n\ne0$, where
$\alpha=\sum\alpha_n$, $\alpha_n\in\A^n$).
Since $E^2\alpha=0$, the vector $E\alpha$ is a linear combination
of highest weight vectors, therefore
$\Supp\,\alpha\subset\{-1,0,1,2,\ldots\}$.
This is equivalent to $J^nJ^m\subset J^{n+m-1}$.
In particular, if $v\in A$ is a weight vector
of the weight $n$ then $v^2=0$ as $n>1$. 
Since $A$ is regular, it follows that
$$A=m_0R_0\oplus m_1R_1=A^{-1}\oplus A^0\oplus A^1.$$
Since $A^1=J^1$, $A^1$ is a subalgebra in $A$.
By Theorem~\ref{JHGVJHGFJHGDKHC} (2),
$A^1$ has a one-dimensional subalgebra $U$ spanned by an idempotent $u$
(since $A$ has no quadratic nilpotents).
Let $v\in A^{-1}$ be a unique vector such that $Ev=u$.
Since $E^2=0$, we can apply formula~(\ref{quasidereq1})
with $x=y=v$. We get
$E(v)^2=2E(E(v)v)$, hence $u=2E(uv)$. 
Therefore,
$uv-{1\over2}v\in J^0$. Notice that the operator of left multiplication
by $u$ preserves $J^0$. 
Since this operator has an eigenvector
in $A/J^0$ with eigenvector $1/2$, it has such an eigenvector
in $A$. Therefore, $U$ is a singular subalgebra in $A$.
\endproof

The following Corollary was first proved 
in \cite{An} using combinatorial methods.

\begin{Corollary}[\cite{An}]
Let $A$ be an $n$-dimensional semisimple commutative algebra,
i.e.~$A$ is a direct sum of $n$ copies of $\C$,
that is $A$ is the algebra of diagonal $n\times n$ matrices.
Then $A$ has no nonzero quasiderivations.
\end{Corollary}

\proof
It is sufficient to check that $A\not\in\D_1\cup\D_2$. 
Clearly, $A$ has no nilpotents. 
Suppose that $U\subset A$ is a one-dimensional
subalgebra spanned by an idempotent $e\in U$.
Since the spectrum of the operator of left multiplication
by $e$ is integer-valued (actually consists of $0$ and $1$),
it follows that $U$ is not singular.
\endproof

\section{Discriminants of Anticommutative Algebras}

If the set of certain objects is parametrized
by an algebraic variety $X$ (for example, by a vector space)
then it makes sense to speak about generic objects.
Namely, we say that a generic object satisfies some property
if there exists a dense Zariski-open subset $X_0\subset X$
such that all objects parametrized by points from $X_0$
share this property. However, sometimes it is possible
to find some discriminant-type closed subvariety $Y\subset X$
and then to study properties of `regular' objects parametrized
by points from $X\setminus Y$.
For example, instead of studying generic hypersurfaces
it is usually worthy to study smooth hypersurfaces.
In this section we implement this 
program for the study of some quite non-classical object,
namely anticommutative algebras, with the multiplication
depending on a number of arguments.

\subsubsection*{Generic anticommutative algebras.}
Let $V=\C^n$. 
We fix an integer $k$, $1<k<n-1$. 
Let $\A_{n,k}=\Lambda^kV^*\otimes V$ be the vector space
of $k$-linear anticommutative maps from
$V$ to $V$. We identify the points of $\A_{n,k}$ 
with the corresponding algebras, that is, we assume that
$A\in \A_{n,k}$ is the space $V$ 
equipped with the structure of a $k$-argument
anticommutative algebra.
Subalgebras in generic algebras with $k=2$ 
were studied in \cite{T2}.
The following theorem is a generalization of these results.

\begin{Theorem}[\cite{T1}]\label{KJHFGKHDGJDH}
Let $A\in\A_{n,k}$ be a generic algebra. Then
\begin{enumerate}
\item[\rm(i)] Every $m$-dimensional subspace is a subalgebra if $m<k$.
\item[\rm(ii)] $A$ contains no $m$-dimensional subalgebras
with $k+1<m<n$.
\item[\rm(iii)] The set of $k$-dimensional subalgebras
is a smooth irreducible $(k-1)(n-k)$-dimensional
subvariety in the Grassmanian $\Gr(k,A)$.
\item[\rm(iv)] There are finitely many $(k+1)$-dimensional subalgebras, and their
number is
\begin{eqnarray*}
\sum_{
\matrix{
\scriptstyle n-k-1\ge\mu_1\ge\ldots\ge\mu_{k+1}\ge 0\cr
\scriptstyle n-k-1\ge\lambda_1\ge\ldots\ge\lambda_{k+1}\ge 0\cr
\scriptstyle\mu_1\le\lambda_1,\ldots,\mu_{k+1}\le\lambda_{k+1}
}
}\!\!\!\!\!\!\!\!\!
(-1)^{|\mu|}
\frac{(\lambda_1+k)!(\lambda_2+k-1)!\ldots\lambda_{k+1}!}
{(\mu_1+k)!(\mu_2+k-1)!\ldots\mu_{k+1}!}
(|\lambda|-|\mu|)!\times\cr\times
\left|\frac{1}{(i-j+\lambda_j-\mu_i)!}\right|^2_{i,j=1,\ldots,k+1},
\end{eqnarray*}
where $|\lambda|=\lambda_1+\ldots+\lambda_{k+1}$, 
$|\mu|=\mu_1+\ldots+\mu_{k+1}$, $1/N!=0$ if $N<0$.
\item[\rm(v)] $A$ contains a $(k+1)$-dimensional subalgebra.
\item[\rm(vi)] If $k=n-2$, then the number of $(k+1)$-dimensional subalgebras
is equal to $$2^n-(-1)^n\over 3.$$
\end{enumerate}
\end{Theorem}

\proof
The $\GL_n$-module $\A_{n,k}$ is a sum of two irreducible submodules:
\begin{equation}
\A_{n,k}=\A_{n,k}^0\oplus\tilde\A_{n,k}.\label{JHGFHGFJ(1)}
\end{equation}
Here $\tilde\A_{n,k}$ is isomorphic to $\Lambda^{k-1}V^*$:
we assign to every $(k-1)$-form
$\omega$ the algebra with multiplication
$$[v_1,\ldots,v_k]=\sum_{i=1}^k(-1)^{i-1}
\omega(v_1,\ldots,\hat v_i,\ldots,v_k)v_i.$$

Note that every subspace of this algebra is a subalgebra.
Hence, the lattice of subalgebras of $A\in\A_{n,k}$
coincides with the lattice of subalgebras of $A^0$, 
where $A\mapsto A^0$ is the $\GL_n$-equivariant
projector on the first summand in (\ref{JHGFHGFJ(1)}).
Algebras in $\A_{n,k}^0$ will be called {\it zero trace algebras\/},
since $A\in\A_{n,k}^0$ if and only if the $(k-1)$-form
$\Tr[v_1,\ldots,v_{k-1},\cdot]$ is equal to zero. 
Hence, the theorem will be proved
once we have proved
it for generic algebras in $\A_{n,k}^0$.

We choose a basis $\{e_1,\ldots,e_n\}$ in $V$, identify
$\GL_n$ with the group of matrices, 
consider the standard diagonal maximal torus $T$, 
and take for $B$ and $B_-$ the subgroups of upper- and lower-triangular
matrices. We fix an  $m\ge k$. 
Consider the parabolic subgroup of matrices
$$P=\left(\matrix{
A&0\cr
*&B\cr}\right),$$
where $B$ is an $m\times m$-matrix and $A$ is an $(n-m)\times(n-m)$-matrix.
Then $G/P$ coincides with $\Gr(m,V)$. 
Consider the vector bundle  
$$\L=\Lambda^k\S^*\otimes V/\S$$
on $G/P$, where $\S$ is the tautological bundle and 
$V/\S$ is the factor-tautological bundle. 
Then the assumptions of Theorem~\ref{JHGFJGFDJCHC} are fulfilled, since
$\L=\L_\lambda$, where $\lambda$ is the highest weight of $\A_{n,k}^0$. 
Therefore, 
$$\A_{n,k}^0=H^0(\Gr(m,V), \Lambda^k\S^*\otimes V/\S).$$
Let $A\in\A_{n,k}^0$, and let $s_A$ be the corresponding
global section.
Then the scheme of zeros 
$(Z_{s_A})_{red}$ coincides with the variety
of $m$-dimensional subalgebras of~$A$.

Let us return to the theorem.
Statement (i) is obvious.
(ii) follows from Theorem~\ref{JHGFJGFDJCHC} (i).
Theorem~\ref{JHGFJGFDJCHC} (ii) implies that if every $k$-argument
anticommutative algebra $A$ contains a $k$-dimensional subalgebra,
then the variety of $k$-dimensional subalgebras of a generic algebra
is a smooth unmixed $(k-1)(n-k)$-dimensional subvariety in $\Gr(k,A)$.
We claim that any $(k-1)$-dimensional subspace $U$
can be included in a $k$-dimensional subalgebra. 
The multiplication in the algebra defines a linear map from $V/U$ to $V/U$. 
Let $v+U$ be a non-zero eigenvector.
It is obvious that $\C v\oplus U$ is a $k$-dimensional subalgebra. 

We have only to prove that the variety
of $k$-dimensional subalgebras is irreducible.
Assume that a section $s$ of the bundle $\L=\Lambda^k\S^*\otimes V/\S$ 
over $\Gr(k,V)$ corresponding to $A$ has transversal intersection
with the zero section. Then the Koszul complex
$$0\rightarrow\Lambda^{n-k}{\L}^*\buildrel s\over \rightarrow\ldots
\buildrel s\over \rightarrow\Lambda^2{\L}^*\buildrel s\over \rightarrow
{\L}^*\buildrel s\over \rightarrow{\cal O}
\rightarrow{\cal O}_{Z(s)}\rightarrow0$$
is exact by Theorem~\ref{KoszulGeneral}.

Note that $\Lambda^p{\L}^*$ is isomorphic to the bundle
$S^p\Lambda^k\S\otimes \Lambda^p(V/\S)^*$. 
This is a homogeneous bundle over $G/P$ 
of the form $\L_\mu$, where
$$\mu=-\eps_1-\ldots-\eps_p+\eps_{n-k+1}+\ldots+\eps_n$$ 
and 
$\eps_i$ are the weights of the diagonal torus
in the tautological representation.
Note that the weight $\mu+\rho$
(where $\rho$ is the half-sum of the positive roots)
is singular (belongs to the wall of the Weil chamber)
for any $p$, $1\le p\le n-k$.
By the Borel--Weil--Bott theorem,
$H^*(\Gr(k,V), \Lambda^p{\L}^*)=0$ for $1\le p\le n-k$.
Hence, $H^0(Z(s), {\cal O}_{Z(s)})=H^0(\Gr(k, V), {\cal O})=\C$,
as was to be shown.

We postpone the proof of assertion (v) till the end of this section
and consider assertions (iv) and (vi), that is,
we shall calculate the highest Chern class of the bundle
$\Lambda^k\S^*\otimes V/\S$ over $\Gr(k+1, V)$. 

We use the standard notation, facts, and formulae from the Schubert
calculus (cf.~\cite{Fu1}. 
The letters $\lambda$ and $\mu$ always denote Young
diagrams
in the rectangle with $k+1$ rows
and $n-k-1$ columns.
It is well known that such diagrams
parametrize the basis of the Chow ring of
$\Gr(k+1, V)$. The cycle corresponding to $\lambda$
is denoted by $\sigma_\lambda$, 
$\sigma_\lambda\in A^{|\lambda|}(\Gr(k+1, V))$. (We grade the Chow ring
by the codimension, $|\lambda|=\lambda_1+\ldots+\lambda_{k+1}$, where
$\lambda_i$ is the length of the $i$th row of $\lambda$.)
We need the total Chern class of the bundles $\S$ and $V/\S$:
\begin{eqnarray*}
c(\S)=1-\sigma_1+\sigma_{1,1}-\ldots+(-1)^{k+1}\sigma_{1,\ldots,1},\cr
c(V/\S)=1+\sigma_1+\sigma_2+\ldots+\sigma_{n-k-1}.
\end{eqnarray*}

We begin by calculating the total Chern class of $\S\otimes V/\S$.
The standard formula for the total Chern class of the tensor products of two bundles implies that
\begin{equation}
c(\S\otimes V/\S)=\sum_{\mu\subset\lambda}d_{\lambda\mu}
\Delta_{\tilde\mu}(c(\S))\Delta_{\lambda'}(c(V/\S)),\label{JHGJHGJeqno(4)}
\end{equation}
where
$$d_{\lambda\mu}=\left|{\lambda_i+k+1-i
\choose \mu_j+k+1-j}\right|_{1\le i,j\le k+1},\quad
\Delta_\lambda(c)=|c_{\lambda_i+j-i}|,$$ 
$\lambda'=(n-k-1-\lambda_{k+1},n-k-1-\lambda_k,\ldots,n-k-1-\lambda_1)$,
and $\tilde\mu$ is the diagram obtained from $\mu$ by transposition.
We shall use the fact that
$\Delta_{\tilde\mu}(c(\S))=\Delta_{\mu}(s(\S))$,
where $s(E)$ is the Segre class of $E$.
Another fact is that 
$\Delta_{\lambda}(c(V/\S))=\sigma_{\lambda}$,
$\Delta_{\lambda}(s(\S))=(-1)^{|\lambda|}\sigma_{\lambda}$
(since $s(\S)=1-\sigma_1+\sigma_2+\ldots+(-1)^{n-k-1}\sigma_{n-k-1}$).
Therefore, (\ref{JHGJHGJeqno(4)}) can be written as
\begin{equation}
c(\S\otimes V/\S)=\sum_{\mu\subset\lambda}d_{\lambda\mu}(-1)^{|\mu|}
\sigma_\mu\sigma_{\lambda'}.\label{JHGFJGJeqno(5)}
\end{equation}
To calculate the highest Chern class of
$\L=\Lambda^k\S^*\otimes V/\S$, we use the formula
$\L=\Lambda^{k+1}\S^*\otimes(\S\otimes V/\S)$. 
The total Chern class of the first factor is equal to
$c(\Lambda^{k+1}\S^*)=1+\sigma_1$. 
The total Chern class of the second is given by (\ref{JHGFJGJeqno(5)}). 
Hence, the highest Chern class of
$\L$ is
\begin{equation}
c_{top}(\L)=\sum_{\mu\subset\lambda}d_{\lambda\mu}(-1)^{|\mu|}
\sigma_\mu\sigma_{\lambda'}\sigma_1^{|\lambda|-|\mu|}.\label{JHJHGJeqno(6)}
\end{equation}
The last result of the Schubert calculus that we need is the exact formula for the
degree of a product of two cycles.
In our case this can be written as
\begin{equation}
\sigma_\mu\sigma_{\lambda'}\sigma_1^{|\lambda|-|\mu|}=
\deg(\sigma_\mu\sigma_{\lambda'})=(|\lambda|-|\mu|)!
\left|\frac{1}{(i-j+\lambda_j-\mu_i)!}\right|_{i,j=1,\ldots,n},\label{JHGHGeqno(7)}
\end{equation}
where $1/N!=0$ if $N<0$.
The formula in assertion (iv) of the Theorem
can be obtained from (\ref{JHJHGJeqno(6)}) and (\ref{JHGHGeqno(7)}) by a slight modification
of the determinant
in the formula for $d_{\lambda\mu}$.

It remains to verify the formula in the assertion (vi)
of the Theorem. Let
$k=n-2$. Then the formula in the assertion (iv) can be written as
\begin{equation}
\sum_{0\le i\le j\le k+1}
(-1)^i\frac{(k+1)!k!\ldots(k-j+2)!(k-j)!\ldots0!}
{(k+1)!k!\ldots(k-i+2)!(k-i)!\ldots0!}
(j-i)!\,{\det}^2A,\label{HGVHGFeqno(8)}
\end{equation}
where
$$A=\left(\matrix{
X&0&0\cr
*&Y&0\cr
*&*&Z\cr
}\right),$$
$X$ is an $i\times i$-matrix, 
$Y$ is a $(j-i)\times (j-i)$-matrix, and
the format of $Z$ is a $(k+1-j)\times (k+1-j)$.
$X$ and $Z$ are lower-triangular matrices with $1$s on the diagonal and $Y$ is given by
$$Y=\left(\matrix{
1   &1   &0   &0&\cdots&0\cr
1/2!&1   &1   &0&\cdots&0\cr
1/3!&1/2!&1   &1&\cdots&0\cr
1/4!&1/3!&1/2!&1&\cdots&0\cr
\vdots&\vdots&\vdots&\vdots&\ddots&\vdots\cr
1/(j-i)!&1/(j-i-1)!&1/(j-i-2)!&1/(j-3)!&\cdots&1\cr
}\right).$$
It is easy to verify that $\det Y=1/(j-i)!$,
which enables us to rewrite (\ref{HGVHGFeqno(8)}) as
\begin{eqnarray*}
&\displaystyle\sum_{i\le j}
(-1)^i\frac
{(1+(k-i))!\ldots(1+(k-j+1))!}
{(k-i)!\ldots(k-j+1)!}
/(j-i)!=\cr
&\displaystyle\sum_{i\le j}
(-1)^i{k+1-i\choose j-i}=
\sum_i(-1)^i2^{k+1-i}={2^{k+2}-(-1)^{k+2}\over 3}=
{2^n-(-1)^n\over 3},\cr
\end{eqnarray*}
as was to be shown.

It remains to prove assertion (v) of the Theorem.
We have to prove that every
$A\in\A_{n,k}$ has a $(k+1)$-dimensional subalgebra.
We fix a basis $\{e_1,\ldots,e_n\}$ in $V$
and consider the subspace
$U=\langle e_{n-k}, \ldots, e_n\rangle$. 
Let $M\subset\A_{n,k}$ be the subspace that consists of the algebras for which
$U$ is a $(k+1)$-dimensional subalgebra.
We claim that $\A_{n,k}=\GL_n\cdot M$. 
It is sufficient to prove that the differential of the canonical morphism
$\phi:\,\GL_n\times M\to\A_{n,k}$ 
is surjective at a point $(e, A)$. 
Consider the algebra $A\in M$ in which
$[e_{n-k},\ldots,\hat e_i,\ldots, e_n]=e_i$
for all $n-k\le i\le n$ and the other products are zero.
We claim that $d\phi$ is surjective at $(e, A)$. 
Consider the map $\pi:\,\A_{n,k}\to\A_{n,k}/M $.
It is sufficient to verify that
$\pi\circ d\phi_{(e, A)}(\gl_n, 0)=\A_{n,k}/M\simeq\Lambda^kU^*\otimes V/U$.
But this is obvious, since multiplication in
$d\phi_{(e, A)}(E_{ji}, 0)$  for
$n-k\le i\le n$, $1\le j\le n-k-1$
($E_{ji}$ is the matrix identity)
is given by
$[e_{n-k},\ldots,\hat e_i,\ldots, e_n]=e_j$
with the other products equal to zero.
This completes the proof of the theorem.
\endproof

\subsubsection*{$D$-regular algebras.}
An essential drawback of Theorem~\ref{KJHFGKHDGJDH} 
is the fact that
it does not enable us to study the structure of subalgebras of any particular algebra.
To correct this, we shall need to introduce some explicit class
of `regular' algebras instead of implicit class of generic algebras.
The natural way to remove degeneracies is to consider 
discriminants.

Let $\dualA$ be the $\GL_n$-module dual to $\A_{n,k}^0$,
and let $S_D$ be the closure of the orbit of the highest
vector, $S_D\subset\dualA$. Let $PS_D\subset P\dualA$
be its projectivization, let 
$P\D\subset P\A_{n,k}^0$ be the subvariety projectively dual
to the subvariety $PS_D$, 
and let $\D\subset \A_{n,k}^0$ be the cone over it.
We shall call $\D$ the {\it $D$-discriminant subvariety}.
The algebras $A\in\D$ are said to be {\it $D$-singular}.
The algebras $A\not\in\D$ are said to be {\it $D$-regular}.

\begin{Theorem}[\cite{T1}]$\ $
\begin{enumerate}
\item[{\rm(i)}] $\D$ is a hypersurface.
\item[{\rm(ii)}] Let $A$ be a $D$-regular algebra. Then the set of 
$k$-dimensional subalgebras of $A$ is a smooth irreducible 
$(k-1)(n-k)$-dimensional subvariety in $\Gr(k,A)$.
\item[{\rm(iii)}] Let $k=n-2$. 
Then the degree of $\D$ is equal to
\begin{equation}
{(3n^2-5n)2^n-4n(-1)^n}\over{18}.\label{KJHGVJGJ(3)}
\end{equation}
\end{enumerate}
\end{Theorem}

Hence, the $D$-singularity of $A$ is determined
by the vanishing of the $SL_n$-invariant polynomial $D$
that defines~$\D$. This polynomial is called the $D$-discriminant.

\medskip

\proof
The fact that $\D$ is a hypersurface follows immediately
from the results of Section~\ref{JKHGFJKHGDKH}.
Assertion (ii) can be deduced from the corresponding assertion
of Theorem~\ref{KJHFGKHDGJDH} by an easy calculation with differentials.
Namely, let $A$ be a $D$-regular algebra.
We use the arguments
in the proof of Theorem~\ref{JHGFJGFDJCHC} (ii) 
and Theorem~\ref{KJHFGKHDGJDH} (iii).
According to these calculations, it is sufficient
to verify that if $U$ is a $k$-dimensional subalgebra of $A$,
then the map
$$\psi:\,\gl_n\to\Lambda^kU\otimes A/U,$$
$$\psi(g)(v_1\wedge\ldots\wedge v_k)=
g[v_1,\ldots,v_k]-[gv_1,\ldots,v_k]-\ldots-[v_1,v_2,\ldots,gv_k]+U$$
is surjective.
Assume the contrary.
Then there is a hyperplane $H\supset U$ 
such that the image of $\psi$
lies in $\Lambda^kU\otimes H/U$.
Consider a non-zero algebra $\tilde A$ in $S_D\subset\dualA$
such that
$[U^\perp,V^*,\ldots,V^*]=0$, $[V^*,\ldots,V^*]\subset H^\perp$,
where $U^\perp$ and $H^\perp$ are the annihilators of $U$ and $H$ in
$\dualA$. (Such an algebra is unique up to a scalar.)
Then $\tilde A$ annihilates $[\gl_n, A]$, which is equivalent
to the fact that $A$ annihilates $[\gl_n, \tilde A]$, that is,
the tangent space to $S_D$ at $\tilde A$. 
This means that $A$ lies in $\D$, that is, it is a $D$-singular algebra.

Finally, assertion (iii) follows from Theorem~\ref{JHGFJGHFJ}.
\endproof

\subsubsection*{$E$-regular algebras.}
We define the $E$-discriminant and $E$-regularity
only for $(n-2)$-argument $n$-dimensional anticommutative algebras. 
Let $\A=\A_{n,n-2}^0$.
Consider the projection $\pi:\,\Gr(n-1,V)\times P\A\to P\A$
on the second summand and the incidence subvariety
$Z\subset\Gr(n-1,V)\times P\A$ that consists of pairs $S\subset PA$,
where $S$ is a subalgebra in $A$.
Let $\tilde\pi=\pi|_Z$.
By Theorem~\ref{KJHFGKHDGJDH}, we have $\tilde\pi(Z)=\A$. 
Let $\tilde\E\subset Z$ be the set of critical points of $\tilde\pi$,
let $P\E=\tilde\pi(\tilde\E)$ be the set of critical values of
$\tilde\pi$, and let $\E\subset\A$ be the cone over $P\E$.
Then $\E$ is called the {\it $E$-discriminant subvariety}.
The algebras $A\in\E$ are said to be {\it $E$-singular}.
The algebras $A\not\in\E$ are said to be {\it $E$-regular}.

\begin{Theorem}[\cite{T1}]$\ $
\begin{enumerate}
\item[{\rm(i)}] $\E$ is an irreducible hypersurface.
\item[{\rm(ii)}] Let $A$ be an $E$-regular algebra. 
Then $A$ has precisely
$$2^n-(-1)^n\over 3$$ 
$(n-1)$-dimensional subalgebras.
\item[{\rm(iii)}] The map $\tilde\pi:\,\tilde\E\to P\E$ is birational.
\end{enumerate}
\end{Theorem}

Hence, the $E$-singularity of $A$ is determined by the vanishing of the 
$SL_n$-invariant polynomial that defines~$\E$. 
This polynomial is called the {\it $E$-discriminant}. 
Assertion (iii) can be formulated as follows: a generic
$E$-singular algebra has precisely one ``critical'' 
$(n-1)$-dimensional subalgebra.
\medskip

\proof
The proof of assertion (ii) is similar to the proof 
of Theorem~\ref{KJHFGKHDGJDH} (iv).
We claim that statement (i) follows from statement (iii).
We choose in $V^*$ a basis $\{f_1,\ldots,f_n\}$
dual to the basis $\{e_1,\ldots,e_n\}$ in $V$.
Let $U\in\Gr(n-1,V)$ be the hyperplane $f_1=0$.
It is clear that $\tilde\E=\GL_n\cdot M_0$, where
$M_0=\tilde\E\cap (U, P\A)$.
Moreover, $M=Z\cap (U, P\A)$ is the linear subspace of the algebras
for which
$U$ is an $(n-1)$-dimensional subalgebra.
Let $P$ be the parabolic subgroup
of matrices
$\left(\matrix{ A&0\cr *&B\cr}\right)$,
and let $\goth u$ be the Lie algebra of matrices
$\left(\matrix{0&X\cr 0&0\cr}\right)$, where
$B$ is an $(n-1)\times(n-1)$-matrix and
$X$ is an $1\times(n-1)$-matrix.
Every algebra $A\in M$ defines a linear map
$\goth u\to\Lambda^{n-1}U^*\otimes V/U$.
Since $A\in M_0$ if and only if this map is degenerate,
we have $\codim_MM_0=1$.
It is obvious that $M_0$ is irreducible, since
$M_0$ is the spreading of the subspace
$$M_0^1=\{A\in M_0\,|\,\text{the algebra}\ E_{12}A
\ \text{has a subalgebra}\  U\}$$
by the group $P$.
Hence, $\tilde\E$ is an irreducible divisor in 
$Z$, and assertion (i) follows from assertion (iii).

We prove (iii).
We shall say that an $(n-1)$-dimensional subalgebra $U'$
of $A\in\E$ is {\it critical\/} if $(U', A)$ lies in $\tilde\E$.
In this case there is an $(n-2)$-dimensional
subspace $W'\subset U'$ such that if 
$V=U'\oplus\Bbb Ce$ and $v\in\gl_n$ is a non-zero linear operator
such that  $v(V)\subset \Bbb Ce$ and $v(W')=0$, 
then $U'$ is a subalgebra of $vA$.
To prove assertion (iii) it is sufficient to prove
that in generic algebras in
$M_0^1$ the subalgebra $U$ is the unique
critical subalgebra. Let
$N\subset M_0^1$ be the subvariety of all algebras that have another 
critical subalgebra.

Note that $M_0^1$ is normalized by the parabolic subgroup 
$$Q=\left(\matrix{
*&0&0&\ldots&0\cr
*&*&0&\ldots&0\cr
*&*&*&\ldots&*\cr
\vdots&\vdots&\vdots&\ddots&\vdots\cr
*&*&*&\ldots&*\cr}\right).$$
Then $N$ is the spreading of the subvarieties
$N_2$ and $N_3$ by the group $Q$, where $A\in N_i$ 
if and only if the hyperplane $f_i=0$
is a critical subalgebra.
In turn, $N_2$ is the spreading of the vector spaces
$N_2^1$ and $N_2^3$, and
$N_3$ is the spreading of the subspaces $N_3^1$, $N_3^2$, and $N_3^4$, 
where $N_i^j\subset N_i$ is the subspace of all algebras
such that the $(n-2)$-dimensional subspace $W'$ (mentioned above)
is given by $f_i=f_j=0$. 
Let $Q_i^j\subset Q$ be the subgroup that normalizes the flag
$f_i\subset\langle f_i,f_j\rangle$.
It is easy to verify that
$$\codim_QQ_2^1=1,\ \codim_QQ_2^3=n-1,\ \codim_QQ_3^1=n-1,$$
$$\codim_QQ_3^2=n,\ \codim_QQ_3^4=2n-3.$$

On the other hand,
$$\codim_{M_0^1}N_2^1=n,$$
$$\codim_{M_0^1}N_2^3=\codim_{M_0^1}N_3^1=\codim_{M_0^1}N_3^2
=\codim_{M_0^1}N_3^4=2n-2.$$
Hence, $\codim_{M_0^1}QN_i^j\ge1$ in all cases,
which completes the proof of the theorem.
\endproof

\subsubsection*{Regular $4$-dimensional anticommutative algebras.}
An $(n-2)$-argu\-ment $n$-dimensional anticommutative algebra is said to be 
{\it regular} if it is $D$-regular and $E$-regular.
In this subsection we consider $2$-argument
$4$-dimensional algebras.
The corresponding generic algebras
were studied in \cite{T2}. 

\begin{Theorem}[\cite{T1}]\label{KJHHJKGFJHCJHCK}
Let $A$ be a $4$-dimensional regular anticommutative algebra.
Then  
\begin{enumerate}
\item[{\rm(i)}] $A$ has precisely five $3$-dimensional subalgebras.
The set of these subalgebras
is a generic configuration of five hyperplanes.
In particular, $A$ has a pentahedral normal form,
that is, it can be reduced by a transformation
that belongs to  $\GL_4$ to an algebra such that the set of its five subalgebras is
a Sylvester pentahedron
$x_1=0$, $x_2=0$, $x_3=0$, 
$x_4=0$, $x_1+x_2+x_3+x_4=0$.
\item[{\rm(ii)}] $A$ has neither one- nor two-dimensional ideals.
\item[{\rm(iii)}] The set of two-dimensional
subalgebras of $A$ is a del Pezzo surface of degree~$5$ 
{\rm(}a blowing up of $\Bbb P^2$ at four generic points{\rm)}.
\item[{\rm(iv)}] $A$ has precisely $10$ fans, that is,
flags $V_1\subset V_3$ of
$1$-dimensional and $3$-dimensional
subspaces such that every intermediate subspace
$U$, $V_1\subset U\subset V_3$,
is a two-dimensional subalgebra.
\end{enumerate}
\end{Theorem}

\proof
We begin with assertion (i).
We have to prove that if $A$ is a $4$-dimensional
regular anticommutative algebra with zero trace,
then the set of its three-dimensional subalgebras $S_1,S_2,\ldots,S_5$ 
is a generic configuration of hyperplanes, that is,
the intersection of any three of them
is one-dimensional
and the intersection
of any four of them
is zero-dimensional.
Indeed, assume, for example,
that $U=S_1\cap S_2\cap S_3$ is two-dimensional.
Let $v\in U$ and $v\ne0$. Then $[v,\cdot]$ induces a linear operator on
$A/U$, since $U$ is a subalgebra.
$S_1/U$, $S_2/U$, and $S_3/U$ are one-dimensional eigenspaces.
Since $\dim A/U=2$, the operator is a dilation.
Since this is true for any $v\in U$, any three-dimensional
subspace that contains $U$ is a three-dimensional subalgebra,
which contradicts the fact
that there are precisely five such subalgebras.

Now assume that $U=S_1\cap S_2\cap S_3\cap S_4$ is
one-dimensional, and let $v\in U$, $v\ne0$.
Then the operator $[v,\cdot]$ induces an operator on $A/U$. 
This operator has four two-dimensional eigenspaces
$S_1/U,\ldots,S_4/U$ of which any three have zero intersection.
Hence, this operator is a dilation.
Let $W\supset U$ be an arbitrary two-dimensional subspace,
and let  $w\in W$ be a vector that is not proportional to $v$. 
Then the operator $[w,\cdot]$ induces a linear operator on $A/W$. 
Let $z$ be a non-zero eigenvector.
Then  $\langle v,w,z\rangle$ is a three-dimensional subalgebra.
Therefore, every vector can be included in a three-dimensional subalgebra,
which contradicts the fact
that there are only five such subalgebras.

This argument also shows that $A$ has no one-dimensional ideals.
Since every three-dimensional subspace that contains a two-dimensional ideal
is a subalgebra, there are no two-dimensional ideals, which completes the proof
of assertion (ii).

To prove assertion (iii), we consider the subvariety
$X$ of two-dimensional subalgebras in $A$. 
Then $X\subset\Gr(2,4)$.
Consider the Pl\"ucker embedding 
$\Gr(2,4)\subset\Bbb P^5=P(\Lambda^2\Bbb C^4)$.
First we claim that the embedding $X\subset\Bbb P^5$ 
is non-degenerate, that is, the image is contained in no hyperplane.
Let $S_1,S_2,S_3,S_4$ be four three-dimensional subalgebras.
Since it is a generic configuration, we can choose
a basis
$\{e_1,e_2,e_3,e_4\}$ in $A$ such that
$S_i=\langle e_1,\ldots,\hat e_i,\ldots,e_4\rangle$. 
Since the intersection of three-dimensional subalgebras
is a two-dimensional subalgebra,
$A$ has six subalgebras
$\langle e_i, e_j\rangle$, $i\ne j$. 
The set of corresponding bivectors
 $e_i\wedge e_j$ is a basis in $\Lambda^2\Bbb C^4$.
Therefore, they can lie in no hyperplane.
Simple calculation with Koszul complexes (cf.~Theorem~\ref{KoszulGeneral})
shows that $H^0(X,{\cal O}_X(1))^*=\Lambda^2\Bbb C^4$.

To prove that  $X$ is a del Pezzo surface
of degree five, we have only to verify
that ${\cal O}_X(1)$ coincides with the anticanonical sheaf (see~\cite{Ma}).
Since $Y=\Gr(2,4)$ is a quadric in $\Bbb P^5$,
we have 
$$\omega_Y={\cal O}_Y(2-5-1)={\cal O}_Y(-4).$$
The set $X$~is a non-singular subvariety of codimension~$2$
in~$Y$. Therefore, $\omega_X=\omega_Y\otimes\Lambda^2{\cal N}_{X/Y}$,
where ${\cal N}_{X/Y}$ is the normal sheaf.
Further, $X$ is the scheme of zeros
of a regular section of the vector bundle
${\cal L}=\Lambda^2\S^*\otimes V/\S$, whence
${\cal N}_{X/Y}={\cal L}|_Y$ and
$\Lambda^2{\cal N}_{X/Y}={\cal O}_X(3)$, since
$c_1({\cal L})=3H$. We obtain that
$\omega_X={\cal O}_X(-4)\otimes{\cal O}_X(3)={\cal O}_X(-1)$, 
as was to be shown.

It remains to prove assertion (iv).
Since $X$ is a del Pezzo surface of degree five,
it contains ten straight lines.
Since the embedding $X\subset P(\Lambda^2\Bbb C^4)$ is anticanonical,
these straight lines
are ordinary straight lines in
$P(\Lambda^2\Bbb C^4)$ that lie in $\Gr(2,4)$. 
It remains to establish a bijection between these straight lines
and fans.
If $b\in\Lambda^2\Bbb C^4$, then
$b$ belongs to the cone over $\Gr(2,4)$ if and only if
$b\wedge b=0$. If $b_1$ and $b_2$ belong to this cone,
then the straight line that joins them
belongs to this cone
if and only if $b_1\wedge b_2=0$,
which coincides with the fan condition.
\endproof

\subsubsection*{Dodecahedral section.}
Let us start with some definitions.
Let $X$ be an irreducible $G$-variety
(a variety with an action of algebraic group $G$),
$S\subset X$ be an irreducible subvariety.
Then $S$ is called a {\it section\/} of $X$
if $\overline{G\cdot S}=X$.
The section $S$ is called a relative section if
the following condition holds:
there exists a dense Zariski-open
subset $U\subset S$ such that
if $x\in U$ and $gx\in S$ then $g\in H$, where
$H=N_G(S)=\{g\in G\,|\,gS\subset S\}$ is the normalizer
of $S$ in~$G$ (see~\cite{PV}).
In this case for any invariant function
$f\in\Bbb C(X)^G$ the restriction $f|_S$
is well-defined and the map
$$\Bbb C(X)^G\to\Bbb C(S)^H,\quad f\mapsto f|_S,$$
is an isomorphism.
Any relative section defines a $G$-equivariant rational map
$\psi:\,X\to G/H$: if $g^{-1}x\in S$ then $x\mapsto gH$.
Conversely, any $G$-equivariant rational map
$\psi:\,X\to G/H$ with irreducible fibers
defines the relative section
$\overline{\psi^{-1}(eH)}$.

We are going to apply Theorem~\ref{KJHHJKGFJHCJHCK}
and to construct a relative section in the
$\SL_4$-module ${\cal A}_0$ (the module of $4$-dimensional
anticommutative algebras with zero trace).
The action of $\SL_4$ on `Sylvester pentahedrons'
is transitive with finite stabilizer $H$
(which is the central extension of the permutation group $S_5$).
In the sequel the Sylvester pentahedron will always mean
the standard configuration formed by the hyperplanes
$$x_1=0,\ x_2=0,\  x_3=0,\ x_4=0,\ x_1+x_2+x_3+x_4=0.$$
Let $S\subset{\cal A}_0$ be a linear subspace 
formed by all algebras such that
the hyperplanes of Sylvester pentahedron are their
subalgebras.
Then Theorem~\ref{KJHHJKGFJHCJHCK} implies that $S$
is a $5$-dimensional linear relative section
of $\SL_4$-module ${\cal A}_0$.

It is easy to see
that the multiplication in algebras from $S$
is given by formulas
$$[e_i\mathop{,}e_j]=a_{ij}e_i+b_{ij}e_j\quad (1\le i<j\le4),$$
where $a_{ij}$ and $b_{ij}$ satisfy the certain set of linear conditions.
Consider $6$ algebras $A_1,\ldots,A_6$
with the following structure constants

\bigskip
{\hfill\vbox{\offinterlineskip
\hrule
\halign
{&\vrule#&
  \strut\quad\hfil#\quad\hfil&        \vrule#&
  \strut\quad\hfil#\quad\hfil&        \vrule#&
  \strut\quad\hfil#\quad\hfil&        \vrule#&
  \strut\quad\hfil#\quad\hfil&        \vrule#&
  \strut\quad\hfil#\quad\hfil&        \vrule#&
  \strut\quad\hfil#\quad\hfil&        \vrule#&
  \strut\quad\hfil#\quad\hfil&        \vrule#
 \cr
&           && $A_1$ && $A_2$ && $A_3$ && $A_4$ && $A_5$ && $A_6$&\cr
\noalign{\hrule}\cr
& $a_{12}$  && $0$   && $1$   && $-1$   && $1$  && $0$    && $-1$&\cr
& $b_{12}$  && $1$   && $0$   && $1$  && $-1$   && $-1$   && $0$&\cr
& $a_{13}$  && $1$   && $1$   && $-1$   && $0$  && $-1$   && $0$&\cr
& $b_{13}$  && $0$   && $-1$  && $0$   && $1$   && $1$    && $-1$&\cr
& $a_{14}$  && $1$   && $0$   && $0$   && $1$   && $-1$   && $-1$&\cr
& $b_{14}$  && $-1$  && $1$   && $-1$   && $0$  && $0$    && $1$&\cr
& $a_{23}$  && $-1$  && $-1$  && $0$   && $1$   && $0$    && $1$&\cr
& $b_{23}$  && $1$   && $0$   && $-1$   && $0$  && $1$    && $-1$&\cr
& $a_{24}$  && $0$   && $-1$  && $-1$   && $0$  && $1$    && $1$&\cr
& $b_{24}$  && $-1$  && $1$   && $0$   && $1$   && $-1$   && $0$&\cr
& $a_{34}$  && $-1$  && $1$   && $1$  && $-1$   && $0$    && $0$&\cr
& $b_{34}$  && $0$   && $0$   && $-1$   && $1$  && $-1$   && $1$&\cr
\cr}
\hrule}\hfill}
\bigskip

Then it is easy to see that $A_i\in S$ for any $i$.
Moreover, algebras
$A_i$ satisfy the unique linear relation
$A_1+\ldots+A_6=0$. It follows that any
$A\in S$ can be written uniquely in the form
$\alpha_1A_1+\ldots+\alpha_6A_6$, 
where $\alpha_1+\ldots+\alpha_6=0$.
The coordinates $\alpha_i$ are called dodecahedral coordinates
and $S$ is called the dodecahedral section
(this name will be clear later).

The stabilizer of the standard Sylvester pentahedron
in $\PGL_4$ is isomorphic to $\Bbb S_5$
represented by permutations of its hyperplanes.
The group $S_5$ is generated by the transposition $(12)$ 
and the cycle $(12345)$.
The preimages of these elements in $\GL_4$ 
are given by matrices
$$\sigma=\left(\matrix{
0&1&0&0\cr
1&0&0&0\cr
0&0&1&0\cr
0&0&0&1\cr
}\right)\qquad\text{¨}\quad
\tau=\left(\matrix{
-1&-1&-1&-1\cr
1&0&0&0\cr
0&1&0&0\cr
0&0&1&0\cr
}\right).
$$
The preimage of $S_5$ in $\SL_4$ is the group $H$
of $480$ elements.
The representation of $H$ in $S$ 
induces the projective representation of
$S_5$ in $\P^4$. We have the following

\begin{Proposition}
This projective representation is the projectivization
of the $5$-dimensional irreducible representation of $S_5$
(any of two possible).
\end{Proposition}

\proof
Recall certain `folklore' facts about the representation theory of  $S_5$.
It is well-known that $S_5$ admits exactly two
embeddings in $S_6$ up to conjugacy. 
One is standard via permutations of the first five elements
of the six-element set permuted by $S_6$.
The second one can be obtained from the first one
by taking the composition
with the unique (up to conjugacy) outer involution of~$S_6$.
$S_5$ has exactly two irreducible $5$-dimensional
representations, which have the same projectivizations.
One of these representations has the following model.
One takes the tautological 
$5$-dimensional irreducible representation of $S_6$
and consider its composition with the non-standard
embedding $S_5\subset S_6$.

Now let us return to our projective representation.
The action of $\sigma$ and $\tau$ on algebras $A_i$
is given by formulas
\begin{equation}
\left(\matrix{
\sigma A_1\cr
\sigma A_2\cr
\sigma A_3\cr
\sigma A_4\cr
\sigma A_5\cr
\sigma A_6\cr
}\right)=
\left(\matrix{
-A_2\cr
-A_1\cr
-A_4\cr
-A_3\cr
-A_6\cr
-A_5\cr
}\right)\qquad\text{and}\qquad
\left(\matrix{
\tau A_1\cr
\tau A_2\cr
\tau A_3\cr
\tau A_4\cr
\tau A_5\cr
\tau A_6\cr
}\right)=
\left(\matrix{
A_1\cr
A_6\cr
A_2\cr
A_3\cr
A_4\cr
A_5\cr
}\right).\label{KJHGJHGKJGK1}
\end{equation}

Therefore, $S_5$ permutes the lines spanned by $A_i$.
Moreover, the induced embedding is clearly non-standard:
transposition in $S_5$ maps to the composition
of $3$ independent transpositions in $S_6$.
It is clear that
that the corresponding 
projective representation of $S_5$
is isomorphic to the projectivization
of the $5$-dimensional irreducible representation 
in the model described above.
\endproof

\begin{Remark}\rm
The section $S$ is called dodecahedral by the following reason.
Though the surjection $H\to S_5$ does not split,
the alternating group $A_5$ can be embedded in $H$.
The induced representation of $A_5$ in $S$ has the following description.
$\Bbb A_5$ can be realised as a group of rotations 
of dodecahedron.
Let $\{\Gamma_1,\ldots,\Gamma_6\}$ be the set of pairs of opposite
faces of dodecahedron.
Consider the vector space of functions
$$f:\{\Gamma_1,\ldots,\Gamma_6\}\rightarrow\C,\quad
\sum\limits_{i=1}^6f(\Gamma_i)=0.$$
Then this vector space is an $A_5$-module.
It is easy to see that this module is isomorphic to $S$
via the identification $A_i\mapsto f_i$, where
$f_i(\Gamma_i)=5$, $f_i(\Gamma_j)=-1$, $j\ne i$.
\end{Remark}

The following proposition follows
from the discussion above

\begin{Theorem}
The restriction of invariants induces an isomorphism
of invariant fields
$$\C({\cal A}_0)^{GL_4}\simeq\C(\C^5)^{\C^*\times \Bbb S_5},$$
where
$\C^*$ acts on $\C^5$ by homotheties and
$\Bbb S_5$ acts via any of two $5$-dimensional
irreducible representations.
\end{Theorem}

\begin{Remark}\rm
The Sylvester pentahedron
also naturally arises in the theory of cubic surfaces.
The $\SL_4$-module of cubic forms $S^3(\C^4)^*$
admits the relative section (the so-called
Sylvester section, or Sylvester normal form).
Namely, a generic cubic form in a suitable
system of homogeneous coordinates $x_1,\ldots,x_5$,
$x_1+\ldots+x_5=0$, can be written as a sum of $5$ cubes
$x_1^3+\ldots+x_5^3$.
The Sylvester pentahedron can be recovered from a generic cubic
form $f$ in a very interesting way: its $10$ vertices coincide
with $10$ singular points of a quartic surface $\det\Hes(f)$.
The Sylvester section has the same normalizer $H$ as our dodecahedral section.
It can be proved \cite{Bek} that in this case
the restriction of invariants induces an isomorphism
$$\C(S^3(\C^4)^*)^{GL_4}\simeq\C(\C^5)^{\C^*\times \Bbb S_5},$$
where $\C^*$ acts via homotheties and
$\Bbb S_5$ via permutations of coordinates
(i.e.~via the {\it reducible} $5$-dimensional representation).
Other applications of the Sylvester pentahedron
to moduli varieties can be found in~\cite{Bar}.
\end{Remark}

These results were used in~\cite{T2} in order to prove that
the field of invariant functions of the $5$-dimensional
irreducible representation of $S_5$ is rational
(is isomorphic to the field of invariant functions
of a vector space).
From this result it is easy to deduce that in fact
the field of invariant functions of any
representation of $S_5$ is rational (see~\cite{T5}).

\section{Self--dual Varieties}

\subsection{Self--dual Polarized Flag Varieties}

A projective variety $X\subset\P^n$ is called\index{self-dual variety}
self-dual if $X$ is isomorphic to $\dual X$ as 
an embedded projective variety, i.e.~there exists
an isomorphism $f:\,\P^N\to\dual{\P^n}$ such that 
$f(X)=\dual X$.
A lot of examples is provided by the Pyasetskii pairing~\ref{PyasetskyTh}.
Indeed, if $G$ is a connected reductive algebraic group acting on
a vector space $V$ then the dual
action $G: V^*$ can be described as follows.
Having passed to a finite covering of $G$,
one may assume that there is an involution $\theta\in\Aut(G)$
(a so-called involution of maximal rank) such that
$\theta(t)=t^{-1}$ for any $t\in T$, where $T$
is the fixed maximal torus in $G$.
Along with the action of $G$ on $V$ we may
consider the twisted action given by
$g\star v=\omega(g)v$.
Then it is easy to see that the twisted action is 
isomorphic to the dual action $G:V^*$.
Therefore, we may canonically identify $G$-orbits in $V$
and in $V^*$.
In particular, if $G$ acts on $V$ with finitely many orbits
then we may consider the Pyasetskii pairing 
(hence the projective duality of projectivizations)
as pairing between $G$-orbits in $V$.

For example, 
let $L$ be a simple algebraic group and $P$ a parabolic subgroup
with Abelian unipotent radical. In this case $\l=\Lie L$ admits
a short $\Z$-grading with only three non-zero parts:
$$\l=\l_{-1}\oplus\l_0\oplus\l_1.$$
Here $\l_0\oplus\l_1=\Lie P$
and $\exp(\l_1)$ is the Abelian unipotent radical of $P$.
Let $G\subset L$ be a reductive subgroup with Lie algebra
$\l_0$. 
Recall that by Theorem~\ref{AURAOrbitDecomp}
$G$ has finitely many orbits in $\l_1$ naturally 
labelled by integers
from the segment $[0,r]$ such that
$\O_k$ corresponds to $\O_{r-k}$ via Pyasetskii pairing
(here we use the above identifications of $G$-orbits in $\l_1$
and $\l_{-1}$).
Therefore, if $r$ is even then the closure of the projectivization of
$\O_{r\over2}$ is a self-dual projective variety.

However, it is well-known that the closure of the projectivization of
$\O_i$ is smooth if and only if $i=1$ or $i=r$
(in the last case $\overline{\P(\O_r)}=\P(\l_1)$),
see e.g.~\cite{Pan} for the equivariant resolution
of singularities of $\overline{\O_i}$.
Therefore, the construction above gives a smooth self-dual
projective variety if and only if $r=2$.
It is worthy to write down all arising cases.
We use the notation from Section~\ref{ParabolicsAura}.

\begin{Example}\rm
Consider the short gradings of $\l=\sl_{n+2}$.
Then $\O_1$ is the variety of $n\times 2$-matrices of rank $1$.
The projectivization of $\O_1$ is identified with 
$X=\P^{n-1}\times\P^1$ in the Segre embedding.
\end{Example}

\begin{Example}\rm
Consider the short grading of $\l=\so_{n+2}$ that corresponds to $\beta$
being the first simple root.
The projectivization of $\O_1$ is identified with 
the quadric hypersurface in $\P^{n-1}$.
The quadric hypersurface in $\P^5$ also arises from the short
grading of $\l=\so_{8}=D_{4}$ that corresponds
to $\beta=\alpha_1$, $\beta=\alpha_3$, or $\beta=\alpha_4$
(Pl\"ucker quadric).
\end{Example}

\begin{Example}\rm
Consider the short grading of $\l=\so_{10}=D_{5}$ that corresponds
to $\beta=\alpha_4$ or $\beta=\alpha_5$.
The projectivization of $\O_1$ is identified with $X=\Gr(2,5)$
in the Pl\"ucker embedding.
\end{Example}

\begin{Example}\rm
The short grading of $E_6$ gives the 
spinor variety $\SS_5$.
\end{Example}

Remarkably, there are no other known examples
of smooth non-linear self-dual varieties.
Moreover, it is widely expected that this list is complete.
At least, this is true for polarized flag varieties:

\begin{Theorem}[\cite{Sn}]\label{LGHLJHFKJV}
Let $X=G/P\subset\P^N$ be a non-linear polarized flag variety.
If $\dim X=\dim\dual X$ the $X$ is one of the following
\begin{itemize}
\item A quadric hypersurface.
\item The Segre embedding of $\P^n\times\P^1$.
\item The Pl\"ucker embedding of $\Gr(2,5)$.
\item The $10$-dimensional spinor variety $\SS_5$.
\end{itemize}
\end{Theorem}

\proof
If $\dim X=\dim\dual X$ then either $X$ is a hypersurface
or $\defect X>0$. 
In the first case $X$ is necessarily a quadric hypersurface.
Indeed, let $Y\subset\dual{\P^N}$ be the projectivization of the highest
weight vector orbit of the dual representation.
Then $Y$ is isomorphic to $X$ as an embedded projective variety.
The intersection $\dual Y\cap X$ is non-empty and $G$-invariant,
therefore, $\dual Y=X$ and $X$ is a self-dual smooth hypersurface.
Then $X$ is a quadric hypersurface by Example~\ref{JHGFJHGFJGHFJH}.
An alternative proof follows from the fact that
the closure of the highest weight vector orbit 
is always cut out by quadrics~\cite{Li}.

If $\defect X>0$ then the claim easily follows from
Theorem~\ref{ThirdClassificationFlag}.
\endproof

Theorem~\ref{LGHLJHFKJV} was used in \cite{Sn}
to recover the classification of homogeneous real hypersurface\index{homogeneous real hypersurface}
in a complex projective space due to~\cite{Ta}.

\begin{Theorem}[\cite{Ta,Sn}]
Let $M$ be a homogeneous complete real hypersurface
embedded equivariantly in $\P^N$. Then $M$ is a tube
over a linear projective space or one of the $4$ self-dual
homogeneous spaces $X\subset \P^N$ listed in Theorem~\ref{LGHLJHFKJV}.
\end{Theorem}

\sketch
Let $M=K/L$, where $K$ is a compact Lie group.
The main idea is to use the fact that 
$M$ is necessarily a tube over a {\it complex} submanifold $X\subset\P^N$
called a focal submanifold, see e.g.~\cite{CR}.
Then $X$ is easily seen to be homogeneous, 
therefore $X$ is a flag variety, $X=G/P$, where $G$
is the complexification of $K$.
Moreover, $G$ acts transitively on the normal directions to $X$.
It follows that the  conormal variety
$\P(N_X\P^N)$ is also homogeneous, therefore the dual variety
$\dual X$ is also homogeneous being the image of the conormal variety.
Hence $\dual X$ is smooth and therefore $\dim X=\dim \dual X$ by 
Theorem~\ref{ZakThTangencies}.
Now we can apply Theorem~\ref{LGHLJHFKJV}.
\endproof

\subsection{Around Hartshorne Conjecture}
In \cite{Ha2} R.~Hartshorne has suggested a number of conjectures
related to the geometry of projective varieties of small codimension.
This work has stimulated a serie of remarkable researches,
for example the Zak's proofs of Theorem~\ref{ZakThLinNorm}
on linear normality and Theorem~\ref{ZakThOnSeveriVarieties}
on Severi varieties. Undoubtfully, the most famous conjecture
from this paper is the so-called 
Hartshorne conjecture on complete intersections:\index{Hartshorne conjecture}

\begin{Conjecture}[\cite{Ha2}]
If $X$ is a smooth $n$-dimensional projective variety in $\P^N$
and $\codim X<N/3$, then $X$ is a complete intersection.
\end{Conjecture}

This conjecture is still very far from being solved.
Only partial results are known, for example
the following Landsberg's theorem on
complete intersections. The proof is based on the technique
of Section~\ref{ProjectiveIISection}.
\begin{Theorem}[\cite{Lan1}]
Let $X\subset\P(V)=\P^N$ be an irreducible projective variety
cut out by quadrics, that is, the homogeneous ideal $I_{\Cone(X)}$ in 
$\C[V]$
is generated by its second homogeneous component $I_{\Cone(X)}^2$.
Let $b=\dim\Sing X$
{\rm(}set $b=-1$ if $X$ is smooth{\rm)}. 
If $\codim X<\bigover{N-(b+1)+3}{4}$, then $X$ is a complete intersection.
\end{Theorem}

If Hartshorne conjecture is true then any smooth projective variety $X$
in $\P^N$ such that $\codim X<N/3$
is a complete intersection. Therefore, $\dual X$ should be
a hypersurface by Theorem~\ref{DualCompleteIntersection}.
It is not known whether it is true or not.
In particular, if $X$ is a smooth self-dual variety
then, up to Hartshorne Conjecture, either $X$ is a quadric
hypersurface or $\codim X\ge N/3$.
In particular, the following theorem should give a complete
list of smooth self-dual varieties:

\begin{Theorem}[\cite{E1}]\label{JHGFJHGJGFJ}
Let $X$ be a nonlinear smooth projective variety
in $\P^N$. Assume that $\codim X\ge N/3$.
Suppose that $\dim X=\dim \dual X$. Then $X$ is one of the following varieties:
\begin{enumerate}
\item[\rm(i)] $X$ is a hypersurface in $\P^2$ or $\P^3$.
\item[\rm(ii)] $X$ is the Segre embedding of $\P^1\times\P^{n-1}$ in $\P^{2n-1}$.
\item[\rm(iii)] $X$ is the Pl\"ucker embedding of $\Gr(2,5)$ in $\P^9$.
\item[\rm(iv)] $X$ is the $10$-dimensional spinor variety $\SS_5$ in 
$\P^{15}$.
\end{enumerate}
\end{Theorem}

\sketch
Let $\dim X=n$.
We may assume that $n\ge3$ by Example~\ref{DualOfTheCurve} and Example~\ref{DualOfASurface}.
Now $\defect X=N-1-n$. Since $\defect X\le n-2$ by Example~\ref{FGGDCGFD},
we conclude that $n\ge\bigover{N+1}2$. If 
$n=\bigover{N+1}2$, then $\defect X=n-2$ and, therefore, we have case (ii) by
Example~\ref{FHGFGHJMFMGHC}.
So we may assume that $n\ge N/2+1$. Then $K_X=\O_X(\bigover{-N-1}2)$
by Theorem~\ref{EINKX} (d).
We conclude that $\defect X=N-1-n\le\bigover{n-2}2$ by
Theorem~\ref{ApplicationBeilinson} (a).
Hence $n\ge 2N/3$. By our assumption $n\le 2N/3$. Therefore, $n=2N/3$.
Now $\defect X=N-1-n={1\over2}n-1$. Therefore $n$ is even.
Since $\defect X\equiv n\mod2$ by 
Theorem~\ref{PairityTheorem}, we conclude that $n=4m+2$
and $\defect X=2m$. It follows that $m\le2$ by 
Theorem~\ref{ApplicationBeilinson} (b).
This leaves two cases: $X$ is either a $6$-dimensional variety in $\P^9$
or a $10$-dimensional variety in $\P^{15}$.
A very delicate treatment of both cases can be found in \cite{E1}.
\endproof

\subsection{Self--dual Nilpotent Orbits}

If $X\subset\P^n$ is a smooth projective variety
then $X$ is almost never self-dual. Moreover, up to Hartshorne
conjecture Theorem~\ref{JHGFJHGJGFJ} provides a complete list of them.
However, there are a lot of non-smooth self-dual varieties.
Many equivariant self-dual varieties 
are provided by Pyasetskii Theorem~\ref{PyasetskyTh}.
Perhaps, the most interesting examples
of self-dual varieties are
the Kummer surface in $\P^3$ \cite{GH1} and
the Coble quartic in $\P^7$ \cite{Pau}.

Another interesting examples were found in~\cite{Pop}.
Let $G$ be a connected reductive group acting
linearly on a vector space $V$. Suppose that
there exists
a $G$-invariant non-degenerate scalar product $(\cdot ,\cdot )$
on $V$,
in particular we have a natural isomorphism $V\simeq V^*$.
Finally, suppose that $G$ acts on the null-cone 
$${\cal R}(V)=\{v\in V\,|\,\overline{G\cdot v}\ni0\}$$
with finitely many orbits.
In particular, any $G$-orbit $\O\subset{\cal R}(V)$ is conical.

\begin{Theorem}[\cite{Pop}]
Let $\O=G\cdot v\subset{\cal R}(V)$ be a non-zero orbit, 
$X=\overline{\P(\O)}$
be the closure of its projectivization.
Then $X=\dual X$ if and only if $(\g\cdot v)^\perp\subset{\cal R}(V)$.
\end{Theorem}

\proof
Since $\O$ is conical, $v\in\g\cdot v$. Therefore,
$\g\cdot v$ is the affine cone $\Cone(\hat T_{[v]}X)$ over the embedded
tangent space $\hat T_{[v]}X$. 
Hence, $\dual X$ is equal to $\overline{\P(G\cdot(\g\cdot v)^\perp)}$.
Since $(\cdot,\cdot)$ is $G$-invariant, we have $v\in(\g\cdot v)^\perp$,
therefore, $X\subset\dual X$.
If $X=\dual X$ then, obviously, $(\g\cdot v)^\perp\subset{\cal R}(V)$.
Suppose that $(\g\cdot v)^\perp\subset{\cal R}(V)$.
Then $\dual X\subset\P({\cal R}(V))$.
Therefore, $\dual X=\overline{\P(\O')}$ 
for some orbit $\O'\subset{\cal R}(V)$.
Applying the same arguments we get that ${\dual X}\subset\dual{\dual X}$.
By Reflexivity Theorem, we finally obtain $X=\dual X$.
\endproof

Let $G$ be semisimple, $V=\g$ be the adjoint representation.
Then the Killing form $(\cdot,\cdot)$ is $G$-invariant, ${\cal R}(\g)$
is the cone of nilpotent elements. Since there are finitely many
nilpotent orbits, we can apply the previous Theorem.
For any $x\in\g$ the orthogonal complement $(\ad(\g)x)^\perp$
is identified with the centralizer $\g_x$ of $x$ in $\g$.
Therefore, a nilpotent orbit $\O=\Ad(G)x$ has a self-dual
projectivization if and only if the centralizer $\g_x$
belongs to ${\cal R}(\g)$, i.e.~has no semi-simple elements.
These orbits are called the distinguished nilpotent orbits.
For example, the regular nilpotent orbit (the orbit dense in ${\cal R}(\g)$)
is distinguished, hence the null-cone ${\cal R}(\g)$
itself has a self-dual projectivization.

The distinguished nilpotent orbits are, in a sense, the building
blocks for the set of all nilpotent orbits. Namely, 
due to the Bala--Carter
correspondence~\cite{BCa1,BCa2} the set of all nilpotent orbits
in a semisimple Lie algebra $\g$
is in the natural bijection with the set of isoclasses of pairs $(\l,\O)$,
where $\l\subset\g$ is a Levi subalgebra and $\O\subset\l$
is a distinguished nilpotent orbit in $[\l,\l]$.

\section{Linear Systems of Quadrics of Constant Rank}

\index{linear systems of quadrics}%
\begin{Theorem}[\cite{IL}]\label{LandsbergThe1}
Let $X\subset\P^N$ be a smooth $n$-dimensional variety.
If $H\in\dual{X}_{sm}$ then the projective second fundamental
form $|II|_H$ of $\dual X$ at $H$ is a system of quadrics
of projective dimension $\defect X$ and constant rank $n-\defect X$.
\end{Theorem}

This Theorem (analogous to Theorem~\ref{IntermediateEinTh})
leads to the examination of linear systems of quadrics of constant
rank. 

More generally, let $V=\C^m$, $W=\C^n$ and let $A\subset\Hom(V,W)$
(resp.~$A\subset S^2V^*$, resp.~$A\subset \Lambda^2V^*$)
be a linear subspace. 

\index{linear systems of bounded rank}%
\index{linear systems of constant rank}%
\index{linear systems of rank bounded below}%
\begin{Definition}\rm
$A$ is said to be of bounded rank $r$ if for all
$f\in A$, $\rank f\le r$.
$A$ is said to be of constant rank $r$ if for all non-zero
$f\in A$, $\rank f=r$.
$A$ is said to be of rank bounded below by $r$ if for all non-zero
$f\in A$, $\rank f\ge r$.

We define the numbers
$$l(r,m,n)=
\max\{\dim A\,|\,A\subset\Hom(V,W)\hbox{\rm \ is of constant rank $r$}\}$$
$$\underline l(r,m,n)=
\max\{\dim A\,|\,A\subset\Hom(V,W)\hbox{\rm \ is of bounded rank $r$}\}$$
$$\overline l(r,m,n)=
\max\{\dim A\,|\,A\subset\Hom(V,W)\hbox{\rm \ is of rank bounded below by $r$}\}$$
Similarly, define $c(r,m)$, $\underline c(r,m)$, and $\overline c(r,m)$
in the symmetric case and
$\lambda(r,m)$, $\underline\lambda(r,m)$, and $\overline\lambda(r,m)$
in the skew-symmetric case.
Recall that any skew-symmetric matrix has even rank.
\end{Definition}

\begin{Proposition}[\cite{IL}]
The upper bounds for dimensions of linear subspaces
of matrices of rank bounded below by $r$ are as follows.
\begin{enumerate}
\item[\rm (a)] $\ov l(r,m,n)=(m-r)(n-r)$;
\item[\rm (b)] $\ov c(r,m)={m-r+1\choose 2}$;
\item[\rm (c)] $\ov \lambda(r,m)={m-r\choose 2}$, $r$ even.
\end{enumerate}
\end{Proposition}

\proof
Let 
$$X_r=\{f\in\P(\Hom(V,W))\,|\,\rank f\le r\}$$
and similarly for $X_r^c\subset\P(S^2V)$
and $X_r^\lambda\subset\P(\Lambda^2 V)$.
Clearly, $\ov l(r,m,n)=\codim X_r$, $\ov c(r,m)=\codim X_r^c$,
and $\ov \lambda(r,m)=\codim X_r^\lambda$.
Now the theorem follows by an easy dimension count.
\endproof

Given a linear subspace $A\subset\Hom(V,W)$ we can also
consider $A$ as a linear subspace $A^c$ of $S^2(V^*\oplus W)$
or a linear subspace $A^\lambda$ of $\Lambda^2(V^*\oplus W)$. 
In matrix form,
for any matrix $a\in A$ we associate a symmetric matrix $a^c\in A^c$
or skew-symmetric matrix $a^\lambda\in A^\lambda$ of the form
$$a^c=\left(\matrix{0&a\cr a^t&0}\right),\quad
a^\lambda=\left(\matrix{0&a\cr -a^t&0}\right).$$
If $A$ has constant rank $r$ (resp.~bounded rank $r$, resp.~rank bounded below by $r$)
then $A^c$ and $A^\lambda$  have
constant rank $2r$ (resp.~bounded rank $2r$, resp.~rank bounded below by $2r$).
$A^c$ and $A^\lambda$ are called
doublings of $A$. \index{doublings of linear systems}%

\begin{Proposition}[\cite{IL}]\label{JGHFJGCDJNFGCNC}
The lower bounds for dimensions of linear subspaces
of matrices of constant rank are as follows.
\begin{enumerate}
\item[\rm (a)] If $0<r\le m\le n$ then $l(r,m,n)\ge n-r+1$.
\item[\rm (b)] If $r\ge2$ is even then $c(r,m)\ge m-r+1$ and $\lambda(r,m)\ge m-r+1$.
\end{enumerate}
\end{Proposition}

\proof
$X_{r-1}\subset\P(\Hom(\C^r,\C^n)$ has codimension $n-r+1$.
Thus we can find an $(n-r+1)$-dimensional subspace
$A\subset\Hom(\C^r,\C^n)$ of constant rank $r$ and (a) follows.
The doublings of $A$ produce $(n-r+1)$-dimensional
subspaces $A^c\subset S^2\C^{r+n}$ and $A^\lambda\subset \Lambda^2\C^{r+n}$
of constant rank $2r$. This implies (b).
\endproof

Any linear map $\psi:\,A\to\Hom(V,W)$ can be viewed as an element of 
\begin{eqnarray*}
A^*\otimes\Hom(V,W)=H^0(\P(A),\O_{\P(A)}(1))\otimes\Hom(V,W)\\
=H^0(\P(A),\O_{\P(A)}(1)\otimes\Hom(V,W)).
\end{eqnarray*}
Thus to give a linear map $A\to\Hom(V,W)$ (say, embedding)
is the same as to give a vector bundle map
$$\ov\psi:\,V\otimes\O_{\P(A)}\to W\otimes\O_{\P(A)}(1)$$
on $\P(A)$.
Let $K$, $N$, and $E$ be the kernel, the cokernel, and the image of $\ov\psi$.
Then it is easy to see that $\psi$ is an embedding of constant rank $r$
if and only if all $K$, $N$, and $E$ are vector bundles of ranks
$\dim V-r$, $\dim W-r$, and~$r$ if and only
if any of $K$, $N$, and $E$ is a vector bundle of rank as above.
Therefore, all embeddings $A\subset\Hom(V,W)$ of constant rank are
in one-one correspondence with exact sequences of vector
bundles
$$0\to K\to \O_{\P(A)}\otimes V\to E\to 0\hbox{\rm\ and\ }
0\to E\to \O_{\P(A)}(1)\otimes W\to N\to 0.$$
Similarly, 
all embeddings $A\subset\S^2 V$ (resp.~$A\subset\Lambda^2V$)
of constant rank are
in one-one correspondence with exact sequences of vector
bundles
$$0\to K\to \O_{\P(A)}\otimes V^*\to E\to 0,$$
where $E$ is a rank $r$ vector bundle such that $E\simeq E^*(1)$
and this isomorphism is symmetric (resp.~skew-symmetric).

\begin{Theorem}[\cite{IL}]\label{JKGFTDRJKKJGH}$\ $
\begin{enumerate}
\item[\rm(a)] Let $A\subset S^2 V$ or $A\subset \Lambda^2V$ be a constant rank $r$
subspace, $\dim A\ge2$. Then $E$ is a uniform vector bundle
of splitting type $\O^{r/2}_{\P^1}\oplus\O^{r/2}_{\P^1}(1)$.
In particular, $r$ is even.
\item[\rm(b)] If $r$ is odd then $c(r,m)=\lambda(r,m)=1$.
\end{enumerate}
\end{Theorem}

\index{uniform vector bundle}%

\proof
(a) The proof is by the same argument as in the proof of Theorem~\ref{SecondEin} (b).

(b) This follows from (a). However, in the symmetric case
it was also known classically as a consequence of the Kronecker--Weierstrass
theory giving a normal form for pencils of symmetric matrices
of bounded rank, see~\cite{Gan,HP,Me}.
\endproof

\index{Kronecker--Weierstrass theory}%
\index{pencils of symmetric matrices}%

\begin{Theorem}[\cite{E1,IL}]\label{JHGDFJGDSG}
A smooth $n$-dimensional projective variety $X\subset\P^N$ with positive defect $\defect X$
determines a $(\defect X+1)$-dimensional linear subspace
$A\subset S^2\C^{N-1-\defect X}$ with associated short exact sequence
$$0\to N_X^*H(1)|_L\to\O_L\otimes\C^{N-1-\defect X}\to N_LX\to0,$$
where $H\in\dual X_{sm}$ and $L=\Sing X\cap H$ is the contact locus.
\end{Theorem}

\proof
Recall that $L$ is a projective subspace of dimension $\defect X$.
We have the following exact sequence
$$0\to N_LX\to N_LH\to N_XH|_L\to0$$
which holds because $H$ is tangent to $X$ along $L$.
Dualising and twisting by $\O_L(1)$ we get
$$0\to N_X^*H(1)|_L\to N_L^*H(1)\to N_L^*X(1)\to0.$$
It remains to use the fact that $N_LH\simeq\O_L(1)\otimes\C^{N-1-\defect X}$
and, most importantly, that $N_LX\simeq N_L^*X(1)$
by Theorem~\ref{SecondEin} (a).
\endproof

The following theorem is a consequence
of the general theory of uniform vector bundles~\cite{Sa}.

\begin{Theorem}[\cite{E1,IL}]$\ $
\begin{enumerate}
\item[\rm(a)] Suppose that $A\subset S^2V$ is a linear $(l+1)$-dimensional
subspace of constant rank $r$. Let
$0\to K\to \O_{\P(A)}\otimes V^*\to E\to 0$
be the corresponding exact sequence of vector bundles.
If $r\le l$ then $E\simeq\O_{\P(A)}^{r/2}\oplus\O_{\P(A)}^{r/2}(1)$
unless $r=l=2$ and $E=T\P^2(-1)$.
\item[\rm(b)] Suppose that $X\subset\P^N$ is a smooth projective variety.
Let $L$ be a contact locus with a generic tangent hyperplane.
If $\defect X\ge N/2$, then 
$$N_LX=\O_L^{(\dim X-\defect X)/2}\oplus\O_L^{(\dim X-\defect X)/2}(1).$$
\end{enumerate}
\end{Theorem}

\proof
(a) Indeed, $E$ is uniform by~Theorem~\ref{JKGFTDRJKKJGH} (a).
Therefore, if $r\le l$ then $E$ either splits into a direct sum 
of line bundles or $E$ is isomorphic to $T\P^l(a)$ or $T^*\P^l(a)$
by~\cite{Sa}. If $E$ splits then 
$E\simeq\O_{\P(A)}^{r/2}\oplus\O_{\P(A)}^{r/2}(1)$
by~Theorem~\ref{JKGFTDRJKKJGH} (a).
Suppose now that $E$ is isomorphic to $T\P^l(a)$ or $T^*\P^l(a)$.
Since $E\simeq E^*(1)$, the first Chern class $c_1(E)$
is equal to $(l/2)H$, where $H$ is the class of the hyperplane section.
Taking first Chern classes we conclude that $E$ is isomorphic either to
$T\P^2(-1)$ or to $T^*\P^2(2)$, but these last two bundles are isomorphic.

(b) Let $\dim X=n$ and $\defect X=k$.
By (a) and by~Theorem\ref{JHGDFJGDSG} we see that
$N_LX$ is isomorphic to either
$\O_L^{(n-k)/2}\oplus\O_L^{(n-k)/2}(1)$
or $T\P^2(-1)$. However, in the last case $n=4$, $k=2$,
therefore $X$ is a projective bundle $\P_C(F)$ over a curve $C$
by Example~\ref{FGGDCGFD}.
Then $L$ is embedded as a $2$-plane in a fibre $f$ of $\P_C(F)$.
Consider the exact sequence
$$0\to N_Lf\to N_LX\to N_fX|_L\to 0.$$
Since $N_Lf=\O_L(1)$ and $N_fX=\O_f$, this exact sequence reduces to
$$0\to\O_L(1)\to N_LX\to\O_L\to0.$$
It follows that $N_LX=\O_L\oplus\O_L(1)$.
\endproof

The following theorem was proved in \cite{W}.
Weaker results (but similar arguments) first appeared in~\cite{Sy}.
\begin{Theorem}[\cite{W}]\label{JHGFJHCDHFDHFC}
Suppose $2\le r\le m\le n$. Then
\begin{enumerate}
\item[\rm(a)] $l(r,m,n)\le m+n-2r+1$.
\item[\rm(b)] $l(r,m,n)=n-r+1$ if $n-r+1$ does not divide $(m-1)!/(r-1)!$.
\item[\rm(c)] $l(r,r+1,2r-1)=r+1$.
\end{enumerate}
\end{Theorem}

\sketch
Let $A\subset\Hom(V,W)$ be an $l$-dimensional linear subspace
of constant rank $r$, $\dim V=m$, $\dim W=n$, $m\ge n$.
Let 
$$0\to K\to \O_{\P(A)}\otimes V\to E\to 0\hbox{\rm\ and\ }
0\to E\to \O_{\P(A)}(1)\otimes W\to N\to 0$$
be the corresponding exact sequences of vector bundles.
Then we have the following equations of Chern classes
$c(K)c(E)=1$ and $c(E)c(N)=(1+H)^n$, where $H$ is the class
of the hyperplane section.
Thus, $c(K)(1+H)^n=c(N)$. 
If $n-r+1\le i\le l-1$ then $c_i(N)=0$
and looking at the coefficient of $H^i$
we get
$$\sum_{j=0}^{m-r}{n\choose i-j}k_j=0,$$
where $c_j(K)=k_jH^j$ and we use the convention
that ${n\choose i}=0$ if $i<0$ or $i>n$.
The coefficient matrix of this system of linear equations
is 
$$\left({n\choose i-j}\right)_{0\le j\le m-r,\ n-r+1\le i\le l-1}.$$
If $l=m+n-2r+2$ then this is square invertible matrix with determinant
$\prod_{j=0}^{m-r}j!$.
Thus $k_0=0$ which is a contradiction since $k_0=1$.
This proves (a).

Proofs of (b) and (c) can be found in \cite{W}. Notice that by 
Proposition~\ref{JGHFJGCDJNFGCNC} (a) we have
$l(r,m,n)\ge n-r+1$, so only the opposite inequality
should be proved.
\endproof

Another possible way of proving Theorem~\ref{JHGFJHCDHFDHFC} (a)
is to use the following theorem of R.~Lazarsfeld
(its application is immediate):

\begin{Theorem}[\cite{La1}]\label{JHGFJHFDHFGCDHG}
Let $X$ be $l$-dimensional projective variety. Let $E$ and $F$
be vector bundles of ranks $m$ and $n$, respectively.
Suppose that $E^*\otimes F$ is ample and there is a constant rank $r$
vector bundle map $E\to F$. Then $l\le m+n-2r$.
\end{Theorem}

The (partial) symmetric analogue of this theorem was found in \cite{IL}.
Theorems~\ref{JHGFJHFDHFGCDHG} and~\ref{JKHGFJHGFCK} (and their proofs)
are analogous to results 
related to the problem of the non-emptiness and the connectedness
of degeneracy loci, see \cite{La1,FL1,HT,HT1,Tu,IL}. \index{degeneracy locus}

\begin{Theorem}[\cite{IL}]\label{JKHGFJHGFCK}
Let $X$ be a smooth simply connected $m$-dimensional variety.
Let $E$ be a rank $n$ vector bundle on $X$, and $L$ a line bundle on $X$.
Suppose that $S^2E^*\otimes L$ is an ample vector bundle
and that there is a constant even rank $r$ symmetric bundle map $E\to E^*\otimes L$.
Then $m\le n-r$.
\end{Theorem}

\proof
Consider the projective bundle $\pi:\,\P(E)\to X$.
The symmetric bundle map $E\to E^*\otimes L$ defines
a section $s\in H^0(X,S^2E^*\otimes L)$ which pulls back to give a section
$\pi^*s\in H^0(\P(E),\pi^*(S^2E^*\otimes L))$.
There is a natural map $\pi^*\pi_*\O_{\P(E)}(2)\to\O_{\P(E)}(2)$.
Since $\pi_*\O_{\P(E)}(2)\simeq S^2E^*$ this gives after tensoring
with $\pi^*L$ a map $u:\,\pi^*(S^2E^*\otimes L)\to\O_{\P(E)}(2)\otimes\pi^*L$.
Let $t=u\circ\pi^*s\in H^0(\P(E),\O_{\P(E)}(2)\otimes\pi^*L)$.
Let $Y\subset\P(E)$ be the zero locus of~ $t$.

Let $x\in X$. The section $t$ at $[v]\in\P(E(x))$ is the linear map
$\lambda\to s(x)(v,v)\lambda$ and so vanishes if and only if
$s(x)(v,v)=0$. But by the hypothesis, this defines
a rank $r$ quadric hypersurface in $\P(E(x))$.
Therefore, the fiber of $\pi_Y:\,Y\to X$ over each $x\in X$
is a rank $r$ quadric hypersurface in $\P(E(x))$.

We claim that $\O_{\P(E)}(2)\otimes\pi^*(L)$ is ample on $\P(E)$.
Indeed, since $S^2E^*\otimes L$ is ample on $X$,
$\O_{\P(S^2E^*\otimes L)}(1)$ is ample on $\P(S^2E^*\otimes L)$.
Clearly, $\P(S^2E^*\otimes L)\simeq\P(S^2E^*)$
and $\O_{\P(S^2E^*\otimes L)}(1)\simeq\O_{\P(S^2E^*)}(1)\otimes\sigma^*L$,
where $\sigma:\,\P(S^2E^*)\to X$ is the natural projection.
Therefore, 
$\O_{\P(S^2E^*)}(1)\otimes\sigma^*L$ is ample on $\P(S^2E^*)$.
The second Veronese map gives an inclusion $i:\,\P(E)\subset\P(S^2E^*)$
such that $\pi=i\circ \sigma$. Then $i^*\O_{\P(S^2E^*)}(1)=\O_{\P(E)}(2)$.
So $i^*(\O_{\P(S^2E^*)}(1)\otimes\sigma^*L)\simeq\O_{\P(E)}(2)\otimes\pi^*L$
is ample on $\P(E)$ as required.
In particular, $\P(E)\setminus Y$ is an affine variety
and its homology vanishes above its middle dimension: 
$H_i(\P(E)\setminus Y,\C)=0$
for $i\ge m+n$.

Let $K$ denote the kernel and $F$ denote the image
of the map $E\to E^*\otimes L$.
Then since the map has constant rank $r$, $K$ and $F$
are vector bundles on $X$ of ranks $n-r$ and $r$ respectively.
Since the map is symmetric, there is a symmetric isomorphism $F\simeq F^*\otimes L$.
We have a natural map $p:\,\P(E)\setminus\P(K)\to\P(F)$
given on the fibers by the linear projection centered at $\P(K(x))\subset\P(E(x))$.
Then $p$ is a $\C^{n-r}$-fiber bundle.

The isomorphism $F\simeq F\otimes L$ determines a hypersurface $Z\subset\P(F)$
such that if $\rho:\,\P(F)\to X$ is the natural projection
then the fiber of $\rho|_Z:\,Z\to X$ over each $x\in X$
is a smooth quadric in $\P(F(x))$. Now $\P(K)\subset Y$ and so $p$
restricts to a $\C^{n-r}$-fiber bundle $\P(E)\setminus Y\to\P(F)\setminus Z$.
Thus, $H_i(\P(E)\setminus Y,\C)\simeq H_i(\P(F)\setminus Z,\C)$
and by Lefschetz duality this is isomorphic to $H^{2(m+r-1)-i}(\P(F),Z,\C)$.
Hence $H^i(\P(F),Z,\C)=0$ for $i\le 2r+m-n-2$
and then using the long exact sequence of the pair $(\P(F),Z)$
we conclude that $H^i(Z,\C)\simeq H^i(\P(F),\C)$ for $i\le 2r+m-n-3$.
Both bundles $\rho:\,\P(F)\to X$ and $\rho|_Z:\,Z\to X$
have smooth fibers and are defined over smooth simply-connected
variety. Therefore, by Deligne' theorem \cite{GH1},
the Leray spectral sequence for $\rho$ and $\rho|_Z$ degenerates
at the $E_2$ term.
It follows that $H^*(\P(F),\C)=H^*(X,\C)\otimes H^*(\P^{r-1},\C)$
and $H^*(Z,\C)=H^*(X,\C)\otimes H^*(Q,\C)$,
where $Q\subset\P^{r-1}$ is a smooth quadric.
Let $b_j(B)=\dim H^j(B,\C$ be the Betti numbers.
Then it is well-known that 
$$b_i(\P^{r-1})=\cases{
1&\hbox{\rm for even $i$, $0\le i\le 2(r-1)$}\cr
0&\hbox{\rm otherwise.}}$$
Since $r$ is even, we have
$$b_i(Q)=\cases{
1&\hbox{\rm for even $i$, $0\le i\le 2(r-2)$}\cr
2&\hbox{\rm $i=r-2$}\cr
0&\hbox{\rm otherwise.}}$$
It follows that $b_{r-2}(\P(F)\ne b_{r-2}(Z)$.
Since $b_i(\P(F)=b_i(Z)$ for 
$i\le 2r+m-n-3$, we finally obtain
$r-2>2r+m-n-3$. Equivalently, $m\le n-r$.
\endproof

As an easy application of Theorem~\ref{JKHGFJHGFCK} we have
the following theorem.

\begin{Theorem}[\cite{IL}]\label{LandsbergThe2}
Let $V$ be an $m$-dimensional vector space, 
let $r>0$ be an even number
$$\max\{\dim A\,|\,A\subset\P(S^2V^*)\hbox{\rm \ is of constant rank $r$}\}=
m-r.$$
\end{Theorem}

\begin{Example}\rm
Theorems~\ref{LandsbergThe1} 
and~\ref{LandsbergThe2} combined
furnish a new proof of 
the formula $\dim\dual X\ge\dim X$ (Theorem~\ref{ZakThTangencies} (b)).
\end{Example}

\cleardoublepage

\addcontentsline{toc}{chapter}{Index}
\flushbottom\printindex
\end{document}